\documentclass[fleqn,10pt]{article}
\usepackage{latexsym, graphicx, epsfig, amsmath, amssymb,amsfonts}
\usepackage{natbib,amsthm,version}
\usepackage{amsbsy,bm,multirow,enumerate}
\usepackage[titletoc,page]{appendix}
\usepackage[mathscr]{eucal}
\usepackage{mathtools}
\usepackage{color}
\usepackage[utf8]{inputenc}
\usepackage[english]{babel}
\usepackage{amsthm}
\usepackage{enumerate}
\usepackage[hidelinks]{hyperref}
\usepackage{url}
\usepackage{subfigure}
\usepackage[lined,boxed,linesnumbered,ruled]{algorithm2e}
\usepackage{float}

\DeclareMathOperator{\sech}{sech}

\newcommand{\mbs}[1]{\mathbf{#1}}

\newtheorem{remark}{Remark}[section]

\theoremstyle{definition}

\title{
  Local Extreme Learning Machines and Domain
  Decomposition for Solving Linear and Nonlinear Partial Differential Equations
} 
\author{
  Suchuan Dong$^1$\thanks{Author of correspondence.
    Email: sdong@purdue.edu}, \ Zongwei Li$^2$
  \\
  $^1$Center for Computational and Applied Mathematics \\
  Department of Mathematics, 
  Purdue University \\
  West Lafayette, Indiana, USA \\
  $^2$Department of Mathematics,
  Purdue University\\
  Fort Wayne, Indiana, USA
 } 

\date{(December 4, 2020)}

\begin{document}
\maketitle


\begin{abstract}

  We present a neural network-based method for solving linear and
  nonlinear partial differential equations, by combining 
  the ideas of extreme learning machines (ELM), domain decomposition
  and local neural networks.
  The field solution on each sub-domain is represented by a local feed-forward
  neural network, and $C^k$ continuity with an appropriate
  integer $k$ is imposed on the sub-domain boundaries.
  Each local neural network consists
  of a small number (one or more) of hidden layers, while its last hidden layer can
  be wide. The weight/bias coefficients in all the hidden layers of the
  local neural networks are
  pre-set to random values and fixed throughout the computation,
  and only the weight coefficients in the output layers of the local
  neural networks are adjustable training parameters.
  The overall neural network is trained by a linear or
  nonlinear least squares computation, not by the back-propagation type
  algorithms. We introduce a block time-marching scheme together with
  the presented method for long-time simulations of time-dependent
  linear/nonlinear partial differential equations.
  The current method exhibits a clear sense of convergence with respect to
  the degrees of freedom in the neural network.
  Its numerical errors typically decrease exponentially or nearly exponentially
  as the number of degrees of freedom (e.g.~the number of training
  parameters, number of training data points, number of sub-domains)
  in the system increases. Extensive numerical experiments have been
  performed to demonstrate the computational performance of the current method
  and to study the effects of the simulation parameters. 
  We also present results to demonstrate its capability
  for long-time dynamic simulations with certain test problems.
  We compare the presented method with the deep Galerkin method (DGM)
  and the physics-informed neural network (PINN) method in terms of
  the accuracy and computational cost. The current method exhibits a clear
  superiority, with its numerical errors and network training time
  considerably smaller (typically by orders of magnitude) than those of DGM and PINN.
  We also compare the current method with the classical finite element method (FEM).
  The computational performance of the current method is on par with,
  and often exceeds, the FEM performance in terms of the accuracy and
  computational cost. To achieve the same accuracy, the network training time of
  the current method is comparable to, and oftentimes less than,
  the FEM computation time. Under the same
  computational cost (training/computation time),
  the numerical errors of the current method are comparable to, and oftentimes
  markedly smaller than, the FEM errors.

\end{abstract}


\vspace{0.05cm}

\tableofcontents


\section{Introduction}
\label{sec:intro}

%
%
%


Neural network based numerical methods,
especially those based on deep learning~\cite{GoodfellowBC2016},
have attracted a significant amount of research in the past few years
for simulating the governing partial differential equations (PDE) of physical
phenomena. These methods provide a new way for approximating the field solutions,
in the form of deep neural networks (DNN), which is different from the ansatz space
with traditional numerical methods such as finite difference or finite element techniques.
This can be a promising approach, potentially more effective and more efficient
than the traditional methods, for solving the governing PDEs
of scientific and engineering importance.
DNN-based methods solve the PDE by transforming the solution finding problem
into an optimization problem.
They typically parameterize the PDE solution by the
training parameters in a deep neural network, in light of the universal approximation
property of DNNs~\cite{HornikSW1989,HornikSW1990,Cotter1990,Li1996}.
Then these methods attempt to minimize a loss function that consists of
the residual norms of the governing equations and also the associated boundary
and initial conditions,
typically by some flavor of gradient descent type techniques (i.e.~back
propagation algorithm~\cite{Werbos1974,Haykin1999}).
This process constitutes the predominant computations in the DNN-based PDE
solvers, commonly known as the training of the neural network.
Upon convergence of the training process,
the solution is represented by the neural network,
with the training parameters set according to their converged values.
Several successful DNN-based PDE solvers have emerged 
in the past years,
such as the deep Galerkin method (DGM)~\cite{SirignanoS2018},
the physics-informed neural network (PINN)~\cite{RaissiPK2019}, and
related approaches
(see e.g.~\cite{LagarisLF1998,LagarisLP2000,RuddF2015,EY2018,HeX2019,ZangBYZ2020,Samaniegoetal2020,Xu2020},
among others).
Neural network-based PDE solutions are smooth analytical functions,
depending on the activation functions used therein.
The solution and its derivatives can be computed exactly,
by evaluation of the neural network or by auto-differentiation~\cite{BaydinPRS2018}.


While their computational performance is promising,
DNN-based PDE solvers, in their current
state, suffer from a number of limitations that make them
numerically less than satisfactory and
computationally uncompetitive.
The first limitation is the solution accuracy of DNN-based
methods~\cite{JagtapKK2020}.
A survey of related literature indicates that the absolute error of
the current DNN-based methods is generally on, and rarely goes below,
the level of $10^{-3}\sim 10^{-4}$.
Increasing the resolution or the number of training epochs/iterations
does not notably improve this error level.
The accuracy of such levels is less than satisfactory for scientific computing,
especially considering that the classical numerical methods
can achieve the machine accuracy given sufficient mesh resolution and
computation time.
Perhaps because of such limited accuracy levels, a sense of convergence with
a certain convergence rate is generally lacking with the DNN-based PDE solvers.
For example, when the number of layers, or the number of nodes within
the layers, or the number of training data points
is varied systematically, one can hardly observe
a consistent improvement in 
the accuracy of the obtained simulation results.
Another limitation concerns the computational cost.
The computational cost of DNN-based PDE solvers is extremely high.
The neural network of these
solvers takes a considerable amount of time to train, in order to reach a reasonable
level of accuracy.
For example, a DNN-based PDE solver can take hours to train
to reach a certain accuracy, while with a traditional numerical method
such as the finite element method it may take only a few seconds
to produce a solution with the same or better accuracy.
Because of their limited accuracy and large computational cost,
there seems to be a general sense that the DNN-based PDE solvers, at least
in their current state, cannot
compete with classical numerical methods, except perhaps for certain  problems
such as high-dimensional
PDEs which can be challenging to classical methods due to
the so-called curse of dimensionality.


In the current work we concentrate on the accuracy and the computational cost
of neural network-based numerical methods.
We introduce a
neural network-based method for solving linear and nonlinear
PDEs that exhibits a disparate computational performance 
from the above DNN-based PDE solvers.
The current method exhibits a clear sense of convergence with respect to
the degrees of freedom in the system.
Its numerical errors typically decrease exponentially or nearly
exponentially as the
number of degrees of freedom (e.g.~the number of training
parameters, number of training data points) in the network increases.
In terms of the accuracy and computational cost,
it exhibits a clear superiority to the often-used DNN-based PDE solvers.
Extensive comparisons with the deep Galerkin method~\cite{SirignanoS2018}
and the physics-informed neural network~\cite{RaissiPK2019} are presented
in this paper.
The numerical errors, and the network training time,
of the current method are typically orders of magnitude
smaller than those of DGM and PINN.
The computational performance of the current method
is competitive compared with traditional numerical methods.
Extensive comparisons with the classical finite element method (FEM)
are provided. The performance of the current method is on par with,
and often exceeds, the performance of FEM with regard to the
accuracy and computational cost.
For example,
to achieve the same accuracy, the network training time of the current
method is comparable to, and oftentimes smaller than, the FEM computation
time. With the same computational cost (training/computation time),
the numerical errors of 
the current method are comparable to, and oftentimes
markedly smaller than, those of the FEM.


The superior computational performance of the current
method  can be
attributed to several of its algorithmic characteristics:
\begin{itemize}

\item
  Network architecture and training parameters. The current method
  is based on shallow feed-forward neural networks. Here ``shallow''
  refers to the configuration that the network contains only
  a small number (e.g.~one, two or three) of hidden layers, while the
  last hidden layer can be wide.
  The weight/bias coefficients in all the hidden layers
  are pre-set to random values and are fixed, and they are not
  training parameters. 
  The training parameters consist of the weight coefficients
  of the output layer.

\item
  Training method. The network is trained and the values for the training
  parameters are determined by a least squares computation, not
  by the back propagation (gradient descent-type) algorithm.
  For linear PDEs, training the neural network involves a
  linear least squares computation.
  For nonlinear PDEs, the network training involves a nonlinear
  least squares computation.

\item
  Domain decomposition and local neural networks.
  We partition the overall domain into sub-domains, and represent the
  solution on each sub-domain
  locally by a shallow feed-forward neural network. $C^k$ continuity
  conditions, where $k\geqslant 0$ is an integer related
  to the PDE order, are enforced across sub-domain boundaries.
  The local neural networks collectively form a
  multi-input multi-output logical network model, and are trained
  in a coupled way with the linear or nonlinear least squares computation.

\item
  Block time marching. For long-time simulations of
  time-dependent PDEs, the current method adopts a block time-marching strategy.
  The overall spatial-temporal domain is first
  divided into a number of windows in time, referred to as
  time blocks. The PDE is then solved on the spatial-temporal
  domain of each time block, individually and successively.
  Block time marching is crucial to long-time simulations, especially for
  nonlinear time-dependent PDEs.

\end{itemize}


The idea of random weight/bias coefficients in the network
and the use of linear least
squares method for network training stem from the so-called
extreme learning machines (ELM)~\cite{HuangZS2006,HuangWL2011}.
ELM was developed for single-hidden layer feed-forward neural networks (SLFN),
and for linear problems. It transforms the linear
classification or regression
problem into a system of linear algebraic equations, which is then
solved by a linear least squares method or by
using the pseudo-inverse (Moore-Penrose inverse)
of the coefficient matrix~\cite{GolubL1996}.
ELM is one example of the so-called randomized neural
networks (see e.g.~\cite{PaoPS1994,IgelnikP1995,MaassM2004,JaegerLPS2007,ZhangS2016}),
which can be traced to Turing's unorganized machine and Rosenblatt's
perceptron~\cite{Webster2012,Rosenblatt1958} and have witnessed a revival
in neuro-computations in recent years.
%
The application of ELM to function approximations and linear differential
equations have been considered in several recent
works~\cite{BalasundaramK2011,YangHL2018,Sunetal2019,PanghalK2020,LiuXWL2020,DwivediS2020}.
Domain decomposition has found widespread applications in classical numerical
methods~\cite{SmithBG1996,ToselliW2005,Dong2010,DongS2015,Dong2018}. Its use in neural network-based methods, however,
has been very limited and is very recent
(see e.g.~\cite{LiTWL2020,JagtapKK2020,DwivediS2020}).


The contribution of the current work lies in several aspects.
A main contribution of this work is the introduction of an ELM-like method
for nonlinear differential equations, based on domain decomposition and local neural networks.
In contrast, existing ELM-based methods for differential equations have been confined to
linear problems, and the neural network is limited to a single hidden layer.
For nonlinear problems, to solve the resultant nonlinear
algebraic system about the training parameters, we have adopted two methods:
(i) a nonlinear least squares method with perturbations (referred to as NLSQ-perturb),
and (ii) a combined Newton/linear least squares method (referred to as Newton-LLSQ).
We find that the random perturbation in the NLSQ-perturb method
is crucial to preventing the method from being trapped to local minima with cost values
exceeding some given tolerance, especially in under-resolved cases and in long-time simulations.
We present an algorithm for effective generation of the random perturbations
for the nonlinear least squares method.

Another contribution of the current work is the afore-mentioned
block time-marching scheme for long-time simulations of time-dependent
linear/nonlinear PDEs. When the temporal dimension of the spatial-temporal domain
is large, if the PDE is solved on the entire domain all at once,
we find that the neural network becomes very hard to train
with the ELM algorithm (and also with the back propagation-based algorithms),
in the sense that the obtained solution can contain pronounced errors, especially toward
later time instants in the spatial-temporal domain.
On the other hand, by using the block time-marching strategy
and with a moderate time block size,
the problem becomes much easier to solve and the neural network is
much easier to train with the ELM algorithm.
Accurate results can be attained with
the block time-marching scheme for very long-time simulations.
The block time marching strategy is often crucial to the simulations of
nonlinear time-dependent PDEs when the temporal dimension becomes
even moderately large.

We would also like to emphasize that, with the current method,
each local neural network is not limited to a single hidden layer,
which is another notable difference from existing ELM-type methods.
Up to three hidden layers in the local neural networks have been tested
in the current paper. We observe that with one or a small number (more than one) of
hidden layers in the local neural networks, the current method can
produce accurate simulation results.


Since the current method is a combination of the ideas of ELM, domain decomposition,
and local neural networks, we refer to this method as locELM (local extreme
learning machines) in the current paper.

We have performed extensive numerical experiments with linear and nonlinear, stationary
and time-dependent,
partial differential equations
to test the performance of the locELM method,
and to study the effects of the simulation parameters involved therein.
For certain test problems (e.g.~the advection equation) we present very long-time
simulations to demonstrate the capability and accuracy of the locELM method
together with the block time-marching scheme.
We compare extensively the current locELM method with
the deep Galerkin method~\cite{SirignanoS2018} and the physics-informed neural
network method~\cite{RaissiPK2019}, and demonstrate the superiority of
the current method in terms of both accuracy and the computational cost.
We also compare the current method with the classical finite element method,
and show that the computational performance of the locELM method is comparable to,
and often exceeds, the FEM performance.
The current locELM method, DGM and PINN have all been implemented in Python,
using the Tensorflow (www.tensorflow.org) and Keras (keras.io) libraries.
The finite element method is also implemented in Python,
by using the FEniCS library (fenicsproject.org).

The rest of this paper is structured as follows.
In Section \ref{sec:method} we outline the locELM representation of field functions
based on domain decomposition and local extreme learning machines, and then
discuss how to solve linear and nonlinear differential equations
using the locELM representation 
and how to
train the overall neural network by the linear or nonlinear
least squares method. For nonlinear differential equations
we present the NLSQ-perturb method and the Newton-LLSQ method for solving
the resultant nonlinear algebraic system.
For time-dependent PDEs,
we present the block time-marching scheme, and discuss how to employ
the locELM method together with block time marching for long-time
simulations. We primarily use second-order differential equations in
two spatial dimensions,
and also plus time if the problem is time-dependent,
as examples in the presentation of the locELM method.
In Section \ref{sec:tests} we present extensive numerical
experiments with the linear and nonlinear Helmholtz equations,
the advection equation, the diffusion equation, nonlinear spring equation,
and the viscous Burger's equation
to test the performance of the locELM
method.
We compare the locELM method with DGM and PINN, and demonstrate the superiority
of locELM in terms of the accuracy and computational cost.
We also compare locELM with the classical finite element method,
and show that locELM is on par with, and often exceeds,
the FEM in computational performance.
Section \ref{sec:summary} concludes the main presentation with
a number of comments on the characteristics and properties of the current method.
The Appendix summarizes some additional numerical tests
not included in the main text.

\section{Domain Decomposition and Local Extreme Learning Machines}
\label{sec:method}

\subsection{Local Extreme Learning Machines (locELM)
  for Representing Functions }
\label{sec:loc_elm}

%
%

Consider the domain $\Omega$ in $d$ ($d=1$, $2$ or $3$) dimensions,
where one of the dimensions may denote time and so $\Omega$ in general
can be a spatial-temporal domain.
We consider a function $f(\mbs x)$ ($\mbs x\in\Omega$) defined on this domain,
and would like to represent this function using neural networks.

We partition $\Omega$ into $N_e$ ($N_e\geqslant 1$) non-overlapping
sub-domains, 
\begin{equation*}
  \Omega = \Omega_1\cup\Omega_2\cup \dots \cup\Omega_{N_e},
\end{equation*}
where $\Omega_i$ denotes the $i$-th sub-domain.
If $\Omega_i$ and $\Omega_j$ ($1\leqslant i,j\leqslant N_e$)
share a common boundary, we will denote this common boundary
by $\Gamma_{ij}$.

We will represent $f(\mbs x)$, in a spirit analogous
to the finite elements or
spectral elements~\cite{KarniadakisS2005,ZhengD2011,DongY2009,DongS2012},
locally on the sub-domains by local neural networks.
More specifically, on each sub-domain $\Omega_i$ ($1\leqslant i\leqslant N_e$)
we represent $f(\mbs x)$ by a shallow feed-forward neural network~\cite{GoodfellowBC2016}.
Here ``shallow'' refers to the configuration that
each local neural network has only a small number 
(e.g.~one, two or perhaps three) of hidden layers,
apart from the input layer (representing $\mbs x$) and the
output layer (representing $f(\mbs x)$, restricted to $\Omega_i$).

Let $f_i(\mbs x)$ ($1\leqslant i\leqslant N_e$) denote the function
$f(\mbs x)$ restricted to $\Omega_i$.
On any common boundary $\Gamma_{ij}$ between $\Omega_i$
and $\Omega_j$ (for all $1\leqslant i,j\leqslant N_e$),
we impose the requirement that $f_i(\mbs x)$ and
$f_j(\mbs x)$ satisfy the $C^{\mbs k}$ continuity conditions
with an appropriate $\mbs k=(k_1,k_2,\dots,k_d)$.
In other words, their function values and partial derivatives
up to the order $k_s$ ($1\leqslant s\leqslant d$) should be continuous across the
sub-domain boundary in the $s$-th direction.
The order $\mbs k$ in the $C^{\mbs k}$ continuity is a user-defined parameter.
When solving differential equations, one can determine
$\mbs k$ for a specific coordinate direction
based on the order of the differential equation along that
direction. For example, if the highest derivative with
respect to the coordinate $x_s$ ($1\leqslant s\leqslant d$)
involved in the equation is $m$, one would typically
impose $C^{m-1}$ continuity to the solution 
on the sub-domain boundary along the $s$-th direction.
Thanks to these $C^{\mbs k}$ continuity conditions,
the local  neural networks for the sub-domains,
while physically separated, are coupled with one another
logically, and need to be trained together in a coupled fashion.
The local neural networks collectively constitute the representation
of the function $f(\mbs x)$ on the overall domain $\Omega$.


We impose further requirements on 
the local neural networks.
Suppose a particular layer in the local neural network
contains $n$ nodes, and the previous layer contains $m$ nodes.
Let $\phi_i(\mbs x)$ ($1\leqslant i\leqslant m$) denote the output
of the previous layer,
and $\varphi_i(\mbs x)$ ($1\leqslant i\leqslant n$)
denote the output of this layer.
Then the logic of this layer is represented 
by~\cite{GoodfellowBC2016},
\begin{equation}
  \varphi_i(\mbs x) = \sigma\left(\sum_{j=1}^m \phi_j(\mbs x)w_{ji} + b_i\right),
  \quad 1\leqslant i\leqslant n,
\end{equation}
where the constants $w_{ji}$ and $b_i$ ($1\leqslant i\leqslant n$,
$1\leqslant j\leqslant m$) are the weight and bias coefficients
associated with this layer, and $\sigma(\cdot)$ is
the activation function of this layer and is in general
nonlinear.
We assume the following for the local neural networks:
\begin{itemize}

\item
  The weight and bias coefficients for all the hidden layers are pre-set
  to uniform random values generated on the interval $[-R_m,R_m]$,
  where $R_m>0$ is a user-defined constant parameter.
  Once these coefficients are set randomly,
  they are fixed throughout the training and computation. These weight/bias
  coefficients are not adjustable, and they are not
  training parameters of the neural network.
  We hereafter refer to $R_m$ as the maximum magnitude of the
  random coefficients of the neural network.

\item
  The last hidden layer, i.e.~the layer before the output layer, can be wide.
  In other words, this layer may contain a large number of nodes.
  We use $M$ to denote the number of nodes in the last hidden layer
  of each local neural network.

\item
  The output layer contains no bias (i.e.~$b_i=0$)
  and no activation function. In other words, the output layer
  is linear, i.e.~$\sigma(x)=x$.
  The weight coefficients in the output layers of the local neural networks
  are adjustable.
  The collection of these weight coefficients
  constitutes the  training parameters of the overall neural network.
  Therefore, the number of training parameters in each local neural network
  equals $M$, the number of nodes in the last hidden layer of
  the local neural network.

\item
  The set of training parameters for the overall neural network
  is to be determined and set by a linear or nonlinear least squares
  computation, not by the back propagation-type algorithm.
  
\end{itemize}

\begin{remark}\label{rem_aa}
When a subset of the above requirements is imposed on a
single global neural network,
containing a single hidden layer, for the entire domain, the resultant
network, when trained with a linear least squares method,
is known as an extreme learning
machine (ELM)~\cite{HuangZS2006}.
In the current work we follow this terminology, and will refer to the local
neural networks presented here as local
extreme learning machines (or locELM).

\end{remark}

Let $N$ ($N\geqslant 1$) denote the number of nodes in
the output layer of the local neural networks.
Based on the above
assumptions, on the sub-domain $\Omega_s$ ($1\leqslant s\leqslant N_e$)
we have the relation,
\begin{equation}
  u_{i}^{s}(\mbs x) = \sum_{j=1}^M V_j^{s}(\mbs x) w^{s}_{ji},
  \quad \mbs x \in \Omega_s, \ \
  1\leqslant i\leqslant N,
  \label{equ_a}
\end{equation}
where $V_j^s(\mbs x)$ ($1\leqslant j\leqslant M$) denote
the output of the last hidden layer, $u_i^{s}(\mbs x)$
denote the the components of output function
of the network, $w_{ji}^s$
are the training parameters on $\Omega_s$, and $M$ denotes the number of
nodes in the last hidden layer.
The function
\begin{equation}
f_s(\mbs x) = (u_1^s, u_2^s, \dots, u_N^s)
\end{equation}
is the local representation of $f(\mbs x)$ on the sub-domain $\Omega_s$.

It should be noted that the set of output
functions of the last hidden layer, $V_j^s(\mbs x)$ ($1\leqslant j\leqslant M$),
are known functions and they are fixed throughout the computation.
%
%
Since the weight/bias coefficients in the hidden layers are pre-set
to random values on $[-R_m,R_m]$ and are fixed,
$V_j^s(\mbs x)$ can be pre-computed
by a forward evaluation of the local neural network (up to
the last hidden layer)
against the input $\mbs x$ data.
The first, second, and higher-order derivatives of $V_j^s(\mbs x)$
with respect to the input $\mbs x$ can then be computed
by auto-differentiations.

The collection of local representations $f_s(\mbs x)$ ($1\leqslant s\leqslant N_e$),
with $C^{\mbs k}$ continuity imposed on the sub-domain boundaries and
with $w_{ij}^s$ ($1\leqslant i\leqslant M$, $1\leqslant j\leqslant N$,
$1\leqslant s\leqslant N_e$) as the  training parameters,
form the set of trial functions for representing the function
$f(\mbs x)$. Hereafter, we will refer to this representation as
the locELM representation of a function. 
Once the data for $f(\mbs x)$ or the data
for the governing equations that describe $f(\mbs x)$ are given,
the adjustable parameters $w_{ij}^s$ can be trained and determined by a
linear or nonlinear least squares computation.

\begin{remark}
\label{rem:rem_1}
In the locELM representation, the hyper-parameters for the local
neural networks associated with different sub-domains
(e.g.~depths, widths and activation
functions of the hidden layers) can in principle assume different values.
This can allow one to place more degrees of freedom locally in
regions where the field function may be more complicated and thus require
more resolution.
For simplicity of implementation, however, in the current work
we will employ the same hyper-parameters for all the local
neural networks for different sub-domains.

\end{remark}

In the following sub-sections we focus on how to use local
extreme learning machines to represent the solutions to
ordinary or partial differential equations (ODE/PDE), and discuss how to train
the overall neural network by least squares computations.
We consider two cases: (i) linear differential equations,
and (ii) nonlinear differential equations, and
discuss how to treat them individually.
Apart from the basic algorithm,
we develop a block time-marching scheme
for long-time simulations of time-dependent linear/nonlinear PDEs.
In the presentations we use two spatial dimensions, and 
plus time if the problem is time-dependent, as examples.
The formulations can be reduced to one spatial dimension
or extended to higher spatial dimensions in a straightforward fashion.
For simplicity we concentrate on rectangular spatial-temporal domains
in the current work.



\subsection{Linear Differential Equations}

\subsubsection{Time-Independent Linear Differential Equations}
\label{sec:steady}


Let us first consider
the boundary value problem involving
linear partial differential equations together
with Dirichlet boundary conditions, and discuss how to
solve the problem by using the locELM representation for
the solution.
To make the discussion concrete, we concentrate on
two dimensions ($d=2$, with the coordinates $x$ and $y$),
and consider second-order partial differential equations with respect to
both $x$ and $y$ (i.e.~highest partial derivatives with respect to
$x$ and to $y$ are both two).
The procedure outlined below can be extended to higher dimensions
or to higher-order differential equations, with appropriate boundary
conditions and $C^{\mbs k}$ continuity conditions taken into account.

Let us consider the following generic
second-order linear partial differential equation
\begin{subequations}
  \begin{align}
    &
    L u = f(x,y), \label{equ_1} \\
    &
    u(x,y) = g(x,y), \quad \text{on} \ \partial\Omega,
    \label{equ_2}
  \end{align}
\end{subequations}
where $L$ is a  linear second-order operator with respect to
both $x$ and $y$, 
$u(x,y)$ is the scalar unknown field function to be
solved for, $f(x,y)$ and $g(x,y)$ are prescribed source
terms for the equation and the Dirichlet boundary condition,
and $\partial\Omega$ denotes the boundary of $\Omega$.
We assume that
this boundary value problem is well-posed.
Our goal here is to illustrate the procedure for
numerically solving this problem
by approximating its solution using local extreme learning
machines.


Here is the general idea for the solution process.
We partition the overall domain
into a number of sub-domains, and represent the field solution
using the locELM representation described in Section \ref{sec:loc_elm}.
We next choose a set of
points (collocation points) within each sub-domain,
which can have a regular or
random distribution.
We enforce the governing equations on the collocation points
within each sub-domain, and enforce the boundary conditions
on those collocation points in those sub-domains that reside
on $\partial\Omega$. 
We further enforce the $C^{\mbs k}$ continuity conditions on
those collocation points that reside on the sub-domain
boundaries.
Auto-differentiations are employed to compute
the first or higher-order derivatives involved in the above operations.
These operations result in a system of algebraic equations,
which may be linear or nonlinear depending on the boundary value problem,
about the training parameters in the locELM representation.
We seek a least squares solution to this algebraic system,
and compute the solution by either a linear least squares method
or a nonlinear least squares method.
The training parameters of the local neural networks are then
determined by the least squares computation.


For simplicity of implementation, we concentrate on the case
with $\Omega$ being a rectangular domain,
i.e.~$\Omega=[a_1,b_1]\times [a_2,b_2]$.
Let $N_x$ ($N_x\geqslant 1$) and $N_y$ ($N_y\geqslant 1$)
denote the number of sub-domains along the $x$ and $y$ directions,
respectively, with a total number of $N_e = N_xN_y$ sub-domains
in $\Omega$.
Let the two vectors $[X_0, X_1, \dots, X_{N_x}]$ and
$[Y_0, Y_1, \dots, Y_{N_y}]$ denote the coordinates of the sub-domain
boundaries along the $x$ and $y$ directions, where
$(X_0,Y_0)=(a_1,a_2)$ and $(X_{N_x},Y_{N_y})=(b_1,b_2)$.
Let $\Omega_{e_{mn}}=[X_m,X_{m+1}]\times [Y_n,Y_{n+1}]$  denote
the region occupied by the sub-domain $e_{mn}$, for
$0\leqslant m\leqslant N_x-1$ and 
$0\leqslant n\leqslant N_y-1$.
Here $e_{mn}$ represents the linear index of the sub-domain associated
with the 2D index $(m,n)$, with $e_{mn}=mN_y+n+1$,
and so $1\leqslant e_{mn}\leqslant N_e$.


We approximate the unknown field function $u(x,y)$ using
the locELM representation as discussed in Section \ref{sec:loc_elm}.
On each sub-domain $e_{mn}$
we represent the solution by a shallow neural network,
which consists of an input layer with two nodes (representing
the coordinates $x$ and $y$), one or a small number of hidden layers,
and an output layer with one node (representing the solution $u^{e_{mn}}$).
Let $V_{j}^{e_{mn}}(x,y)$ ($1\leqslant j\leqslant M$) denote the output
of the last hidden layer, where $M$ is the number of nodes in this layer.
Then equation \eqref{equ_a} becomes
\begin{equation}\label{equ_b}
  u^{e_{mn}}(x,y) = \sum_{j=1}^M V_j^{e_{mn}}(x,y) w_{j}^{e_{mn}}, \quad
  (x,y)\in\Omega_{e_{mn}}, \quad
  0\leqslant m\leqslant N_x-1, \ \ 0\leqslant n\leqslant N_y-1,
\end{equation}
where $w_j^{e_{mn}}$ ($1\leqslant j\leqslant M$) are the training parameters
in the sub-domain $e_{mn}$.
Again note that $V_j^{e_{mn}}(x,y)$ is known, once the weight/bias coefficients
in the hidden layers have been pre-set to random values on $[-R_m,R_m]$.

\begin{remark}\label{rem_bb}
Apart from the above logical operations,
in the implementation we incorporate an additional normalization
layer immediately behind the input layer in each of the local
neural networks.
For each sub-domain $e_{mn}$,
the normalization layer performs an affine mapping
and normalizes the input data,
$(x,y)\in\Omega_{e_{mn}}= [X_m,X_{m+1}]\times [Y_n, Y_{n+1}]$, such that
the output data of the normalization layer fall into
the domain $[-1,1] \times [-1,1]$.
This extra normalization layer contains
no adjustable (training) parameters.
  
\end{remark}

On the sub-domain $e_{mn}$ ($0\leqslant m\leqslant N_x-1$,
$0\leqslant n\leqslant N_y-1$), let $(x_{p}^{e_{mn}},y_q^{e_{mn}})$
($0\leqslant p\leqslant Q_x-1$,
$0\leqslant q\leqslant Q_y-1$)
denote a set of distinct collocation points,
where $x_p^{e_{mn}}$ ($0\leqslant p\leqslant Q_x-1$) denote a set of $Q_x$
collocation points on the interval $[X_m, X_{m+1}]$
and $y_q^{e_{mn}}$ denote a set of $Q_y$ collocation points on
the interval $[Y_n, Y_{n+1}]$.
The total number of collocation points is $Q=Q_xQ_y$
within each sub-domain $e_{mn}$. 
In the current work we primarily consider the following uniform distribution
for the collocation points:
\begin{itemize}
\item
  Uniform distribution: $x_{p}^{e_{mn}}$ forms a set of $Q_x$ uniform
  grid points on $[X_m,X_{m+1}]$, with both end points included,
  i.e.~$x_0^{e_{mn}}=X_m$ and $x_{Q_x-1}^{e_{mn}}=X_{m+1}$.
  $y_q^{e_{mn}}$ forms a set of $Q_y$ uniform grid points on $[Y_n,Y_{n+1}]$,
  with both end points included,
  i.e.~$y_0^{e_{mn}}=Y_n$ and $y_{Q_y-1}^{e_{mn}}=Y_{n+1}$.

  
\end{itemize}

\begin{remark}\label{rem_cc}
  Besides the uniform distribution, we also consider a quadrature-point distribution
  and a random distribution for the collocation points.
  With the quadrature-point distribution,
  $x_{p}^{e_{mn}}$ are taken to be a set of $Q_x$ Gauss-Lobatto-Legendre
  quadrature points on the interval $[X_m,X_{m+1}]$, and
  $y_q^{e_{mn}}$ are taken to be a set of $Q_y$ Gauss-Lobatto-Legendre quadrature points
  on the interval $[Y_n,Y_{n+1}]$.
  With the random distribution, the collocation points in the sub-domain $e_{mn}$
  are taken to be uniformly generated
  random points $(x_l^{e_{mn}},y_l^{e_{mn}})\in \Omega_{e_{mn}}$ ($0\leqslant l\leqslant Q-1$),
  where $Q$ is the total number of collocation points in the sub-domain,
  among which a certain number of points are generated 
  on the sub-domain boundaries and the rest are located inside the sub-domain.  
  Numerical experiments indicate that, with the same number of collocation points,
  the result with the quadrature-point distribution is generally more accurate than that with
  the uniform distribution, which in turn is more accurate than that with
  the random distribution of collocation points.
  The quadrature-point distribution however poses some practical issues in
  the current implementation. When the number of quadrature points exceeds $100$,
  the library on which the current implementation is based cannot compute the
  Gaussian quadrature points accurately. This is the reason why in the current work
  we predominantly employ the uniform distribution
  of collocation points in the numerical tests of
  Section \ref{sec:tests}.
  
\end{remark}


With the above setup, we solve the boundary value problem consisting of
equations \eqref{equ_1} and \eqref{equ_2} as follows.
On each sub-domain $e_{mn}$
we enforce the equation \eqref{equ_1} on
all the collocation points $(x_p^{e_{mn}},y_q^{e_{mn}})$,
\begin{equation}
  \begin{split}
    &
    \sum_{j=1}^M \left[LV_j^{e_{mn}}\left(x_p^{e_{mn}},y_q^{e_{mn}} \right) \right] w_j^{e_{mn}}
    = f(x_p^{e_{mn}}, y_q^{e_{mn}}), \\
    &
    \text{for} \ 0\leqslant m\leqslant N_x-1, \ 0\leqslant n\leqslant N_y-1,
    \ 0\leqslant p\leqslant Q_x-1, \ 0\leqslant q\leqslant Q_y-1,
  \end{split}
  \label{equ_3}
\end{equation}
where we have used equation \eqref{equ_b}.
We enforce equation \eqref{equ_2} on the four boundaries of
the domain $\Omega$,
\begin{subequations}
  \begin{align}
    &
    \sum_{j=1}^M V_{j}^{e_{0n}}\left(a_1,y_q^{e_{0n}} \right)w_j^{e_{0n}} =
    g\left(a_1,y_q^{e_{0n}} \right),
    \ \ 0\leqslant n\leqslant N_y-1, \ 0\leqslant q\leqslant Q_y-1;
    \label{equ_4} \\
    &
    \sum_{j=1}^M V_{j}^{e_{mn}}\left(b_1,y_q^{e_{mn}} \right)w_j^{e_{mn}} =
    g\left(b_1,y_q^{e_{mn}} \right),
    \ \ m=N_x-1, \ 0\leqslant n\leqslant N_y-1, \ 0\leqslant q\leqslant Q_y-1;
    \label{equ_5} \\
    &
    \sum_{j=1}^M V_{j}^{e_{m0}}\left(x_p^{e_{m0}},a_2 \right)w_j^{e_{m0}} =
    g\left(x_p^{e_{m0}},a_2 \right),
    \ \ 0\leqslant m\leqslant N_x-1, \ 0\leqslant p\leqslant Q_x-1;
    \label{equ_6} \\
    &
    \sum_{j=1}^M V_{j}^{e_{mn}}\left(x_p^{e_{mn}},b_2 \right)w_j^{e_{mn}} =
    g\left(x_p^{e_{mn}},b_2 \right),
    \ \ n=N_y-1, \ 0\leqslant m\leqslant N_x-1, \ 0\leqslant p\leqslant Q_x-1,
    \label{equ_7}
  \end{align}
\end{subequations}
where equation \eqref{equ_b} has again been used.


The local representations of the field solution are coupled together
by the $C^{\mbs k}$ continuity conditions. Since the equation \eqref{equ_1}
is assumed to be of second order with respect to both $x$ and $y$,
we impose $C^1$ continuity conditions across the sub-domain
boundaries in both the $x$ and $y$ directions.
On the vertical sub-domain boundaries $x=X_{m+1}$ ($0\leqslant m\leqslant N_x-2$),
the $C^1$ conditions are reduced to,
\begin{subequations}
  \begin{align}
    &
    \sum_{j=1}^M V_j^{e_{mn}}\left(X_{m+1},y_q^{e_{mn}} \right) w_j^{e_{mn}}
    - \sum_{j=1}^M V_j^{e_{m+1,n}}\left(X_{m+1},y_q^{e_{m+1,n}} \right) w_j^{e_{m+1,n}}
    = 0, \label{equ_8a} \\
    &
    \sum_{j=1}^M \left.\frac{\partial V_j^{e_{mn}}}{\partial x}\right|_{\left(X_{m+1},y_q^{e_{mn}} \right)} w_j^{e_{mn}}
    - \sum_{j=1}^M \left.\frac{\partial V_j^{e_{m+1,n}}}{\partial x}\right|_{\left(X_{m+1},y_q^{e_{m+1,n}} \right)} w_j^{e_{m+1,n}}
    = 0, \label{equ_8b} \\
    &
    \text{for}\ 0\leqslant m\leqslant N_x-2, \
    0\leqslant n\leqslant N_y-1, \ 0\leqslant q\leqslant Q_y-1, \nonumber
  \end{align}
\end{subequations}
where it should be noted that $y_q^{e_{mn}}=y_q^{e_{m+1,n}}$.
On the horizontal sub-domain boundaries $y=Y_{n+1}$ ($0\leqslant n\leqslant N_y-2$),
the $C^1$ continuity conditions are reduced to,
\begin{subequations}
  \begin{align}
    &
    \sum_{j=1}^M V_j^{e_{mn}}\left(x_p^{e_{mn}},Y_{n+1} \right) w_j^{e_{mn}}
    - \sum_{j=1}^M V_j^{e_{m,n+1}}\left(x_p^{e_{m,n+1}},Y_{n+1} \right) w_j^{e_{m,n+1}}
    = 0, \label{equ_9a} \\
    &
    \sum_{j=1}^M \left.\frac{\partial V_j^{e_{mn}}}{\partial y}\right|_{\left(x_p^{e_{mn}},Y_{n+1} \right)} w_j^{e_{mn}}
    - \sum_{j=1}^M \left.\frac{\partial V_j^{e_{m,n+1}}}{\partial y}\right|_{\left(x_p^{e_{mn+1}},Y_{n+1} \right)} w_j^{e_{m,n+1}}
    = 0, \label{equ_9b} \\
    &
    \text{for} \ 0\leqslant m\leqslant N_x-1, \
    0\leqslant n\leqslant N_y-2, \ 0\leqslant p\leqslant Q_x-1, \nonumber
  \end{align}
\end{subequations}
where it should be noted that
$x_p^{e_{mn}} = x_p^{e_{m,n+1}}$.

The set of equations consisting of \eqref{equ_3}--\eqref{equ_9b}
is a system of linear algebraic equations
about the training parameters
$w_j^{e_{mn}}$ ($0\leqslant m\leqslant N_x-1$, $0\leqslant n\leqslant N_y-1$,
$1\leqslant j\leqslant M$).
%
In these equations, $V_j^{e_{mn}}(x,y)$, $LV_j^{e_{mn}}(x,y)$,
$\frac{\partial V_j^{e_{mn}}}{\partial x}$ and
$\frac{\partial V_j^{e_{mn}}}{\partial y}$
are all known functions, 
once the weight/bias coefficients
in the hidden layers are randomly set.
These functions can be evaluated on the collocation
points, including those on the domain boundaries and
the sub-domain boundaries.
The derivatives involved in these functions can be computed
by auto-differentiation.

This linear algebraic system consists of
$ 
N_xN_y(Q_xQ_y+2Q_x+2Q_y)
$ 
equations, and $N_xN_yM$ unknown variables of $w_j^{e_{mn}}$.
We seek the least squares solution to this system
with the minimum norm. Linear least squares routines are available
in a number of scientific libraries,
and we take advantage of these numerical libraries in our implementation.
In the current
work we employ the linear least squares routine
from LAPACK, available through wrapper functions
in the scipy package in Python.
Therefore, the adjustable parameters $w_j^{e_{mn}}$
in the neural network
are trained by this linear least squares computation.


In the current work we have employed Tensorflow and Keras to implement
the neural network architecture as outlined above.
Each local neural network consists of several ``dense'' Keras layers.
The set of $N_e=N_xN_y$ local neural networks collectively forms
an overall logical neural network, in the form of
a multi-input multi-output Keras model. 
The input data to the model consist of the coordinates of
the collocation points for all sub-domains,
$(x_p^{e_{mn}},y_q^{e_{mn}})$, for $0\leqslant m\leqslant N_x-1$,
$0\leqslant n\leqslant N_y-1$, $0\leqslant p\leqslant Q_x-1$
and $0\leqslant q\leqslant Q_y-1$.
The output of the Keras model consists
of the solution $u^{e_{mn}}(x,y)$ on the collocation points
for all the sub-domains.
The output of the last hidden layer of each sub-domain, $V_j^{e_{mn}}(x,y)$,
are obtained by creating a Keras sub-model using the Keras functional
APIs (application programming interface). The derivatives
of $V_j^{e_{mn}}(x,y)$, and those involved in $LV_j^{e_{mn}}(x,y)$,
are computed using auto-differentiation with
these Keras sub-models.
After the parameters $w_j^{e_{mn}}$ are obtained by
the linear least squares computation, the weight coefficients in
the output layer of the Keras model are then set based on
these parameter values.


\begin{remark}
  \label{rem_2}
  We observe from numerical experiments that
  the simulation result obtained using the current method is considerably more accurate,
  typically by orders of magnitude,
  than those obtained using DNN-based PDE solvers,
  trained using gradient descent-type algorithms. 
  Furthermore, the current method is computationally fast. Its computational
  cost is essentially the cost of the linear least squares computation.
  We observe that 
  the network training time of the current method
  is considerably lower, typically by orders of magnitude,
  than those of the DNN-based PDE solvers trained with gradient descent-type
  algorithms.
  These points will be demonstrated by extensive numerical experiments 
  in Section \ref{sec:tests}, in which we compare the current method with
  the deep Galerkin method~\cite{SirignanoS2018}
  and the Physics-Informed Neural Network~\cite{RaissiPK2019}.
  
\end{remark}

\begin{remark} \label{rem_2a}
  The computational performance of the current locELM method,
  in terms of the accuracy and the computational cost,
  is comparable to, and oftentimes exceeds, that of the classical
  finite element method.
  These points will be demonstrated by extensive numerical experiments in
  Section \ref{sec:tests} with time-independent and
  time-dependent problems.
  We observe that, with the same training/computation time, the accuracy of
  the current method is comparable, and oftentimes considerably superior,
  to that of the finite element method. To achieve the same accuracy,
  the training time of the current method is comparable to,
  and oftentimes markedly smaller than, the computation time
  of the classical finite element method.
  
\end{remark}



\subsubsection{Time-Dependent Linear Differential Equations}
\label{sec:unsteady}

We next consider initial-boundary value problems involving time-dependent
linear differential equations together with Dirichlet boundary conditions,
and discuss how to solve such problems using the locELM method.
We again concentrate on
two spatial dimensions (with coordinates $x$ and $y$) plus time ($t$), 
and assume second spatial orders
in the differential equation with respect to both $x$ and $y$.

\paragraph{Basic Method}
\label{sec:basic}

We consider the following generic
time-dependent second-order linear PDE, together with
the Dirichlet boundary condition and the initial condition,
\begin{subequations}
  \begin{align}
    &
    \frac{\partial u}{\partial t} = Lu + f(x,y,t),
    \label{equ_10a} \\
    &
    u(x,y,t) = g(x,y,t), \quad \text{for} \ (x,y) \
    \text{on spatial domain boundary},
    \label{equ_10b} \\
    &
    u(x,y,0) = h(x,y), \label{equ_10c}
  \end{align}
\end{subequations}
where 
$u(x,y,t)$ is the unknown field function to
be solved for, $L$ is a second-order linear differential operator
with respect to both $x$ and $y$, $f(x,y,t)$ is a prescribed source
term, $g(x,y,t)$ is the Dirichlet boundary data, and $h(x,y)$
denotes the initial field distribution.
We assume that
this initial-boundary
value problem is well posed, and would like to solve this problem
by approximating $u(x,y,t)$ using the locELM representation.

We seek the solution on a rectangular spatial-temporal
domain,
$
\Omega = \{
(x,y,t)\ |\ x\in[a_1,b_1], \ y\in [a_2,b_2], \
t\in [0, \Gamma]
\},
$
where $a_i$, $b_i$ ($i=1,2$) and $\Gamma$ are prescribed constants.
The solution procedure is analogous to that discussed in Section \ref{sec:steady}.
We partition $\Omega$ into $N_x$ ($N_x\geqslant 1$) sub-domains along the
$x$ direction, $N_y$ ($N_y\geqslant 1$) sub-domains along the $y$ direction,
and $N_t$ ($N_t\geqslant 1$) sub-domains in time, leading to a total of
$N_{e}=N_xN_yN_t$  sub-domains in $\Omega$.
Let the vectors $[X_0, X_1, \dots, X_{N_x}]$,
$[Y_0, Y_1, \dots, Y_{N_y}]$
and $[T_0, T_1, \dots, T_{N_t}]$ denote the coordinates of
the sub-domain boundaries along the $x$, $y$ and temporal directions,
respectively,
where $(X_0,Y_0,T_0)=(a_1,a_2,0)$
and $(X_{N_x},Y_{N_y},T_{N_t})=(b_1,b_2,\Gamma)$.
We use
$
\Omega_{e_{mnl}} = [X_m,X_{m+1}]\times[Y_n,Y_{n+1}]\times[T_l,T_{l+1}]
$
to denote the spatial-temporal region
occupied by the sub-domain with the index
$
e_{mnl} = mN_yN_t+nN_t+l+1,
$
for $0\leqslant m\leqslant N_x-1$,
$0\leqslant n\leqslant N_y-1$ and $0\leqslant l\leqslant N_t-1$.


We approximate $u(x,y,t)$ using the locELM representation
from Section \ref{sec:loc_elm}. More specifically,
we employ a local shallow feed-forward neural network for
the solution on each sub-domain $e_{mnl}$.
The local neural network consists of an input layer with
three nodes, representing
the coordinates $x$, $y$ and $t$, respectively,
a small number of hidden layers, and
an output layer consisting of one node, representing
the solution $u^{e_{mnl}}(x,y,t)$ on this sub-domain.
The output layer is linear and contains no bias.
The weight/bias coefficients in all the hidden layers are pre-set
to uniform random values generated on $[-R_m,R_m]$ and are fixed,
as discussed in Section \ref{sec:loc_elm}.
Additionally, in the implementation, we incorporate an affine mapping operation
right behind the input layer to normalize the 
input data, $(x,y,t)\in\Omega_{e_{mnl}}$, to
the interval $[-1,1]\times[-1,1]\times[-1,1]$.
Let $V_j^{e_{mnl}}$ ($1\leqslant j\leqslant M$) denote the
output  of the last hidden layer,
where $M$ denotes the number of nodes in this layer.
Then we have, in accordance with equation \eqref{equ_b},
\begin{equation}
  \begin{split}
    &
    u^{e_{mnl}}(x,y,t) = \sum_{j=1}^M V_j^{e_{mnl}}(x,y,t) w_j^{e_{mnl}}, \\
    & 
    \text{for}\ 0\leqslant m\leqslant N_x-1, \
    0\leqslant n\leqslant N_y-1, \ 0\leqslant l\leqslant N_t-1,
  \end{split}
\end{equation}
where the coefficients $w_j^{e_{mnl}}$ ($1\leqslant j\leqslant M$) are the training
parameters of the local neural network.
Note that $V_j^{e_{mnl}}(x,y,t)$ and its
derivatives are all known functions,
since the weight/bias coefficients of all the hidden layers
are pre-set and fixed. 


On each sub-domain $e_{mnl}$,
let $(x_p^{e_{mnl}}, y_q^{e_{mnl}},t_r^{e_{mnl}})$
($0\leqslant p\leqslant Q_x-1$, $0\leqslant q\leqslant Q_y-1$,
and $0\leqslant r\leqslant Q_t-1$) denote a set of
distinct collocation points,
where $x_p^{e_{mnl}}$ ($0\leqslant p\leqslant Q_x-1$) denotes a set of
$Q_x$ collocation points on $[X_m,X_{m+1}]$
with $x_0^{e_{mnl}}=X_m$ and $x_{Q_x-1}^{e_{mnl}}=X_{m+1}$,
$y_q^{e_{mnl}}$ ($0\leqslant q\leqslant Q_y-1$) denotes a set of
$Q_y$ collocation points on $[Y_n,Y_{n+1}]$ with
$y_0^{e_{mnl}}=Y_n$ and $y_{Q_y-1}^{e_{mnl}}=Y_{n+1}$,
and $t_{r}^{e_{mnl}}$ ($0\leqslant r\leqslant Q_t-1$) denotes
a set of $Q_t$ collocation points on $[T_l,T_{l+1}]$
with $t_0^{e_{mnl}}=T_l$ and $t_{Q_t-1}^{e_{mnl}}=T_{l+1}$.
We primarily consider the uniform distribution of regular grid points
as the collocation points, analogous to that in Section \ref{sec:steady}.


With these setup,
we next enforce the equations \eqref{equ_10a}--\eqref{equ_10c}
on the collocation points inside each sub-domain and on the domain
boundaries.
On the sub-domain $e_{mnl}$,
equation \eqref{equ_10a} is reduced to
\begin{equation}\label{equ_11}
  \begin{split}
    &
  \sum_{j=1}^M \left.\left(
  \frac{\partial V_j^{e_{mnl}}}{\partial t} - LV_j^{e_{mnl}}
  \right)\right|_{(x_p^{e_{mnl}},y_q^{e_{mnl}},t_r^{e_{mnl}})} w_j^{e_{mnl}}
  = f\left(x_p^{e_{mnl}},y_q^{e_{mnl}},t_r^{e_{mnl}}\right), \\
  & \quad
  \text{for}\ 0\leqslant m\leqslant N_x-1, \
  0\leqslant n\leqslant N_y-1, \ 0\leqslant l\leqslant N_t-1, \\
  & \quad\quad\ \
  0\leqslant p\leqslant Q_x-1, \ 0\leqslant q\leqslant Q_y-1, \
  0\leqslant r\leqslant Q_t-1,
  \end{split}
\end{equation}
where $\left(x_p^{e_{mnl}},y_q^{e_{mnl}},t_r^{e_{mnl}}\right)$
are the collocation points.
The boundary condition \eqref{equ_10b},
when enforced on the spatial domain boundaries corresponding to
$x=a_1$ or $b_1$ and $y=a_2$ or $b_2$, is reduced to
\begin{subequations}
  \begin{align}
    \begin{split}
    &
    \sum_{j=1}^M V_j^{e_{0nl}}(a_1,y_q^{e_{0nl}},t_r^{e_{0nl}}) w_j^{e_{0nl}}
    - g(a_1,y_q^{e_{0nl}},t_r^{e_{0nl}})=0, \label{equ_13a} \\
    & \qquad\qquad \text{for} \
    0\leqslant n\leqslant N_y-1, \ 0\leqslant l\leqslant N_t-1, \
    0\leqslant q\leqslant Q_y-1, \ 0\leqslant r\leqslant Q_t-1;
    \end{split} \\
    \begin{split}
    &
    \sum_{j=1}^M V_j^{e_{mnl}}(b_1,y_q^{e_{mnl}},t_r^{e_{mnl}}) w_j^{e_{mnl}}
    - g(b_1,y_q^{e_{mnl}},t_r^{e_{mnl}})=0, \\
    & \qquad\qquad \text{for} \
    m = N_x-1, \ 
    0\leqslant n\leqslant N_y-1, \ 0\leqslant l\leqslant N_t-1, \
    0\leqslant q\leqslant Q_y-1, \ 0\leqslant r\leqslant Q_t-1;
    \end{split} \\
    \begin{split}
      &
      \sum_{j=1}^M V_j^{e_{m0l}}(x_p^{e_{m0l}},a_2,t_r^{e_{m0l}})w_j^{e_{m0l}}
      - g(x_p^{e_{m0l}},a_2,t_r^{e_{m0l}})=0, \\
      & \qquad\qquad \text{for} \
      0\leqslant m\leqslant N_x-1, \ 0\leqslant l\leqslant N_t-1, \
      0\leqslant p\leqslant Q_x-1, \ 0\leqslant r\leqslant Q_t-1;
    \end{split} \\
    \begin{split}
      &
      \sum_{j=1}^M V_j^{e_{mnl}}(x_p^{e_{mnl}},b_2,t_r^{e_{mnl}})w_j^{e_{mnl}}
      - g(x_p^{e_{mnl}},b_2,t_r^{e_{mnl}})=0, \\
      & \qquad\qquad \text{for} \
      n = N_y-1, \
      0\leqslant m\leqslant N_x-1, \ 0\leqslant l\leqslant N_t-1, \
      0\leqslant p\leqslant Q_x-1, \ 0\leqslant r\leqslant Q_t-1. \label{equ_13d}
    \end{split}
  \end{align}
\end{subequations}
%
On the boundary $t=0$ of the spatial-temporal domain,
the initial condition \eqref{equ_10c} is reduced to 
\begin{equation}
  \begin{split}
    &
    \sum_{j=1}^M V_j^{e_{mn0}}(x_p^{e_{mn0}},y_q^{e_{mn0}},0) w_j^{e_{mn0}}
    - h(x_p^{e_{mn0}},y_q^{e_{mn0}}) = 0, \\
    & \qquad\qquad \text{for} \
    0\leqslant m\leqslant N_x-1, \ 0\leqslant n\leqslant N_y-1, \
    0\leqslant p\leqslant Q_x-1, \ 0\leqslant q\leqslant Q_y-1. \label{equ_14}
  \end{split}
\end{equation}


Since $L$ is assumed to be a second-order operator with respect to
both $x$ and $y$, we impose $C^1$ continuity conditions across the
sub-domain boundaries in both the $x$ and $y$ directions.
Because equation \eqref{equ_10a} is of first order in time, we impose the
$C^0$ continuity condition across the sub-domain boundaries along
the temporal direction.
On the sub-domain boundaries $x=X_{m+1}$ ($0\leqslant m\leqslant N_x-2$),
the $C^1$ conditions become,
\begin{subequations}
  \begin{align}
    \begin{split}
      &
      \sum_{j=1}^M V_j^{e_{mnl}}(X_{m+1},y_q^{e_{mnl}},t_r^{e_{mnl}}) w_j^{e_{mnl}}
      - \sum_{j=1}^M V_j^{e_{m+1,nl}}(X_{m+1},y_q^{e_{m+1,nl}},t_r^{e_{m+1,nl}}) w_j^{e_{m+1,nl}}
      = 0, \label{equ_15a} \\
      & \qquad
      0\leqslant m\leqslant N_x-2, \ 0\leqslant n\leqslant N_y-1, \
      0\leqslant l\leqslant N_t-1, \ 0\leqslant q\leqslant Q_y-1, \
      0\leqslant r\leqslant Q_t-1; 
    \end{split} \\
    \begin{split}
      &
      \sum_{j=1}^M \left.\frac{\partial V_j^{e_{mnl}}}{\partial x}\right|_{(X_{m+1},y_q^{e_{mnl}},t_r^{e_{mnl}})} w_j^{e_{mnl}}
      - \sum_{j=1}^M \left.\frac{\partial V_j^{e_{m+1,nl}}}{\partial x}\right|_{(X_{m+1},y_q^{e_{m+1,nl}},t_r^{e_{m+1,nl}})} w_j^{e_{m+1,nl}}
      = 0, \\
      & \qquad
      0\leqslant m\leqslant N_x-2, \ 0\leqslant n\leqslant N_y-1, \
      0\leqslant l\leqslant N_t-1, \ 0\leqslant q\leqslant Q_y-1, \
      0\leqslant r\leqslant Q_t-1.
    \end{split}
  \end{align}
\end{subequations}
On the sub-domain boundaries $y=Y_{n+1}$ ($0\leqslant n\leqslant N_y-2$)
the $C^1$ continuity conditions become,
\begin{subequations}
  \begin{align}
    \begin{split}
      &
      \sum_{j=1}^M V_j^{e_{mnl}}(x_p^{e_{mnl}},Y_{n+1},t_r^{e_{mnl}}) w_j^{e_{mnl}}
      - \sum_{j=1}^M V_j^{e_{m,n+1,l}}(x_p^{e_{m,n+1,l}},Y_{n+1},t_r^{e_{m,n+1,l}}) w_j^{e_{m,n+1,l}}
      = 0, \\
      & \qquad
      0\leqslant m\leqslant N_x-1, \ 0\leqslant n\leqslant N_y-2, \
      0\leqslant l\leqslant N_t-1, \ 0\leqslant p\leqslant Q_x-1, \
      0\leqslant r\leqslant Q_t-1; 
    \end{split} \\
    \begin{split}
      &
      \sum_{j=1}^M \left.\frac{\partial V_j^{e_{mnl}}}{\partial y}\right|_{(x_p^{e_{mnl}},Y_{n+1},t_r^{e_{mnl}})} w_j^{e_{mnl}}
      - \sum_{j=1}^M \left.\frac{\partial V_j^{e_{m,n+1,l}}}{\partial y}\right|_{(x_p^{e_{m,n+1,l}},Y_{n+1},t_r^{e_{m,n+1,l}})} w_j^{e_{m,n+1,l}}
      = 0, \\
      & \qquad
      0\leqslant m\leqslant N_x-1, \ 0\leqslant n\leqslant N_y-2, \
      0\leqslant l\leqslant N_t-1, \ 0\leqslant p\leqslant Q_x-1, \
      0\leqslant r\leqslant Q_t-1. \label{equ_16b}
    \end{split}
  \end{align}
\end{subequations}
On the sub-domain boundaries $t=T_{l+1}$ ($0\leqslant l\leqslant N_t-2$),
the $C^0$ continuity conditions become,
\begin{equation}\label{equ_17}
  \begin{split}
    &
    \sum_{j=1}^M V_j^{e_{mnl}}(x_p^{e_{mnl}},y_q^{e_{mnl}},T_{l+1}) w_j^{e_{mnl}}
      - \sum_{j=1}^M V_j^{e_{mn,l+1}}(x_p^{e_{mn,l+1}},y_q^{e_{mn,l+1}},T_{l+1}) w_j^{e_{mn,l+1}}
      = 0, \\
      & \qquad
      0\leqslant m\leqslant N_x-1, \ 0\leqslant n\leqslant N_y-1, \
      0\leqslant l\leqslant N_t-2, \ 0\leqslant p\leqslant Q_x-1, \
      0\leqslant q\leqslant Q_y-1.
  \end{split}
\end{equation}


The equations consisting of \eqref{equ_11}--\eqref{equ_17}
form a system of linear algebraic equations about
the training parameters $w_j^{e_{mnl}}$ ($1\leqslant j\leqslant M$,
$0\leqslant m\leqslant N_x-1$, $0\leqslant n\leqslant N_y-1$
and $0\leqslant l\leqslant N_t-1$).
In these equations, $V_j^{e_{mnl}}$, $\frac{\partial V_j^{e_{mnl}}}{\partial t}$,
$\frac{\partial V_j^{e_{mnl}}}{\partial x}$,
$\frac{\partial V_j^{e_{mnl}}}{\partial y}$ and
$LV_j^{e_{mnl}}$ are all known functions and can be evaluated on the collocation
points by the local neural networks.
In particular, the partial derivatives therein
can be computed based on auto-differentiation.

This linear system consists of
$ 
  N_{equ} = N_xN_yN_t\left[
    Q_xQ_yQ_t + 2(Q_x+Q_y)Q_t + Q_xQ_y
    \right]
$ 
equations, and is about $N_xN_yN_tM$ unknown variables $w_j^{e_{mnl}}$.
We seek a least squares solution to this system with minimum norm, and
compute this solution by the linear least squares method.
In the implementation we employ the linear least squares routine from LAPACK
to compute the least squares solution.
The weight coefficients in the output layers of the local neural networks
are then determined by the least squares solution to the above system.
Training the neural network basically consists of
computing the least squares solution.

\paragraph{Block Time-Marching for Long-Time Simulations}
\label{sec:block}


Since the linear least squares computation, and hence
the neural network training, 
is computationally fast,
longer-time dynamic simulations of time-dependent PDEs
become feasible using
the current method.
With the basic method,
we observe that 
as the temporal dimension
of the spatial-temporal domain (i.e.~$\Gamma$)
increases, the network training 
generally becomes more difficult, in the sense that
the obtained solution tends to become less accurate
corresponding to the later time instants in the domain.
When $\Gamma$ is large, the solution can contain pronounced errors.
Therefore, using a large dimension in time (i.e.~large $\Gamma$)
 with the basic method is generally not advisable.

To perform long-time simulations,
we will employ the following block time-marching strategy.
Given a spatial-temporal domain with a large dimension in time,
we divide the domain into a number of windows, referred to
as time blocks, along the temporal direction, so that
the temporal dimension of each time block has a moderate size.
We then solve the initial-boundary value problem using
the basic method as discussed above
on the spatial-temporal domain of each time block,
individually and successively.
We use the solution from the previous time block evaluated at
the last time instant as the initial condition for the
computations of the current time block.
We start with the first time block, and march
in time block by block, until the last time block is completed.

Specifically, let
$
\Omega = \{
(x,y,t)\ |\ x\in[a_1,b_1], \ y\in [a_2,b_2], \
t\in [0, t_f]
\}
$
 denote the spatial-temporal domain on which the
initial-boundary value problem \eqref{equ_10a}--\eqref{equ_10c}
is to be solved, where $t_f$ can be large.
We divide the domain into $N_b$ ($N_b\geqslant 1$)
uniform blocks in time, with each block the size of
$\Gamma = \frac{t_f}{N_b}$. We choose
$N_b$ such that the block size $\Gamma$
is a moderate value.

On the $k$-th ($0\leqslant k\leqslant N_b-1$) time block,
we introduce a time shift and a new dependent
variable as a function of the shifted time
based on the following transform:
\begin{align}
  &
  \xi = t - k\Gamma, \ \ U(x,y,\xi) = u(x,y,t),
  \quad t\in [k\Gamma, (k+1)\Gamma], \ \
  \xi \in [0, \Gamma], \label{equ_18}
\end{align}
where $\xi$ denotes the shifted time and
$U(x,y,\xi)$ denotes the new dependent variable.
The equations \eqref{equ_10a} and \eqref{equ_10b}
are then transformed into,
\begin{subequations}
  \begin{align}
    &
    \frac{\partial U}{\partial \xi} = LU + f(x,y,\xi+k\Gamma), \label{equ_19a}
    \\
    &
    U(x,y,\xi) = g(x,y,\xi+k\Gamma),
    \quad \text{for} \ (x,y) \ \text{on spatial domain boundary}.
    \label{equ_19b}
  \end{align}
\end{subequations}
This is supplemented by the initial condition,
\begin{equation}\label{equ_20}
  U(x,y,0) = U_0(x,y),
\end{equation}
where $U_0(x,y)$ denotes the initial distribution on the time block $k$,
given by
\begin{equation}\label{equ_21}
  U_0(x,y) = \left\{
  \begin{array}{ll}
    u(x,y,0) = h(x,y), & \text{if} \ k=0, \\
    u(x,y,k\Gamma) \ \text{computed on time block} \ (k-1), & \text{if} \ k>0.
  \end{array}
  \right.
\end{equation}
Note that $h(x,y)$ is the initial condition for the problem.

The initial-boundary value problem on time block $k$ now
consists of equations \eqref{equ_19a}, \eqref{equ_19b} and
\eqref{equ_20}, to be solved on the
spatial-temporal domain
$ 
  \Omega^{st} = \{
  (x,y,\xi) \ |\ x\in[a_1,b_1], \ y\in[a_2,b_2],\
  \xi \in [0,\Gamma]
  \}
$ 
for the function $U(x,y,\xi)$.
This is the same problem we have considered previously,
and it can be solved using the basic method discussed before. 
With $U(x,y,\xi)$ obtained, the function $u(x,y,t)$ on time block $k$
is recovered by the transform \eqref{equ_18}.
By solving the initial-boundary value problem on successive
time blocks, we can attain the solution $u(x,y,t)$
on the entire spatial-temporal domain $\Omega$.
This is the block time-marching scheme 
for potentially long-time simulations of time-dependent linear
PDEs.

\subsection{Nonlinear Differential Equations}
\label{sec:nonlinear}

%



In this section
we look into how to solve the initial/boundary value problems
involving nonlinear differential equations using domain
decomposition and the locELM representation
for the solutions.
The overall procedure is analogous to that for linear
differential equations. The main difference lies in that
here the set of local neural networks
needs to be  trained by a nonlinear least squares computation.

\subsubsection{Time-Independent Nonlinear Differential Equations}
\label{sec:nonl_steady}

%

We first consider the boundary value problems involving nonlinear
differential equations together with Dirichlet boundary conditions, and
discuss how to solve such problems using the locELM method.
We assume that the highest-orer terms in the equation are linear,
and that the nonlinear terms involve
the unknown function and also possibly its derivatives of lower orders.
To make the discussions more concrete, we again
focus on two dimensions (with coordinates
$x$ and $y$), and assume that the highest partial derivatives
with respect to both $x$ and $y$ are of second order in the equation.

Let us consider the following generic second-order nonlinear differential
equation of such a form on domain $\Omega$,
together with the Dirichlet boundary condition on $\partial\Omega$,
\begin{subequations}
  \begin{align}
    &
    Lu + F\left(u,u_x,u_y\right)
    = f(x,y), \label{equ_32a}
    \\
    &
    u(x,y) = g(x,y), \quad \text{on}\ \partial\Omega,
     \label{equ_32b}
  \end{align}
\end{subequations}
where $u(x,y)$ is the field function to be solved for,
$u_x=\frac{\partial u}{\partial x}$,
$u_y=\frac{\partial u}{\partial y}$,
$L$ is a second-order linear differential operator with respect to
both $x$ and $y$, $F$ denotes the nonlinear term,
$f(x,y)$ is a prescribed source term, and
$g(x,y)$ denotes the Dirichlet boundary data.

The overall procedure for 
solving equations~\eqref{equ_32a}--\eqref{equ_32b} using the locELM method
is analogous to that in Section \ref{sec:steady}.
We focus on a rectangular domain,
$
\Omega = \{
(x,y)\ |\ x\in[a_1,b_1], \ y\in[a_2,b_2]
\},
$
and partition this domain into $N_x$ and $N_y$ sub-domains along the
$x$ and $y$ directions, respectively, thus leading to a total of
$N_e=N_xN_y$ sub-domains in $\Omega$.
Following the notation of Section \ref{sec:steady},
we denote the sub-domain boundary coordinates along the $x$ and $y$ directions
by two vectors $[X_0,X_1,\dots,X_{N_x}]$ and $[Y_0,Y_1,\dots,Y_{N_y}]$,
respectively.
Let $\Omega_{e_{mn}}=[X_m,X_{m+1}]\times[Y_n,Y_{n+1}]$ denote
the sub-domain with index $e_{mn}$ for $0\leqslant m\leqslant N_x-1$
and $0\leqslant n\leqslant N_y-1$.
We use $(x_p^{e_{mn}},y_q^{e_{mn}})$ ($0\leqslant p\leqslant Q_x-1$,
$0\leqslant q\leqslant Q_y-1$)
to denote a set of uniform collocation points
in the sub-domain $e_{mn}$, where $Q_x$ and $Q_y$ are
the number of collocation points in the $x$ and $y$ directions
on the sub-domain.
%
The input layer of the local neural network consists of two nodes ($x$ and $y$),
and the output layer consists of one node (representing $u$).
Let $u^{e_{mn}}(x,y)$ denote the output of the local neural network
on the sub-domain $e_{mn}$, and $V_j^{e_{mn}}(x,y)$ ($1\leqslant j\leqslant M$)
denote the output of the last hidden layer of the local neural network,
where $M$ is the number of nodes in the last hidden layer.
We have the following relations,
\begin{equation}\label{equ_33}
  \begin{split}
    &
    u^{e_{mn}}(x,y) = \sum_{j=1}^M V_j^{e_{mn}}(x,y) w_j^{e_{mn}}, \quad
    \frac{\partial u^{e_{mn}} }{\partial x}
    = \sum_{j=1}^M \frac{\partial V_j^{e_{mn}}}{\partial x} w_j^{e_{mn}}, \quad
    \frac{\partial u^{e_{mn}} }{\partial y}
    = \sum_{j=1}^M \frac{\partial V_j^{e_{mn}}}{\partial y} w_j^{e_{mn}}, \\
    & \qquad \text{for} \
    0\leqslant m\leqslant N_x-1, \
    0\leqslant n\leqslant N_y-1,
  \end{split}
\end{equation}
where the constants $w_j^{e_{mn}}$ ($1\leqslant j\leqslant M$) denote
the weight coefficients in the output layer of the local neural network
on sub-domain $e_{mn}$, and they constitute the training parameters
of the neural network.


Enforcing equation \eqref{equ_32a} on the collocation points
$(x_p^{e_{mn}},y_q^{e_{mn}})$ for each sub-domain leads to
\begin{equation}\label{equ_34}
  \begin{split}
    &
  \sum_{j=1}^M \left[LV_j^{e_{mn}}(x_p^{e_{mn}},y_q^{e_{mn}})\right] w_j^{e_{mn}}
  + F\left.\left(u^{e_{mn}},u_x^{e_{mn}},u_y^{e_{mn}}\right)\right|_{(x_p^{e_{mn}},y_q^{e_{mn}})}
  - f(x_p^{e_{mn}},y_q^{e_{mn}}) = 0,
  \\
  & \qquad \text{for}\
  0\leqslant m\leqslant N_x-1, \
  0\leqslant n\leqslant N_y-1, \
  0\leqslant p\leqslant Q_x-1, \
  0\leqslant q\leqslant Q_y-1,
  \end{split}
\end{equation}
where $u^{e_{mn}}$, $u_x^{e_{mn}}$ and $u_y^{e_{mn}}$
are given by \eqref{equ_33} in terms of
the training parameters $w_j^{e_{mn}}$.
Enforcing the boundary condition~\eqref{equ_32b}
on the collocation points of the
four domain boundaries $x=a_1$ or $b_1$
and $y=a_2$ or $b_2$ leads to the equations
\eqref{equ_4}, \eqref{equ_5}, \eqref{equ_6} and \eqref{equ_7}.
Since equation \eqref{equ_32a} is of second-order 
with respect to both $x$ and $y$, we impose $C^1$
continuity conditions across the sub-domain boundaries
along both the $x$ and $y$ directions.
Enforcing the $C^1$ continuity conditions on
the collocation points of the sub-domain boundaries
$x=X_{m+1}$ ($0\leqslant m\leqslant N_x-2$) and
$y=Y_{n+1}$ ($0\leqslant n\leqslant N_y-2$) leads to
the equations \eqref{equ_8a}--\eqref{equ_8b}
and \eqref{equ_9a}--\eqref{equ_9b}.

The set of equations consisting of \eqref{equ_34},
\eqref{equ_4}--\eqref{equ_7}, \eqref{equ_8a}--\eqref{equ_8b}
and \eqref{equ_9a}--\eqref{equ_9b}
is a system of nonlinear algebraic equations about
the training parameters
$w_j^{e_{mn}}$ ($1\leqslant j\leqslant M$,
$0\leqslant m\leqslant N_x-1$, $0\leqslant n\leqslant N_y-1$).
In these equations the functions $V_j^{e_{mn}}(x,y)$ are all known
and their partial derivatives can be computed by auto-differentiation.
%
This nonlinear algebraic system consists of
$ 
  N_xN_y(Q_xQ_y + 2Q_x + 2Q_y)
$ 
equations with $N_xN_yM$ unknowns.

This system is to be solved for the determination of the training parameters.
In this paper we consider two methods for
solving this system.
In the first method
we seek a least squares solution to this system for
the training parameters $w_j^{e_{mn}}$, thus leading to
a nonlinear least squares problem.
In the second method we adopt a simple Newton's method combined with
a linear least squares computation for solving this system.

\begin{algorithm}[tb]
  \DontPrintSemicolon
  \SetKwInOut{Input}{input}\SetKwInOut{Output}{output}

  \Input{constant $\delta>0$, initial guess  $\mbs x_0$}
  \Output{solution vector $\mbs x$, associated cost $c$}
  \BlankLine\BlankLine
  call scipy.optimize.least\_squares routine using $\mbs x_0$ as the initial guess\;  
  set $\mbs x\leftarrow$ returned solution\;
  set $c\leftarrow$ returned cost\;
  \If{c is below a threshold}{return\;}
  \BlankLine\BlankLine
  \For{$i\leftarrow 0$ \KwTo maximum number of sub-iterations}{
    generate a random number $\xi_1$ on the interval $[0,1]$\;
    set $\delta_1 \leftarrow \xi_1\delta$\;
    generate a uniform random vector $\Delta\mbs x$ of the
    same shape as $\mbs x$ on the interval [$-\delta_1$, $\delta_1$]\;
    \BlankLine\BlankLine
    generate a random number $\xi_2$ on the interval $[0,1]$\;
    set $\mbs y_0 \leftarrow \xi_2\mbs x + \Delta\mbs x$\;
    \BlankLine\BlankLine
    call scipy.optimize.least\_squares routine using $\mbs y_0$ as the initial guess\;
    \If{the returned cost is less than $c$}{
      set $\mbs x\leftarrow$ the returned solution\;
      set $c\leftarrow$ the returned cost\;
    }
    \If{the returned cost is below a threshold}{
      return\;
    }
  }
  \caption{NLSQ-perturb (nonlinear least squares with perturbations)}
  \label{alg:alg_1}
\end{algorithm}

With the first method, to solve the nonlinear least squares problem,
we take advantage of the nonlinear least squares implementations
from the scientific libraries.
In the current implementation, we employ the nonlinear least squares
routine ``least\_squares''
from the scipy.optimize package.
This method typically works quite well, and exhibits a
smooth convergence behavior.
However, we observe that
in certain cases, e.g.~when the simulation resolution is not sufficient or
sometimes in longer-time simulations with time-dependent nonlinear equations,
this method at times can be attracted to
and trapped in a local
minimum solution. While the method indicates that the nonlinear iterations
have converged,
the norm of the converged equation residuals can turn out to be quite pronounced
in magnitude. 
In the event this takes place,  the obtained solution
can contain significant errors and the simulation 
loses accuracy from that point onward. This issue is typically encountered
when the resolution of the computation
(e.g.~the number of collocation points in the domain
or the number of training parameters in the neural network) decreases 
to a certain point.
This has been a main issue with the nonlinear least squares
computation using this method.

To alleviate this problem and make the nonlinear least squares
computation more robust,
we find it necessary to incorporate a sub-iteration procedure with random
perturbations to the
initial guess when invoking the nonlinear
least squares routine.
The basic idea is as follows.
If the nonlinear least squares routine converges with
the converged cost (i.e.~norm of the equation residual) exceeding a threshold,
the sub-iteration procedure will be triggered. Within
each sub-iteration a random initial guess for the solution
is generated, based on e.g.~a perturbation to the current
approximation of the solution vector,
and is fed to the nonlinear least squares routine.

Algorithm \ref{alg:alg_1} illustrates
the nonlinear least squares computation combined with the sub-iteration procedure,
which will be referred to as the NLSQ-perturb
(Nonlinear Least SQuares with perturbations) method hereafter.
In this algorithm the parameter $\delta$ controls the maximum range on
which the random perturbation vector is generated.
Numerical experiments indicate that the method works better
if $\delta$ is not large. A typical value is $\delta=0.5$, which is observed
to work well in numerical simulations.
Combined with an appropriate resolution (the number of collocation
points in domain, and the number of training parameters in the neural network)
for a given problem,
the NLSQ-perturb method turns out to be very effective.
The solution can typically be attained with only
a few (e.g.~around $4$ or $5$) sub-iterations if such an iteration
is triggered.
For the numerical tests reported in Section \ref{sec:tests},
we employ a threshold value $10^{-3}$ in the lines $4$ and $18$
of Algorithm \ref{alg:alg_1}. The final converged cost value
is typically on the order $10^{-13}$.

\begin{remark}\label{rem_9}
In Algorithm \ref{alg:alg_1}
the value $\xi_2$ controls around which point the random perturbation 
will be generated. In Algorithm \ref{alg:alg_1}, $\xi_2$
is taken to be a random value from $[0,1]$.
An alternative to this is
to fix this value at $\xi_2=0$ or $\xi_2=1$, which has been observed
to work  well in actual simulations.
By using $\xi_2=0$, one is effectively generating a random perturbation
around the origin and use it as the initial guess.
By using $\xi_2=1$, one is effectively setting the initial guess
as a random perturbation
to the best approximation obtained so far.
  
\end{remark}

The second method for solving the nonlinear algebraic system
is a combination of Newton iterations
with linear least squares computations, which we will refer to
as the Newton-LLSQ (Newton-Linear Least SQuares) method hereafter.
The convergence behavior of this method is not
as regular as the first method, but it appears less likely
to be trapped to local minimum solutions.
To outline the idea of the method, let
\begin{equation}
  {\mbs G}(\mbs W) = 0, \quad
  \text{where}\ \mbs G=(G_1,G_2,\dots,G_m), \
  \mbs W=(w_1,w_2,\dots,w_n)
\end{equation}
denote a system of $m$ nonlinear algebraic equations about
$n$ variables $\mbs W$. Let the superscript in ${\mbs W}^{(k)}$ denote
the approximation of the solution at the $k$-th iteration, and
$\Delta \mbs W$ denote the solution increment.
We update the solution iteratively as follows in a way similar to the
Newton's method,
\begin{align}
  &
  \mbs J(\mbs W^{(k)})\Delta \mbs W = -\mbs G(\mbs W^{(k)}), \label{equ_39}
  \\
  &
  \mbs W^{(k+1)} = \mbs W^{(k)} + \Delta \mbs W,
\end{align}
where $\mbs J(\mbs W^{(k)})$ is the Jacobian matrix given by
$
J_{ij} = \frac{\partial G_i}{\partial w_j}
$
($1\leqslant i\leqslant m$, $1\leqslant j\leqslant n$)
and evaluated at $\mbs W^{(k)}$.
The departure point from the standard Newton method lies in
that the linear algebraic system \eqref{equ_39}
involves a non-square coefficient
matrix (Jacobian matrix).
We seek a least squares solution to the linear system \eqref{equ_39},
and solve this system for
the increment $\Delta \mbs W$ using the linear least squares routine
from LAPACK.


\begin{remark}
  \label{rem_3}
  It is observed that the computational cost of the
  Newton-LLSQ method is typically considerably smaller than that of the
  NLSQ-perturb method in training the locELM neural networks. On the other hand,
  the locELM solutions obtained with the Newton-LLSQ method are
  in general markedly less accurate than those obtained using
  the NLSQ-perturb method.

\end{remark}




In the current work, we implement the local neural networks
for each sub-domain $e_{mn}$
using one or several dense Keras layers, with the collocation
points $(x_p^{e_{mn}},y_q^{e_{mn}})$ as the input data and
$u^{e_{mn}}$ as the output. In the implementation, an affine mapping is
incorporated into each local neural network
behind the input layer to normalize the input
$(x,y)$ data to the interval $[-1,1]\times[-1,1]$ for each sub-domain.
The set of local neural networks logically forms
a multiple-input multiple-output Keras model.
The weight/bias coefficients in all the hidden layers are set to
uniform random values generated on $[-R_m,R_m]$.
The weight coefficients of the output layers ($w_j^{e_{mn}}$) of
the local neural networks are determined and
set by the solution to the nonlinear algebraic system
obtained using the NLSQ-perturb or Newton-LLSQ methods.
The partial derivatives
involved in the formulation are computed by  auto-differentiation
from the Tensorflow package.

\subsubsection{Time-Dependent Nonlinear Differential Equations}
\label{sec:tnleq}

We next consider the initial-boundary value problems
involving time-dependent nonlinear different equations together with Dirichlet
boundary conditions, and discuss how to solve such problems using the locELM method.
We make the same assumptions about the differential equation
as in Section \ref{sec:nonl_steady}: The highest-order terms
are assumed to be linear, and the nonlinear terms may involve
the unknown function or its partial derivatives of lower orders.
We again focus on two spatial dimensions, plus time $t$,
and assume that
the equation is of second order with respect to both
spatial coordinates ($x$ and $y$).

Consider the following generic nonlinear partial
differential equation of such a form on a spatial-temporal domain $\Omega$,
supplemented by the Dirichlet boundary condition and an initial condition,
\begin{subequations}
  \begin{align}
    &
    \frac{\partial u}{\partial t} = Lu + F(u,u_x,u_y) + f(x,y,t),
    \label{equ_35a}
    \\
    &
    u(x,y,t) = g(x,y,t), \quad
    \text{for} \ (x,y)\ \text{on the spatial domain boundary}, \label{equ_35b}
    \\
    &
    u(x,y,0) = h(x,y), \label{equ_35c}
  \end{align}
\end{subequations}
where $u(x,y,t)$ is the unknown field function to be solved for,
$L$ is a second-order linear differential operator with
respect to both $x$ and $y$, $F$ denotes the nonlinear term,
$f(x,y,t)$ is a prescribed source term, $g(x,y,t)$ denotes
the Dirichlet boundary data, and $h(x,y)$ is the initial
field distribution.

Our discussion below largely parallels that of Section \ref{sec:unsteady}.
We first discuss the basic method on a spatial-temporal domain,
and then develop the block time-marching idea for longer-time simulations
of the nonlinear partial differential equations.

\paragraph{Basic Method}

We focus on a rectangular spatial-temporal domain
$
\Omega = \{
(x,y,t)\ |\ x\in[a_1,b_1], \ y\in[a_2,b_2], \ t\in[0,\Gamma]
\},
$
and solve the initial-boundary value problem consisting of
equations \eqref{equ_35a}--\eqref{equ_35c} on this domain.

Following the notation of Section~\ref{sec:basic}, we
use $N_x$, $N_y$ and $N_t$ to denote the number of sub-domains
along the $x$, $y$ and $t$ directions,
where the locations of the sub-domain boundaries along the three directions
are given by
the vectors $[X_0,X_1,\dots,X_{N_x}]$, $[Y_0,Y_1,\dots,Y_{N_y}]$
and $[T_0,T_1,\dots,T_{N_t}]$,
respectively. A sub-domain with the index $e_{mnl}$ corresponds to
the spatial-temporal region
$
\Omega_{e_{mnl}}=[X_m,X_{m+1}]\times[Y_n,Y_{n+1}]\times[T_l,T_{l+1}],
$
for $0\leqslant m\leqslant N_x-1$, $0\leqslant n\leqslant N_y-1$
and $0\leqslant l\leqslant N_t-1$.
Let $(x_p^{e_{mnl}},y_q^{e_{mnl}},t_r^{e_{mnl}})$
($0\leqslant p\leqslant Q_x-1$, $0\leqslant q\leqslant Q_y-1$,
$0\leqslant r\leqslant Q_t-1$)
denote the set of
$Q=Q_xQ_yQ_t$ collocation points on each sub-domain $e_{mnl}$.
%
Let $u^{e_{mnl}}(x,y,t)$ denote the output of the local neural network
corresponding to  the sub-domain $e_{mnl}$,
and $V_j^{e_{mnl}}(x,y,t)$ ($1\leqslant j\leqslant M$)
denote the output of the last hidden layer of the local neural network,
where $M$ is the number of nodes in the last hidden layer.
The following relations hold,
\begin{equation}\label{equ_36}
  \left\{
  \begin{split}
    &
    u^{e_{mnl}}(x,y,t) = \sum_{j=1}^M V_j^{e_{mnl}}(x,y,t) w_j^{e_{mnl}}, \quad
    u_x^{e_{mnl}}(x,y,t) = \sum_{j=1}^M \frac{\partial V_j^{e_{mnl}}}{\partial x} w_j^{e_{mnl}},
    \\
    &
    u_y^{e_{mnl}}(x,y,t) = \sum_{j=1}^M \frac{\partial V_j^{e_{mnl}}}{\partial y} w_j^{e_{mnl}},
    \quad
    \frac{\partial u^{e_{mnl}}}{\partial t}
    = \sum_{j=1}^M \frac{\partial V_j^{e_{mnl}}}{\partial t} w_j^{e_{mnl}},
    \\
    &
    \text{for}\ 0\leqslant m\leqslant N_x-1, \
    0\leqslant n\leqslant N_y-1, \
    0\leqslant l\leqslant N_t-1,
  \end{split}
  \right.
\end{equation}
where $w_j^{e_{mnl}}$ denote the weight coefficients
in the output layers of the local neural networks
and they constitute the training parameters of the network.


Enforcing equation \eqref{equ_35a} on the collocation points
$(x_p^{e_{mnl}},y_q^{e_{mnl}},t_r^{e_{mnl}})$ of each sub-domain $e_{mnl}$
leads to
\begin{equation}\label{equ_37}
  \begin{split}
    &
    \sum_{j=1}^M\left.\left[
      \frac{\partial V_j^{e_{mnl}}}{\partial t}
      -L V_j^{e_{mnl}}
      \right]\right|_{(x_p^{e_{mnl}},y_q^{e_{mnl}},t_r^{e_{mnl}})}w_j^{e_{mnl}}
    - \left.F(u^{e_{mnl}},u_x^{e_{mnl}},u_y^{e_{mnl}}) \right|_{(x_p^{e_{mnl}},y_q^{e_{mnl}},t_r^{e_{mnl}})}
    \\
    & \qquad
    - f(x_p^{e_{mnl}},y_q^{e_{mnl}},t_r^{e_{mnl}}) = 0, \\
    &
    \text{for} \
    0\leqslant m\leqslant N_x-1, \
    0\leqslant n\leqslant N_y-1, \
    0\leqslant l\leqslant N_t-1, \
    0\leqslant p\leqslant Q_x-1, \
    0\leqslant q\leqslant Q_y-1, \\
    & \quad \ \
    0\leqslant r\leqslant Q_t-1,
  \end{split}
\end{equation}
where $u^{e_{mnl}}$, $u_x^{e_{mnl}}$ and $u_y^{e_{mnl}}$
are given by \eqref{equ_36} in terms of the known function
$V_j^{e_{mnl}}$ and its partial derivatives.
This is a set of nonlinear algebraic equations about
the training parameters $w_j^{e_{mnl}}$.
Enforcing the boundary condition \eqref{equ_35b}
on the collocation points of the four spatial
boundaries at $x=a_1$ or $b_1$ and $y=a_2$ or $b_2$
leads to the equations \eqref{equ_13a}--\eqref{equ_13d}.
Enforcing the initial condition \eqref{equ_35c}
on the spatial collocation points at $t=0$
results in equation \eqref{equ_14}.
We impose the $C^1$ continuity conditions on the unknown field 
$u(x,y,t)$ across the sub-domain boundaries
along the $x$ and $y$ directions, since $L$ is assumed to be
a second-order operator with respect to both $x$ and $y$.
We impose the $C^0$ continuity condition across the
sub-domain boundaries in
the temporal direction, since equation \eqref{equ_35a} is
first-order with respect to time.
Enforcing the $C^1$ continuity conditions
on the collocation points on the sub-domain boundaries
$x=X_{m+1}$ ($0\leqslant m\leqslant N_x-2$) and
$y=Y_{n+1}$ ($0\leqslant n\leqslant N_y-2$) leads to
the equations \eqref{equ_15a}--\eqref{equ_16b}.
Enforcing the $C^0$ continuity condition on
the collocation points on the sub-domain boundaries
$t=T_{l+1}$ ($0\leqslant l\leqslant N_t-2$) leads
to the equation \eqref{equ_17}.

The set of equations consisting of \eqref{equ_37} and
\eqref{equ_13a}--\eqref{equ_17} is a nonlinear algebraic system
of equations about the training parameters $w_j^{e_{mnl}}$.
%
%
This system consists of
$
N_xN_yN_t[Q_xQ_yQ_t+2(Q_x+Q_y)Q_t + Q_xQ_y]
$
coupled nonlinear algebraic equations with
$N_xN_yN_tM$ unknowns.
This system can be solved using the NLSQ-perturb or
Newton-LLSQ methods from Section \ref{sec:nonl_steady}
to determine the training parameters $w_j^{e_{mnl}}$.


\paragraph{Block Time-Marching}


For longer-time simulations of time-dependent nonlinear
differential equations, we employ a block time-marching
strategy analogous to that of Section \ref{sec:block}.
Let
$
\Omega = \{
(x,y,t) | x\in[a_1,b_1],\ y\in [a_2,b_2],\ t\in [0,t_f]
\}
$
denote the spatial-temporal domain on which the problem is to be solved,
where $t_f$ can be large.
We divide the temporal dimension into $N_b$ uniform time blocks,
with the block size $\Gamma = \frac{t_f}{N_b}$ being a moderate value,
and solve the problem on each time block separately and successively.
On the $k$-th ($0\leqslant k\leqslant N_b-1$) time block, we introduce
a shifted time $\xi$ and a new dependent variable $U(x,y,\xi)$ as given by
equation \eqref{equ_18}.
Then equation \eqref{equ_35a} is transformed into
  \begin{align}
    &
    \frac{\partial U}{\partial \xi} = LU + F(U, U_x, U_y)+ f(x,y,\xi+k\Gamma),
    \label{equ_41}
  \end{align}
where $U_x=\frac{\partial U}{\partial x}$ and $U_y=\frac{\partial U}{\partial y}$.
Equation~\eqref{equ_35b} is transformed into \eqref{equ_19b}.
The initial condition for time block $k$ is given by \eqref{equ_20},
in which the initial distribution data is given by \eqref{equ_21}.

The initial-boundary value problem consisting of equations
\eqref{equ_41}, \eqref{equ_19b} and \eqref{equ_20},
on the spatial-temporal domain
$
\Omega^{st}=[a_1,b_1]\times[a_2,b_2]\times[0,\Gamma],
$
is the same problem we have considered before,
and can be solved for $U(x,y,\xi)$ using the basic method.
The solution $u(x,y,t)$ on time block $k$
can then be recovered by the transform \eqref{equ_18}.


Starting with the first time block, we can solve the initial-boundary
value problem on each time block successively.
After the problem on the $k$-th block is solved,
the obtained solution can be evaluated at $t=(k+1)\Gamma$ 
and used as the initial condition for the computation on
the subsequent time block.


\begin{remark}
  \label{rem_6}
  We observe from numerical experiments that the time block size $\Gamma$
  can play a crucial role
  in long-time simulations of time-dependent
  nonlinear differential equations.
  In general, reducing $\Gamma$ can improve the
    convergence of the nonlinear iterations on the time blocks.
    If $\Gamma$ is too large, the nonlinear iterations
    can become hard to converge.
    With the other simulation parameters (such as the number of collocation
    points in the time block and the number of training parameters in
    the neural network) fixed, reducing the time block size effectively amounts to
    an increase in the resolution of the data on each time block.

  

\end{remark}





\begin{remark}\label{rem_dd}
  We will present numerical experiments with nonlinear 
  PDEs in Section \ref{sec:tests} to compare
  the current locELM method with the deep Galerkin method (DGM) and
  the physics-informed neural network (PINN), and also compare
  the current method with the classical finite element method (FEM).
  We observe that for these problems the locELM method is considerably superior
  to DGM and PINN,
  with regard to both the accuracy and the computational cost.
  In terms of the computational performance, 
  the locELM method is on par with the finite element method,
  and oftentimes the locELM performance  exceeds the FEM performance.

\end{remark}


\section{Numerical Examples}
\label{sec:tests}

%


In the forthcoming section we provide a number of numerical
examples to test the locELM method developed here.
These examples pertain to stationary and time-dependent,
linear and nonlinear differential equations.
They are in general one- or
two-dimensional (1D/2D) in space, and also plus time if time-dependent.
For certain problems (e.g.~the advection equation) we provide results
from long-time simulations, to demonstrate the capability
of the locELM method combined with the block time-marching scheme.
We employ $\tanh$ as the activation function in all the local
neural networks of this section.

In our discussion we focus on the
accuracy and the computational cost.
For locELM, the computational cost here refers to the total training time
of the overall neural network, which includes the computation time for
the output functions of the last hidden layer and its derivatives
(e.g.~$V_j^{e_{mnl}}$, $\frac{\partial V_j^{e_{mnl}}}{\partial x}$, etc),
the computation time for the coefficient matrix and the right hand side
of the least squares problem, and the solution time for the linear/nonlinear
squares problem. It does not include, after the training is over,
the evaluation of the neural network on a set of given points
for the output of the solution data.
The timing data is collected using the ``timeit'' module in Python.

We compare the current locELM method with
the deep Galerkin method (DGM)~\cite{SirignanoS2018}
and the physics-informed neural network (PINN) method~\cite{RaissiPK2019},
which are both based on deep neural
networks (DNN), in terms of the accuracy and the neural-network training time.
The DGM and PINN are trained using both the Adam~\cite{KingmaB2014}
and the L-BFGS~\cite{NocedalW2006} optimizers. For L-BFGS, we have employed
the routine available from the Tensorflow-Probability library
(www.tensorflow.org/probability).
For DGM and PINN, the training time refers to the time interval
between the start and the end of the Adam or L-BFGS training loop
for a given number of epochs/iterations.
The locELM, the DGM and the PINN methods
are all implemented in Python
with the Tensorflow (www.tensorflow.org) and Keras (keras.io) libraries.


Additionally, we compare the  locELM method with the classical
finite element method (linear elements, second-order), in terms of
the accuracy and computational cost.
For the numerical tests reported below,
the finite element method (FEM) is implemented also in Python,
using the FEniCS library (fenicsproject.org).
When the FEM code is run for the first time, the
FEniCS library uses Just-In-Time (JIT) compilers to compile
certain key finite element operations in the Python code into C++ code,
which is in turn compiled by the C++ compiler and then cached. 
This is done only once. So the FEM code is slower as JIT compilation
occurs when run for the first time, but it
is much faster in subsequent runs.
For FEM, the computational cost here refers to the computation time
collected using the ``timeit'' module
after the code has been compiled by the JIT compilers.
All the timing data with the locELM, DGM, PINN and FEM methods
is collected on a MAC computer ($3.2$GHz Intel Core i5
CPU, $24$GB memory) at the authors' institution.





\subsection{One-Dimensional Helmholtz Equation}
\label{sec:helm1d}


In the first test we consider the boundary value problem with
the one-dimensional (1D) Helmholtz equation on the domain $x\in[a,b]$,
\begin{subequations}
\begin{align}
  &
  \frac{d^2u}{dx^2} - \lambda u = f(x), \label{equ_t1} \\
  &
  u(a) = h_1, \\
  &
  u(b) = h_2, \label{equ_t3}
\end{align}
\end{subequations}
where $u(x)$ is the field function to be solved for,
$f(x)$ is a prescribed source term, $h_1$ and $h_2$ are the boundary
values, and the other constants in the above equations and the domain
specification are
\begin{equation*}
  \lambda = 10, \quad
  a = 0, \quad
  b = 8.
\end{equation*}
We choose the source term $f(x)$ such that the equation~\eqref{equ_t1}
has the following solution,
\begin{equation}\label{equ_t2}
  u(x) = \sin \left(3\pi x+\frac{3\pi}{20}\right) \cos \left(2\pi x+\frac{\pi}{10}\right) + 2.
\end{equation}
We choose $h_1$ and $h_2$ according to this
analytic solution by setting $x=a$ and $x=b$ in \eqref{equ_t2},
respectively.
Under these settings
the boundary value problem 
\eqref{equ_t1}--\eqref{equ_t3} has the analytic solution \eqref{equ_t2}.

\begin{figure}
  \centerline{
    \includegraphics[width=1.5in]{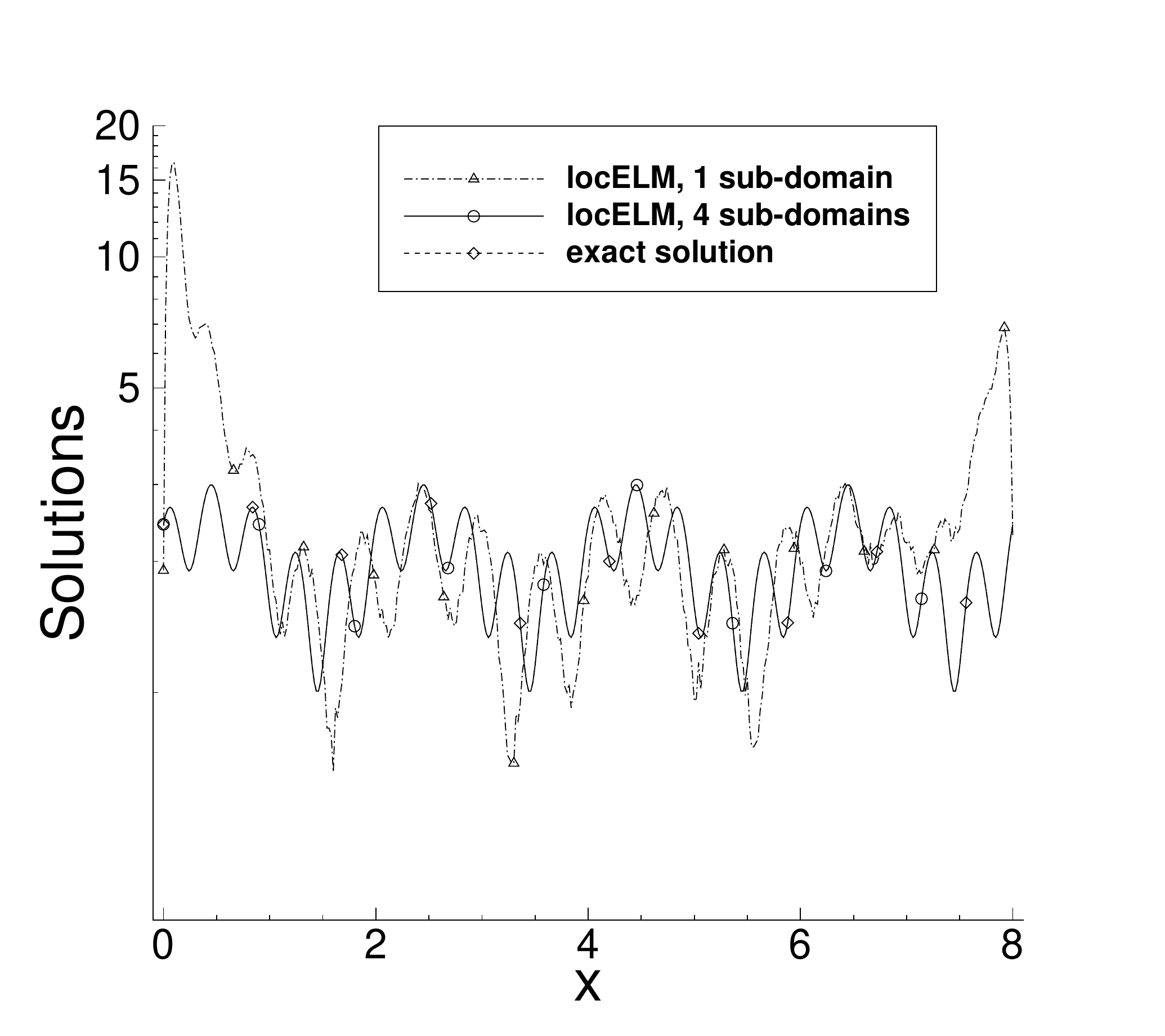}(a)
    \includegraphics[width=1.5in]{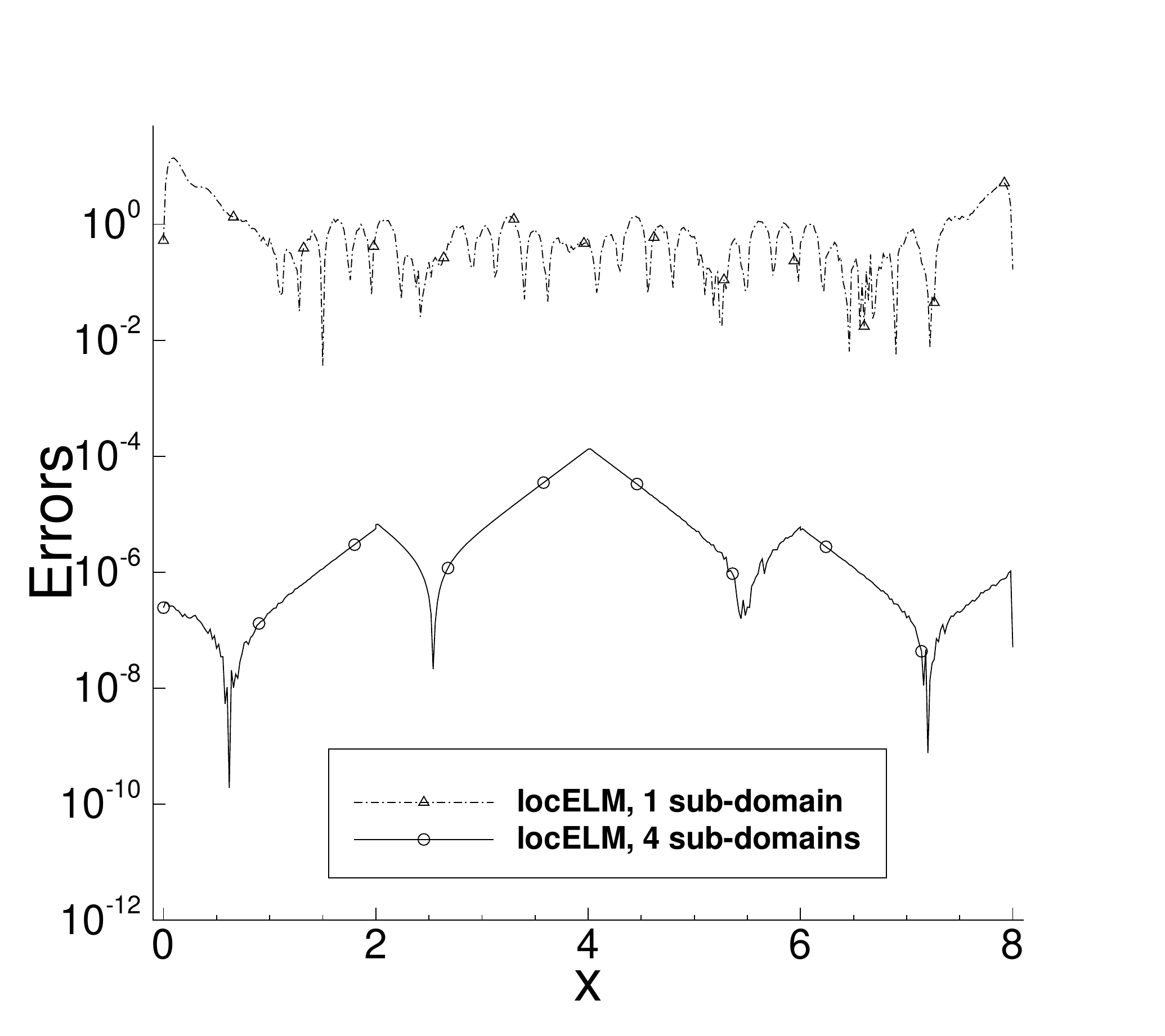}(b)
    \includegraphics[width=1.5in]{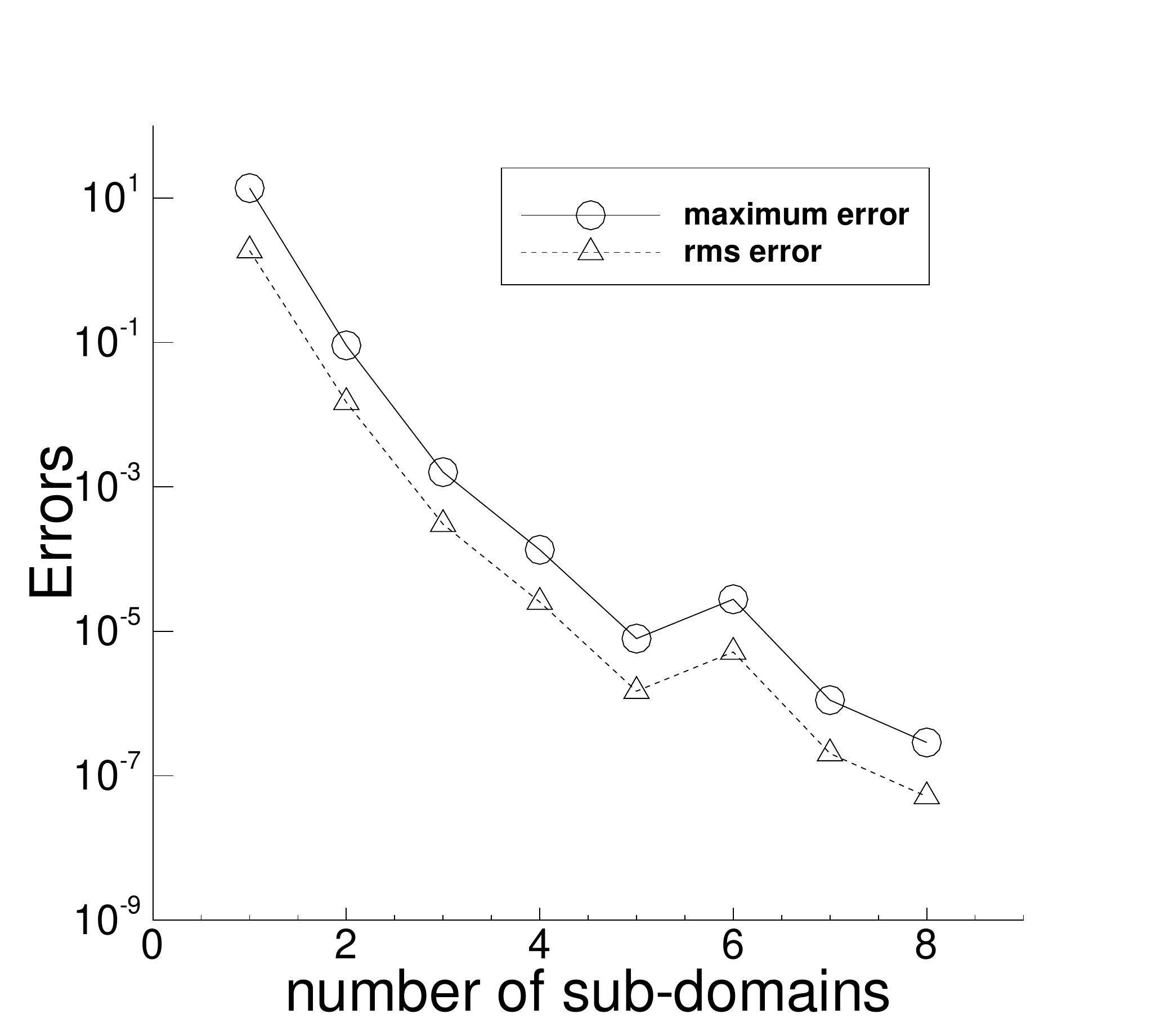}(c)
    \includegraphics[width=1.5in]{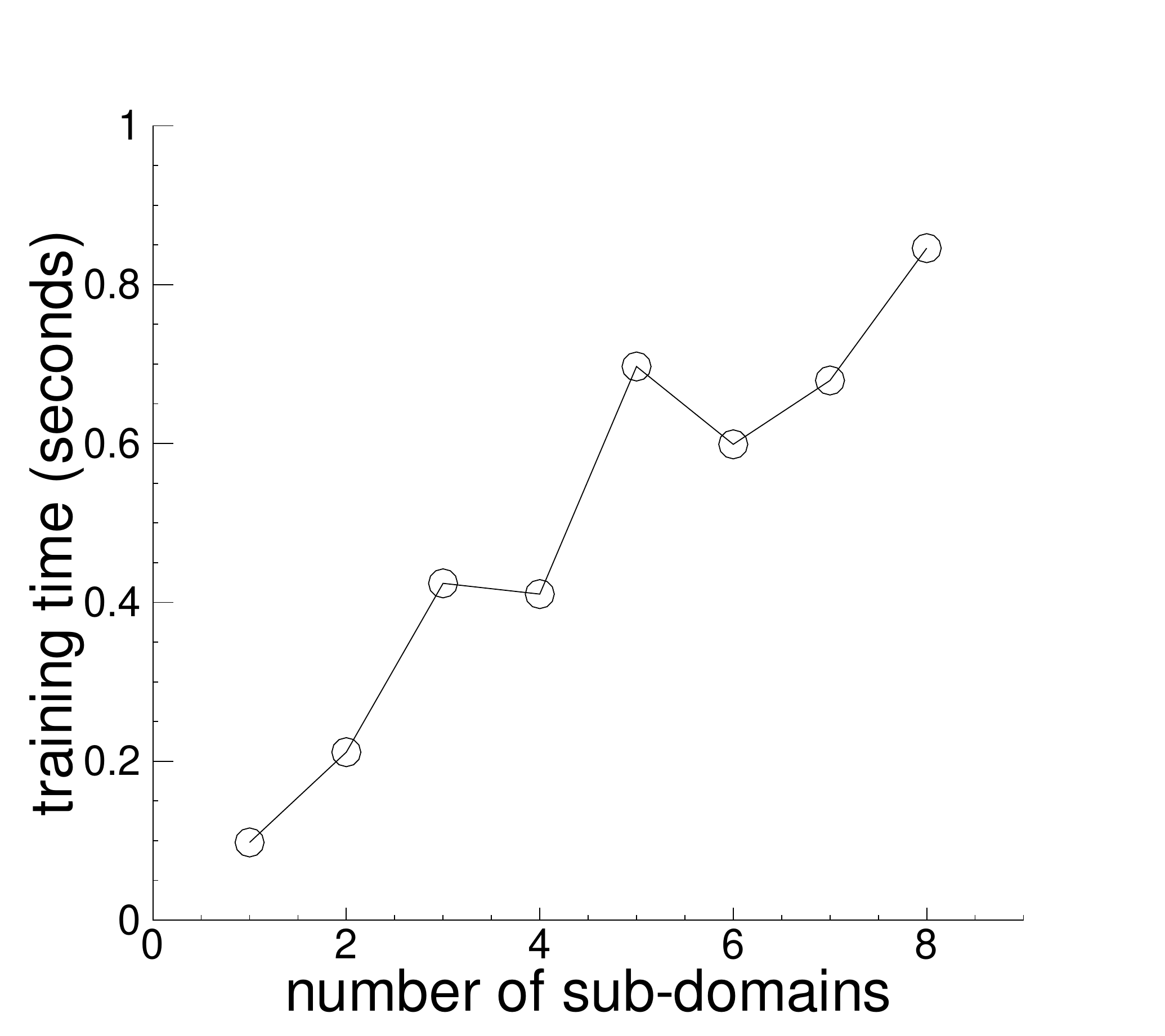}(d)
  }
  \caption{
    Effect of the number of sub-domains, with fixed degrees of freedom
    per sub-domain (1D Helmholtz equation): Profiles of
    (a) the locELM solutions
    and (b) their absolute errors, computed using one sub-domain and four sub-domains.
    (c) The maximum and rms errors in the domain, and 
    (d) the neural-network training time,
    as a function of the number of sub-domains.
  }
  \label{fig:helm1d_1}
\end{figure}


We solve this problem using the locELM method
presented in Section \ref{sec:steady}, by restricting the scheme to
one spatial dimension. We partition $[a,b]$ into $N_e$ uniform
sub-domains (sub-intervals), and impose the $C^1$
continuity conditions across the sub-domain boundaries.
Let $Q$ denote the number of
collocation points within each sub-domain, and consider three types
of collocation points: uniform grid points, the Gauss-Lobatto-Legendre
quadrature points, and random points. The majority of tests reported below
are performed with uniform collocation points in each sub-domain.

For the majority of tests in this subsection,
each local neural network consists of
an input layer with one node (representing $x$),
an output layer with one node
(representing the solution $u$), and one hidden layer in between.
We have also considered local neural networks with two or three
hidden layers between the input
and the output layers.
We employ $\tanh$ as the activation function for all the hidden layers.
The output layer contains no bias and no activation function, as
discussed in Section \ref{sec:loc_elm}.
Additionally, an affine mapping  operation
that normalizes the input $x$ data on each sub-domain to the interval $[-1,1]$ is
incorporated into the local neural networks right behind the input layer.
This operation is 
implemented using the ``lambda'' layer in Keras, which contains
no adjustable parameters and we do not count it
toward the number of  hidden layers.
Following Section \ref{sec:method}, let $M$  denote the number of
nodes in the last hidden layer, which is also the number of training parameters
for each sub-domain.
As discussed in Section \ref{sec:loc_elm},
the weight and bias coefficients in the hidden layers are pre-set
to uniform random values generated on the interval $[-R_m,R_m]$ and
are fixed in the computation.

The main simulation parameters with locELM include
the number of sub-domains ($N_e$), the number of collocation points
per sub-domain ($Q$), the number of training parameters per
sub-domain ($M$), the maximum magnitude of the random coefficients ($R_m$),
the number of hidden layers in the local neural network,
and the type of collocation points in
each sub-domain. We will use the total number of
collocation points ($N_eQ$) and the total number of training parameters ($N_eM$)
to characterize the total degrees of freedom in the simulation.
The effects of the above  parameters on
the simulation results will be investigated.
To make the numerical tests repeatable, all the random numbers are
generated by the Tensorflow library, and we employ a fixed seed value $1$ for
the random number generator with all the tests in this sub-section.

Figure \ref{fig:helm1d_1} illustrates the effect of the number of
sub-domains in the locELM simulation, with the degrees
of freedom per sub-domain (i.e.~the number of collocation points
and the number of training parameters per sub-domain) fixed.
Figures \ref{fig:helm1d_1}(a) and (b) show the solution and error profiles
obtained with one sub-domain and $4$ sub-domains in the locELM simulation.
Figure \ref{fig:helm1d_1}(c) shows the maximum ($L^{\infty}$)
and the rms ($L^2$) errors of the locELM solution in the overall domain
as a function of the number of sub-domains.
Figure \ref{fig:helm1d_1}(d) shows the training time of the overall neural network
as a function of the number of sub-domains.
Here the error refers to the absolute value of the difference between
the locELM solution and the exact solution give by equation \eqref{equ_t2}.
As discussed before,
the training time refers to the total computation time of the locELM method,
and includes the time for computing the output of the last hidden layer
$V_j^{s}(x)$ ($1\leqslant s\leqslant N_e$, $1\leqslant j\leqslant M$)
and its derivatives, the coefficient matrix and the right hand side,
and for solving the linear least squares problem.
In this set of tests, we have employed $Q=50$ uniform collocation points per sub-domain
and $M=50$ training parameters per sub-domain.
Each local neural network contains a single hidden layer,
and we have employed $R_m=3.0$ when generating the random weight/bias coefficients
for the hidden layers of the local neural networks.
It can be observed that the locELM method produces dramatically (nearly exponentially)
more accurate
 results with increasing number of sub-domains, with
the maximum error in the domain reduced from around $10^1$ for a single sub-domain
to about $10^{-7}$ for $8$ sub-domains.
The training time for the neural network, on the other hand,
increases approximately linearly with increasing sub-domains, 
with the training time from about $0.1$ seconds for a single sub-domain
to about $0.8$ seconds for 8 sub-domains.

\begin{figure}
  \centerline{
    \includegraphics[width=2.2in]{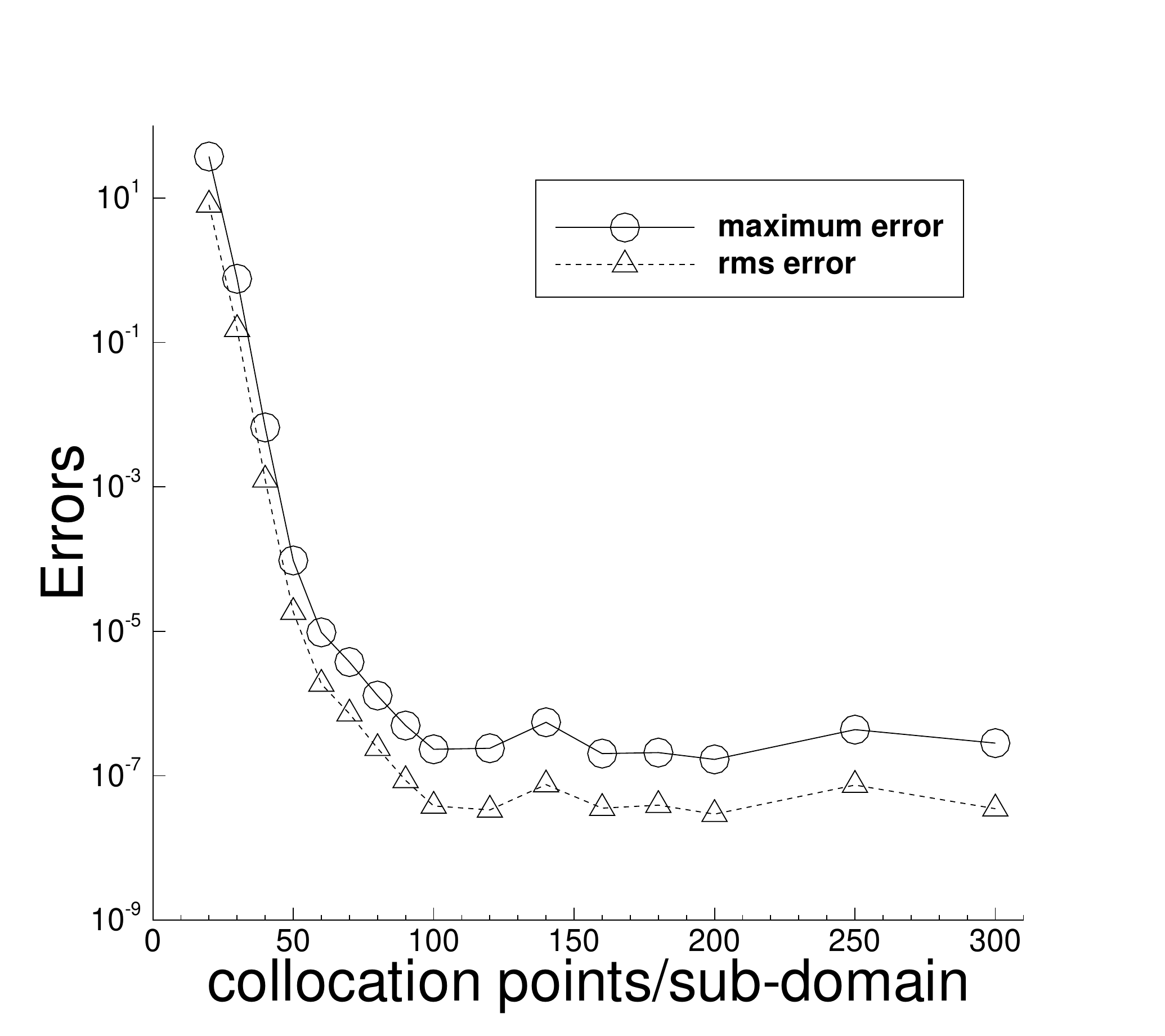}(a)
    \includegraphics[width=2.2in]{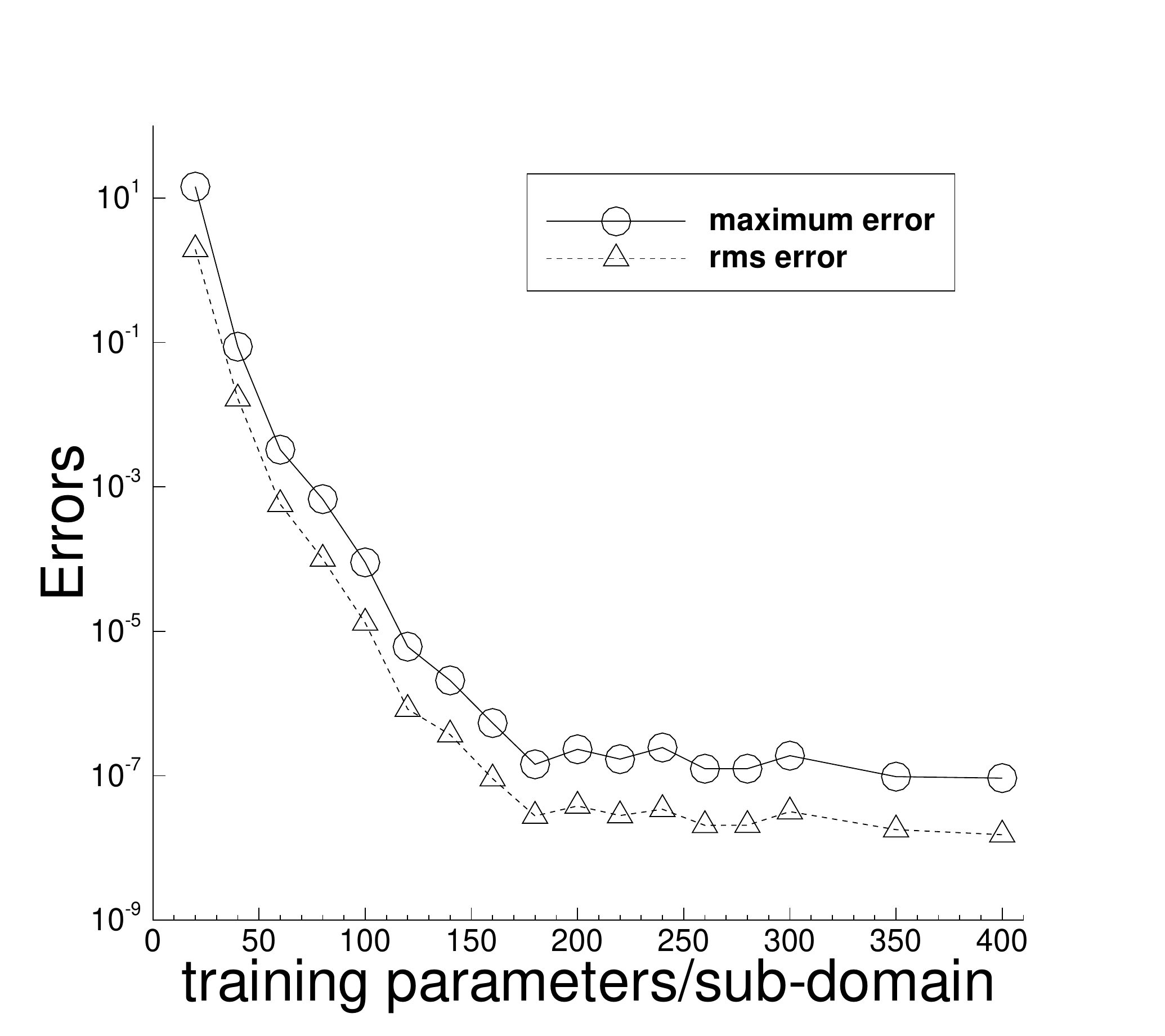}(b)
  }
  \caption{Effect of the number of collocation points and training parameters
    (1D Helmholtz equation):
    the maximum and rms errors as a function of (a) the number of collocation
    points/sub-domain, and (b) the number of training parameters/sub-domain.
    Two uniform sub-domains  are used.
  }
  \label{fig:helm1d_2}
\end{figure}


Figure \ref{fig:helm1d_2} illustrates the effects of the number of
collocation points and the number of training parameters per sub-domain
on the simulation accuracy.
Figure \ref{fig:helm1d_2}(a) depicts the maximum and rms errors in the domain
versus the number of collocation points/sub-domain.
Figure \ref{fig:helm1d_2}(b) depicts the maximum and rms errors in the domain
versus the number of training parameters/sub-domain.
In these tests we have employed $N_e=2$ uniform sub-domains,
uniform collocation points in each sub-domain,
one hidden layer in each local neural network,
and $R_m=3.0$ when generating the random weight/bias coefficients
for the hidden layer.
For the tests in plot (a) the number of training parameters/sub-domain
is fixed at $M=200$, and for the tests in plot (b) the number of
collocation points/sub-domain is fixed at $Q=100$.
Increasing the collocation points per sub-domain
causes an exponential decrease in the numerical errors initially.
The errors then stagnate as the number of collocation points/sub-domain
exceeds a certain point ($Q\sim 100$ in this case).
The error stagnation is due to the fixed number of training
parameters/sub-domain ($M=200$) here.
The number of training parameters/sub-domain appears to have
a similar effect on the errors.
Increasing the training parameters per sub-domain also causes
a nearly exponential decrease in the errors initially.
The errors then stagnate as the number of training parameters increases
beyond a certain point ($M\sim 175$ in this case).

The results in Figures \ref{fig:helm1d_1} and \ref{fig:helm1d_2} show
that the current locELM method exhibits a clear sense of convergence
with respect to the degrees of freedom.
The numerical errors decrease exponentially or nearly exponentially,
as the number of sub-domains, or the number of collocation points
per sub-domain, or the number of training parameters per sub-domain
increases.

\begin{figure}
  \centerline{
    \includegraphics[width=2.2in]{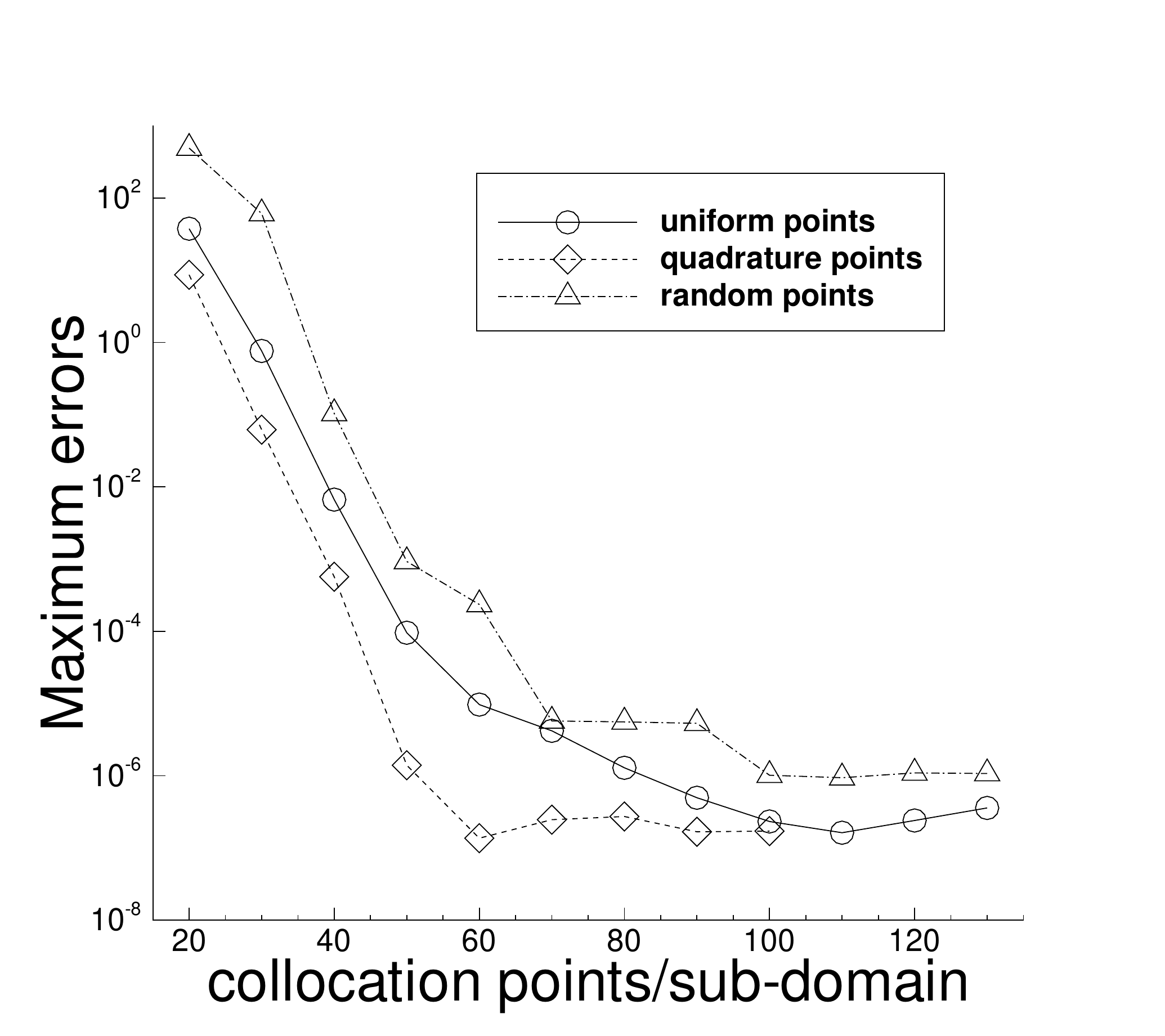}
  }
  \caption{
    Effect of the collocation-point distribution (1D Helmholtz equation):
    the maximum error in the domain versus the number of
    collocation points/sub-domain, obtained with three 
    collocation-point distributions: uniform points, quadrature points,
    and random points.
  }
  \label{fig:helm1d_3}
\end{figure}

Figure \ref{fig:helm1d_3} illustrates the effect of the collocation-point
distribution on the simulation accuracy.
It shows the maximum error in the domain versus
the number of collocation points/sub-domain in
the locELM simulation using three types of
collocation points:
uniform regular points, Gauss-Lobatto-Legendre quadrature points,
and random points (see Remark~\ref{rem_cc}).
In this group of tests we have employed two
sub-domains ($N_e=2$) with $M=200$ training parameters/sub-domain,
and the local neural networks each contains a single hidden layer
with $R_m=3.0$ when generating the random weight/bias coefficients.
With the same number of collocation points, we observe that
the results corresponding to
the random collocation points are the least accurate.
The results obtained with the quadrature points
are the most accurate among the three, whose errors can be orders
of magnitude smaller than those with the random collocation points.
The accuracy corresponding to the uniform regular collocation points
lies between the other two.
With the quadrature points, however, we have encountered practical
difficulties in
our implementation  when the
number of quadrature points becomes larger (above $100$), because the
library our implementation is based on is unable to
compute the quadrature points accurately when the number of quadrature
points exceeds $100$ due to an inherent limitation.
Consequently, we are unable to obtain results with more than $100$ collocation
points/sub-domain when quadrature points are used, which hampers
our ability to perform certain types of tests.
Therefore, the majority of locELM  simulations in
the current work are conducted with uniform
collocation points.


\begin{figure}
  \centerline{
    \includegraphics[width=1.5in]{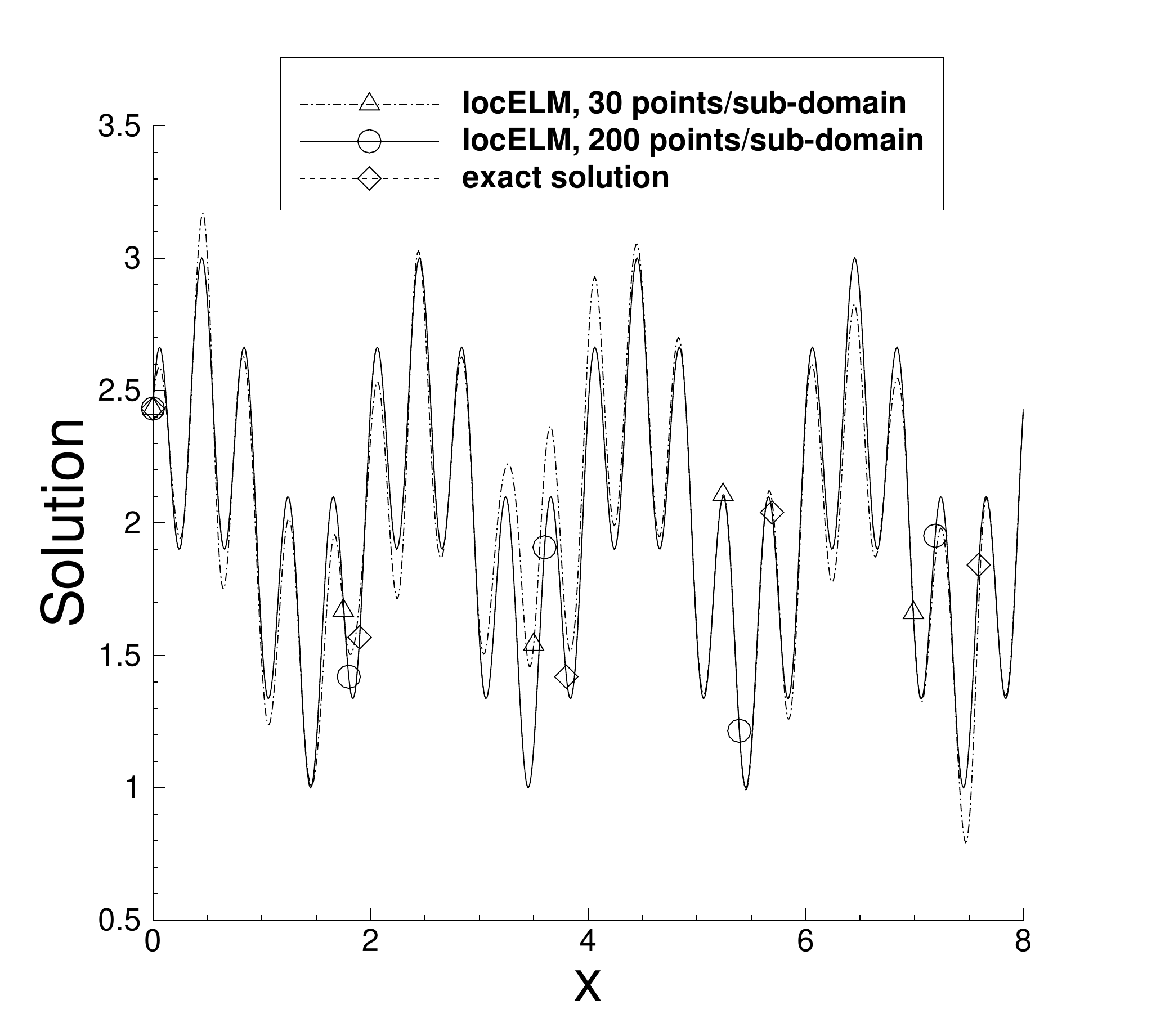}(a)
    \includegraphics[width=1.5in]{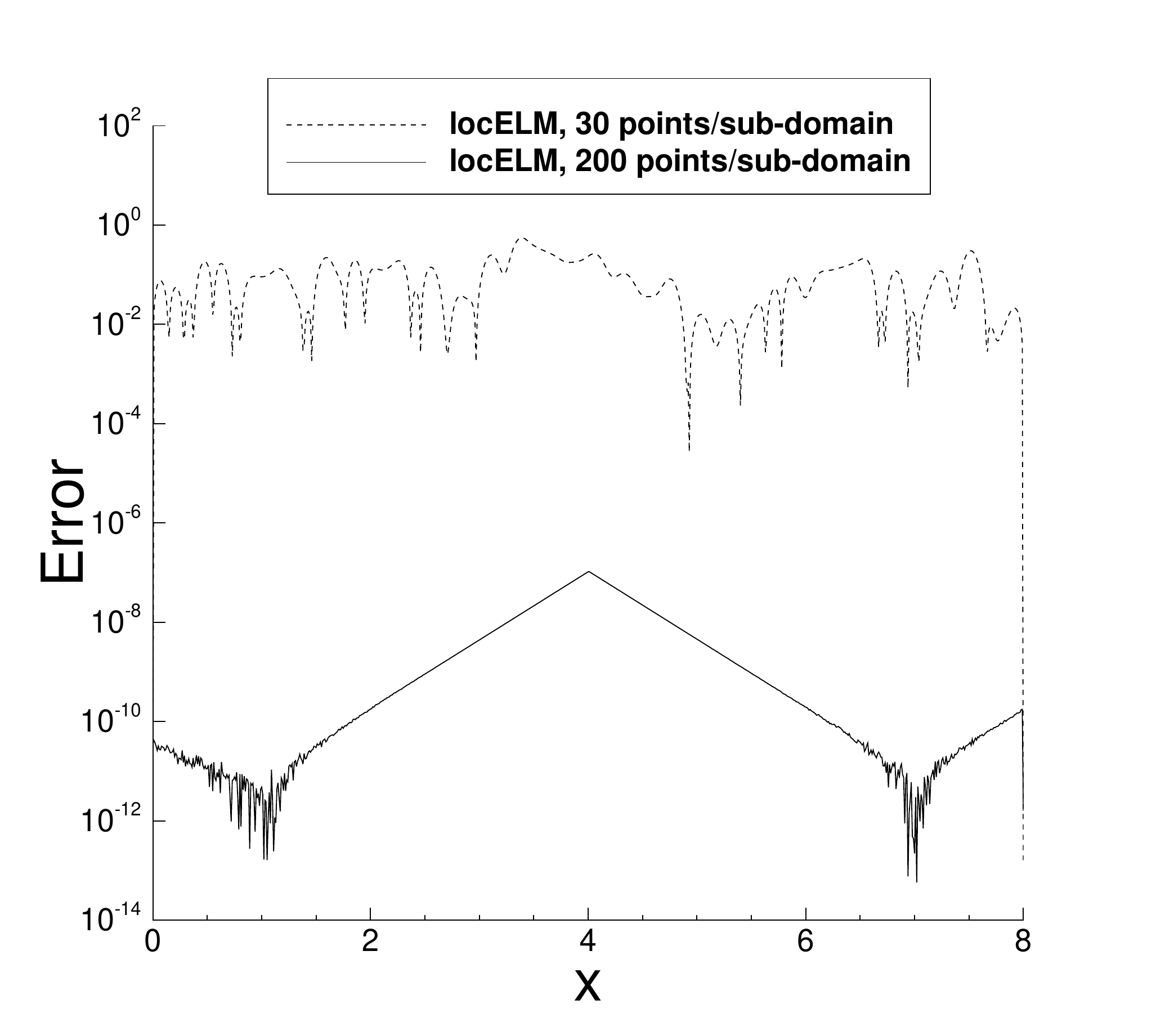}(b)
    \includegraphics[width=1.5in]{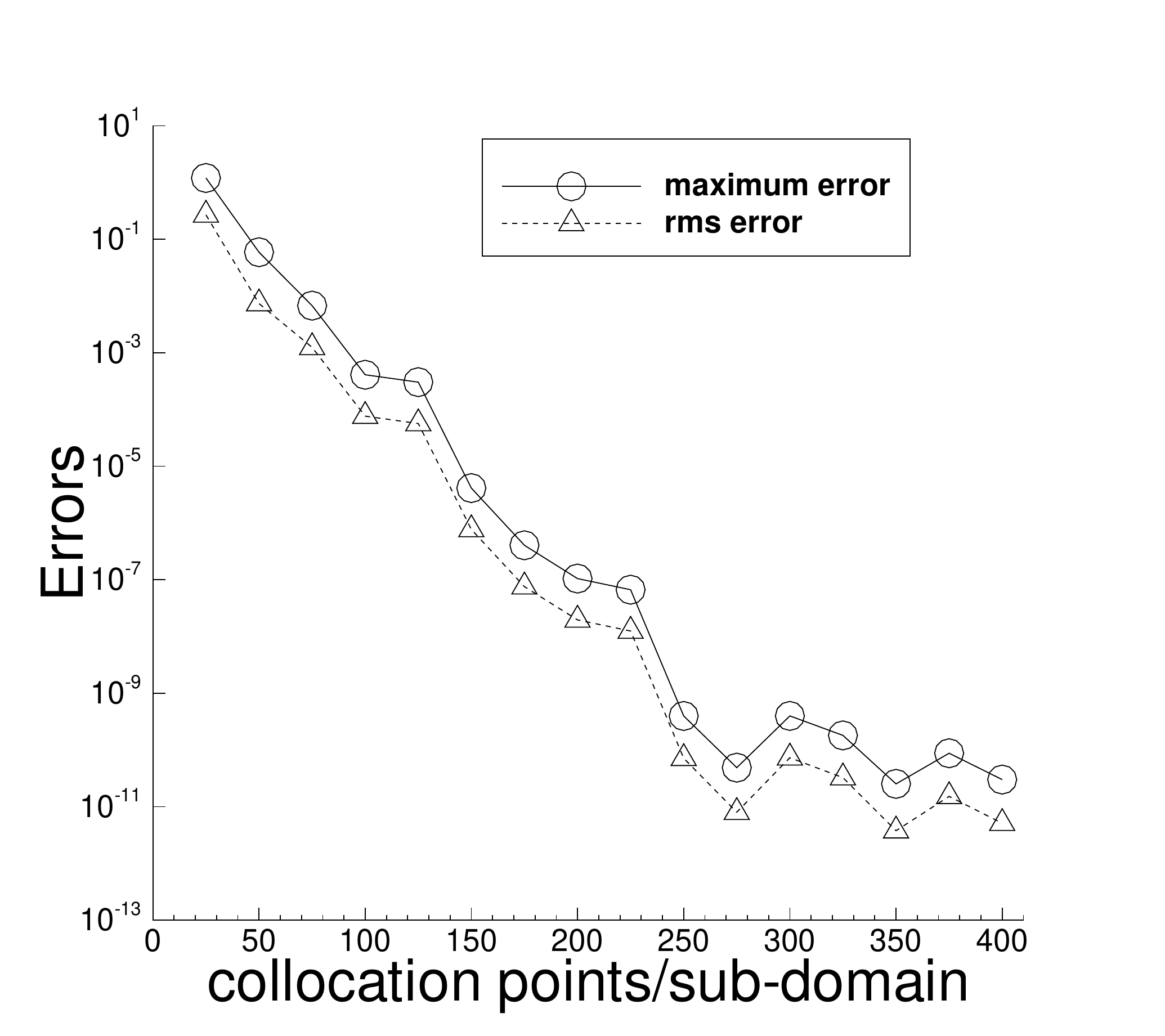}(c)
    \includegraphics[width=1.5in]{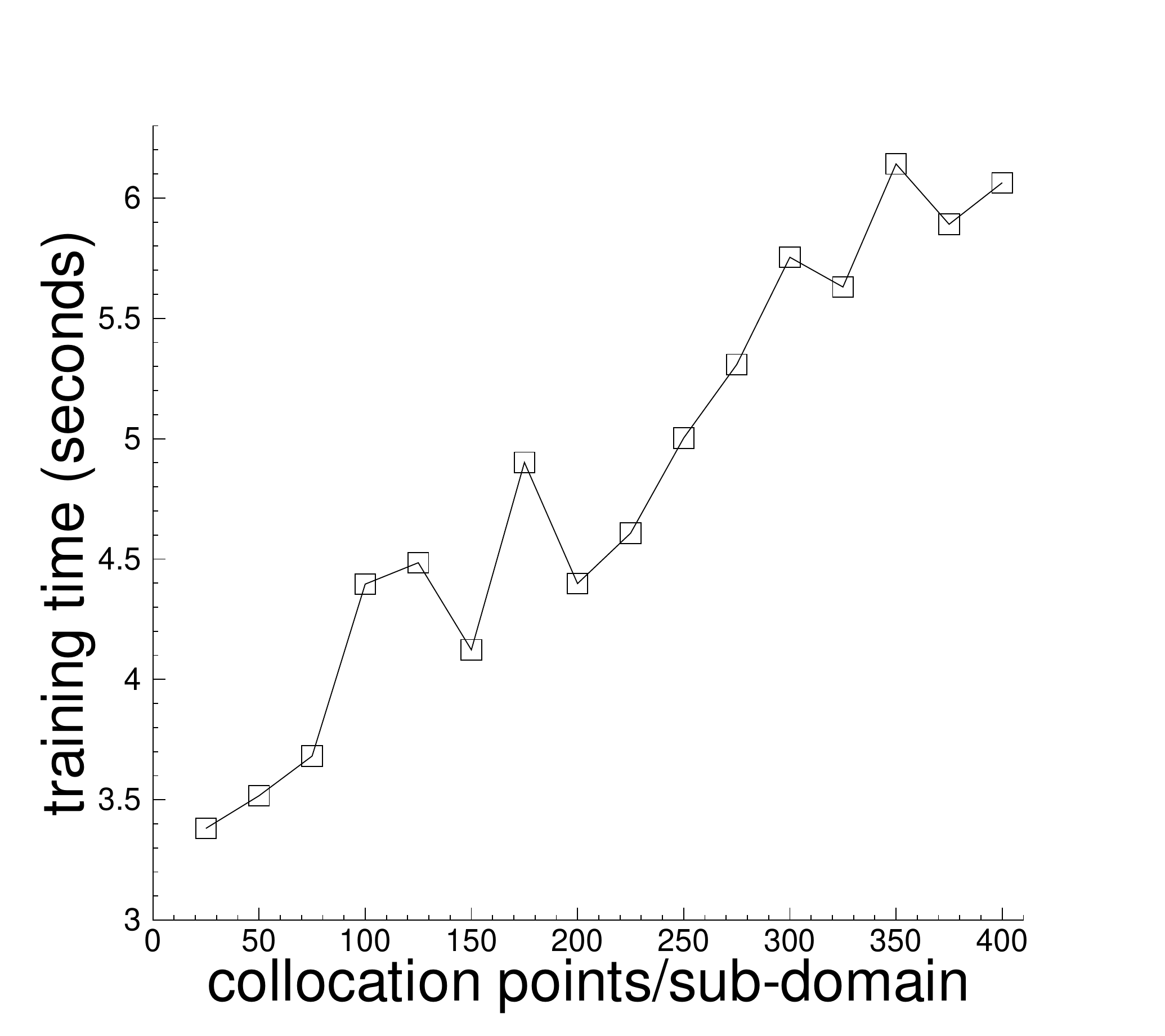}
  }
  \caption{locELM simulations with 2 hidden layers in local neural networks
    (1D Helmholtz equation): profiles of (a) the locELM solutions
    and (b) their absolute errors, computed
    with $30$ and $200$ uniform collocation points per sub-domain.
    (c) the maximum and rms errors in the domain, and (d) the training time,
    as a function of the number of uniform collocation points per sub-domain.
  }
  \label{fig:helm1d_4}
\end{figure}

\begin{figure}
  \centerline{
    \includegraphics[width=1.5in]{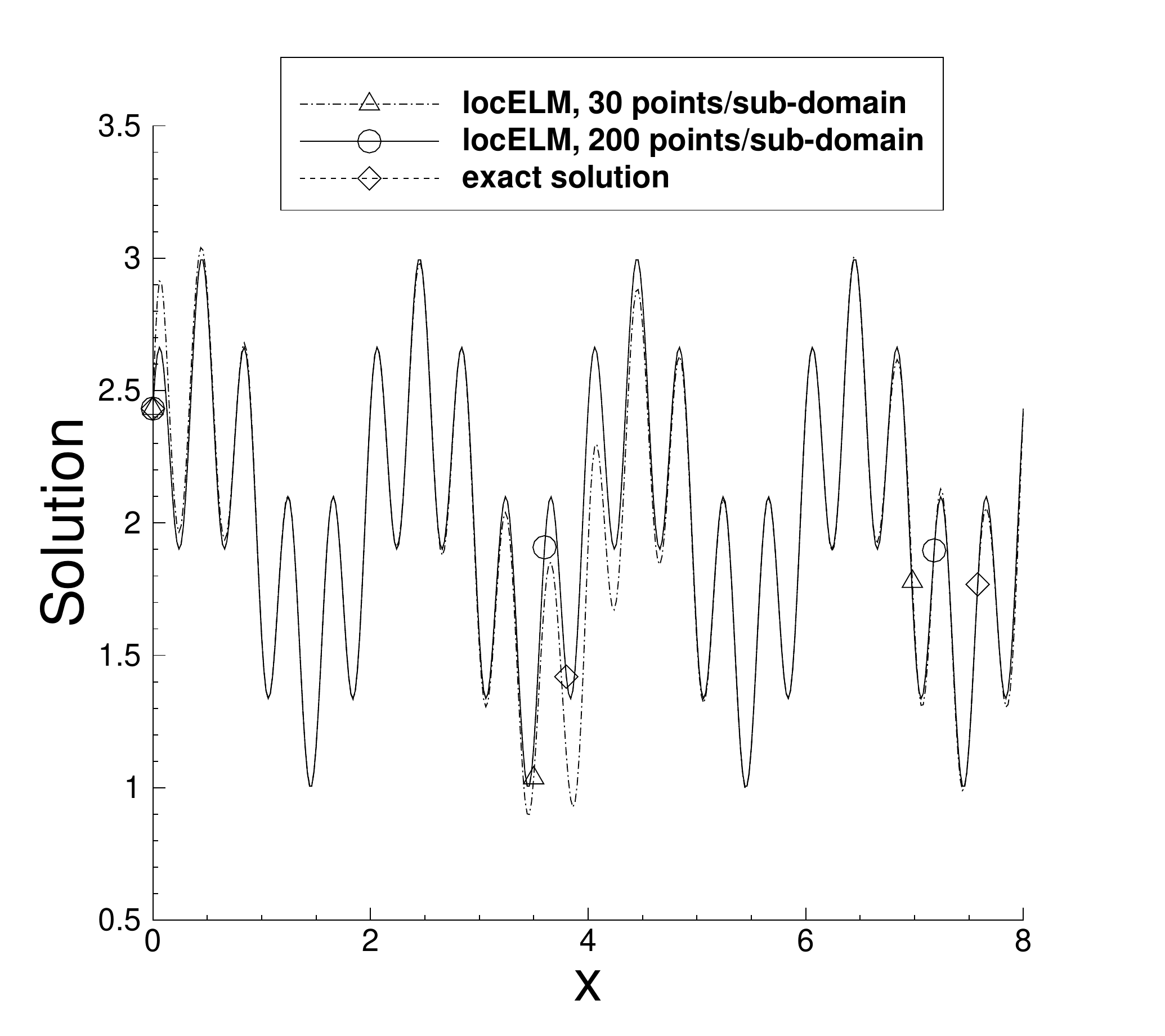}(a)
    \includegraphics[width=1.5in]{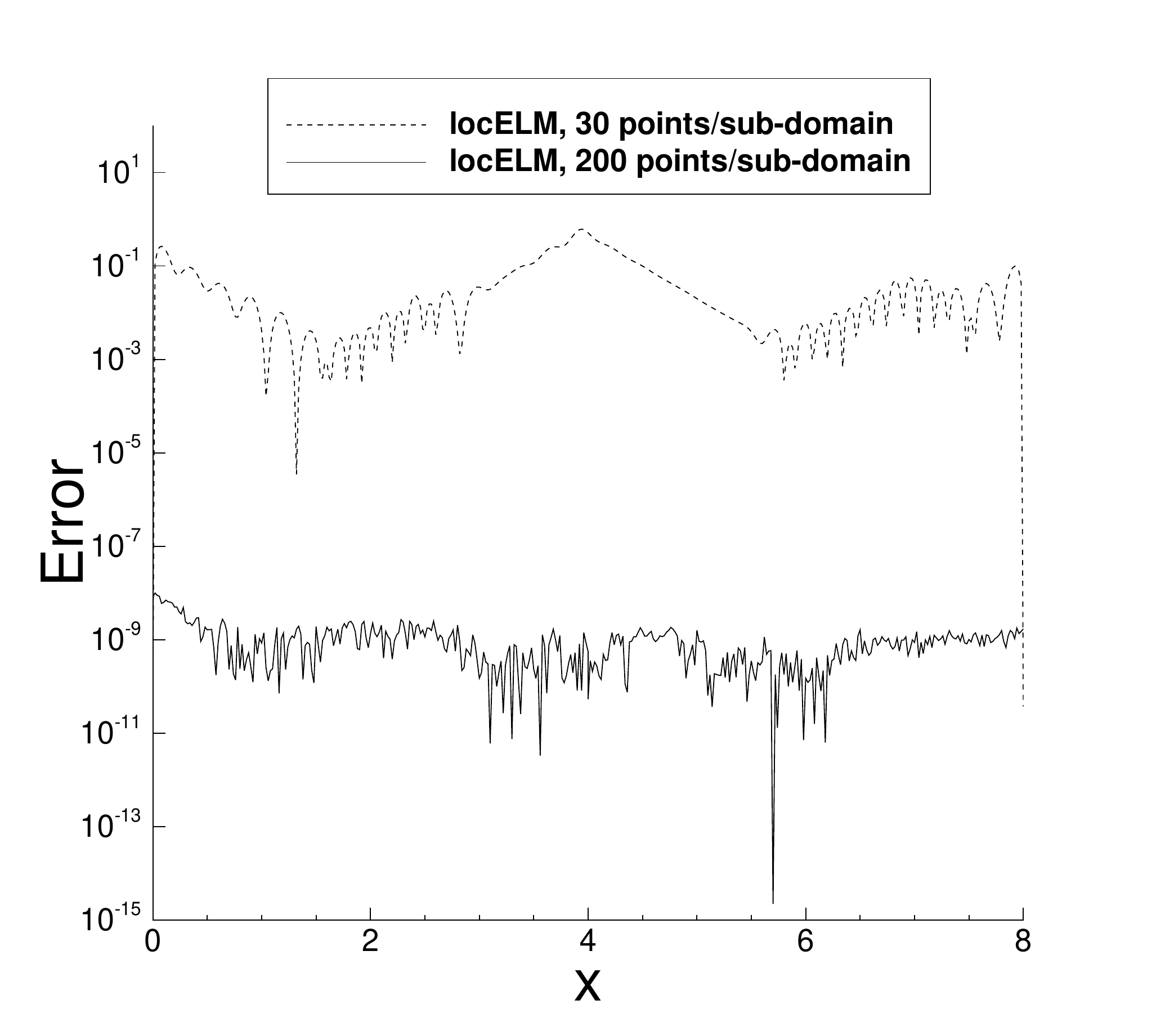}(b)
    \includegraphics[width=1.5in]{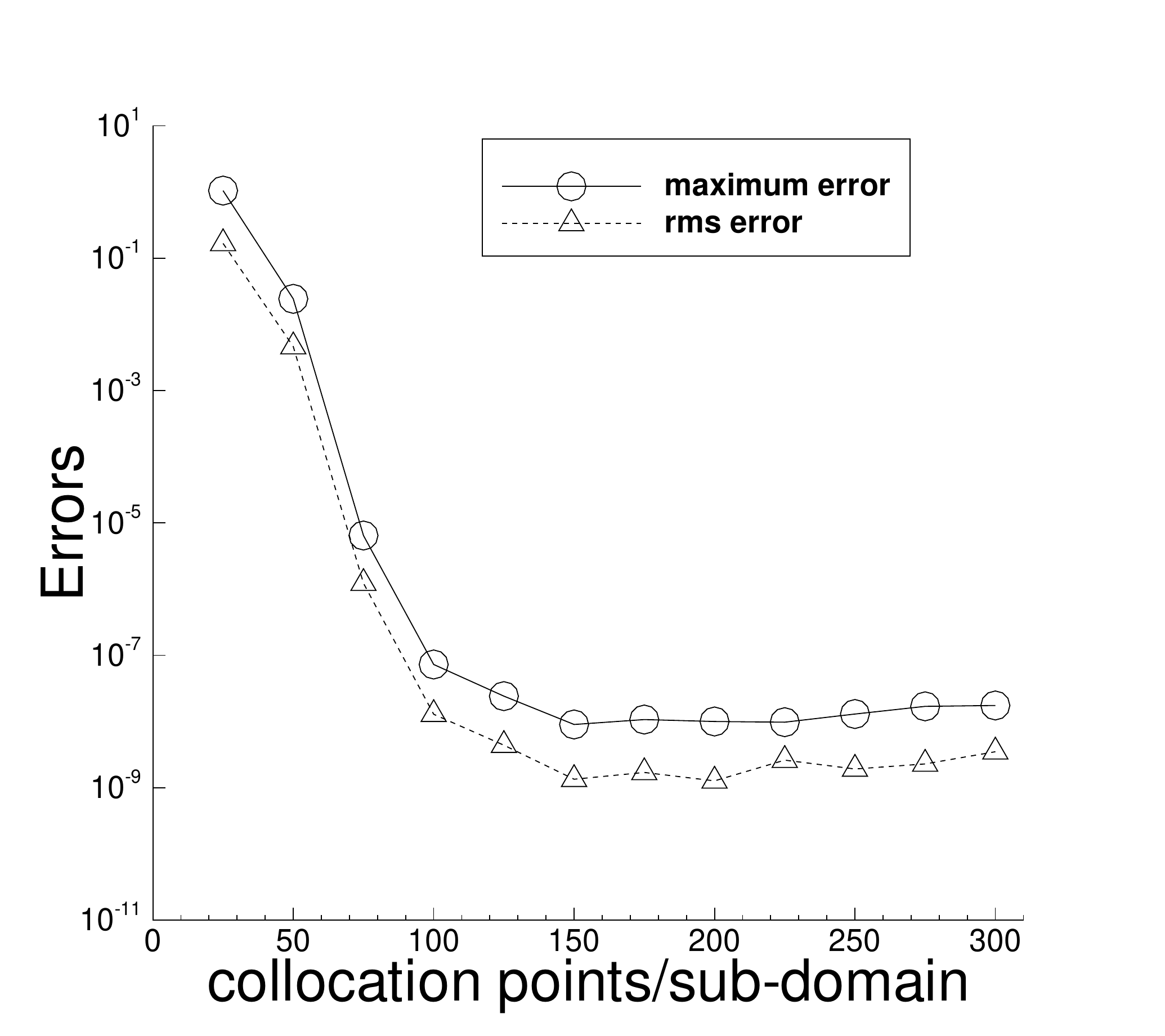}(c)
    \includegraphics[width=1.5in]{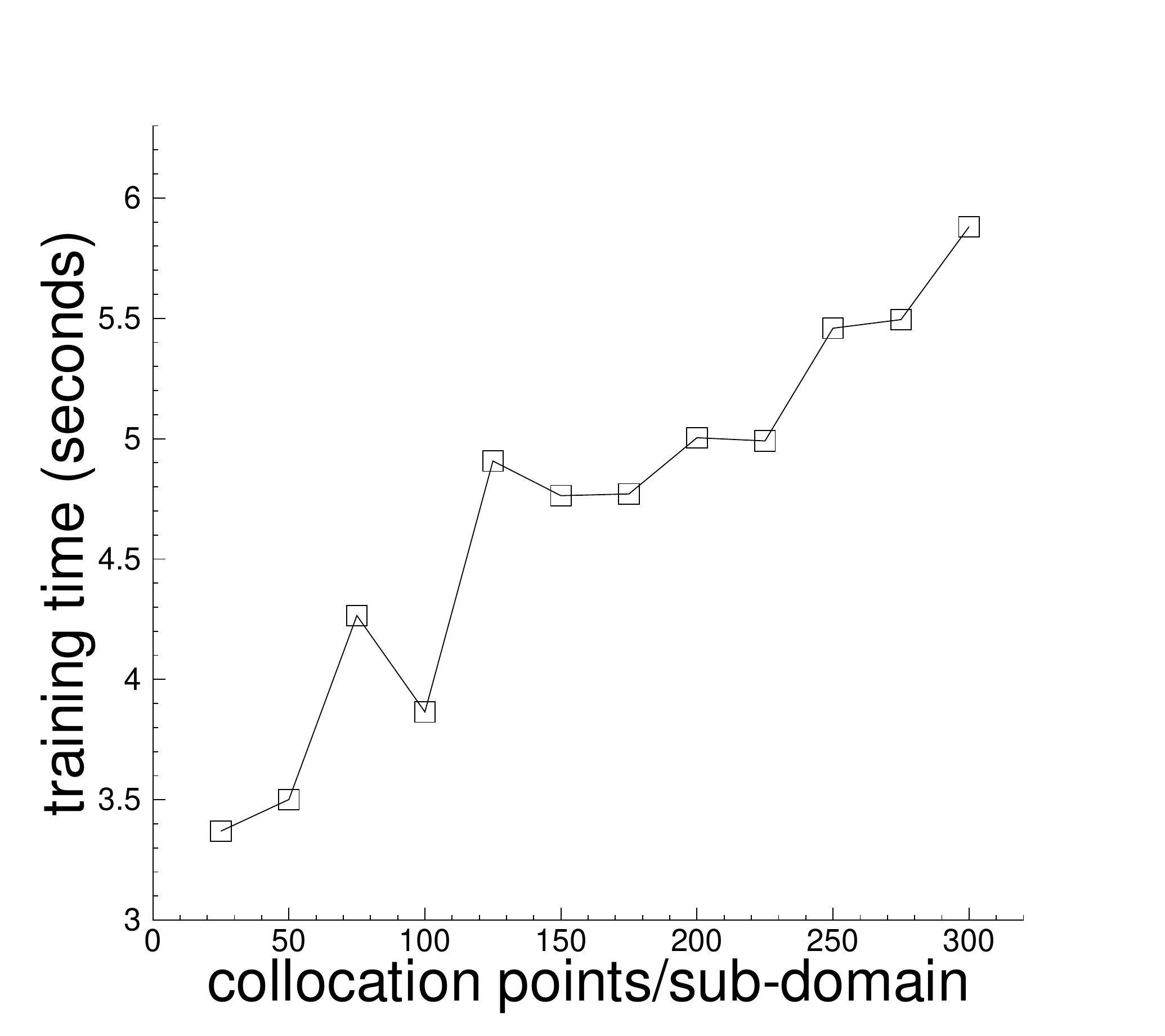}(d)
  }
  \caption{locELM simulations with 3 hidden layers in local neural networks
    (1D Helmholtz equation): profiles of (a) the locELM solutions
    and (b) their absolute errors, computed
    with $30$ and $200$ uniform collocation points per sub-domain.
    (c) the maximum and rms errors in the domain, and (d) the training time,
    as a function of the number of uniform collocation points per sub-domain.
  }
  \label{fig:helm1d_5}
\end{figure}

The test results discussed so far are obtained using a single hidden layer
in the local neural networks.
Traditional studies of global extreme learning machines 
are confined to such a configuration, using  a single hidden
layer in the neural network~\cite{HuangZS2006}.
With the current locELM method,
it is observed that using more than one hidden layer in the local neural networks
one can also obtain accurate results.
This is demonstrated by the results in Figures \ref{fig:helm1d_4}
and \ref{fig:helm1d_5}.
Figure \ref{fig:helm1d_4} shows locELM simulation results obtained with
2 hidden layers in each of the local neural networks, and
Figure \ref{fig:helm1d_5} shows locELM results obtained with 3 hidden layers
in the local neural networks.
In these tests two uniform sub-domains ($N_e=2$) have been used.
The local neural networks corresponding to Figure \ref{fig:helm1d_4} each
contains 2 hidden layers with  $20$ and $300$ nodes, respectively,
and $R_m=3.0$ is employed when the random weight/bias coefficients
for the hidden layers are generated.
The local neural networks corresponding to Figure \ref{fig:helm1d_5} each
contains 3 hidden layers with $20$, $20$ and $300$ nodes, respectively,
and $R_m=1.0$ is employed when the random weight/bias coefficients are generated
for the hidden layers.
The number of training parameters per sub-domain in these tests
is therefore fixed at $M=300$, which corresponds to the number of nodes
in the last hidden layer.
We have used $\tanh$ as the activation function for all the hidden layers.
Uniform collocation points have been used in each sub-domain, and
the number of collocation points is varied in the tests.
In each of these two figures, the plots (a) and (b) are profiles of
the locELM solutions and their absolute errors
computed with $30$ and $200$ uniform collocation
points per sub-domain, respectively.
The plots (c) and (d) show the maximum/rms errors in the domain and the training time
as a function of the number of collocation points per sub-domain, respectively.
It is evident that the numerical errors decrease exponentially
with increasing collocation points/sub-domain, similar to what has been observed
with a single hidden layer from Figure \ref{fig:helm1d_2}(a),
until the errors saturate as the number of collocation points increases beyond
a certain point.
With more than one hidden layer, the locELM method can similarly produce
accurate results with a sufficient number of collocation points per sub-domain. 
The training time is also observed to increase essentially linearly
with respect to the number of collocation points per sub-domain.
Numerical experiments with even more hidden layers in the local neural networks
suggest that the simulation tends to be not as accurate
as those corresponding to one, two or three hidden layers.
It appears to be harder to obtain accurate or
more accurate results with even more hidden layers.

\begin{figure}
  \centerline{
    \includegraphics[width=2.2in]{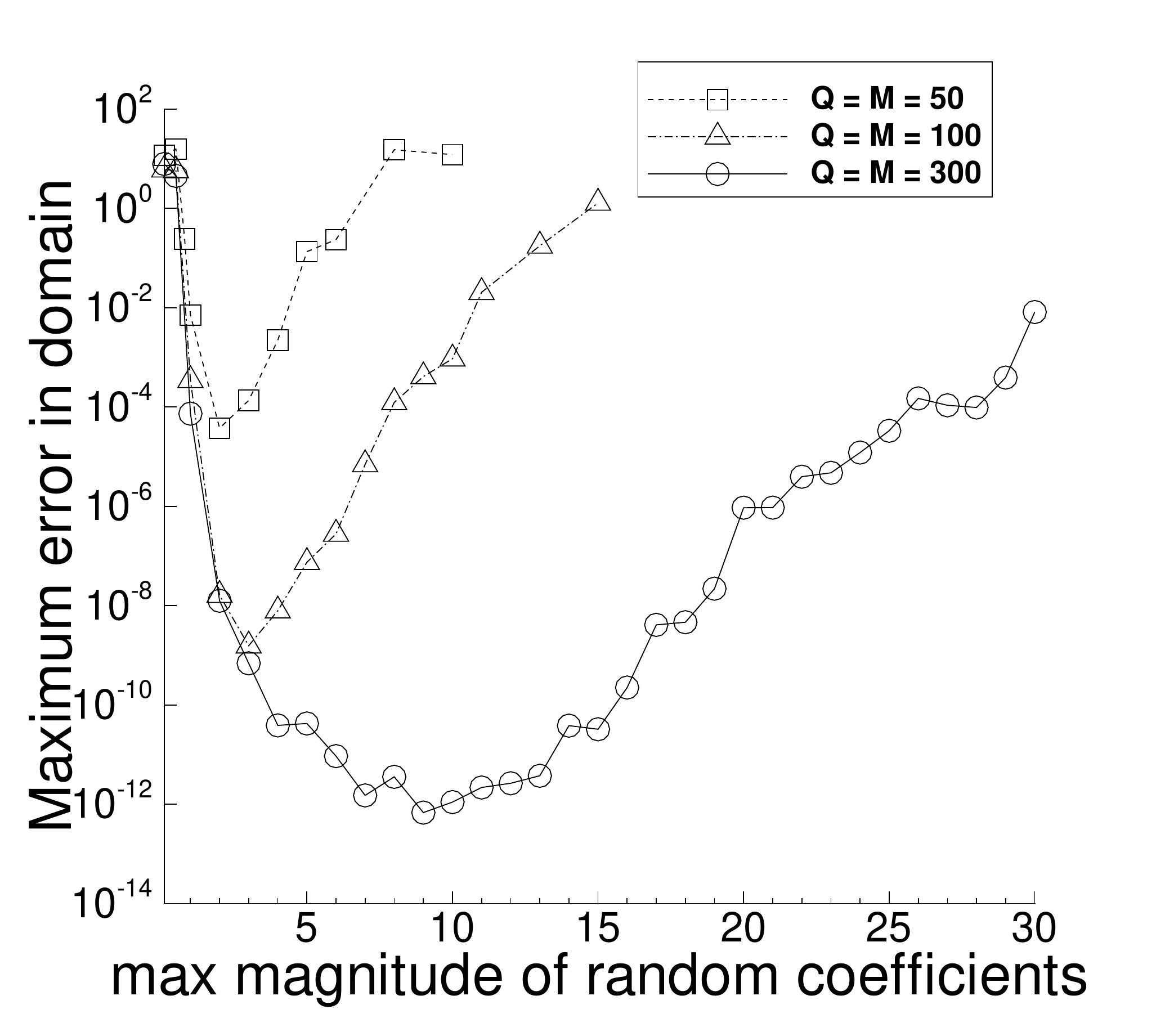}(a)
    \includegraphics[width=2.2in]{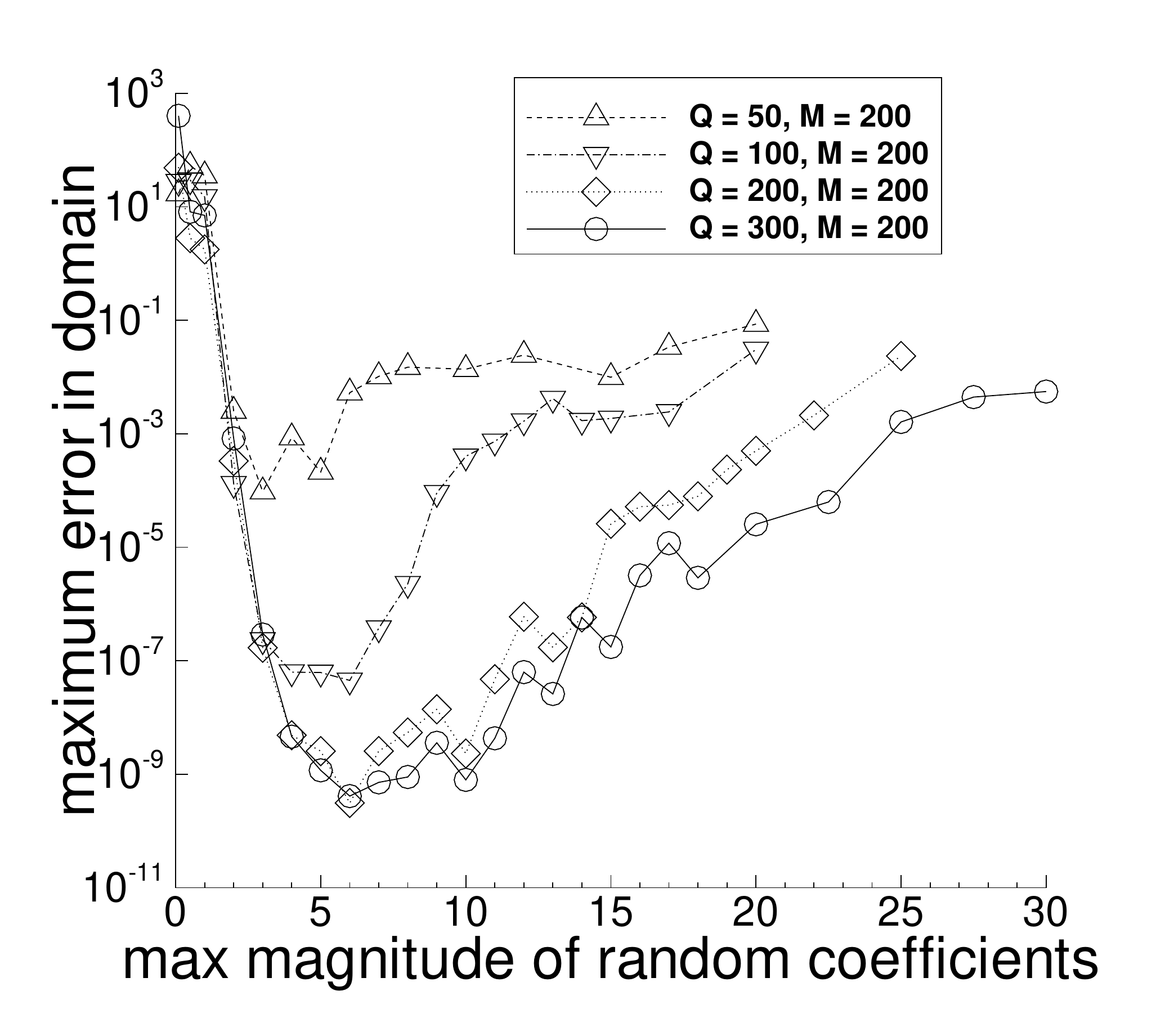}(b)
  }
  \caption{Effect of random weight/bias coefficients in hidden layers
    (1D Helmholtz equation): (a) The maximum error in the domain versus $R_m$,
    for several cases with the number of collocation points/sub-domain ($Q$)
    and the number of training parameters/sub-domain ($M$) kept identical.
    (b) The maximum error in the domain versus $R_m$,
    for several cases with the number of training parameters/sub-domain fixed
    and the number of collocation points/sub-domain varied.
    Four uniform sub-domains are used in (a), and two uniform sub-domains
    are used in (b).
  }
  \label{fig:helm1d_6}
\end{figure}

Apart from the number of collocation points and the number of training parameters
in each sub-domain,
we observe that the random weight/bias coefficients in the hidden layers
can influence the accuracy of the locELM simulation results.
As discussed in Section \ref{sec:loc_elm}, the weight/bias coefficients
in the hidden layers of the local neural networks are pre-set to uniform random values
generated on the interval $[-R_m,R_m]$, and they are
fixed throughout the computation. It is observed that $R_m$,
the maximum magnitude of the random coefficients, can influence significantly
the simulation accuracy.
Figure \ref{fig:helm1d_6} demonstrates this effect with two groups of tests.
In the first group, four uniform sub-domains ($N_e=4$) are used.
The number of (uniform) collocation points per sub-domain (Q)
and the number of training parameters per sub-domain ($M$) are kept to be the same,
and several of these values have been considered ($Q=M=50$, $100$, $300$).
Then for each of these
cases we vary $R_m$ systematically and record the errors of the simulation
results. Figure \ref{fig:helm1d_6}(a) shows the maximum error in the domain
as a function $R_m$ for this group of tests.
In the second group of tests, two uniform sub-domains ($N_e=2$) are used.
The number of training parameters per sub-domain is fixed
at $M=200$, and several values for the number of (uniform) collocation
points are considered ($Q=50$, $100$, $200$, $300$).
For each of these cases, $R_m$ is varied systematically and the corresponding
errors of the simulation results are recorded.
Figure \ref{fig:helm1d_6}(b) shows the maximum error in the domain as a function
of $R_m$ for this group of tests.
In both groups of tests, the local neural networks each contains a single
hidden layer.
These results indicate that, for a fixed simulation resolution (i.e.~fixed
$Q$ and $M$),
the error  tends be worse as $R_m$ becomes very large
or very small. The simulation tends to produce more accurate results for a range
of moderate $R_m$ values, which is typically around $R_m\approx 1 \sim 10$. 
As the simulation resolution increases, the optimal range of $R_m$
values tends to expand and shift rightward (toward larger values) on the $R_m$ axis.
Further tests also suggest that with increasing number of sub-domains
the optimal range of $R_m$ values tends to shift leftward (toward smaller values)
along the $R_m$ axis.

\begin{figure}
  \centerline{
    \includegraphics[width=1.5in]{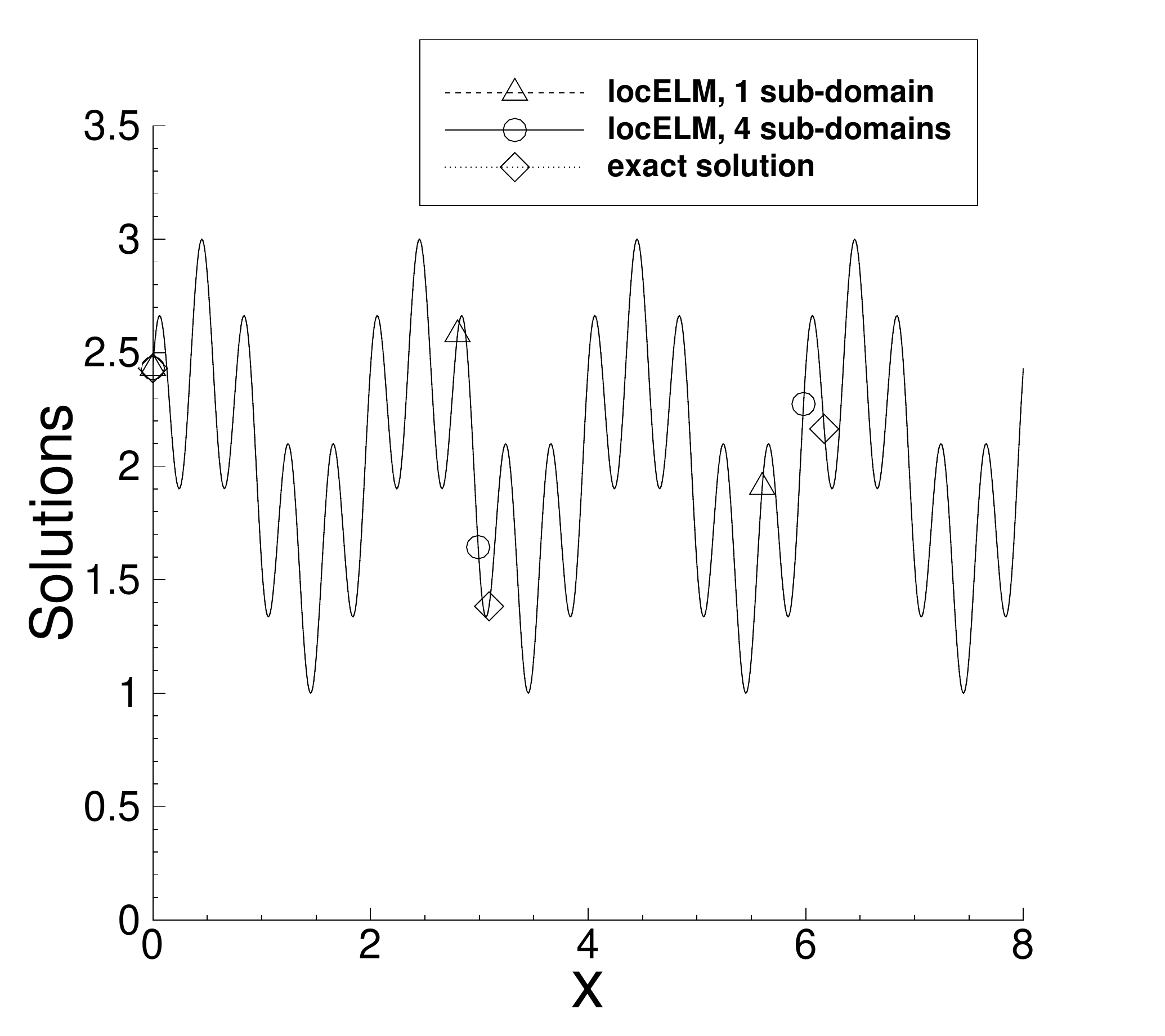}(a)
    \includegraphics[width=1.5in]{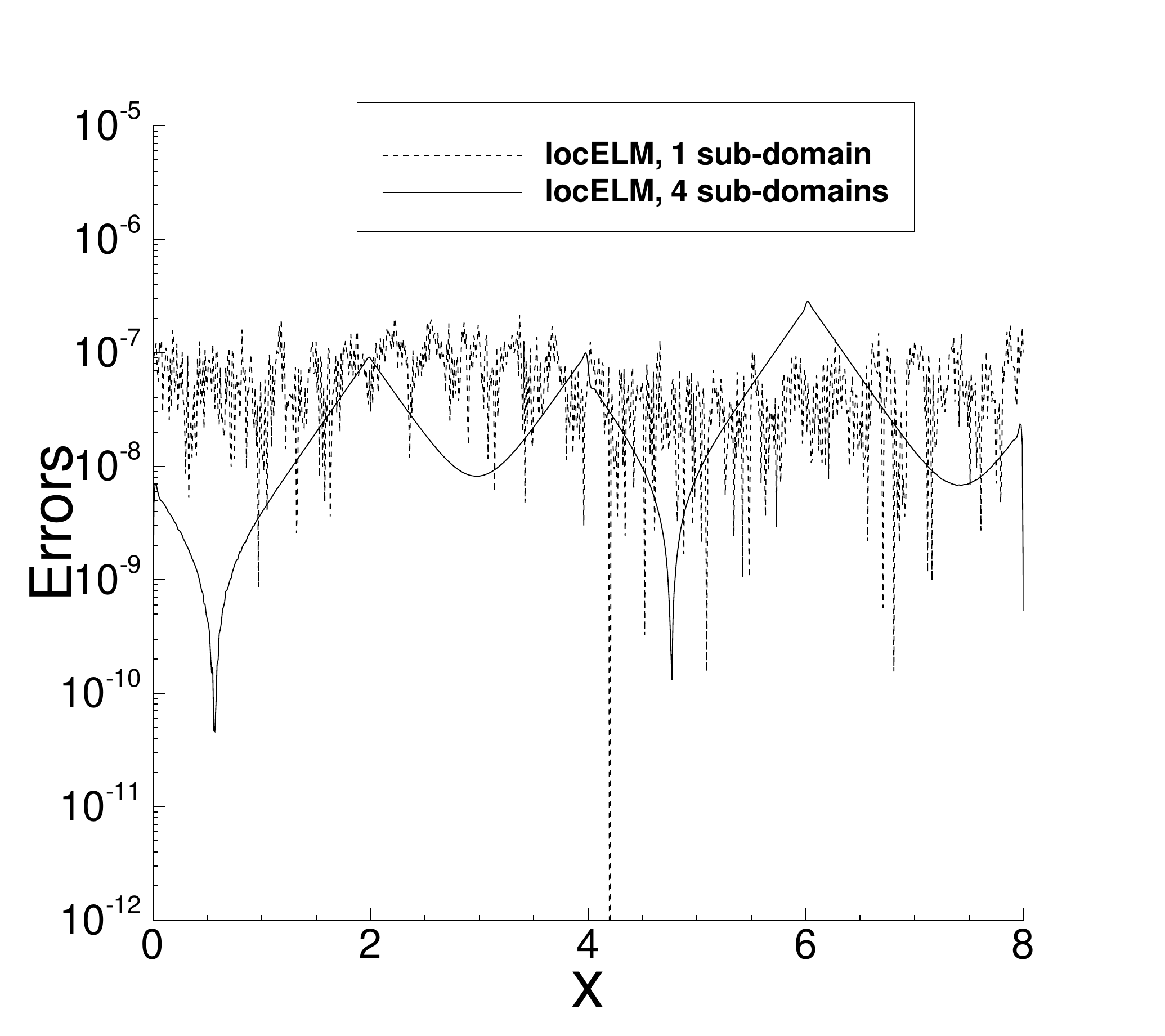}(b)
    \includegraphics[width=1.5in]{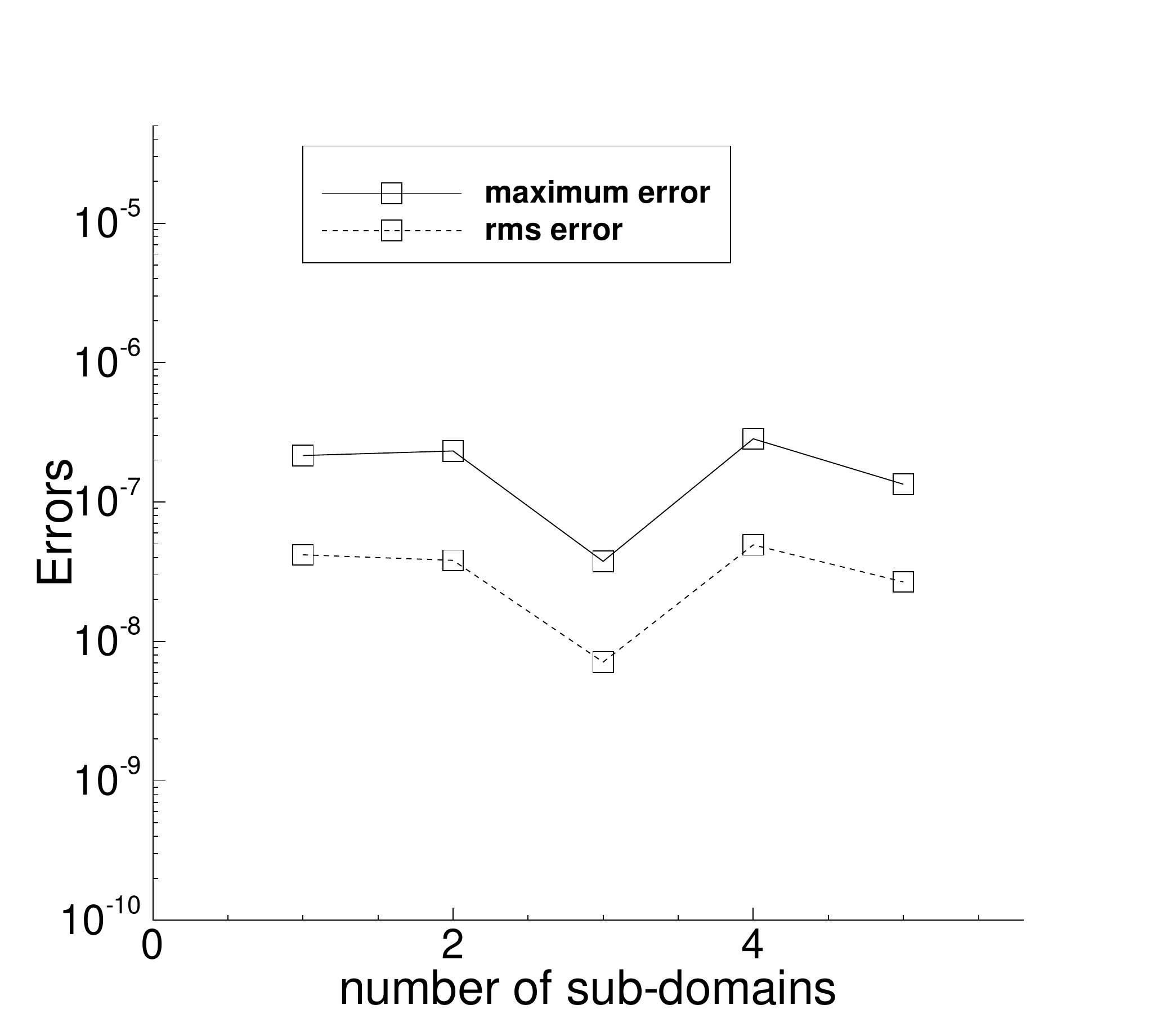}(c)
    \includegraphics[width=1.5in]{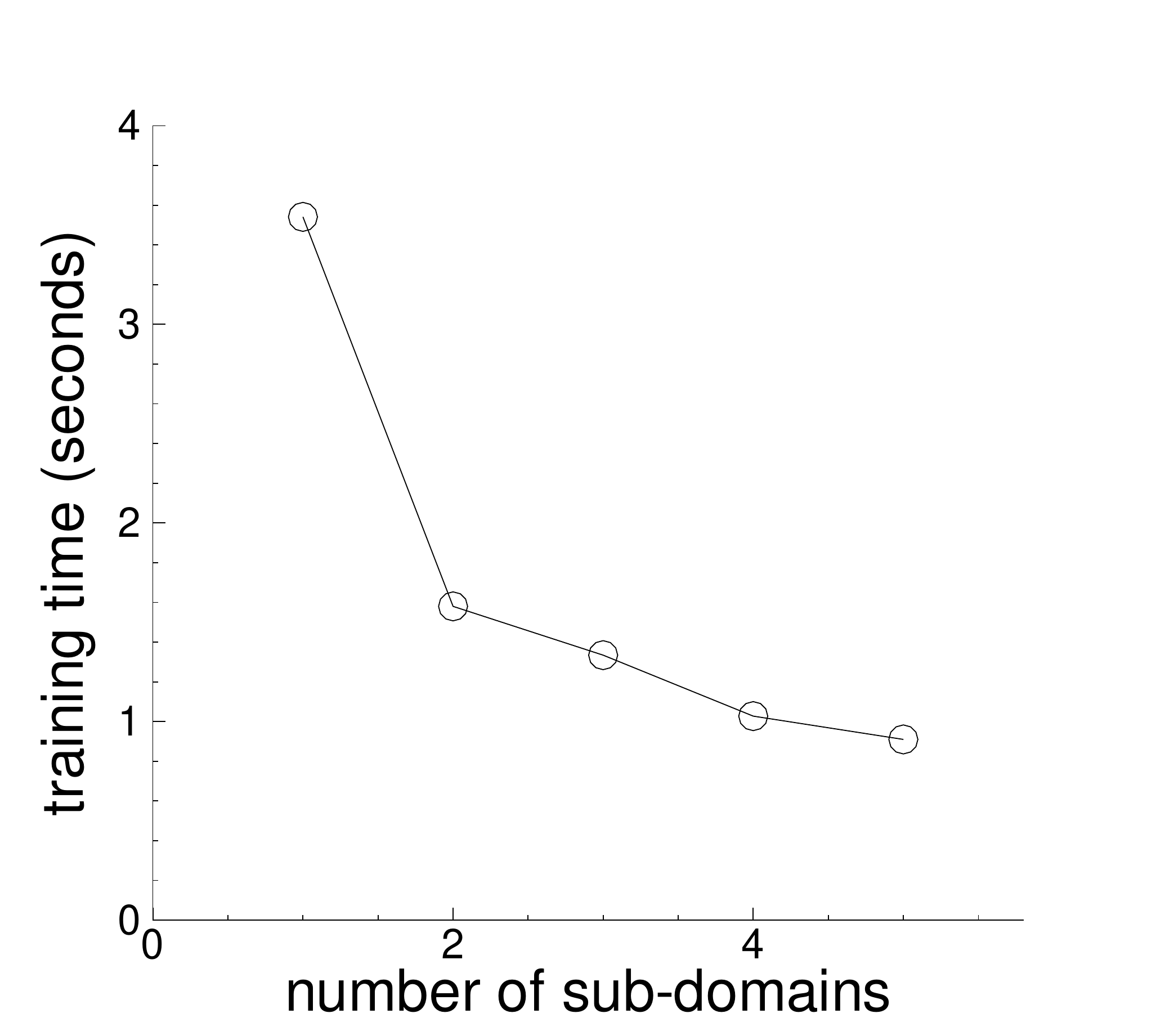}(c)
  }
  \caption{Effect of the number of sub-domains,
    with fixed total degrees of freedom in the domain
    (1D Helmholtz equation):
    profiles of (a) the locELM solutions and (b) their
    absolute errors, computed using
    one and four uniform sub-domains in the simulation.
    (c) The maximum and rms errors
    in the domain, and (d) the training time,
    as a function of the number of uniform sub-domains.
  }
  \label{fig:helm1d_7}
\end{figure}


We observe that
the use of multiple sub-domains and local extreme learning machines
can significantly accelerate the computation and reduce the network training time,
without seriously compromising the accuracy, when compared with 
global extreme learning machines.
This point is demonstrated by Figure \ref{fig:helm1d_7}.
Here we fix the total degrees of freedom in the domain,
i.e.~the total number of collocation points and the total number of training
parameters in the domain, and vary the number of sub-domains in
the locELM simulation. The locELM case with a single sub-domain is equivalent to
a global ELM.
The total number of collocation points in the domain is fixed at $N_eQ=200$, and
the total number of training parameters is fixed at $N_eM=400$.
Uniform sub-domains are employed in these tests, with uniform collocation points
in each sub-domain.
So with multiple
sub-domains the total degrees of freedom
are evenly distributed to different sub-domains and local
neural networks. The local neural networks each contains a single
hidden layer, and the maximum magnitudes of
the random coefficients ($R_m$) employed in the tests here
are approximately in their optimal range of values.
Figures \ref{fig:helm1d_7}(a) and (b) illustrates profiles of the localELM
solutions and their absolute errors obtained using a single sub-domain
($Q=200$, $M=400$, $R_m=6.0$)
and using four sub-domains ($Q=50$, $M=100$, $R_m=3.0$) in the locELM simulations.
Both simulations have produced accurate results, with comparable error levels.
Figure \ref{fig:helm1d_7}(c) shows the maximum and rms errors in the domain
versus the number of sub-domains in the simulations,
and Figure \ref{fig:helm1d_7}(d) shows the training time as a function
of the number of  sub-domains.
It can be observed that the error levels corresponding to multiple sub-domains
are comparable to, or in certain cases maybe slightly
better or worse than,  those of a single sub-domain.
But the training time of the neural network is dramatically reduced
with multiple sub-domains, when compared with a single sub-domain.
The reduction in the training time is due to the fact that,
with multiple sub-domains, the coefficient matrix in the linear
least squares problem becomes very sparse, because only neighboring sub-domains
are coupled through the $C^k$ continuity conditions while those sub-domains
that are not adjacent to each other are not coupled.
On the other hand, with a single sub-domain, all the degrees of freedom
in the domain are coupled with one another, leading to a dense coefficient
matrix in the linear least squares problem and larger computation time.
These results suggest that,
when compared with global ELM,
the use of domain decomposition and local
neural networks can reduce the coupling
among the degrees of freedom
in different sub-domains without seriously compromising the accuracy, and this can
significantly reduce the computation time for the least squares
problem, and hence the network training time.


\begin{figure}
  \centerline{
    \includegraphics[width=2.2in]{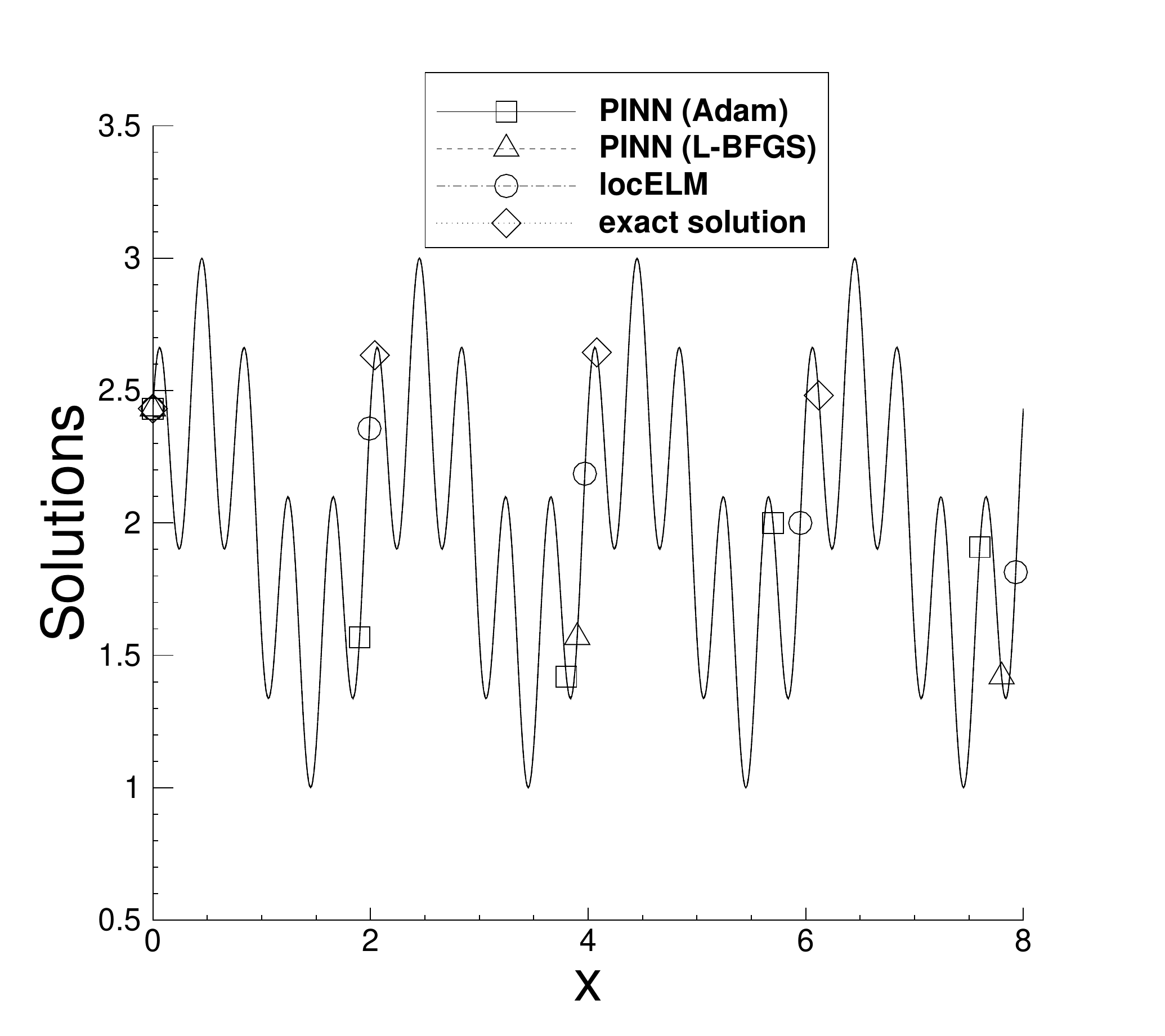}(a)
    \includegraphics[width=2.2in]{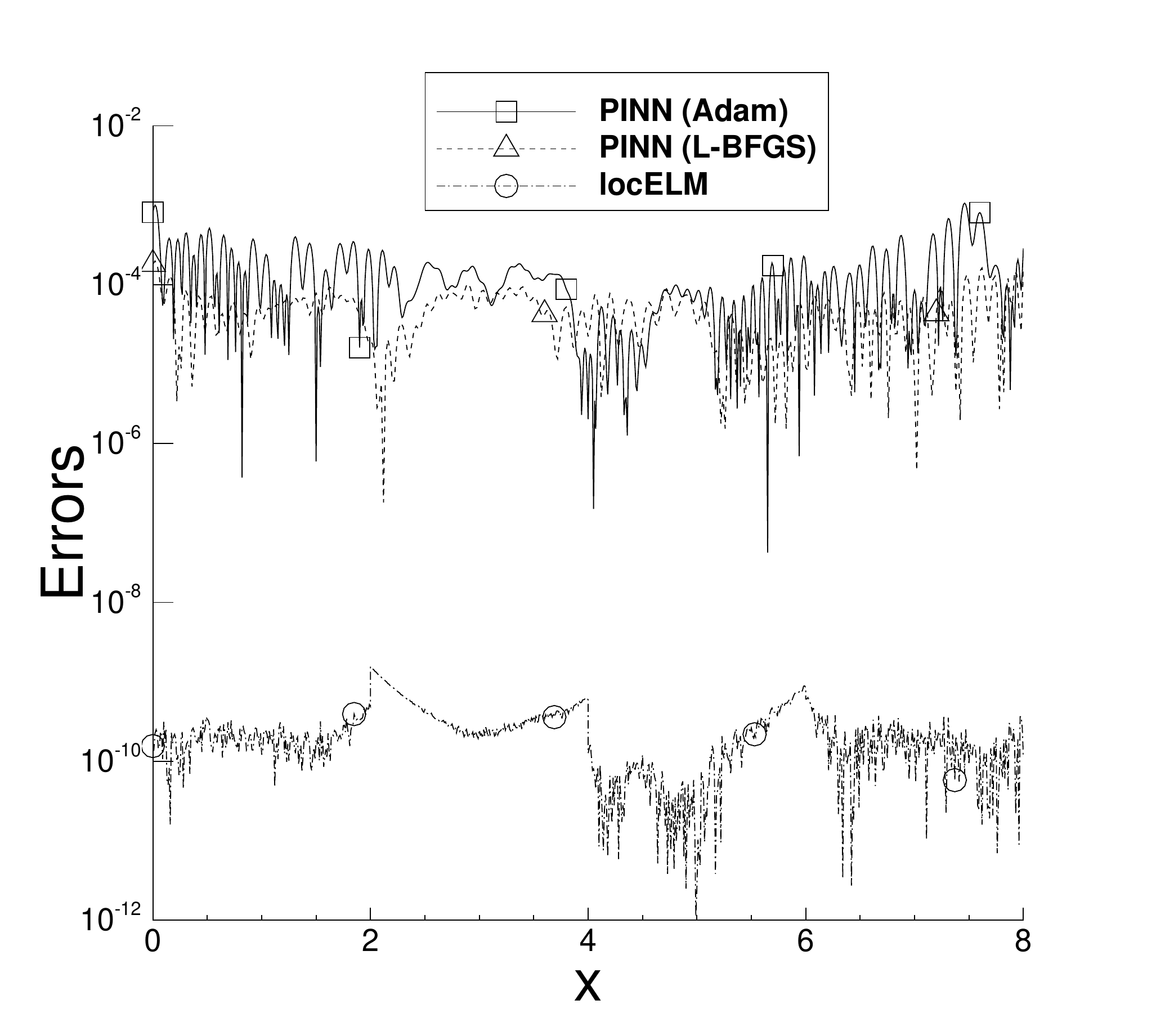}(b)
  }
  \caption{Comparison between locELM and PINN
    (1D Helmholtz equation):
    profiles of (a) the solutions and (b) their absolute errors, obtained using
    PINN~\cite{RaissiPK2019} 
    with the Adam and L-BFGS optimizers,
    and using the current locELM method.
  }
  \label{fig:helm1d_8}
\end{figure}

\begin{table}[tb]
  \centering
  \begin{tabular}{lllll}
    \hline
    method & maximum error & rms error & epochs/iterations & training
    time (seconds)\\
    PINN (Adam) & $1.06e-3$ & $1.57e-4$ & $45,000$ & $507.7$ \\
    PINN (L-BFGS) & $1.98e-4$ & $3.15e-5$ & $22,500$ & $1035.8$ \\
    locELM & $1.56e-9$ & $2.25e-10$ & $0$ & $1.1$ \\
    \hline
  \end{tabular}
  \caption{1D Helmholtz equation: Comparison between the current locELM method and
    PINN,
    in terms of the maximum/rms errors in the domain, the number of
    epochs or iterations in the training of neural networks, and
    the training time.
    The problem settings correspond to those of Figures \ref{fig:helm1d_8}.
  }
  \label{tab:tab_1}
\end{table}

We next compare the the current locELM method with
the physics-informed neural network (PINN)~\cite{RaissiPK2019} method,
an often-used PDE solver based on deep neural networks.
Figure \ref{fig:helm1d_8} compares profiles of the solutions (plot (a))
and their absolute errors (plot (b)) obtained using PINN
with the Adam and the L-BFGS optimizers, and using the current locELM method.
In the PINN simulations, the neural network contains $6$ hidden layers
with $50$ nodes and the $\tanh$
activation function in each layer, and
the output layer contains no activation function.
The input data consist of $300$ uniform
collocation points in the domain.
In the PINN/Adam simulation, the network has been trained on the input data
for $45,000$ epochs, with the learning rate gradually decreasing from $0.001$
at the beginning to $2.5\times 10^{-5}$ at the end.
In the PINN/L-BFGS simulation,
the network has been trained for $22,500$ L-BFGS iterations.
In the locELM simulation, four uniform sub-domains ($N_e=4$) have been used,
with $M=100$ training parameters per sub-domain and $Q=100$ uniform collocation
points per sub-domain.
The four local neural networks each consists of one hidden layer
with $M=100$ nodes and the $\tanh$ activation function,
and we have employed $R_m=3.0$ for generating
random weight/bias coefficients in the hidden layer.
Figure \ref{fig:helm1d_8} shows that both PINN and the current
locELM method have captured the solution quite accurately. But the current
method is considerably more accurate than PINN, by a factor of nearly five
orders of magnitude in terms of the errors. 

Table \ref{tab:tab_1} is a further comparison of PINN and locELM
in terms of the maximum/rms errors in the domain, and the computational
cost (the network training time and the number of epochs or iterations).
The problem setting corresponds to that of Figure \ref{fig:helm1d_8}.
The current method is not only  much more
accurate than PINN, but also considerably cheaper in terms of the
computational cost. The training time with the current locELM method
is on the order of a second.
In contrast, it takes
over $500$ seconds to train PINN with Adam
and over $1000$ seconds to train it with L-BFGS.
We observe a clear superiority of the current locELM
method to the PINN solver in terms of both accuracy and
computational cost.
These observations will be confirmed and reinforced with other problems
in subsequent sections.


\begin{figure}
  \centerline{
    \includegraphics[width=2in]{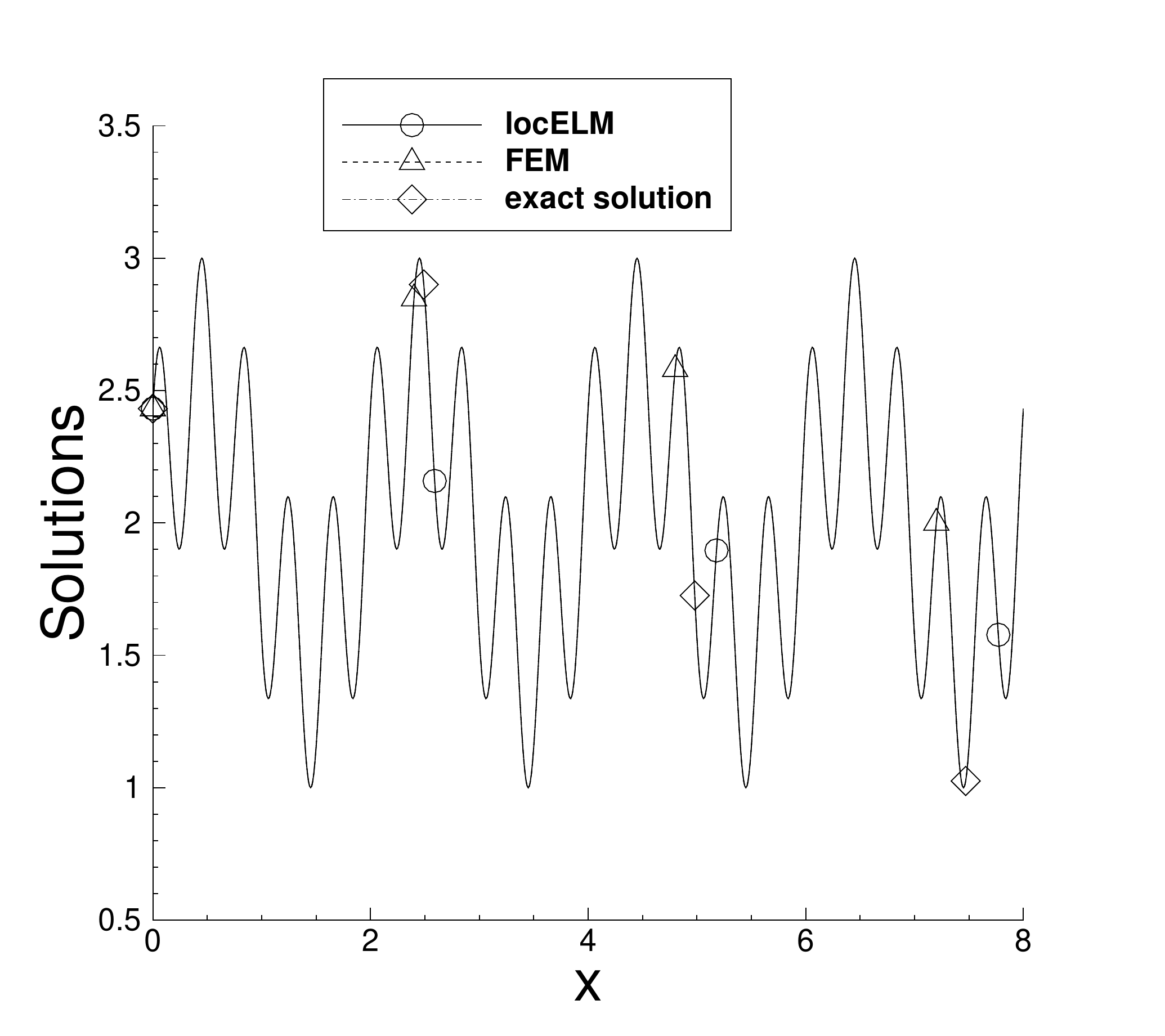}(a)
    \includegraphics[width=2in]{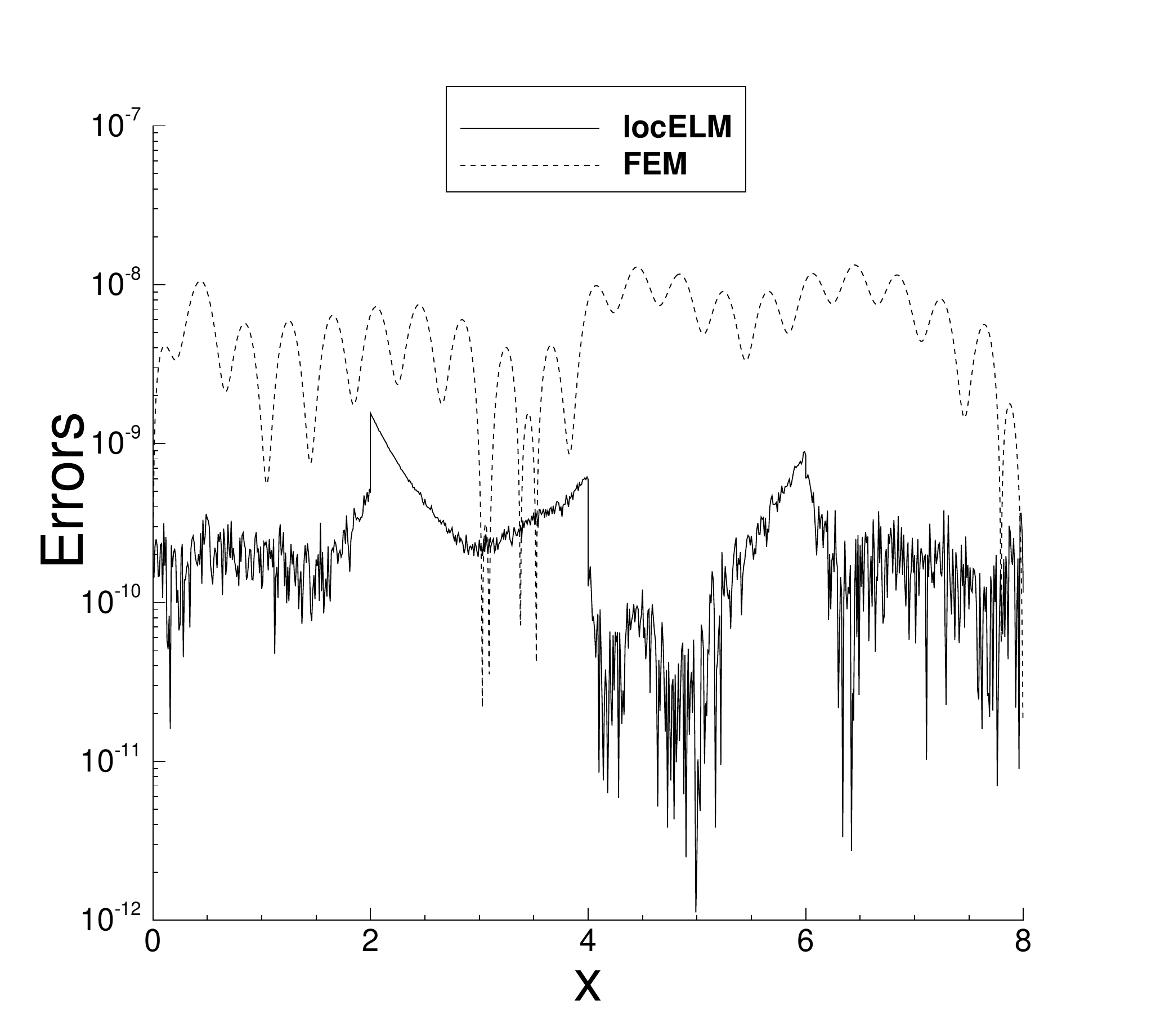}(b)
    \includegraphics[width=2in]{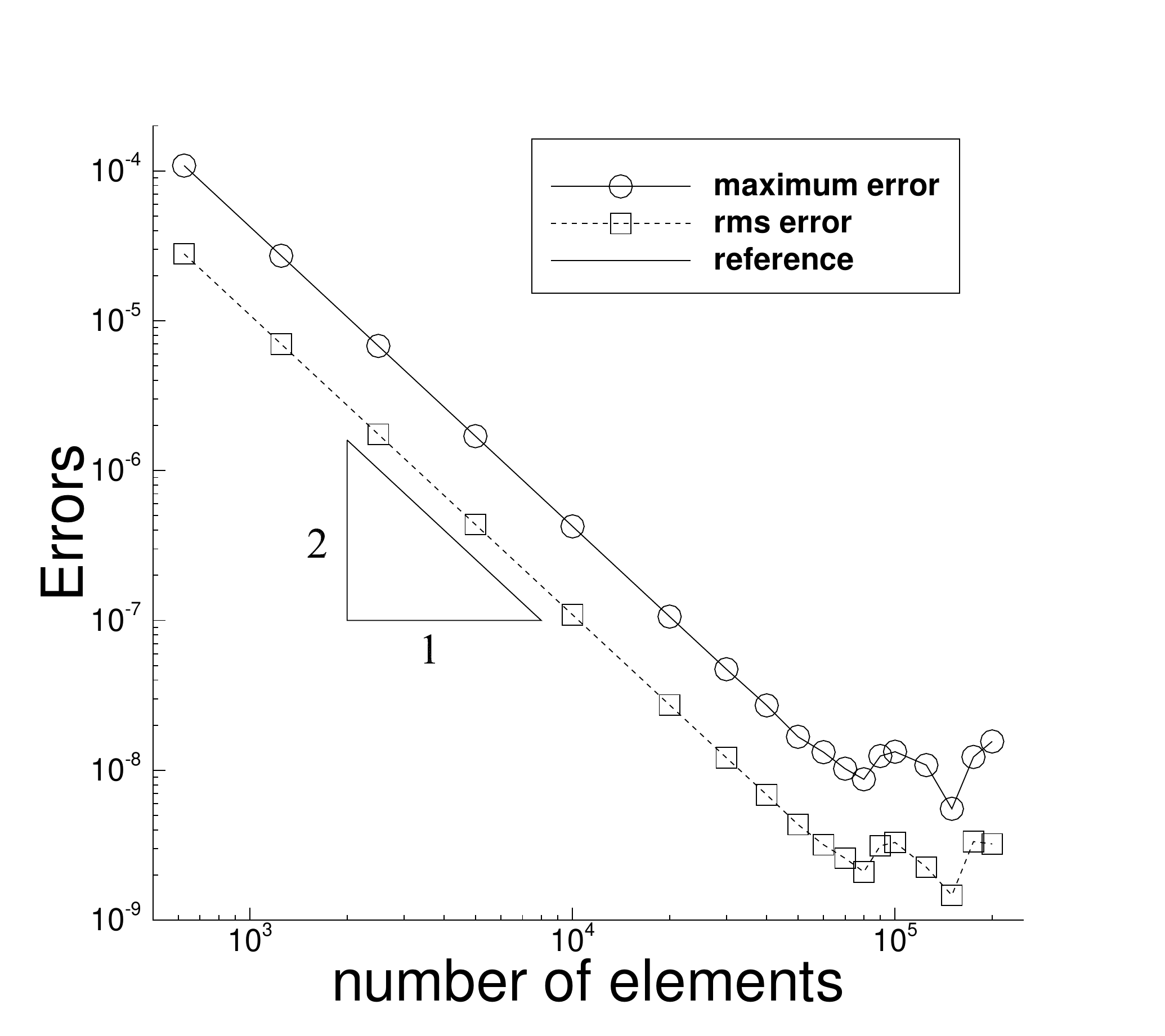}(c)
  }
  \caption{Comparison between locELM and FEM (1D Helmholtz equation):
    Profiles of (a) the solutions and (b) their absolute errors, computed using
    the finite element method (FEM) and the current locELM method.
    (c) The maximum and rms errors in the domain versus 
    the number of elements from the FEM simulations, showing its second-order
    convergence rate.
  }
  \label{fg_helm1d_9}
\end{figure}

\begin{table}[tb]
  \centering
  \begin{tabular}{l|lllllll}
    \hline
    method & elements & sub-domains & $Q$ & $M$ & maximum error & rms error  & wall
    time (seconds)\\ \hline
    locELM
    & -- & $4$ & $100$ & $75$ & $4.02e-8$ & $5.71e-9$ & $0.67$ \\
    & -- & $4$ & $100$ & $100$ & $1.56e-9$ & $2.25e-10$  & $1.1$ \\
    & -- & $4$ & $100$ & $125$ & $1.42e-10$ & $2.55e-11$ & $1.3$ \\
    \hline
    FEM
    & $25,000$ & -- & -- & -- & $6.82e-8$ & $1.74e-8$ & $0.32$ \\
    & $50,000$ & -- & -- & -- & $1.67e-8$ & $4.35e-9$ & $0.62$ \\
    & $100,000$ & -- & -- & -- & $1.33e-8$ & $3.30e-9$ & $1.24$ \\
    \hline
  \end{tabular}
  \caption{1D Helmholtz equation: Comparison between the current locELM
    method and the finite element method (FEM),
    in terms of the maximum/rms errors in the domain and
    the training or computation
    time. The problem settings correspond to those of Figure
    \ref{fg_helm1d_9}.
  }
  \label{tb_helm1d_10}
\end{table}

Finally we compare the current locELM method with the classical
finite element method (FEM).
We observe that the computational
performance of locELM is comparable to that of FEM, and oftentimes
the locELM performance surpasses that of FEM, in terms of the accuracy and
computational cost.
Figures \ref{fg_helm1d_9}(a) and (b) are comparisons of the solution profiles
and the error profiles obtained using locELM and FEM.
Figure \ref{fg_helm1d_9}(c) shows the maximum and rms errors as a function
of the number of elements obtained using FEM, demonstrating its
second-order convergence rate.
As mentioned before,
the finite element method is implemented in Python
using the FEniCS library. In these
tests uniform linear elements have been used.
For the plots (a) and (b), $100,000$ elements are used in the FEM simulation.
In the locELM simulation, we have employed $N_e=4$ uniform sub-domains,
$Q=100$ uniform collocation points per sub-domain, $M=100$ training parameters
per sub-domain, a single hidden layer in the local neural networks,
and $R_m=3.0$ when generating the random coefficients.
It is evident that both FEM and locELM produce accurate
solutions.

Table \ref{tb_helm1d_10} provides a more comprehensive comparison
between locELM and FEM for the 1D Helmholtz equation,
with regard to the accuracy and computational cost.
Here we list the maximum and rms errors in the domain,
and the training or computation time, obtained using
locELM and FEM corresponding
to several numerical resolutions.
The data show that the current locELM method is very
competitive compared with FEM. For example, the locELM case with
$M=75$ training parameters/sub-domain is similar in performance
to the FEM case with $50,000$ elements, with comparable
values for the numerical errors and
the wall time. The locELM cases with $M=100$ and $M=125$ training
parameters/sub-domain have wall time values comparable to the FEM case
with $100,000$ elements, but the numerical errors of these locELM
cases are considerably smaller than those of the FEM case.

%


\subsection{Advection Equation}

We next test the locELM method using the advection equation
in one spatial dimension plus time, and we will
demonstrate the capability of the method, when combined with the block
time-marching strategy, for long-time simulations.
Consider the spatial-temporal domain,
$\Omega=\{(x,t)\ |\ x\in[a_1,b_1], \ t\in[0,t_f]  \}$,
and the initial/boundary-value problem with the advection equation
on this domain,
\begin{subequations}
  \begin{align}
    &
    \frac{\partial u}{\partial t} - c\frac{\partial u}{\partial x} = 0,
    \label{wav1_1} \\
    &
    u(a_1,t) = u(b_1,t), \\
    &
    u(x,0) = h(x) \label{wav1_2}
  \end{align}
\end{subequations}
where $u(x,t)$ is the field function to be solved for,
the constant $c$ denotes the wave speed, and we impose the periodic boundary condition
on the spatial domain boundaries
$x=a_1$ and $b_1$.
$h(x)$ denotes
the initial wave profile given by
\begin{equation}\label{eq_wav1_ic}
  h(x) = 2\sech\left[\frac{3}{\delta_0}\left(x-x_0 \right)\right],
\end{equation}
where $x_0$ is the peak location of the wave and $\delta_0$ is a constant
that controls the width of the wave profile.
The above equations and the domain specification contain several constant parameters,
and we employ the following values in this problem,
\begin{equation}\label{eq_wav1_par}
  a_1 = 0, \quad
  b_1 = 5, \quad
  c = -2, \quad
  \delta_0 = 1, \quad
  x_0 = 2.5, \quad
  t_f=2, \ \text{or}\ 10,\ \text{or}\ 100.
\end{equation}
The temporal domain size $t_f$ is varied in different tests  and will be specified in
the discussions below.
This problem has the following solution
\begin{equation}\label{wav1_3}
  u(x,t) = 2\sech\left[\frac{3}{\delta_0}\left(-\frac{L_1}{2}+ \xi \right)  \right], \quad
  \xi = \bmod\left(x-x_0+ct+\frac{L_1}{2}, L_1\right), \quad
  L_1 = b_1 - a_1,
\end{equation}
where $\bmod$ denotes the modulo operation.

\begin{figure}
  \centerline{
    \includegraphics[height=2.5in]{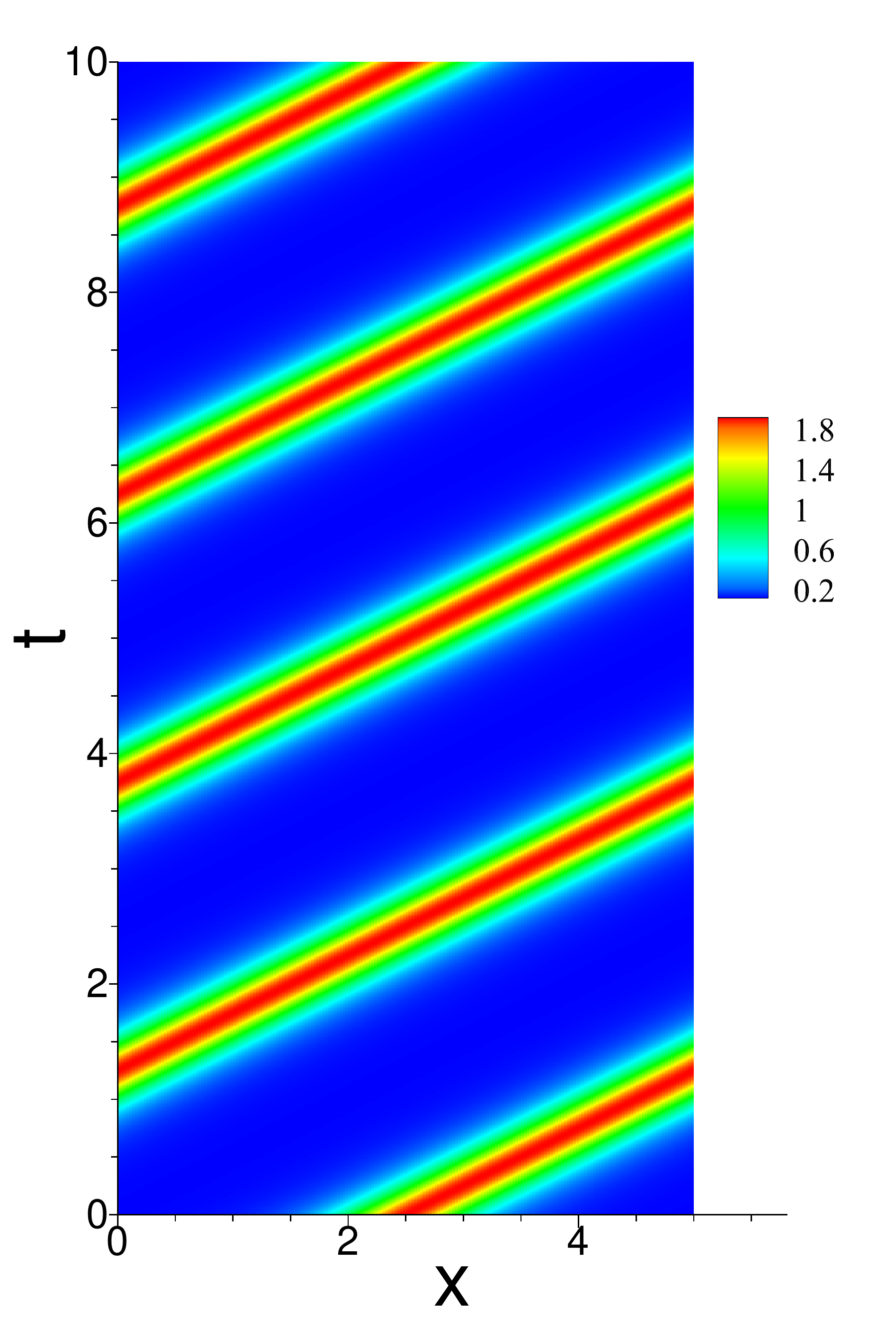}(a)
    \includegraphics[height=2.5in]{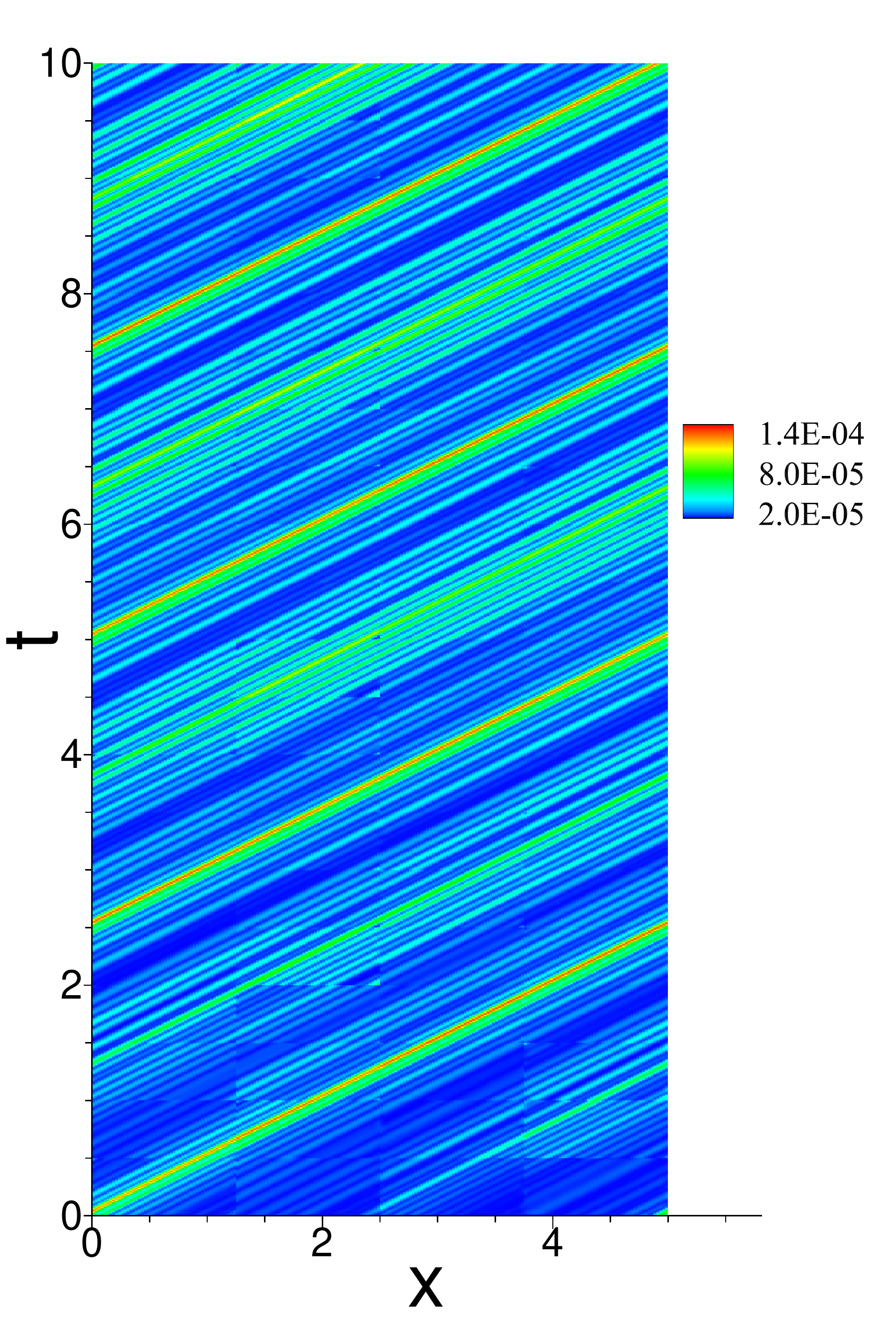}(b)
  }
  \caption{
    Advection equation: Distributions of (a) the locELM solution
    and (b) its absolute error in the spatial-temporal plane.
    The temporal domain size is $t_f=10$, and
    $10$ time blocks are used in the simulation.
  }
  \label{fg_wav1_1}
\end{figure}


We simulate this problem using the locELM method together with the
block time-marching strategy
from Section \ref{sec:unsteady}, by restricting the method to one spatial dimension.
We divide the overall spatial-temporal domain
into $N_b$ uniform blocks along the temporal direction, with a time block size
$\Gamma = \frac{t_f}{N_b}$. The spatial-temporal domain of each time block is then partitioned
into $N_x$ uniform sub-domains along the $x$ direction and $N_t$ uniform sub-domains in time,
leading to $N_e=N_xN_t$ uniform sub-domains in each time block.
$C^0$ continuity is imposed on the sub-domain boundaries in both the $x$ and $t$ directions.
Within each sub-domain, let $Q_x$ denote the number of uniform collocation points along
the $x$ direction and $Q_t$ denote the number of uniform collocation points in time,
leading to $Q=Q_xQ_t$ uniform collocation points in each sub-domain.

The local neural network corresponding to each sub-domain contains
an input layer of two nodes (representing $x$ and $t$),
a single hidden layer
with $M$ nodes and the $\tanh$ activation function,
and an output layer (representing the solution $u$) of a single node.
The output layer is linear and contains no bias.
An additional affine mapping
normalizing the input $x$ and $t$ data to the interval $[-1,1]\times[-1,1]$
has been incorporated into the local neural networks right behind the input layer
for each sub-domain. The number of training parameters per sub-domain
corresponds to $M$, the width of the hidden layer.
The weight and bias coefficients in the hidden layer
are pre-set to uniform random values generated on $[-R_m,R_m]$,
as in the previous section.

The locELM simulation parameters include
the number of sub-domains ($N_x$, $N_t$, $N_e$), the number of collocation
points per sub-domain ($Q_x$, $Q_t$, $Q$), the number of training parameters
per sub-domain ($M$), and the maximum magnitude of the random coefficients ($R_m$).
The degrees of freedom within a sub-domain are characterized by
($Q,M$). The degrees of freedom in each time block are characterized by
($N_eQ$, $N_eM$).
We use a fixed seed value $1$ for the Tensorflow random number generators in all the tests
of this sub-section, so that all the numerical tests here are repeatable.

Figure \ref{fg_wav1_1} illustrates the solution from the locELM simulation.
Plotted here are the distributions of the locELM solution and its absolute error in
the spatial-temporal plane.
In this test, the temporal domain size is $t_f=10$, and we  employ
$10$ uniform time blocks ($N_b=10$) in this domain. Within each time block,
we have employed $N_e=8$ uniform sub-domains (with $N_x=4$ and $N_t=2$),
and $Q=20\times 20$ uniform collocation points ($Q_x=Q_t=20$)
in each sub-domain. We employ $M=300$ training parameters per sub-domain,
and $R_m=1.0$ when generating the random weight/bias coefficients.
It is evident that the current method has captured the wave solution accurately.

\begin{figure}
  \centerline{
    \includegraphics[width=2.in]{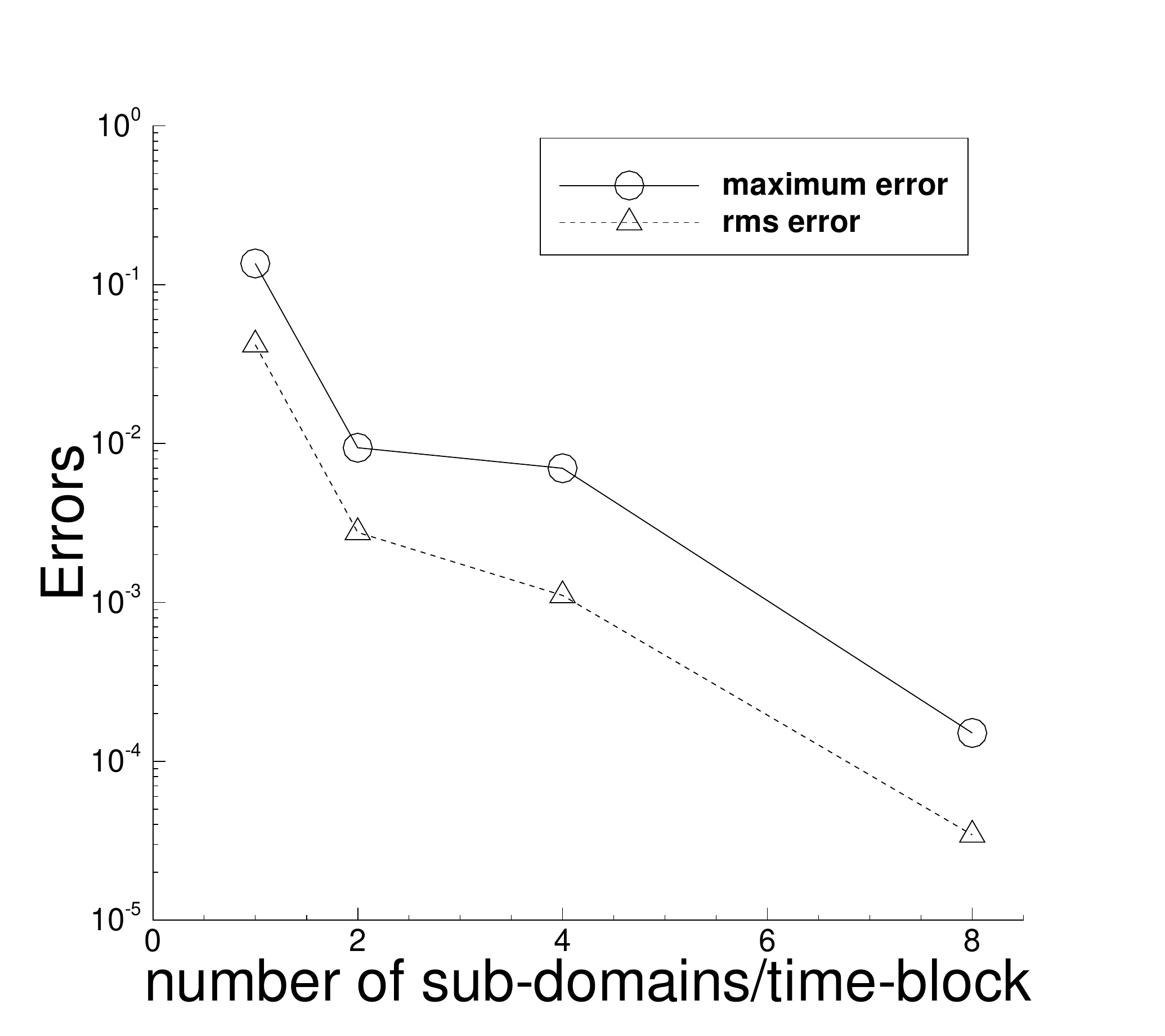}(a)
    \includegraphics[width=2.in]{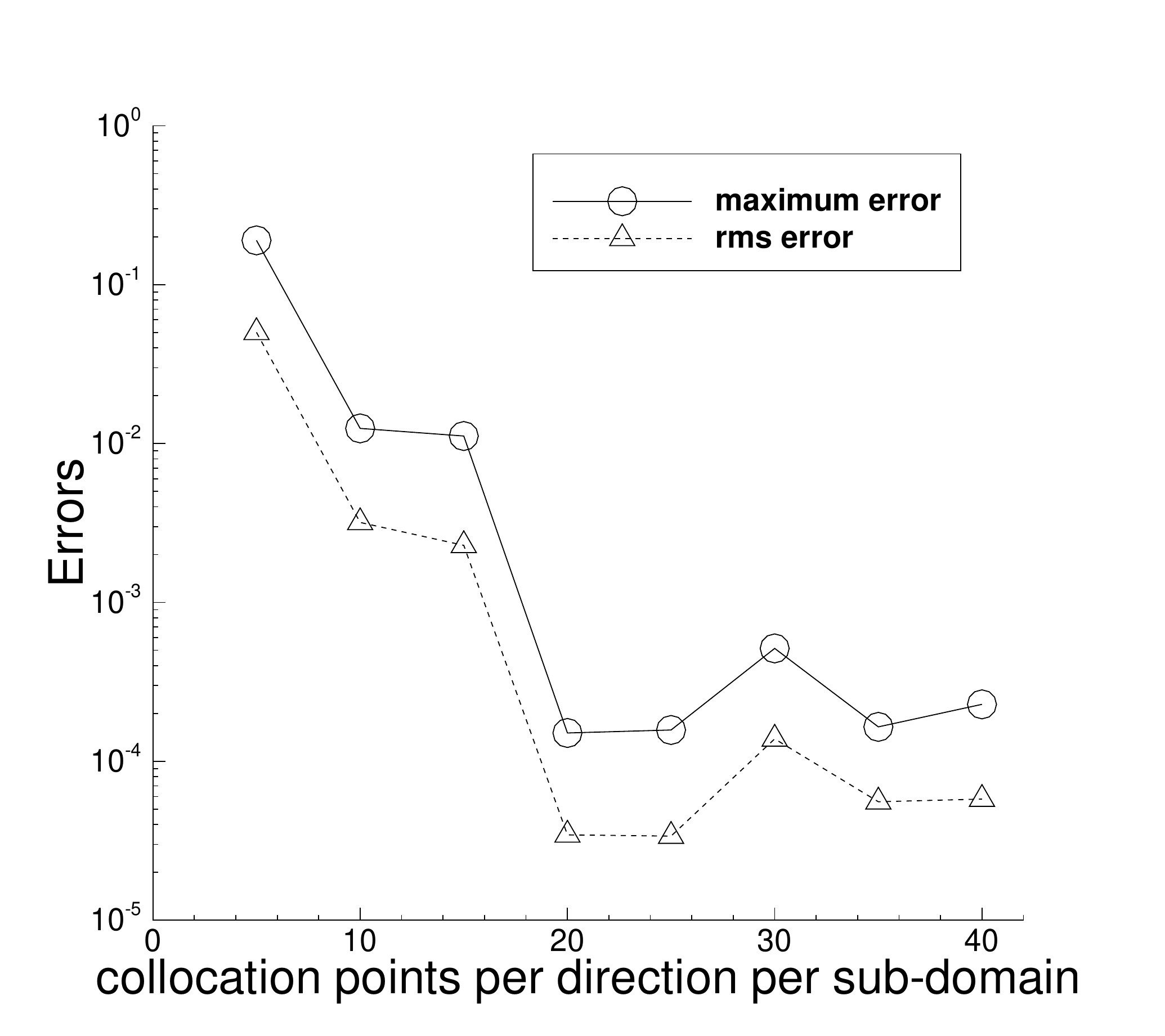}(b)
    \includegraphics[width=2.in]{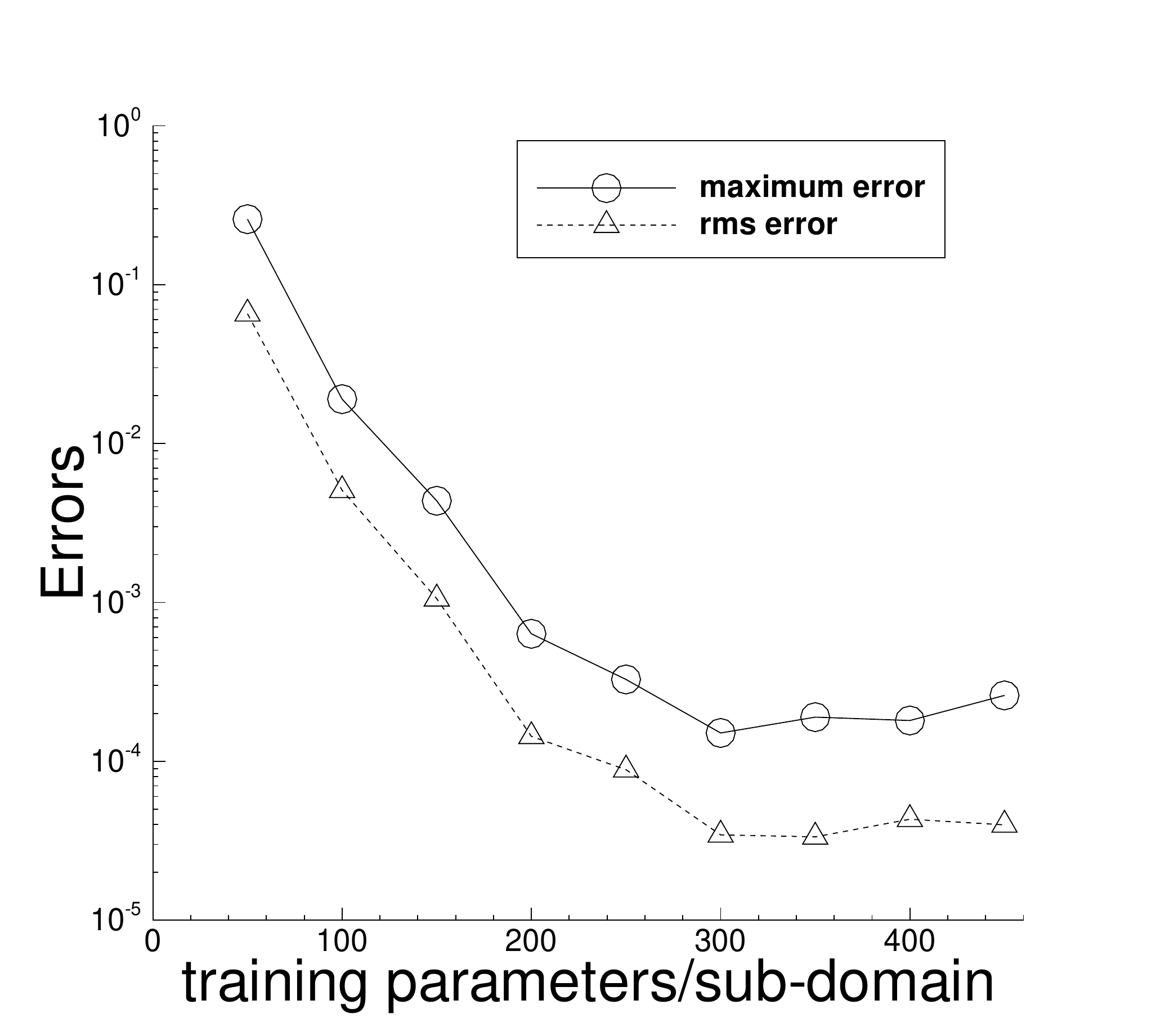}(c)
  }
  \caption{Effect of the degrees of freedom
    (advection equation): the maximum and rms errors in the overall domain
    as a function of (a) the number of sub-domains,
    (b) the number of collocation points in each direction
    per sub-domain,
    and (c) the number of training parameters per sub-domain.
    Temporal domain size is $t_f=10$, and $10$ time blocks have been used.
    In (a), the degrees of freedom per sub-domain are fixed.
    In (b) and (c), $N_e=8$ sub-domains per time block are used.
  }
  \label{fg_wav1_3}
\end{figure}

The effect of the degrees of freedom  on the simulation accuracy 
is illustrated by Figure \ref{fg_wav1_3}.
In this group of tests, the temporal domain size is fixed at $t_f=10$. 
We have employed $N_b=10$ uniform time blocks within the domain,
one hidden layer in each local nueral
network, and $R_m=1.0$
when generating the random weight/bias coefficients for the hidden layers.
Figure \ref{fg_wav1_3}(a) illustrates the effect of the number of sub-domains per time
block, when the degrees of freedom per sub-domain are fixed.
Here the number of sub-domains within each time block is varied systematically.
We employ a fixed set of $Q=20\times 20$ uniform 
collocation points per sub-domain ($Q_x=Q_t=20$),
and fix the number of training parameters per sub-domain
at $M=300$.
This plot shows the maximum and rms errors in the domain
as a function of the number of sub-domains per time block.
Here the case with $N_e=2$ sub-domains/time-block corresponds to $(N_x,N_t)=(2,1)$.
The case with $N_e=4$ sub-domains corresponds to $(N_x,N_t)=(2,2)$,
and the case with $N_e=8$
sub-domains corresponds to $(N_x,N_t)=(4,2)$.
It can be observed that, with increasing sub-domains/time-block, the rate of
reduction in the
errors, while not very regular, is approximately exponential.

Figure \ref{fg_wav1_3}(b) shows the maximum and rms errors in the entire
spatial-temporal domain as a function of the number of collocation
points in each direction (with $Q_x=Q_t$ maintained) in each sub-domain.
Figure \ref{fg_wav1_3}(c) shows the maximum and rms errors in the entire domain
as a function of the number of training parameters per sub-domain.
In these tests, 
we have employed $8$ sub-domains
($N_x=4$, $N_t=2$) per time block.
For those tests of Figure \ref{fg_wav1_3}(b), the number of training parameters/sub-domain
is fixed at $M=300$. For the tests of Figure \ref{fg_wav1_3}(c), the number of
collocation points/sub-domains is fixed at $Q=20\times 20$ ($Q_x=Q_t=20$).
With the increase of the collocation points in each direction,
or the increase of the training parameters per sub-domain, we can observe
an approximately exponential decrease in the maximum and rms
errors.
When the number of collocation points (or training parameters) increases above
a certain point, the errors start to stagnate, apparently because of
the fixed number of training parameters (or the fixed number of collocation points)
in these tests.
The sense of convergence exhibited by the current locELM method
is unmistakable.

\begin{figure}
  \centerline{
    \includegraphics[width=2.2in]{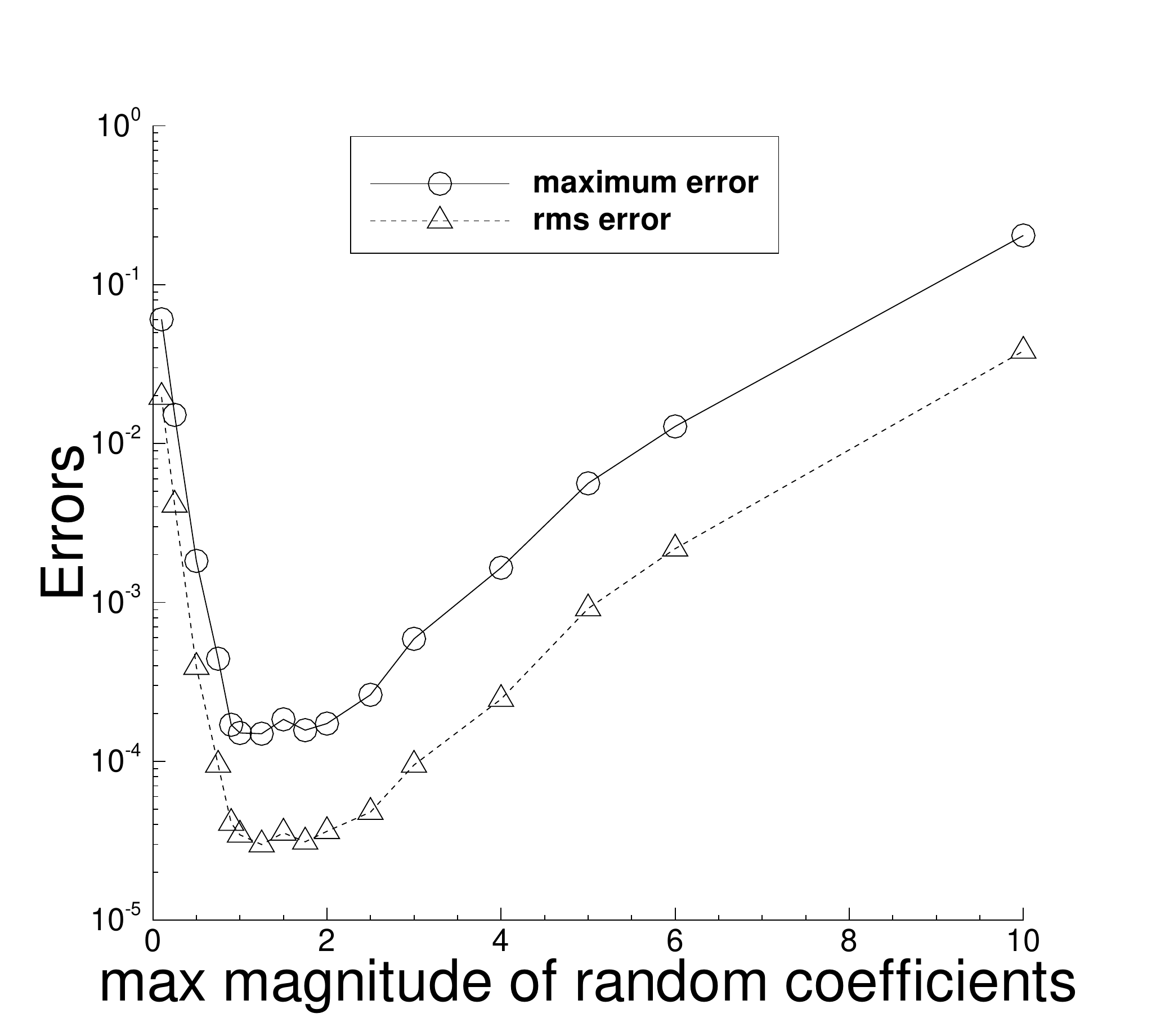}
  }
  \caption{Effect of the random coefficients
    (advection equation): the maximum and rms errors in the domain
    as a function of $R_m$, the maximum magnitude of the random coefficients in
    hidden layer of local neural networks.
  }
  \label{fg_wav1_4}
\end{figure}

The effect of the random coefficients in the hidden layers of the local
neural networks on the simulation accuracy is illustrated in Figure \ref{fg_wav1_4}.
This plot shows the maximum and rms errors in the domain
as a function of $R_m$, the maximum magnitude
of the random weight/bias coefficients.
In this set of experiments, the temporal
domain size is $t_f=10$, and $N_b=10$ time blocks are used in the domain.
We have employed $8$ uniform sub-domains per time block ($N_x=4$, $N_t=2$),
$Q=20\times 20$ uniform collocation points per sub-domain ($Q_x=Q_t=20$),
and $M=300$ training parameters per sub-domain.
The weight/bias coefficients in the hidden layers of the local neural
networks are set to uniform random values generated on $[-R_m,R_m]$,
and $R_m$ is varied systematically in these tests.
Very large or very small values of $R_m$ have an adverse effect on
the simulation accuracy. Better accuracy generally
corresponds to a range of moderate $R_m$ values.

\begin{figure}
  \centerline{
    \includegraphics[width=6.2in]{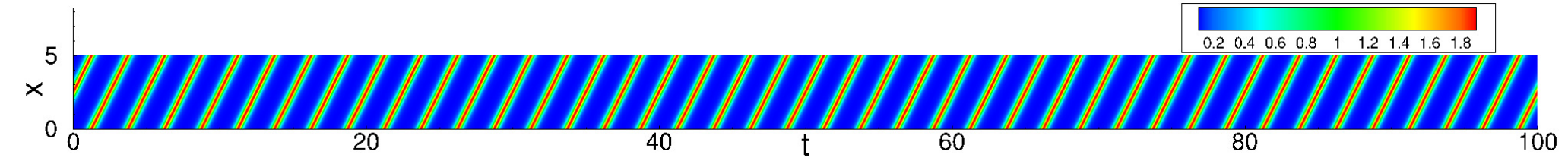}(a)
  }
  \centerline{
    \includegraphics[width=6.2in]{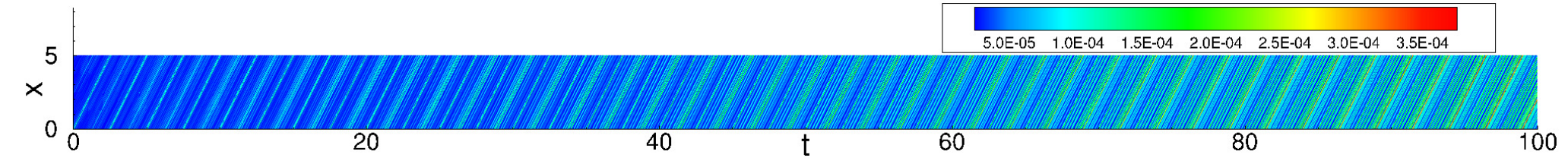}(b)
  }
  \caption{Long-time simulation of the
    advection equation: distributions of (a) the locELM solution
    and (b) its absolute error in the spatial-temporal plane for
    a long-time simulation.
    In these tests $100$ time blocks in the domain and $8$ sub-domains
    per time block are used.
  }
  \label{fg_wav1_5}
\end{figure}

\begin{figure}
  \centerline{
    \includegraphics[width=6in]{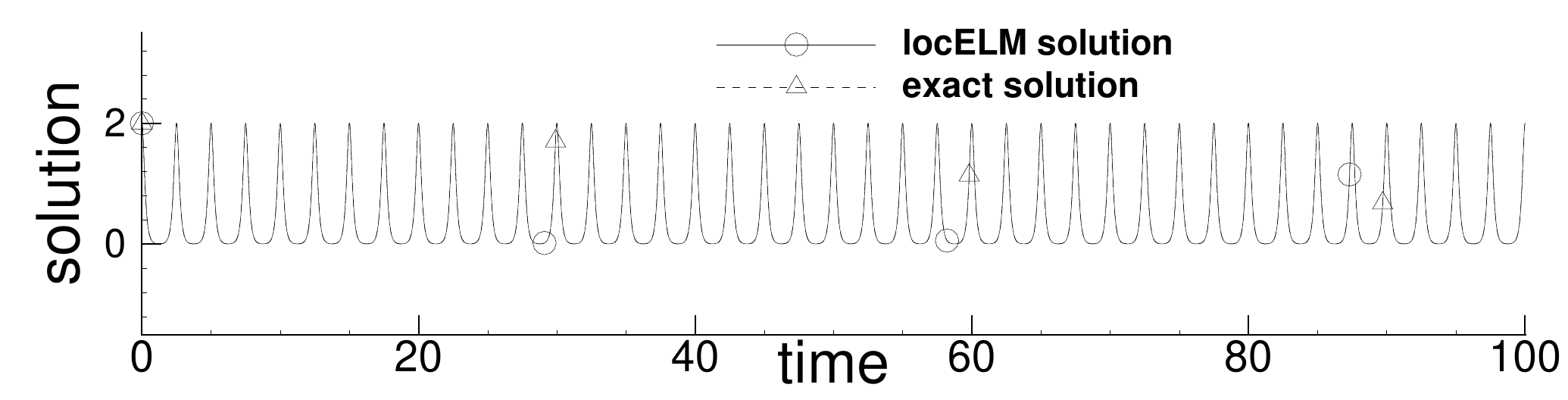}(a)
  }
  \centerline{
    \includegraphics[width=6in]{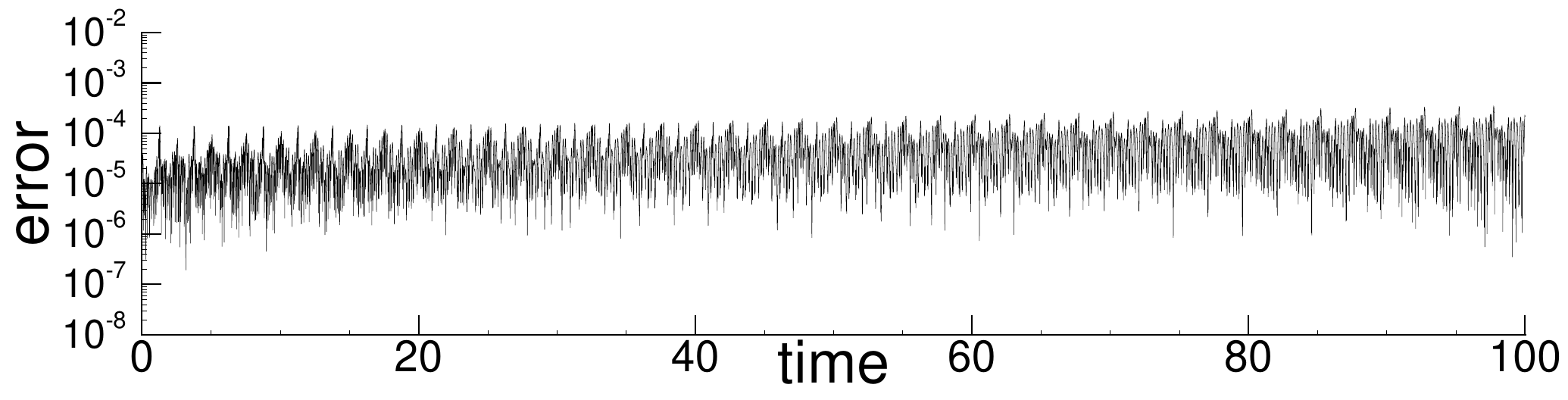}(b)
  }
  \centerline{
    \includegraphics[width=2in]{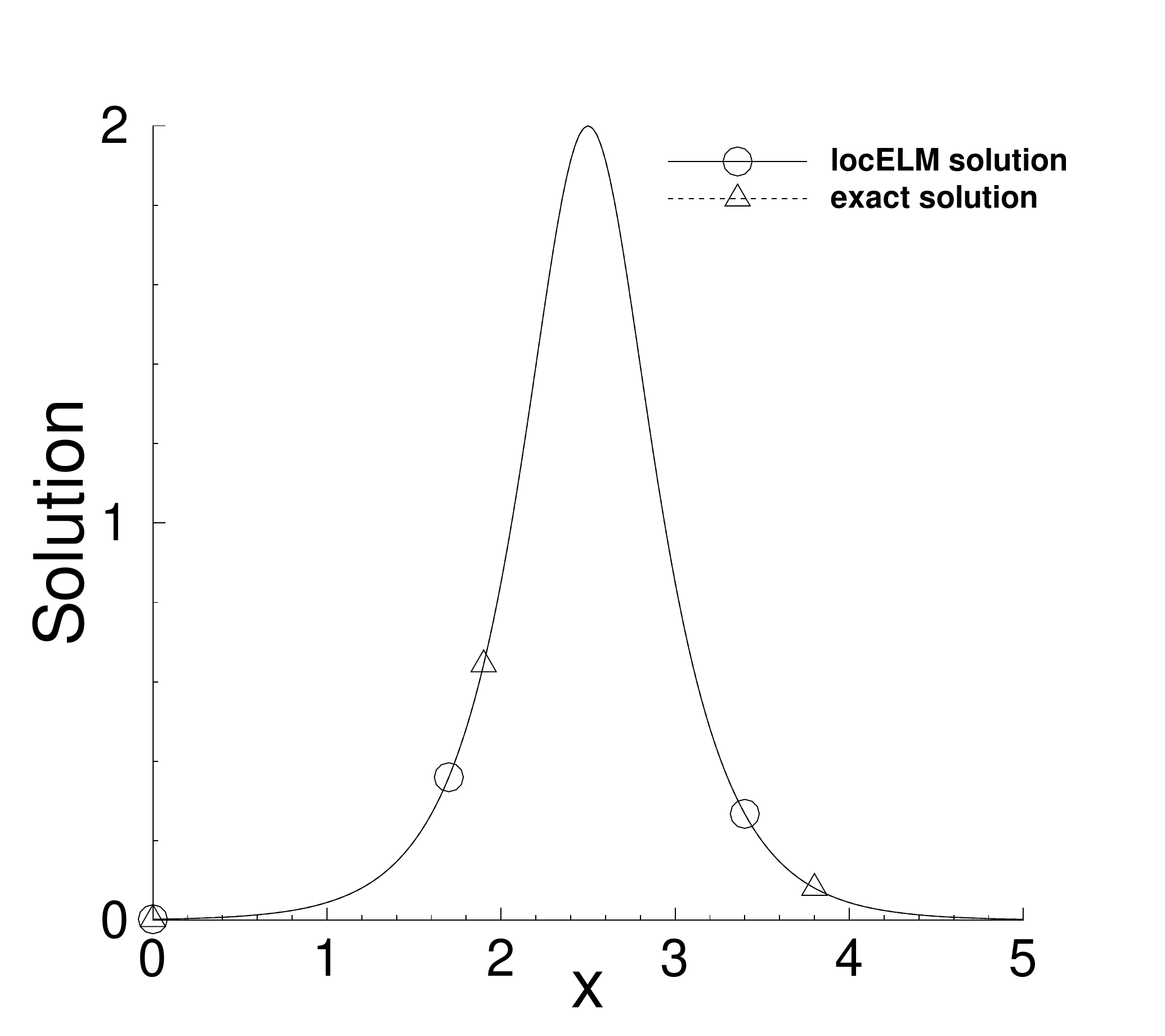}(c)
    \includegraphics[width=2in]{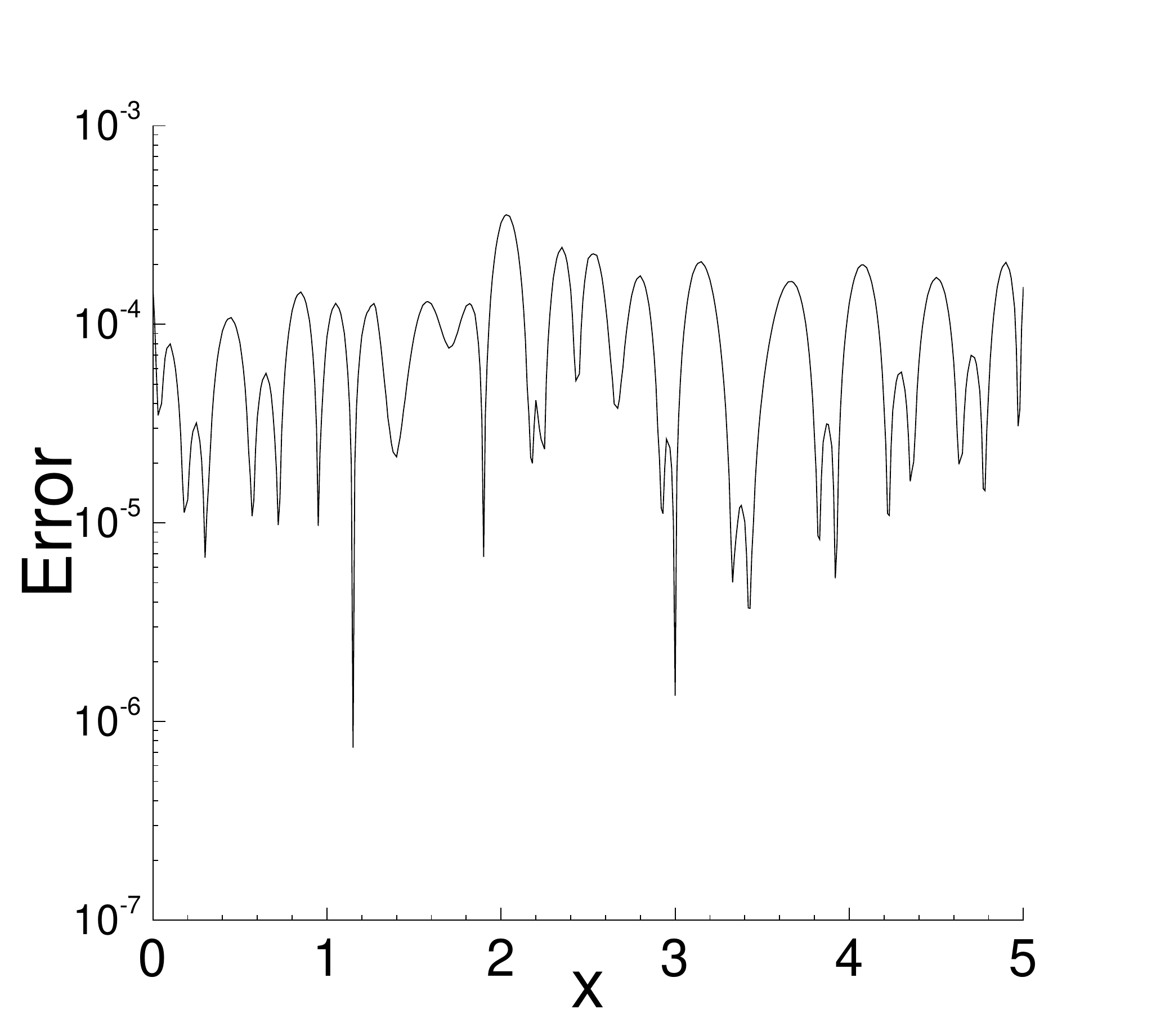}(d)
    \includegraphics[width=2in]{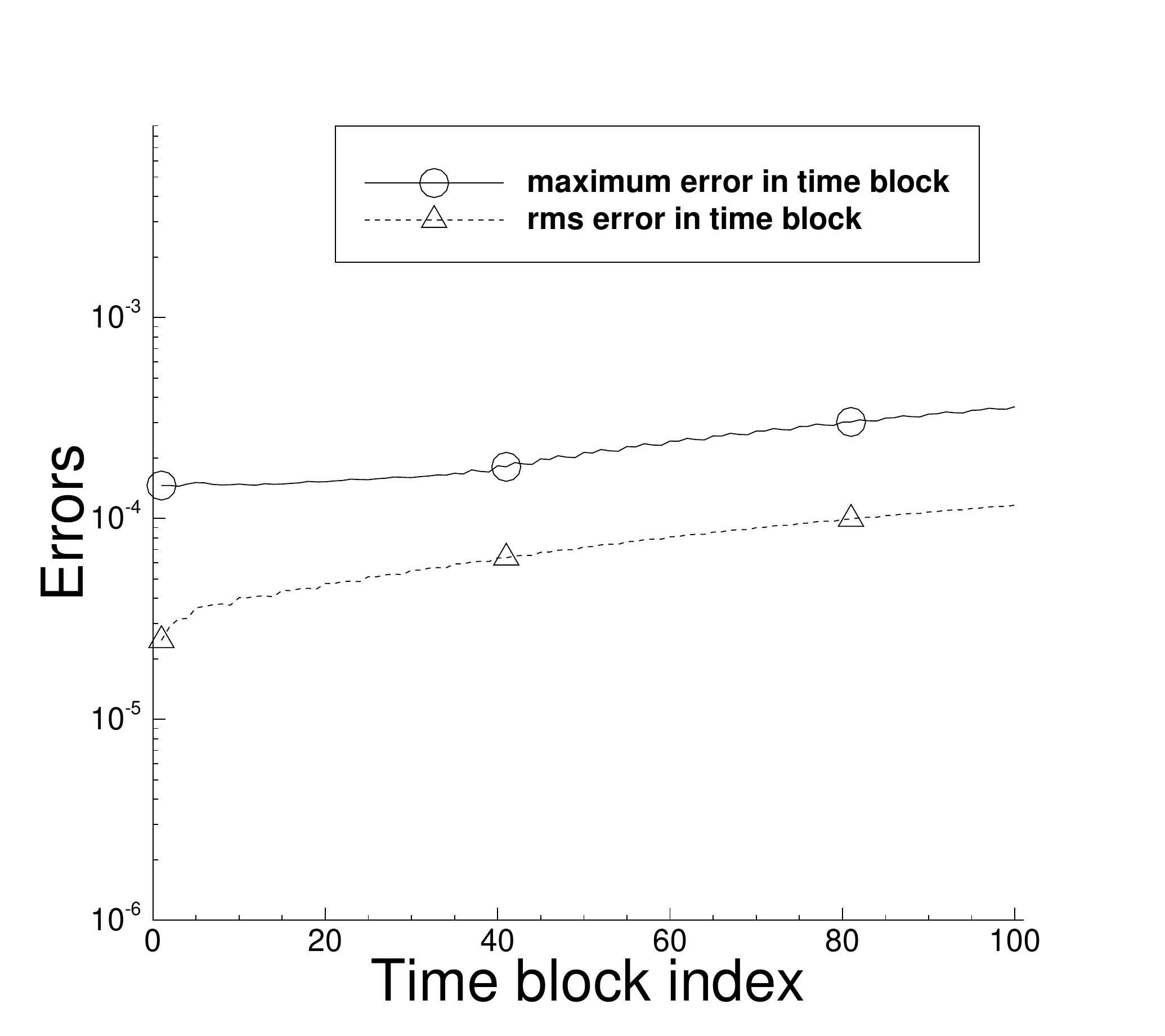}(e)
  }
  \caption{Long-time simulation of the advection equation:
    Time histories of the locELM solution (a) and its absolute error
    against the exact solution (b) at the mid-point ($x=2.5$) of the spatial
    domain.
    Profiles of the locELM solution (c) and its absolute
    error against the exact solution (d)
    at the last time instant $t=100$.
    (e) Time histories of the maximum and rms errors
    in each time block.
    The problem settings correspond to those of Figure \ref{fg_wav1_5}.
  }
  \label{fg_wav1_6}
\end{figure}

Thanks to its accuracy and favorable computational cost,
it is feasible to perform
long-time simulations of time-dependent PDEs
using the current locELM method.
Figures \ref{fg_wav1_5} and \ref{fg_wav1_6} demonstrate
a long-time simulation
of the advection equation with the current method.
In this simulation, the temporal domain size is set to $t_f=100$,
which amounts to approximately $40$ periods of the wave propagation time.
In the simulation
we have employed $100$ uniform time blocks in the domain,
$8$ uniform sub-domains per time block (with $N_x=4$ and $N_t=2$),
$20\times 20$ uniform collocation points per sub-domain (i.e.~$Q_x=Q_t=20$),
$300$ training parameters per sub-domain ($M=300$),
a single hidden layer in each local neural network,
and $R_m=1.0$ when generating the random weight/bias coefficients
for the hidden layers of the local neural networks.
The total network training time for this locELM
computation is about $892$ seconds.
Figure \ref{fg_wav1_5} shows the distributions of the locELM solution
and its absolute error in the spatial-temporal plane. 
Figures \ref{fg_wav1_6}(a) and (b) are the time histories of
the locELM solution and its absolute error at the mid-point ($x=2.5$)
of the spatial domain. The time history of the exact solution
at this point is also shown in Figure \ref{fg_wav1_6}(a),
which can be observed to overlap with that of the locELM solution.
Figures \ref{fg_wav1_6}(c) and (d) show the locELM-computed wave profile and
its absolute-error profile at the last time instant $t=100$.
We have also computed and monitored the maximum and rms errors of the locELM solution
within each time block.
Figure \ref{fg_wav1_6}(e) shows these errors versus the time block index,
which represents essentially the time histories of these block-wise
maximum and rms errors.
All these results show that the current method has captured
the solution to the advection equation quite accurately
in the long-time simulation.
Accurate simulation of the advection equation
in long time integration is challenging, even for classical
numerical methods. 
The results presented here demonstrate the capability and the promise
of the current method in tackling long-time dynamical simulations
of these challenging problems.

\begin{figure}
  \centerline{
    \includegraphics[width=2in]{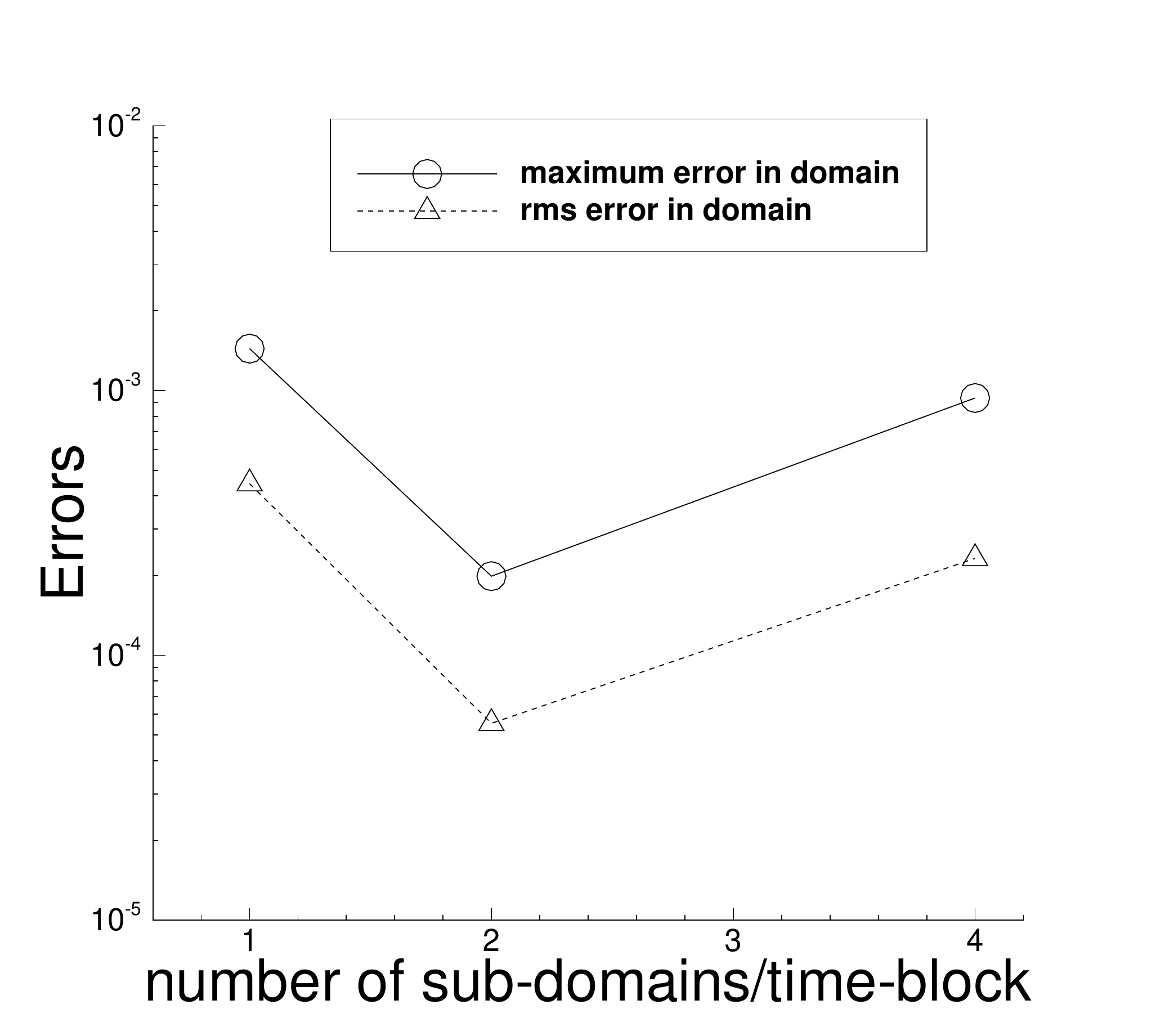}(a)
    \includegraphics[width=2in]{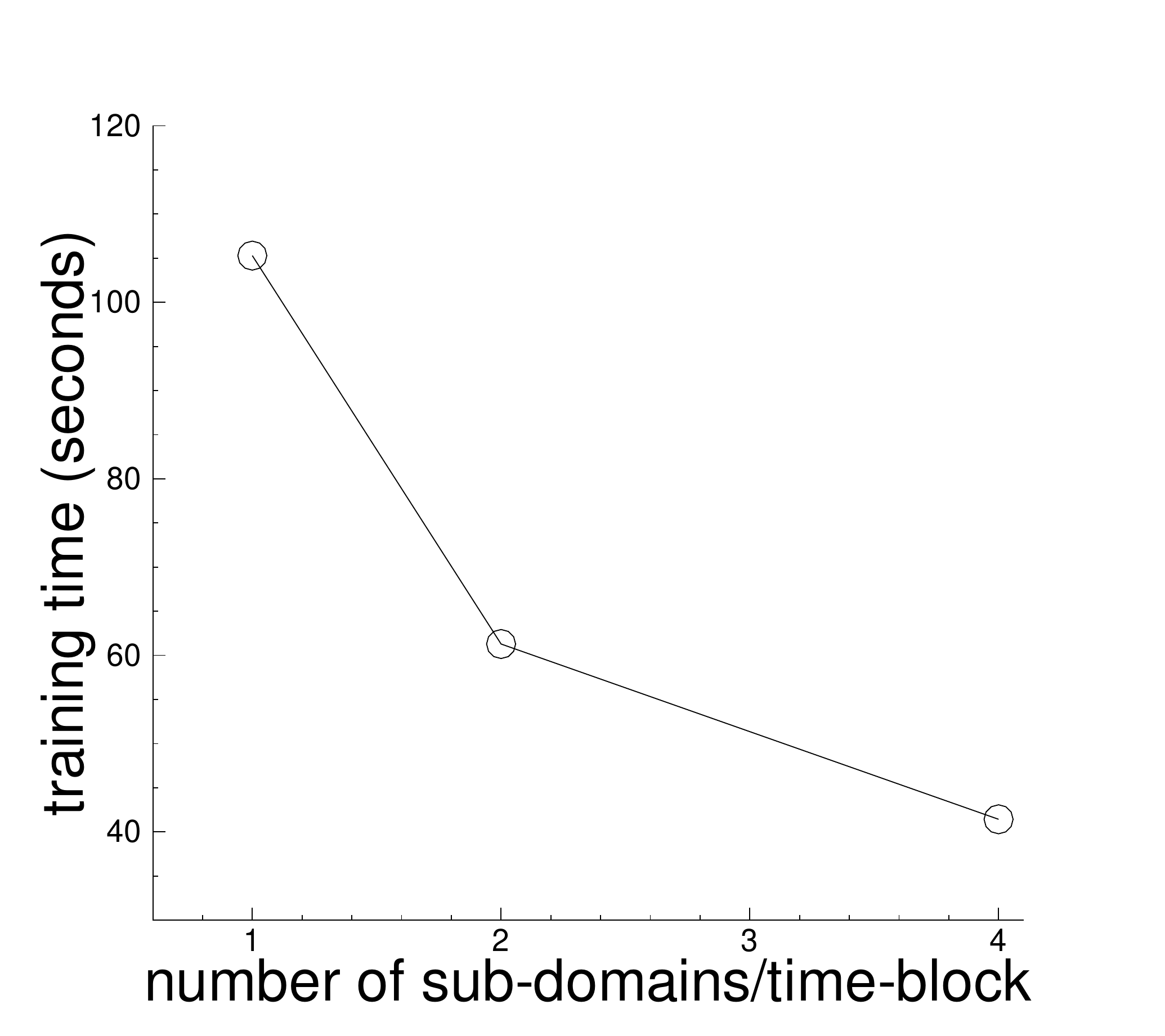}(b)
  }
  \caption{Effect of the number of sub-domains,
    with fixed total degrees of freedom in the domain
    (advection equation): (a) the maximum and rms errors in the overall
     domain,
     and (b) the training time, as a function of the number of sub-domains per
     time-block.
    10 time blocks are used in the domain.
    For all cases,
    the total number of training parameters per time block is fixed at $2500$,
    and the total number of collocation points per time block is
    approximately $2500$.
  }
  \label{fg_wav1_9}
\end{figure}

Figure \ref{fg_wav1_9} provides a comparison between the
locELM and global ELM results.
We fix the total degrees of freedom in each time block (temporal
dimension $t_f=10$, $10$ time blocks,
$1600$ training parameters/time-block, approximately $2500$ collocation
points/time-block), and vary the
number of sub-domains per time block.
Figure \ref{fg_wav1_9}(a) shows the maximum and rms errors  in the overall
spatial-temporal domain as a function of the number of sub-domains per
time block.
For the case with $4$ sub-domains
per time-block, we have employed the configuration of $N_x=4$ and $N_t=1$,
$M=400$ training parameters
per sub-domain, $Q=25\times 25$ uniform collocation points per sub-domain,
and $R_m=2.0$ when generating the random coefficients in the hidden layers
of the local neural networks.
The error levels with one sub-domain and  multiple sub-domains per time block
are observed to be comparable, with the results of $2$ sub-domains
per time block  more accurate than the others.
Figure \ref{fg_wav1_9}(b) compares the network training time corresponding
to different sub-domains.
The use of multiple sub-domains is observed to significantly
reduce the training time of the neural network,
from around $105$ seconds with a single sub-domain
per time block to around $40$ seconds with $4$ sub-domains per
time block.
The results here confirm what has been observed 
in the previous section. With the same total degrees of freedom in the domain,
the use of multiple sub-domains and local neural networks with
the current locELM method can significantly reduce
the training/computation time, while 
producing results with comparable accuracy when
compared with the global ELM method.

\begin{figure}
  \centerline{
    \includegraphics[width=2.1in]{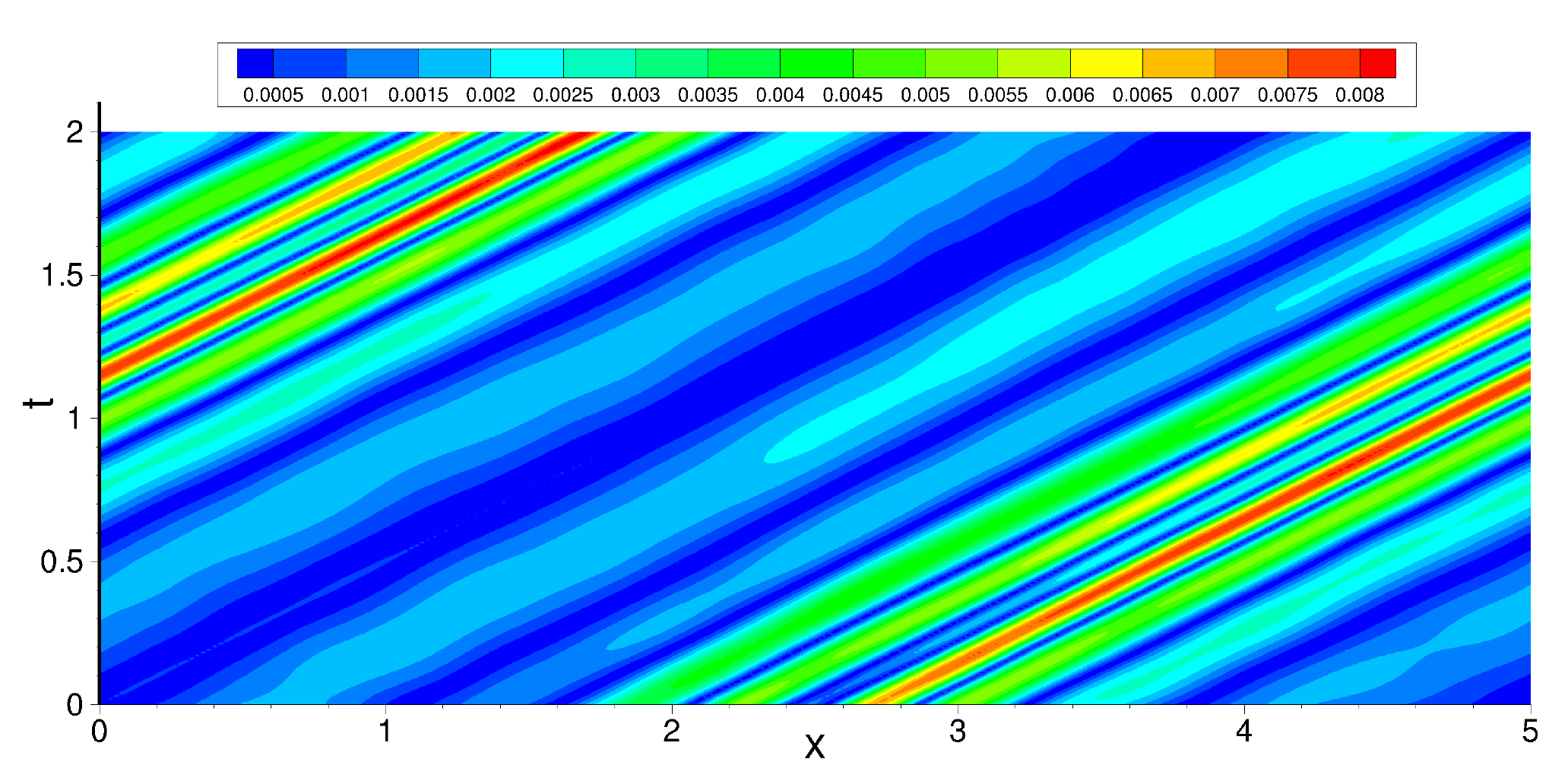}(a)
    \includegraphics[width=2.1in]{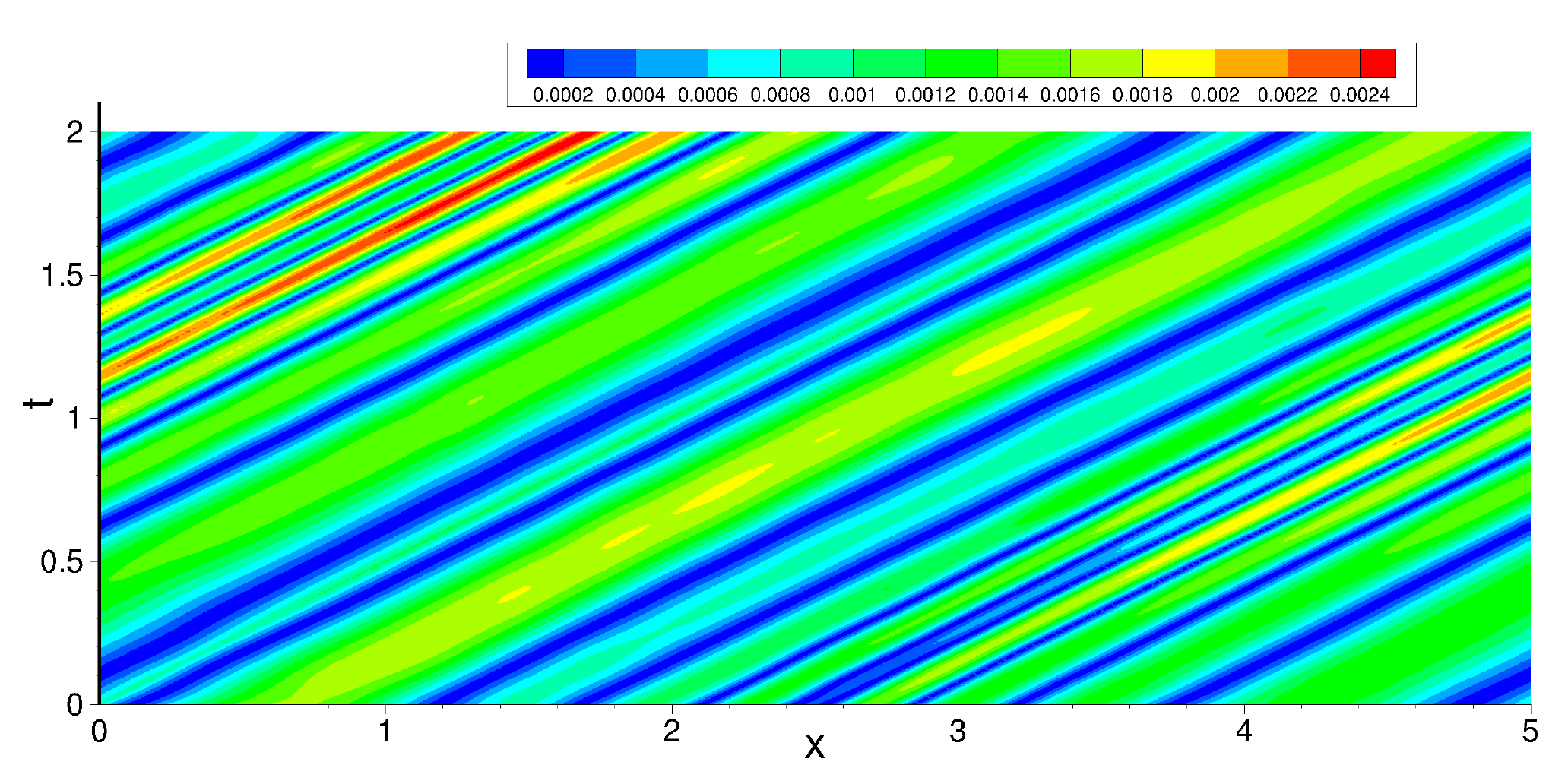}(b)
    \includegraphics[width=2.1in]{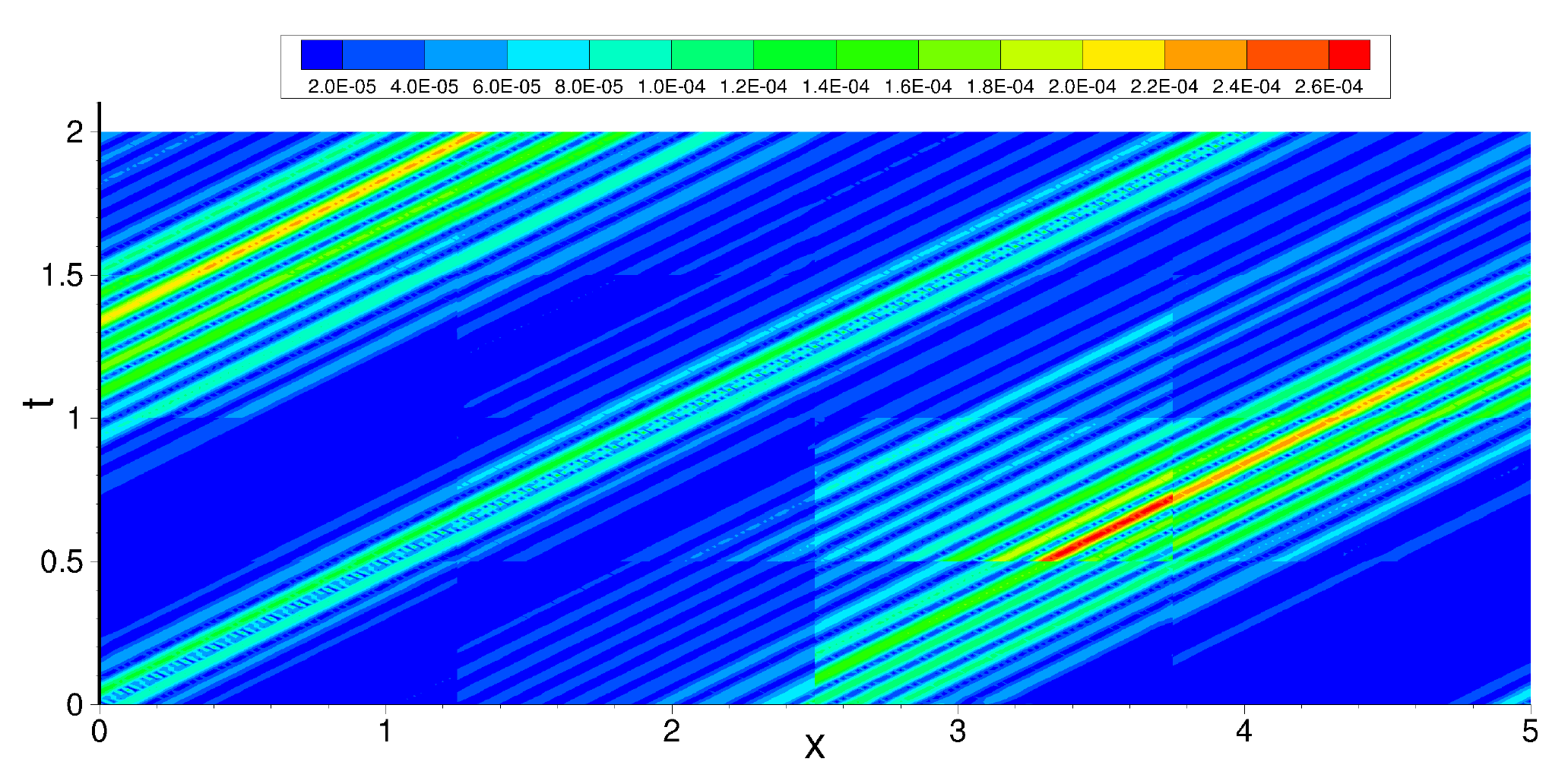}(c)
  }
  \caption{Comparison between locELM and DGM (advection equation):
    distributions of the  absolute errors, computed using the deep Galerkin
    method (DGM)~\cite{SirignanoS2018} with the Adam
    optimizer (a) and the L-BFGS optimizer (b), and computed using
    the current locELM method (c).
  }
  \label{fg_wav1_10}
\end{figure}

\begin{table}
  \centering
  \begin{tabular}{lllll}
    \hline
    method & maximum error & rms error & epochs/iterations & training time (seconds) \\
    DGM (Adam) & $8.37e-3$ & $1.64e-3$ & $60,000$ & $2527.8$ \\
    DGM (L-BFGS) & $2.59e-3$ & $5.37e-4$ & $12,000$ & $1675.9$ \\
    locELM (no block time-marching) & $2.74e-4$ & $6.05e-5$ & $0$ & $43.4$ \\
    locELM (with block time-marching) & $1.83e-4$ & $4.34e-5$ & $0$ & $19.3$ \\
    \hline
  \end{tabular}
  \caption{Advection equation:
    comparison between locELM and DGM.
    The problem settings correspond to those of Figure \ref{fg_wav1_10}.
    The two DGM cases and the locELM case with no block time-marching correspond
    to those of Figure \ref{fg_wav1_10}.
    In the locELM case with block time-marching, two time blocks
    in the domain and $8$ sub-domains per time block are used.
    The total degrees of freedom for this case are identical to
    those of the locELM case with no
    block time marching.
  }
  \label{tab_wav1_11}
\end{table}

Finally we compare the current locELM method with the deep Galerkin method
(DGM)~\cite{SirignanoS2018}, another often-used DNN-based PDE solver,
for solving the advection equation.
Figure \ref{fg_wav1_10} shows distributions of the solutions and their
absolute errors obtained using DGM with the Adam and the
L-BFGS optimizers and using the current locELM method.
The temporal domain size is $t_f=2.0$ in these tests.
With DGM, four hidden layers with a width of $40$ nodes  and
the $\tanh$ activation function in each layer
have been employed in the neural networks.
When computing the residual norms in the DGM loss function, we have partitioned
the domain into $8$ sub-regions ($4$ in $x$ and $2$ in time) and
used $30\times 30$ Gaussian quadrature points in each sub-region
for calculating the integrals.
The periodic boundary condition is enforced exactly using the method
from~\cite{DongN2020} for DGM.
With the Adam optimizer, the neural network has been trained for $60,000$
epochs, with the learning rate decreasing gradually from $0.001$
at the beginning to $2.5\times 10^{-5}$ at the end of training.
With the L-BFGS optimizer, the neural network has been trained for
$12,000$ iterations.
In the locELM simulation, a single time block has been used in the spatial-temporal
domain, i.e.~without block time marching. We employ $16$ sub-domains
(with $4$ sub-domains in $x$ and time) per time block,
$20\times 20$ uniform collocation points in each sub-domain,
$250$ training parameters per sub-domain, a single hidden layer
in each local neural network, and $R_m=2.0$ for
generating the random weight/bias coefficients for the hidden layer
of the local neural networks.
One can observe that the current method produces considerably more
accurate result than DGM for the advection equation.

Table \ref{tab_wav1_11} provides further comparisons between
locELM and DGM. Here we list the maximum and rms
errors in the overall spatial-temporal domain, the number of epochs or iterations
in the network training, and the training time obtained using
DGM (Adam/L-BFGS optimizers) and using locELM without block time marching,
and additionally using locELM together with block time marching.
The problem settings here correspond to those of Figure \ref{fg_wav1_10}, and
the DGM cases and the locELM case without block time marching
correspond to those in Figure \ref{fg_wav1_10}.
In the locELM case with block time marching, we
have employed $2$ uniform time blocks in the spatial-temporal domain,
$8$ sub-domains ($N_x=4$, $N_t=2$)
per time block, $20\times 20$ uniform collocation points per sub-domain,
$250$ training parameters per sub-domain, a single hidden layer in
the local neural networks, and $R_m=2.0$
when generating the random weight/bias coefficients.
So the total degrees of freedom
in this case are identical to
those of the locELM case without block time marching.
The data in the table shows that the current locELM method (with and without
block time marching) is much more accurate than DGM (by an order of magnitude),
and is dramatically faster to train than DGM (by nearly two
orders of magnitude).
In addition, we observe that the locELM method
with block time marching and
a moderate time block size
can significantly reduce
the training time, and simultaneously achieve
the same or better accuracy,
when compared with that without block time marching.




\subsection{Diffusion Equation}



In this subsection we the test the locELM method using the diffusion
equation in one and two spatial dimensions (plus time).
Let us first study the 1D diffusion equation.
We consider the spatial-temporal domain,
$\Omega=\{ (x,t)\ |\ x\in[a_1,b_1], \ t\in[0,t_f] \}$,
and the following initial/boundary-value problem,
\begin{subequations}
  \begin{align}
    &
    \frac{\partial u}{\partial t} - \nu\frac{\partial^2 u}{\partial x^2}
    = f(x,t), \label{eq_diffu_1} \\
    &
    u(a_1,t) = g_1(t), \\
    &
    u(b_1,t) = g_2(t), \\
    &
    u(x,0) = h(x), \label{eq_diffu_2}
  \end{align}
\end{subequations}
where $u(x,t)$ is the field function to be solved for, $f(x,t)$
is a prescribed source term,
the constant $\nu$ is the diffusion coefficient,
$g_1(t)$ and $g_2(t)$ are the boundary conditions, and
$h(x)$ is the initial field distribution.
The values for the constant parameters involved in the above equations
and in the domain specification are,
\begin{equation*}
  a_1=0,\quad
  b_1=5, \quad
  \nu = 0.01, \quad
  t_f = 10 \ \text{or}\ 1.
\end{equation*}
We choose the source term $f$ such that the following function satisfies
equation \eqref{eq_diffu_1},
\begin{equation}\label{eq_diffu_3}
  u(x,t) = \left[2\cos\left(\pi x+\frac{\pi}{5}\right)
    + \frac32\cos\left(2\pi x - \frac{3\pi}{5} \right) \right]
  \left[2\cos\left(\pi t+\frac{\pi}{5}\right)
    + \frac32\cos\left(2\pi t - \frac{3\pi}{5}\right) \right].
\end{equation}
We choose the boundary conditions $g_1(t)$ and $g_2(t)$ and the initial
condition $h(x)$ according to equation \eqref{eq_diffu_3},
by restricting this expression to the corresponding boundaries
of the spatial-temporal domain.
Therefore, the function given by \eqref{eq_diffu_3} solves
the initial/boundary value problem represented by
equations \eqref{eq_diffu_1}--\eqref{eq_diffu_2}.

\begin{figure}
  \centerline{
    \includegraphics[height=2.3in]{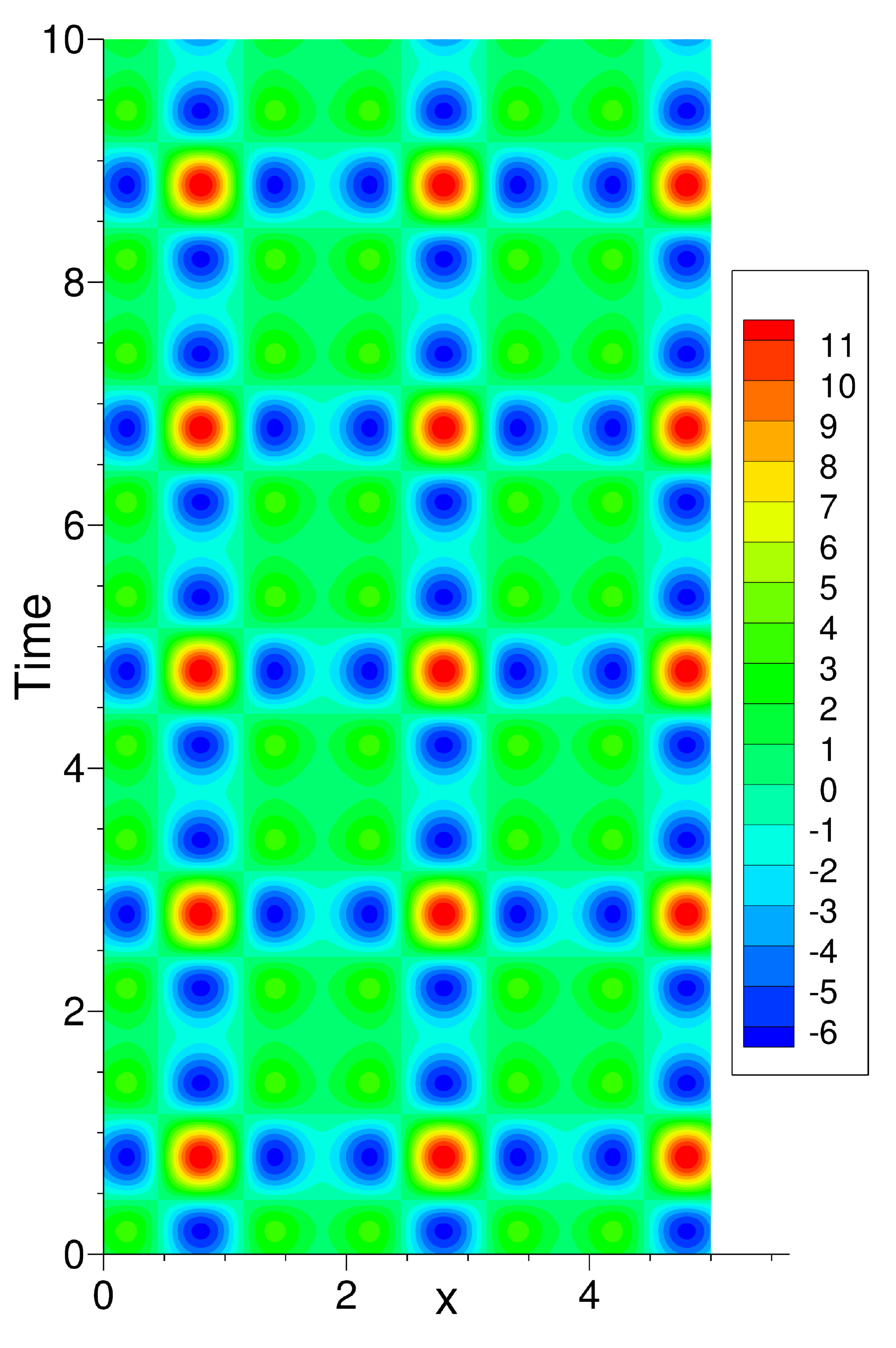}(a)
    \includegraphics[height=2.3in]{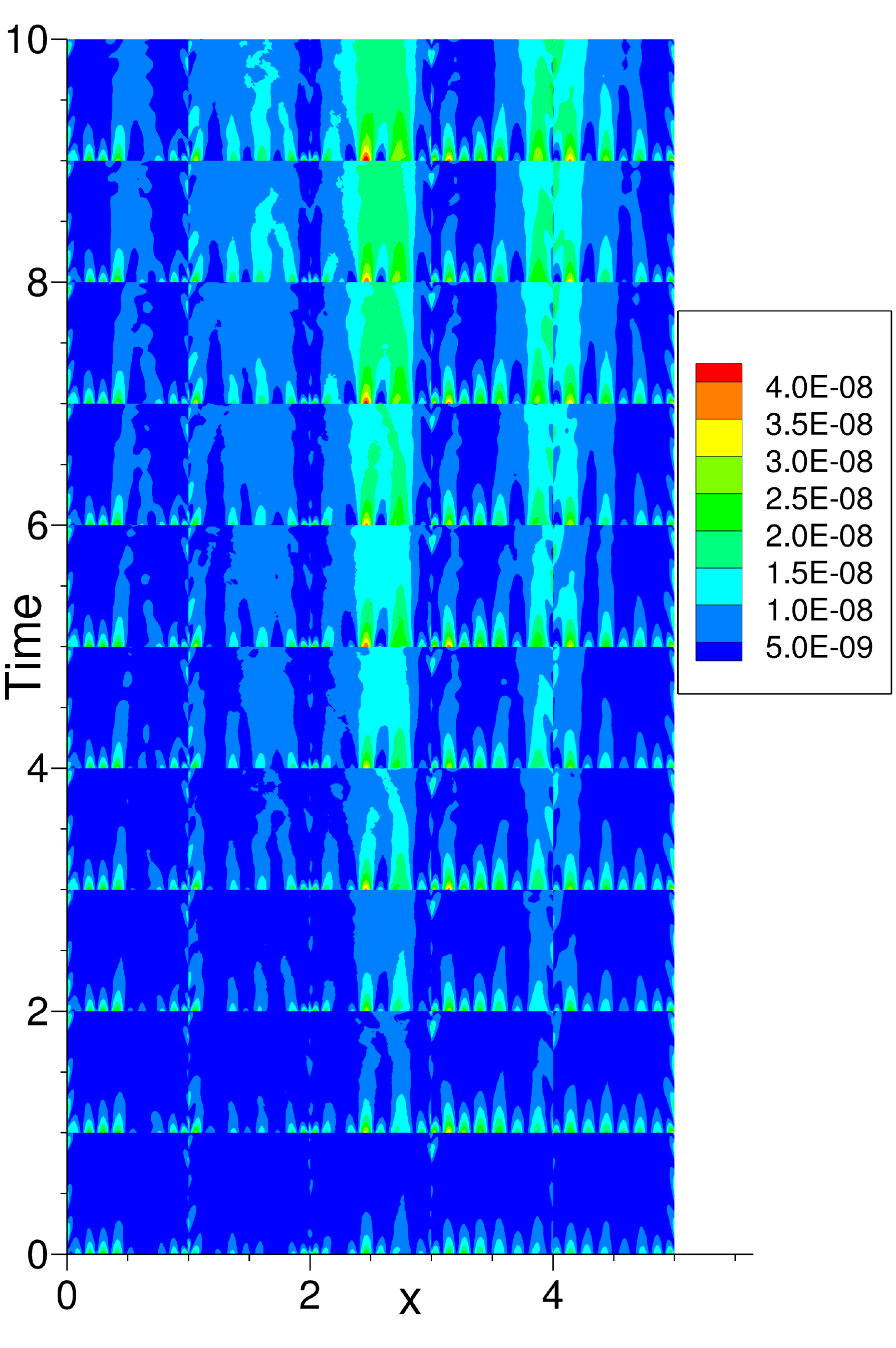}(b)
  }
  \caption{1D diffusion equation: distributions of the solution (a)
    and its absolute error (b) computed using the current locELM method.
    10 time blocks and 5 sub-domains per time block are employed.
  }
  \label{fg_diffu_1}
\end{figure}


We employ the locELM method together with block time marching from Section
\ref{sec:unsteady} to solve this initial/boundary value problem, by restricting
the method to one spatial dimension.
We partition the spatial-temporal domain $\Omega$ in time into $N_b$ uniform blocks,
and compute these time blocks individually and successively.
Within each time block, we further partition its spatial-temporal domain
into $N_x$ uniform sub-domains along the $x$ direction and $N_t$ uniform
sub-domains in time,
leading to $N_e=N_xN_t$ uniform sub-domains per time block.
We impose $C^1$ continuity conditions on the sub-domain boundaries in the $x$ direction
and $C^0$ continuity on the sub-domain boundaries in the temporal direction.
Within each sub-domain we use $Q_x$ uniform collocation points along the
$x$ direction and $Q_t$ uniform collocation points in time as the input training data,
leading to a total of $Q=Q_xQ_t$ uniform collocation points per sub-domain.

We use one local neural network to approximate the solution on each
sub-domain within the time block, thus leading to a total of $N_e$ local
neural networks in the locELM simulation. In the majority of
tests of this subsection the local neural networks each
contains a single hidden layer with $M$ nodes and the
$\tanh$  activation function. We also report results
obtained with the local neural networks containing more than one hidden layer.
The input layer of the local neural networks consists of two nodes,
representing $x$ and $t$. The output layer consists of a single node,
representing the solution $u$, and has no bias coefficient
and no activation function.
As in previous subsections, we incorporate an additional affine mapping operation
right behind the input layer of the local neural network
to normalize the input $x$ and $t$ data to the interval $[-1,1]\times[-1,1]$
in each sub-domain. The weight/bias coefficients in the hidden layer of
each of the local neural networks are set to uniform random values
generated on the interval $[-R_m,R_m]$.
We use a fixed seed value $22$ for the Tensorflow random number generator
for all the tests in this subsection.

The locELM simulation parameters  include
the number of time blocks ($N_b$), the number of sub-domains per
time block ($N_e, N_x, N_t$), the number of training parameters per
sub-domain ($M$), the number of collocation points per sub-domain ($Q_x$, $Q_t$, $Q$),
and the maximum magnitude of the random coefficients ($R_m$).
In accordance with previous subsections, we use $(Q,M)$ to characterize
the degrees of freedom within a sub-domain, and $(N_eQ,N_eM)$ to
characterize the degrees of freedom within a time block.

Figure \ref{fg_diffu_1} shows distributions of the locELM solution
and its absolute error in the spatial-temporal plane.
In this test the temporal domain size is set to $t_f=10$.
We have employed $N_b=10$ uniform time blocks in the simulation,
$N_e=5$ uniform sub-domains per time block (with $N_x=5$ and $N_t=1$),
$Q=30\times 30$ uniform collocation points in each sub-domain
(with $Q_x=Q_t=30$), $M=300$ training parameters per sub-domain,
a single hidden layer in the local neural networks,
and $R_m=1.0$ when generating the random weight/bias coefficients for the
hidden layers of the local neural networks.
It is evident that the locELM method has captured the solution
very accurately, with the absolute error on the order of $10^{-9}\sim 10^{-8}$.

\begin{figure}
  \centerline{
    \includegraphics[width=2in]{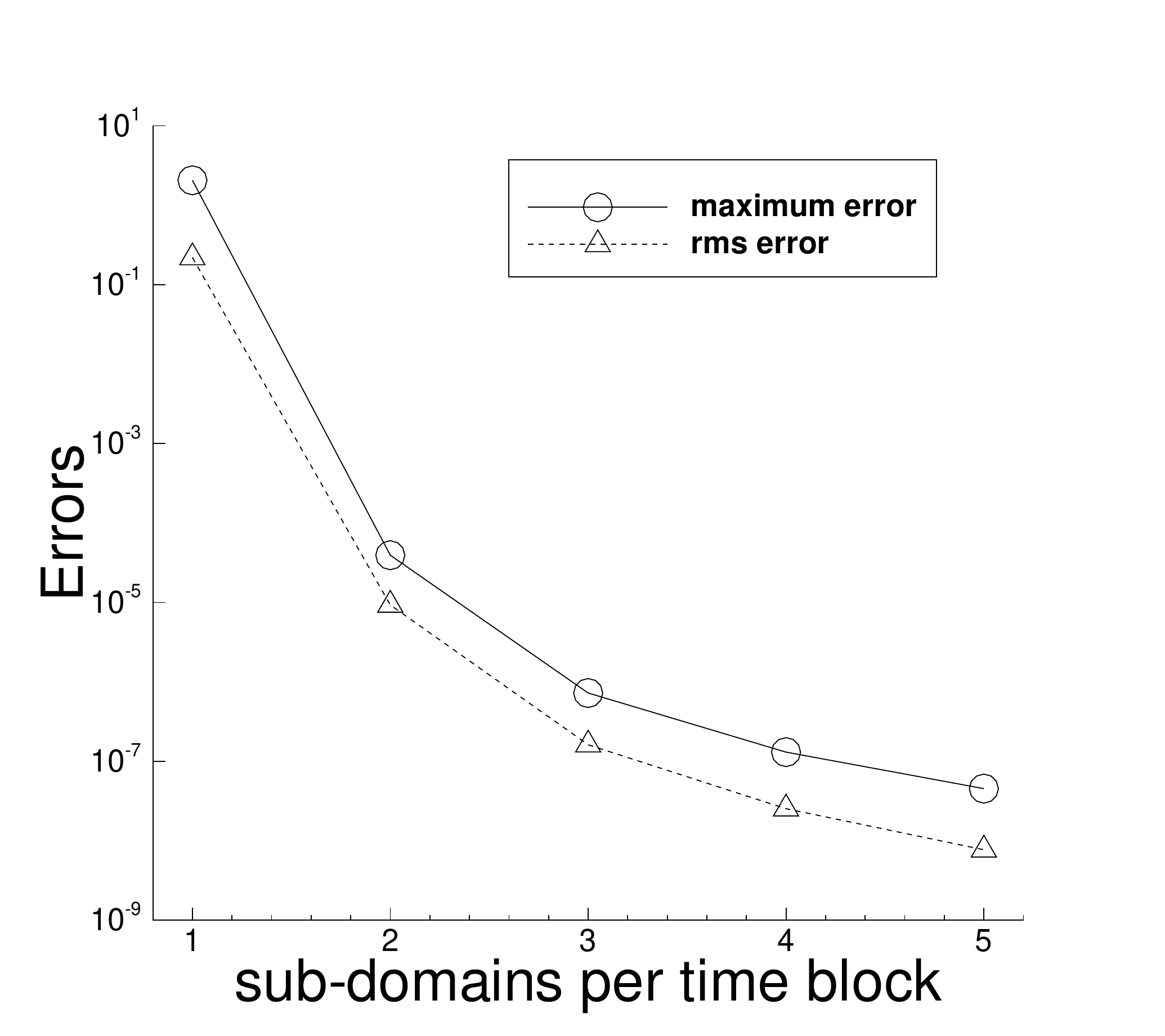}(a)
    \includegraphics[width=2in]{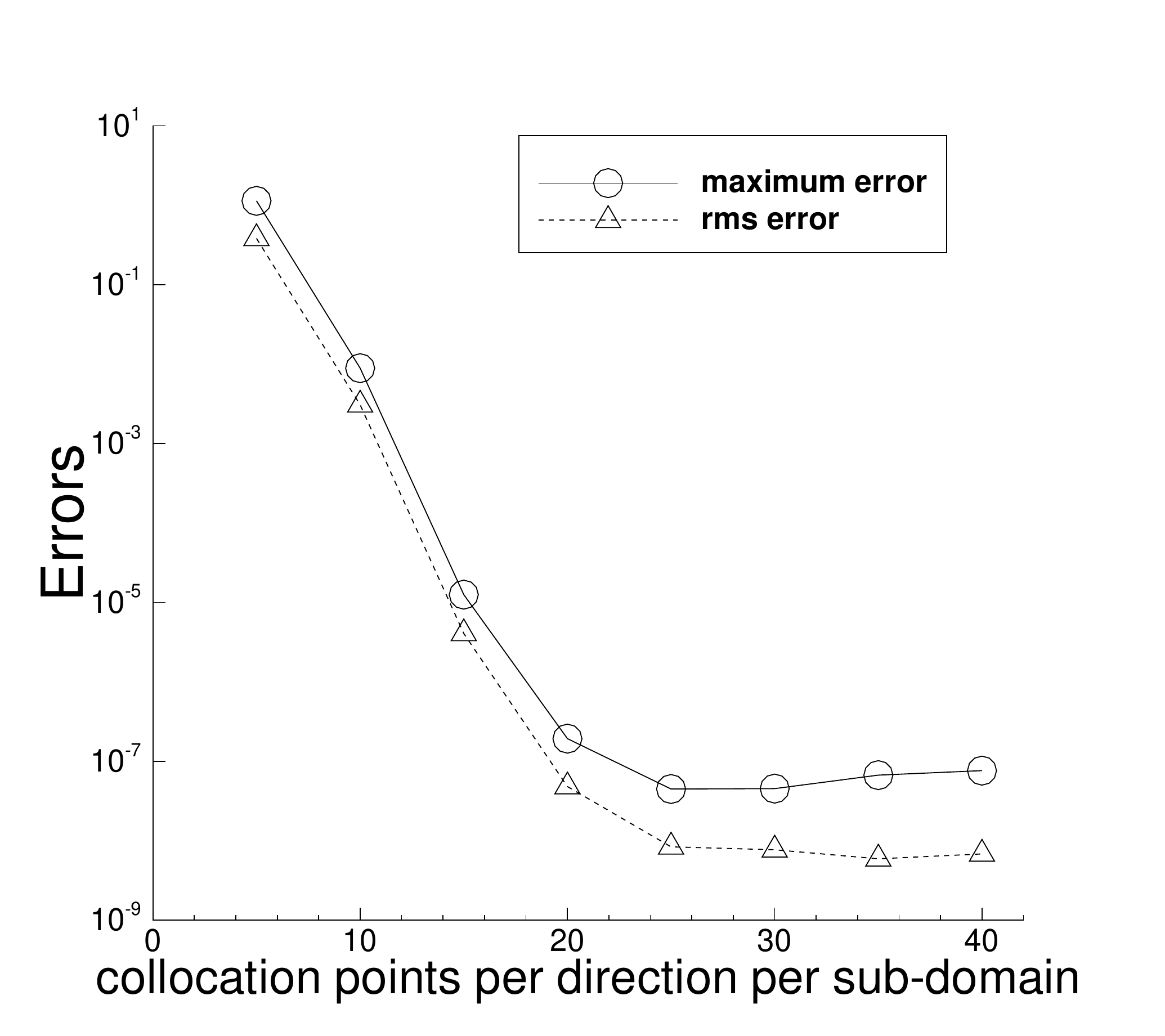}(b)
    \includegraphics[width=2in]{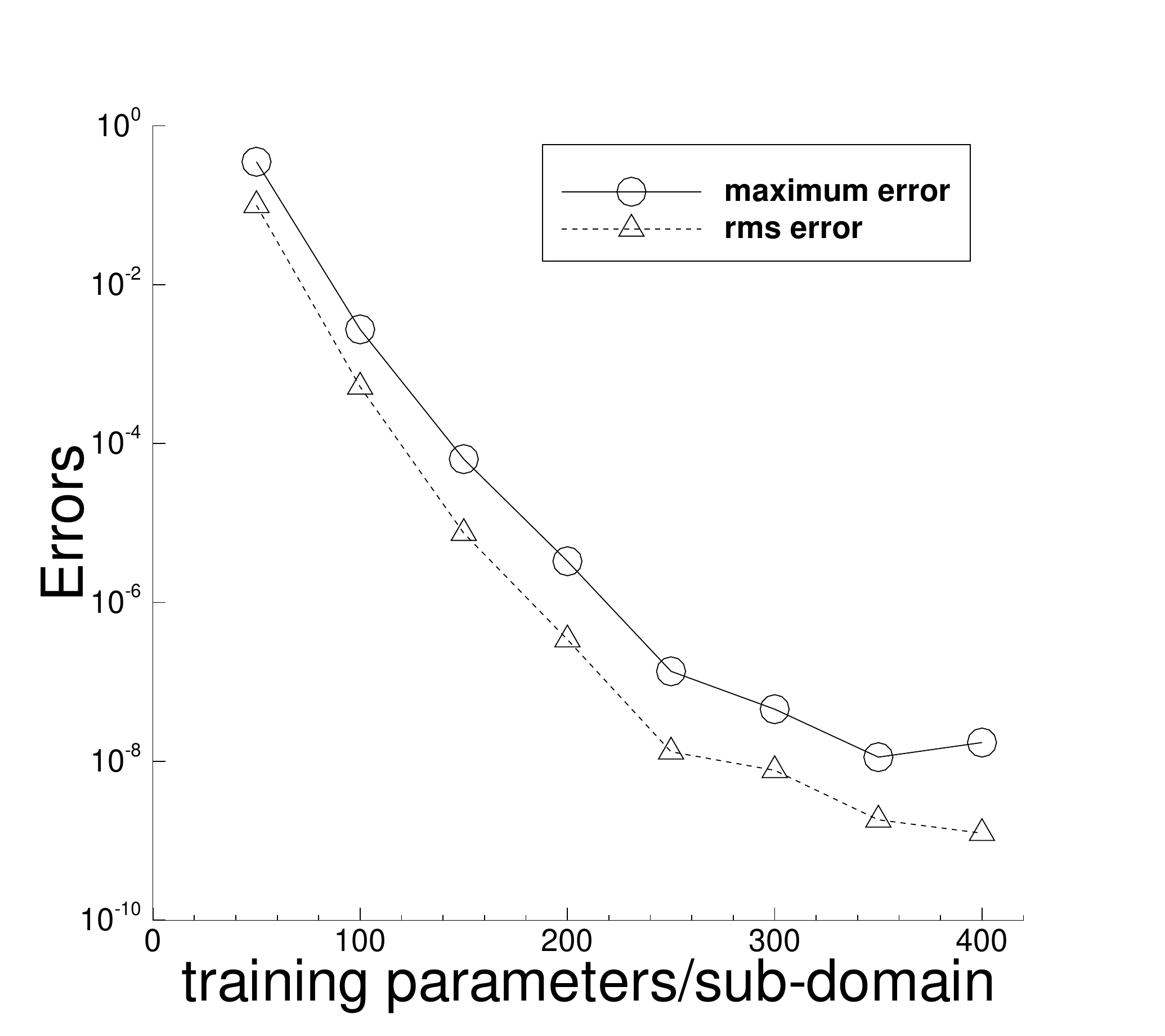}(c)
  }
  \caption{Effect of the degrees of freedom on simulation accuracy
    (1D diffusion equation): the maximum and rms errors in the domain as
    a function of (a) the number of sub-domains in each time block,
    (b) the number of collocation points in each direction
    in each sub-domain, and (c) the number of training parameters in each
    sub-domain.
    Temporal domain size is $t_f=10$ and $10$ uniform time blocks are used.
  }
  \label{fg_diffu_3}
\end{figure}

The effect of the degrees of freedom on the simulation accuracy
is illustrated by Figure \ref{fg_diffu_3}.
In this group of tests, the temporal domain size is set to $t_f=10$,
and we have employed $N_b=10$ uniform time blocks in the spatial-temporal domain,
a single hidden layer in each local neural network, 
and $R_m=1.0$ when generating the random coefficients for the hidden layers
of the local neural networks.
The number of sub-domains in each time block, or the number of collocation
points per sub-domain, or the number of training parameters
per sub-domain has been varied in the tests.

Figure \ref{fg_diffu_3}(a) illustrates the effect of the number of sub-domains
within each time block, while the degrees of freedom per sub-domain are fixed.
Here we fix the number of uniform collocation points per sub-domain at $Q=30\times 30$
($Q_x=Q_t=30$)
and the number of training parameters per sub-domain at $M=300$, and then
vary the number of uniform sub-domains per time block systematically.
This plot shows the maximum and rms errors of the locELM
solution in the overall spatial-temporal domain as a function of the
number of sub-domains per time block in the simulations.
With increasing number of sub-domains, the numerical errors are observed to decrease
dramatically, from around $10^{-1}$ with one sub-domain/time-block
to around $10^{-8}$ with $5$ sub-domains/time-block.

Figure  \ref{fg_diffu_3}(b) illustrates the effect of the number of
collocation points per sub-domain on the simulation accuracy.
Here we use $N_e=5$ uniform sub-domains (with $N_x=5$ and $N_t=1$) in each time block,
fix the the number of training parameters per sub-domain
at $M=300$,
and vary the number of collocation points per sub-domain while maintaining
$Q_x=Q_t$.
This plot shows the maximum and rms errors
in the overall spatial-temporal domain as a function of the number of
collocation points in each direction in each sub-domain.
The numerical errors can be observed to initially
decrease exponentially with increasing
number of collocation points per direction when it is below about $20$,
and then stagnate at a level around $10^{-8}\sim 10^{-7}$ as the number
of collocation points per direction further increases.

Figure \ref{fg_diffu_3}(c) illustrates the effect of the number of training
parameters on the simulation accuracy.
Here we use $N_e=5$ sub-domains (with $N_x=5$ and $N_t=1$) in each time block,
fix the number of collocation points per sub-domain at
$Q=30\times 30$ ($Q_x=Q_t=30$), and vary the number of training parameters per
sub-domain.
The plot shows the maximum/rms errors in the overall domain
as a function of the number of training parameters per sub-domain.
One can observe that the  errors initially decrease exponentially
with increasing number of training parameters/sub-domain when
it is below about $250$, and then the error reduction slows down
as the number of training parameters/sub-domain further increases.
These behaviors are  consistent with what have been observed with
other problems in previous subsections.


\begin{figure}
  \centerline{
    \includegraphics[height=2.4in]{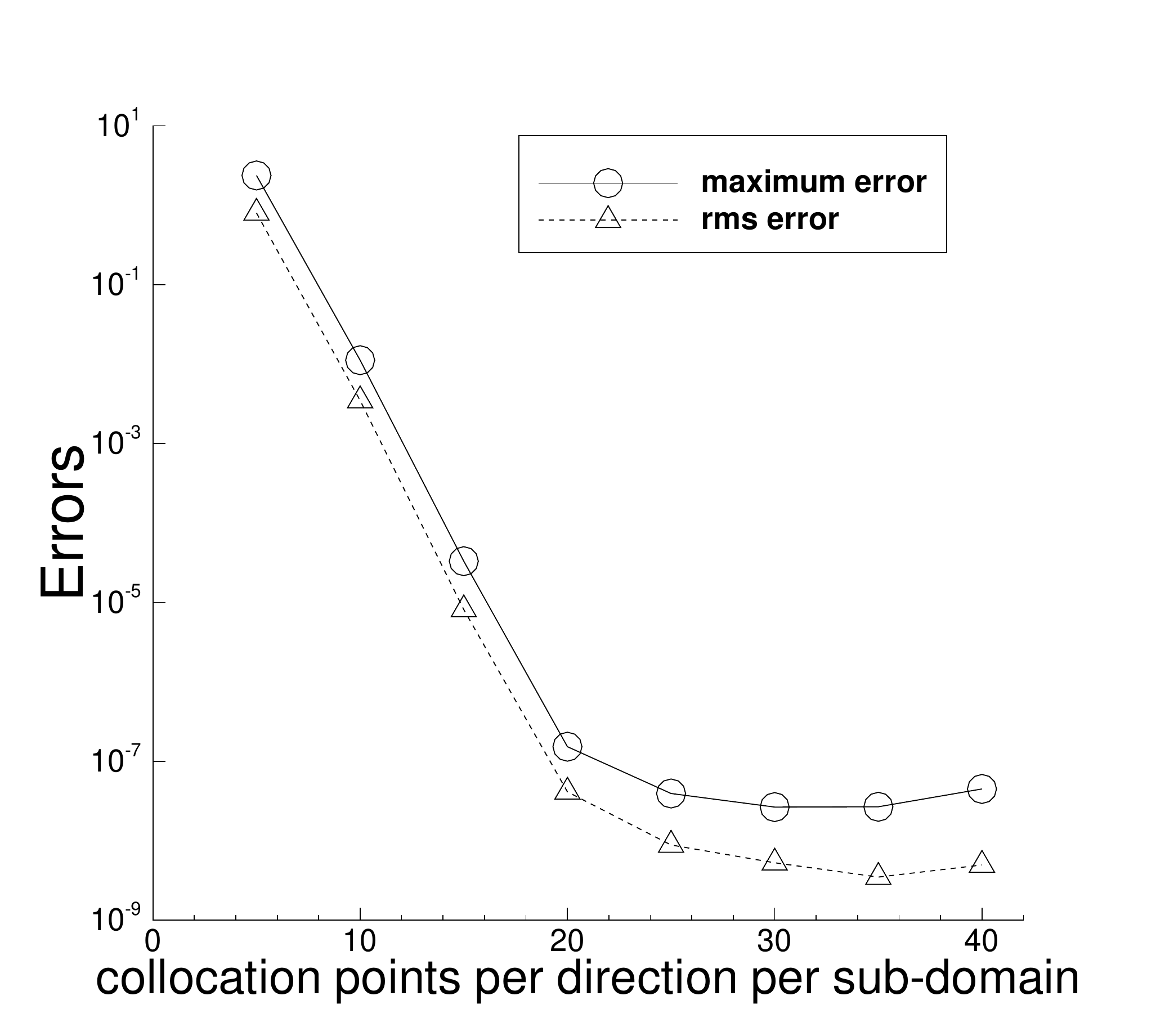}(b)
  }
  \caption{Results obtained with two hidden layers in
    local neural networks (1D diffusion
    equation): 
    The maximum/rms errors in the domain versus
    the number of collocation points in each direction per sub-domain.
  }
  \label{fg_diffu_a}
\end{figure}

In Figure \ref{fg_diffu_a} we show results obtained with local neural networks
containing more than one hidden layer.
In this group of simulations, each local neural network contains two
hidden layers, with $30$ and $300$ nodes in these two layers, respectively.
The activation function is $\tanh$ in both hidden layers.
The temporal domain size is $t_f=10$, and $10$ uniform time blocks have been
used. We employ $N_e=5$ sub-domains per time block ($N_x=5$ and $N_t=1$),
$M=300$ training parameters per sub-domain (width of the last hidden layer),
and $R_m=0.5$ when generating the random weight/bias coefficients for
the hidden layers of the local neural networks.
The number of collocation points per sub-domain ($Q$) is varied systematically
while $Q_x=Q_t$ is maintained.
Figure \ref{fg_diffu_a}(a) shows the distribution of the absolute
error of the locELM solution in the spatial-temporal plane,
obtained with $Q=30\times 30$ uniform collocation points per sub-domain.
This figure can be compared with Figure \ref{fg_diffu_1}(b),
which corresponds to the same simulation resolution but is obtained with
local neural networks containing a single hidden layer.
Figure \ref{fg_diffu_a}(b) shows the maximum and rms errors in the overall domain
as a function of the number of collocation points in each direction
per sub-domain. This figure can be compared with
Figure \ref{fg_diffu_3}(b), which corresponds to a single hidden layer
in the local neural networks.
It is evident that the solution has been captured accurately by the current locELM
method using two hidden layers in the local neural networks.
The results shown here and those results in Section \ref{sec:helm1d}
(see Figures \ref{fig:helm1d_4} and \ref{fig:helm1d_5})
demonstrate that the current locELM method,
using local neural networks with a small number of (more than one) hidden layers,
is able to produce accurate simulation results.


\begin{figure}
  \centerline{
    \includegraphics[width=2.2in]{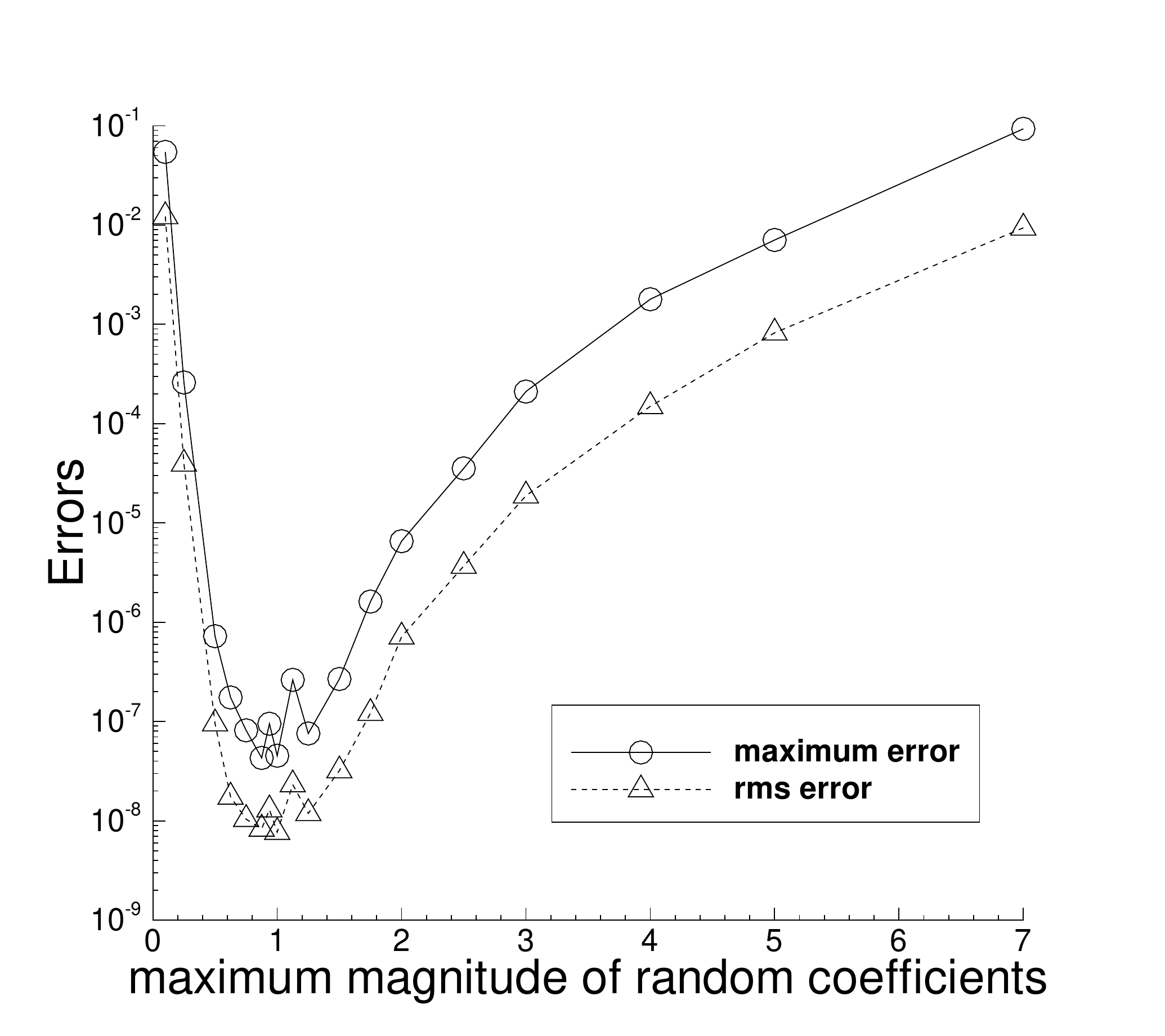}
  }
  \caption{Effect of the random coefficients in local
    neural networks
    (1D diffusion equation): the maximum and rms errors in the domain as a
    function of $R_m$, the maximum magnitude of the random coefficients.
  }
  \label{fg_diffu_4}
\end{figure}

Figure \ref{fg_diffu_4} illustrates  the effect of the random weight/bias
coefficients in the local neural networks on
the simulation accuracy.
It shows the maximum and rms errors of the locELM solution in
the overall domain as a function of $R_m$, the maximum magnitude of the
random coefficients.
In this group of tests, the temporal domain size is $t_f=10$,
and we have employed $N_b=10$ uniform time blocks in the domain,
$N_e=5$ uniform sub-domains per time block ($N_x=5$, $N_t=1$),
$Q=30\times 30$ uniform collocation points per sub-domain,
$M=300$ training parameters per sub-domain, and a single hidden layer
in the local neural networks.
The random coefficients in the hidden layers are generated on $[-R_m,R_m]$,
and $R_m$ is varied systematically in these tests.
We observe a similar behavior to those in previous subsections.
A better accuracy can be attained with a range of moderate $R_m$ values,
while very large or very small $R_m$ values tend to produce less accurate
results.

\begin{figure}
  \centerline{
    \includegraphics[width=2in]{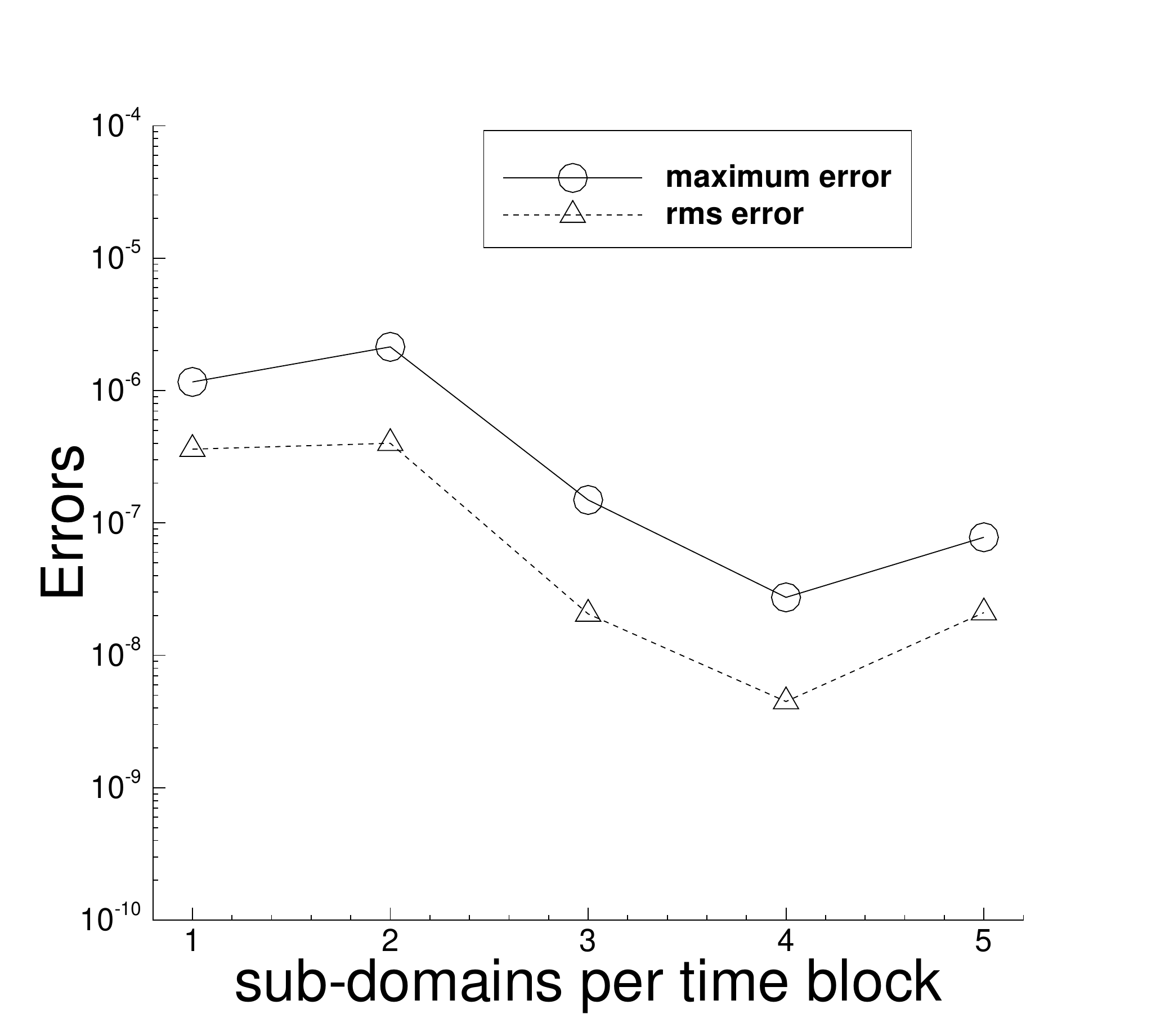}(a)
    \includegraphics[width=2in]{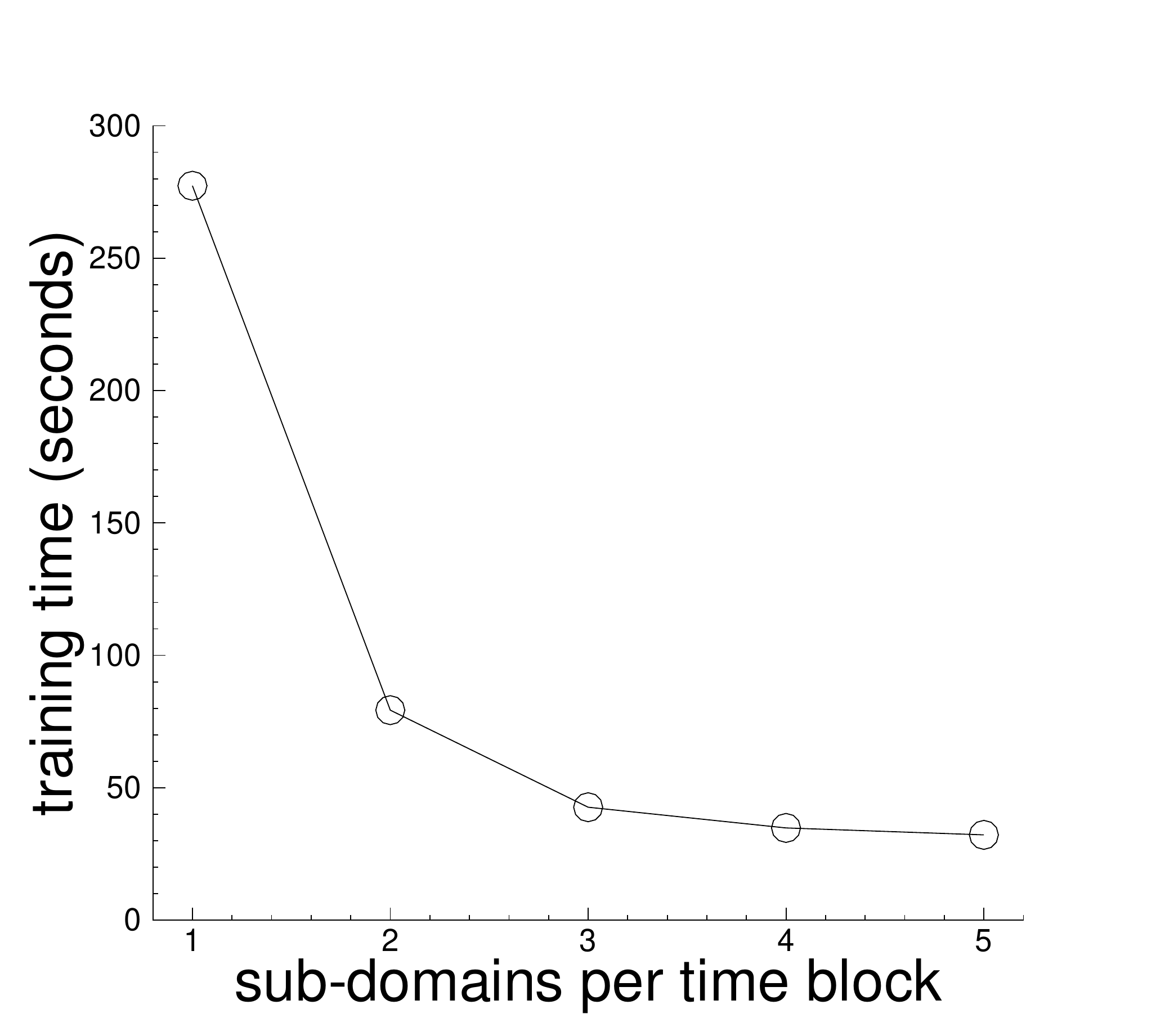}(a)
  }
  \caption{Effect of the number of sub-domains, with fixed total degrees of freedom in
    the domain (1D diffusion equation): (a) the maximum and rms errors in the domain, and
    (b) the training time, as a function of the number of sub-domains per time
    block.
  }
  \label{fg_diffu_5}
\end{figure}

Figure \ref{fg_diffu_5} depicts a study of the effect of the number of sub-domains
in the simulation on the simulation accuracy and on the network training time,
while the total number of degrees of freedom
in the domain is fixed.
In this group of tests, the temporal domain size is $t_f=10$,
and we have used $N_b=10$ time blocks in the overall spatial-temporal domain.
The number of uniform sub-domains per time block is varied systematically
between $N_e=1$ and $N_e=5$ in the simulations, implemented by fixing
$N_t=1$ and varying $N_x$ between $1$ and $5$.
We set the number of training parameters per sub-domain ($M$), and the number of
uniform collocation points per sub-domain ($Q$, with $Q_x=Q_t$),
in a way such that the total number of training parameters per time block is
fixed at $N_eM = 1500$ and the total number of collocation points per
time block is approximately fixed at $N_eQ\approx 2500$.
Specifically, $M$ and $Q$ in different cases  are:
$M=1500$ and $Q=50\times 50$ for $1$ sub-domain per time block,
$M=750$ and $Q=35\times 35$ for $2$ uniform sub-domains per time block,
$M=500$ and $Q=29\times 29$ for $3$ uniform sub-domains per time block,
$M=375$ and $Q=25\times 25$ for $4$ uniform sub-domains per time block,
and $M=300$ and $Q=22\times 22$ for $5$ uniform sub-domains per time block.
We employ $R_m=3.0$ when generating the random coefficients
for the case with one sub-domain per time block,
which is approximately at the optimal range of $R_m$ values for this case.
We employ $R_m=1.0$ when generating the random coefficients
for the rest of the cases with $N_e=2 \sim 5$ sub-domains per time block. 
Note that the case
with one sub-domain per time block is equivalent to the configuration of
a global ELM in the simulation.
Figure \ref{fg_diffu_5}(a) shows a comparison of the maximum and rms errors
in the overall spatial-temporal domain as a function of the number of sub-domains
per time block in the simulations.
It can be observed that the numerical errors with $2$ or more sub-domains
are comparable to or smaller than the errors corresponding to one sub-domain
in the simulations.
Figure \ref{fg_diffu_5}(b) shows the neural-network training time
as a function of the number of sub-domains per time block.
One can observe that the training time decreases significantly
with increasing number of
sub-domains. Compared with the case of one sub-domain per time block,
the training time corresponding to $2$ and more sub-domains in the simulations
has been considerably reduced, e.g.~$277$ seconds with one sub-domain versus $79$
seconds with $2$ sub-domains.
These results confirm and reinforce our observations with the other problems
that, compared with global ELM, the use of domain decomposition and locELM
with multiple sub-domains can significantly reduce the network training time, and
hence the computational cost,
while attaining the same or sometimes
even better accuracy in the simulation results.

\begin{figure}
  \centerline{
    \includegraphics[width=2.5in]{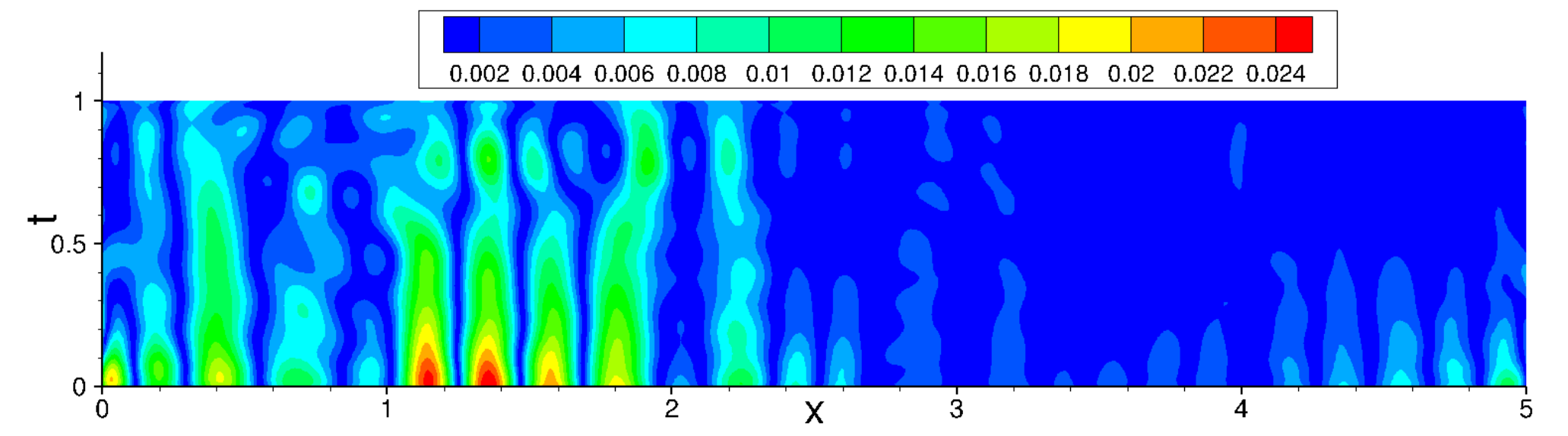}(a)
  }
  \centerline{
    \includegraphics[width=2.5in]{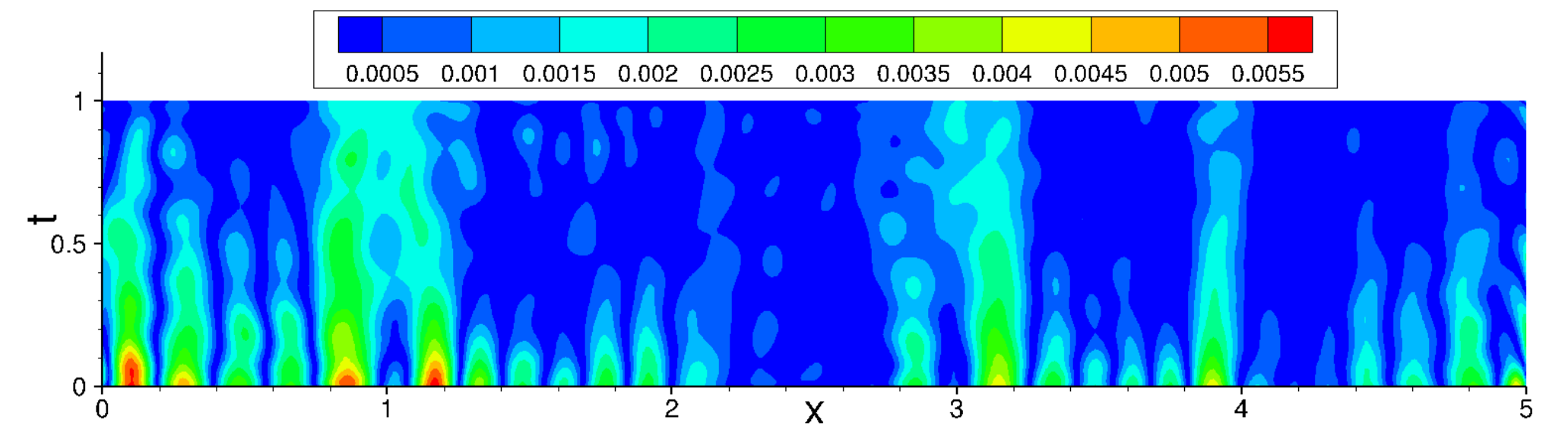}(b)
  }
  \centerline{
    \includegraphics[width=2.5in]{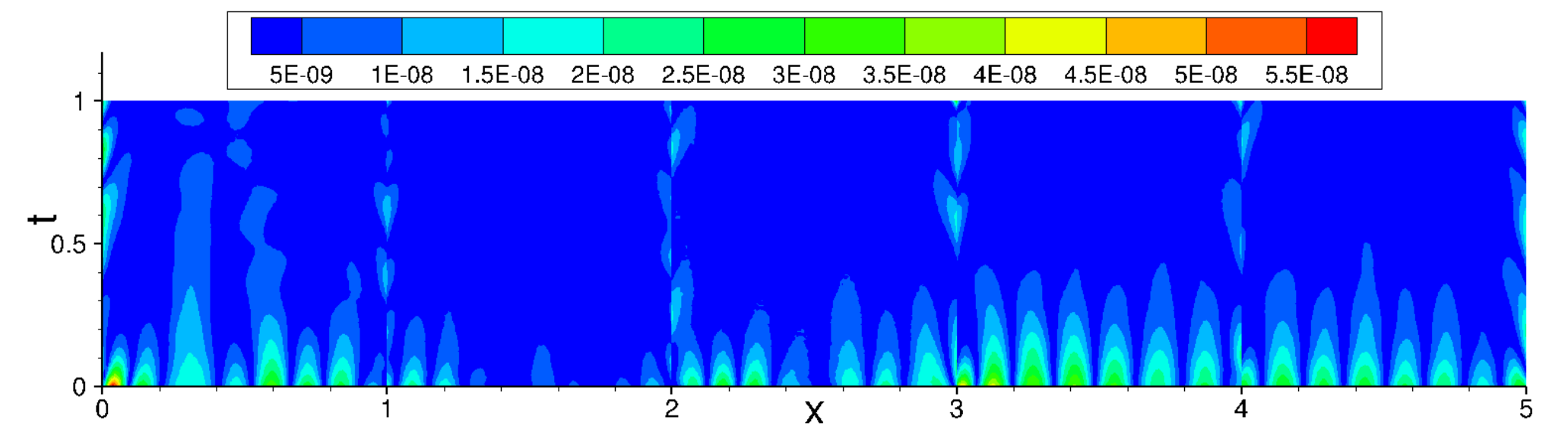}(c)
  }
  \caption{Comparison between locELM and DGM
    (1D diffusion equation): Distributions of the
     absolute errors  computed
    using DGM with the Adam optimizer (a) and the L-BFGS
    optimizer (b) and using the current locELM method (c).
  }
  \label{fg_diffu_6}
\end{figure}

\begin{figure}
  \centerline{
    \includegraphics[width=2.2in]{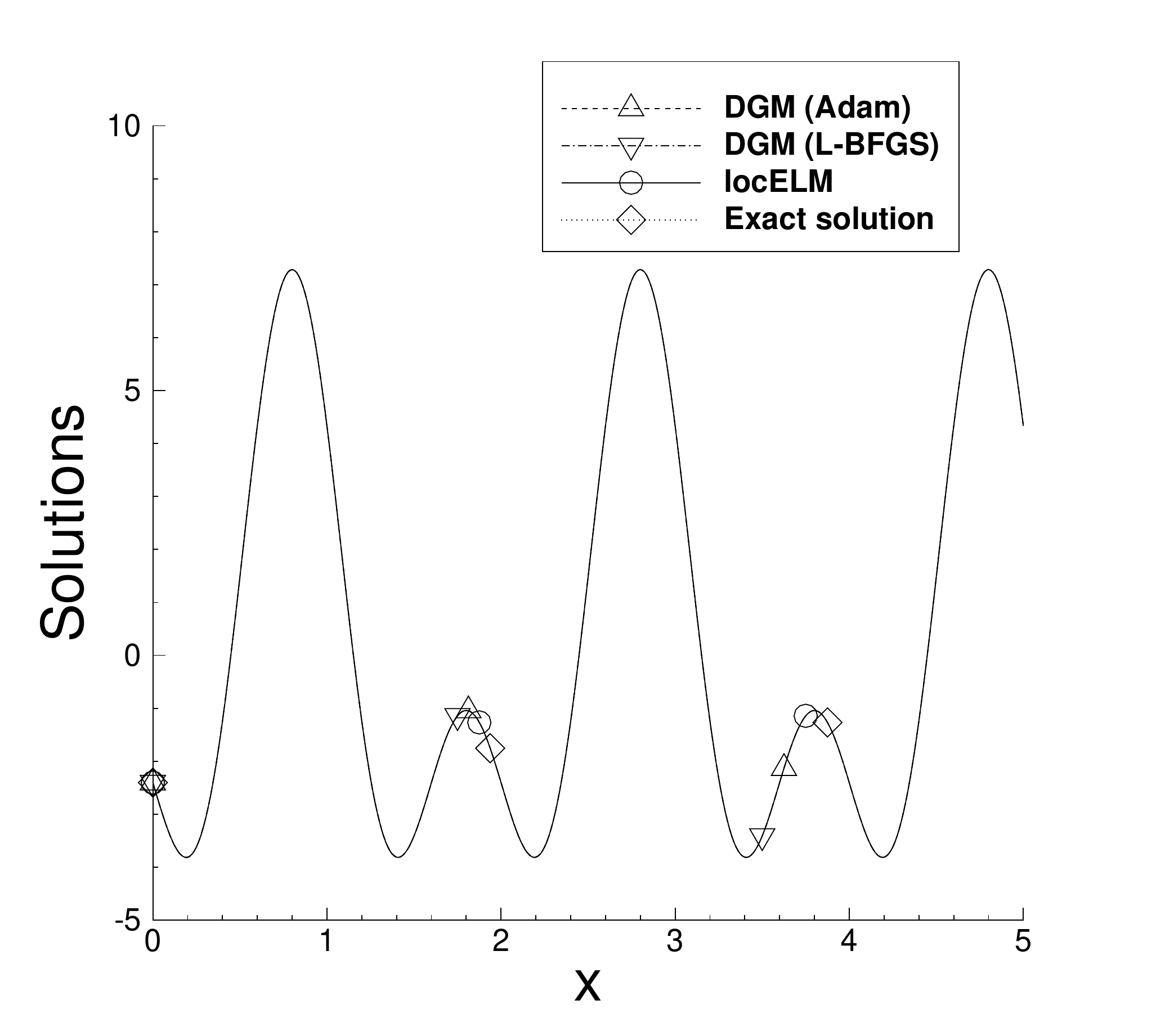}(a)
    \includegraphics[width=2.2in]{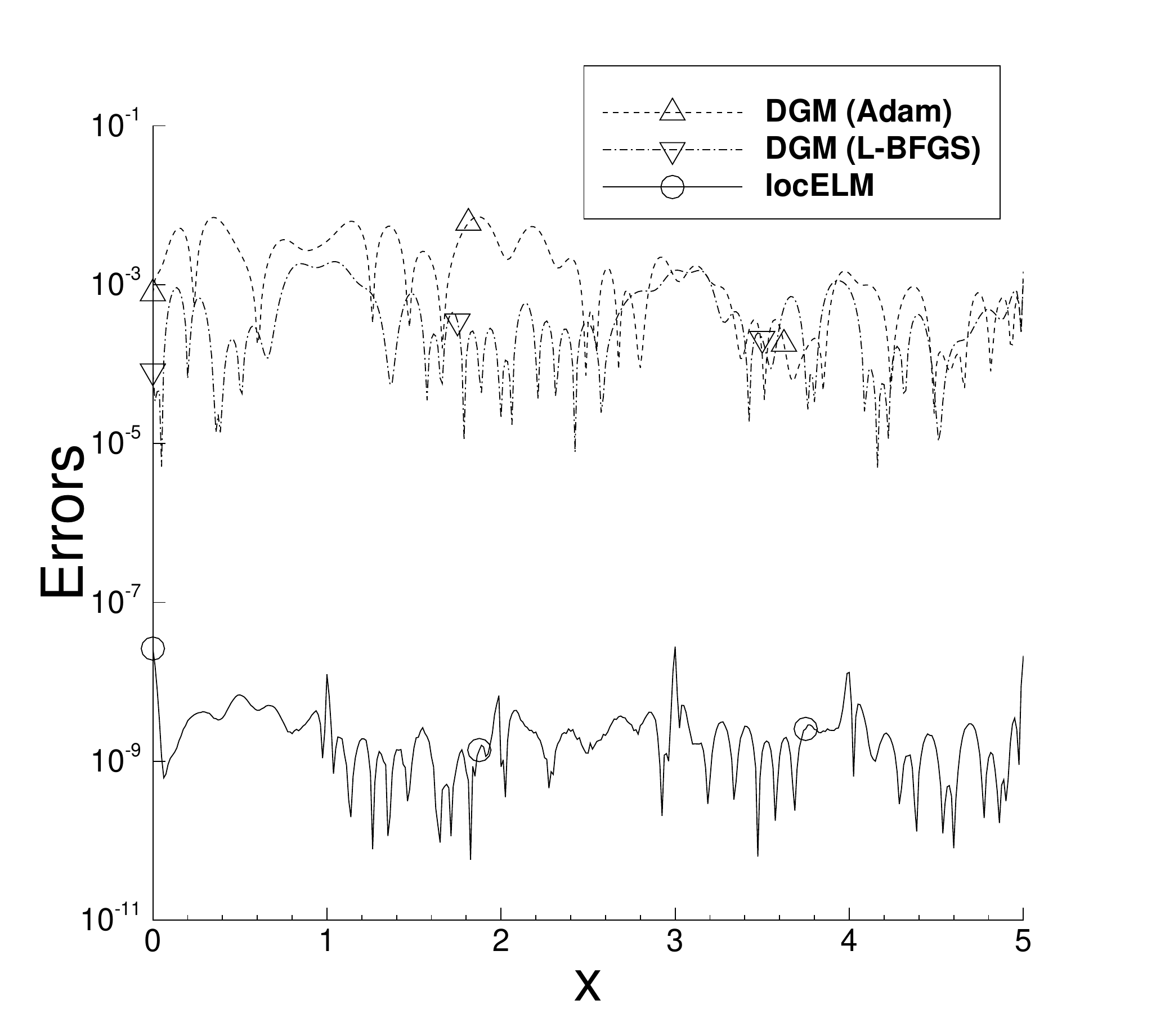}(b)
  }
  \caption{Comparison between locELM and DGM
    (1D diffusion equation): Profiles of the solutions (a)
    and their absolute errors (b) at $t=1.0$
    obtained using DGM (Adam/L-BFGS optimizers)
    and using the current locELM method.
    The problem settings and the simulation parameters
    correspond to those of Figure \ref{fg_diffu_6}.
  }
  \label{fg_diffu_7}
\end{figure}

Let us now compare the current locELM method and the deep Galerkin method (DGM)
for solving the 1D diffusion equation. 
Figure \ref{fg_diffu_6} compares distributions of the solutions and their
absolute errors obtained using DGM with the Adam optimizer and the L-BFGS
optimizer and using the current locELM method.
The temporal domain size is set to $t_f=1$ in these tests.
With DGM, the neural network consists of 4 hidden layers, with a width of
$40$ nodes and the $\tanh$ activation function in each layer.
When computing the loss function of the network,
we have divided the domain  into
$5$ uniform sub-regions along the $x$ direction, and computed the
residual norm integral by the Gaussian quadrature rule
on $20\times 20$ Gauss-Lobatto-Legendre quadrature points in each sub-region.
With the Adam optimizer, the neural network has been trained for
$135,000$ epochs, with the learning rate decreasing gradually from $0.001$
at the beginning to $2.5\times 10^{-6}$ at the end of the training.
With the L-BFGS optimizer, the neural network has been trained for
$36,000$ L-BFGS iterations.
In the simulation with the current locELM method, we
employ $N_b=1$ time block in the spatial-temporal domain,
$N_e=5$ sub-domains (with $N_x=5$, $N_t=1$) per time block,
$Q=30\times 30$ uniform collocation points per sub-domain,
$M=300$ training parameters per sub-domain, $1$ hidden layer in
each of the local neural networks, and
$R_m=1.0$ when generating the random weight/bias coefficients.
The DGM has captured the solution reasonably well. But its error levels
are considerably higher, by about five orders of magnitude,
than that of the current locELM method ($10^{-3}$ versus $10^{-8}$).

A comparison of the solution and the error profiles between locELM and DGM
is provided in Figure \ref{fg_diffu_7}.
Figure \ref{fg_diffu_7}(a) compares the solution profiles at $t=1.0$
obtained using DGM (Adam/L-BFGS optimizers) and using the
current locELM method, together with that of the exact solution.
The settings and the parameters correspond to those of Figure \ref{fg_diffu_6}.
The computed profiles all agree with the exact solution quite well.
Figure \ref{fg_diffu_7}(b) compares profiles of the absolute errors at $t=1.0$
obtained with DGM and the current method.
The numerical error of the current method, which is at a level around $10^{-9}$,
is significantly smaller than those from DGM, which are at
a level around $10^{-4}$.

\begin{table}
  \centering
  \begin{tabular}{lllll}
    \hline
    method & maximum error & rms error & epochs/iterations & training time (seconds)\\
    DGM (Adam) & $2.59e-2$ & $3.84e-3$ & $135,000$ & $4194.5$ \\
    DGM (L-BFGS) & $5.82e-3$ & $8.21e-4$ & $36,000$ & $3201.4$ \\
    locELM & $5.82e-8$ & $6.25e-9$ & $0$ & $28.4$ \\
    \hline
  \end{tabular}
  \caption{1D diffusion equation: comparison between DGM (Adam/L-BFGS
    optimizers) and locELM.
    The settings and parameters correspond to those of Figure \ref{fg_diffu_6}.
  }
  \label{tab_diffu_8}
\end{table}

In Table \ref{tab_diffu_8} we provide some further comparisons between locELM and
DGM  in terms of the accuracy and the computational cost.
Here we list the maximum and rms errors in the overall spatial-temporal domain,
the number of epochs or iterations in network training, and the training time
of DGM with the Adam and L-BFGS optimizers and of the current locELM
method. The problem settings and the simulation parameters correspond to
those of Figure \ref{fg_diffu_6}.
The data demonstrate the clear superiority of locELM to DGM,
with the locELM errors five orders of magnitude smaller and
the training time over two orders of magnitude less.


\begin{figure}
  \centerline{
    \includegraphics[width=2in]{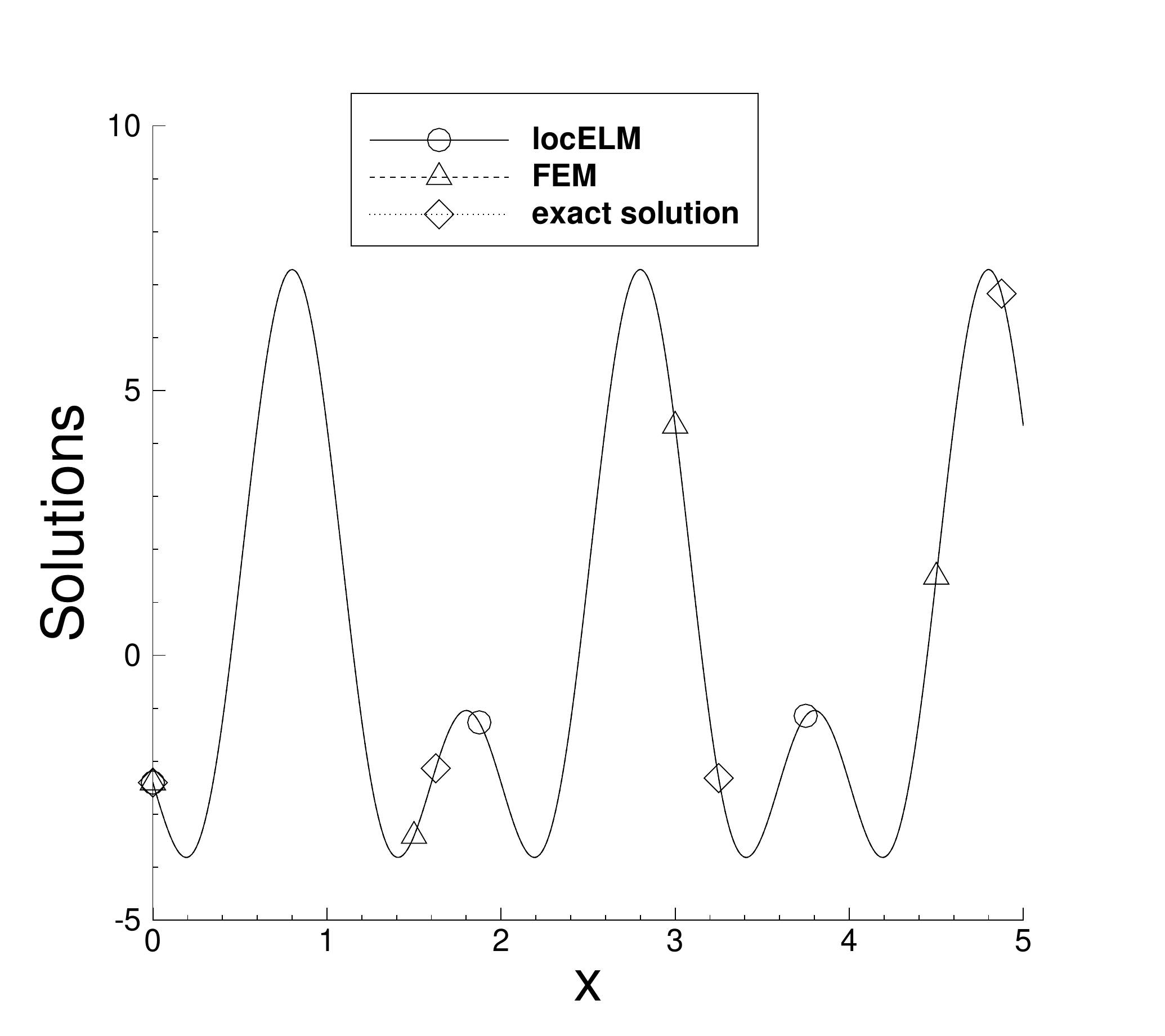}(a)
    \includegraphics[width=2in]{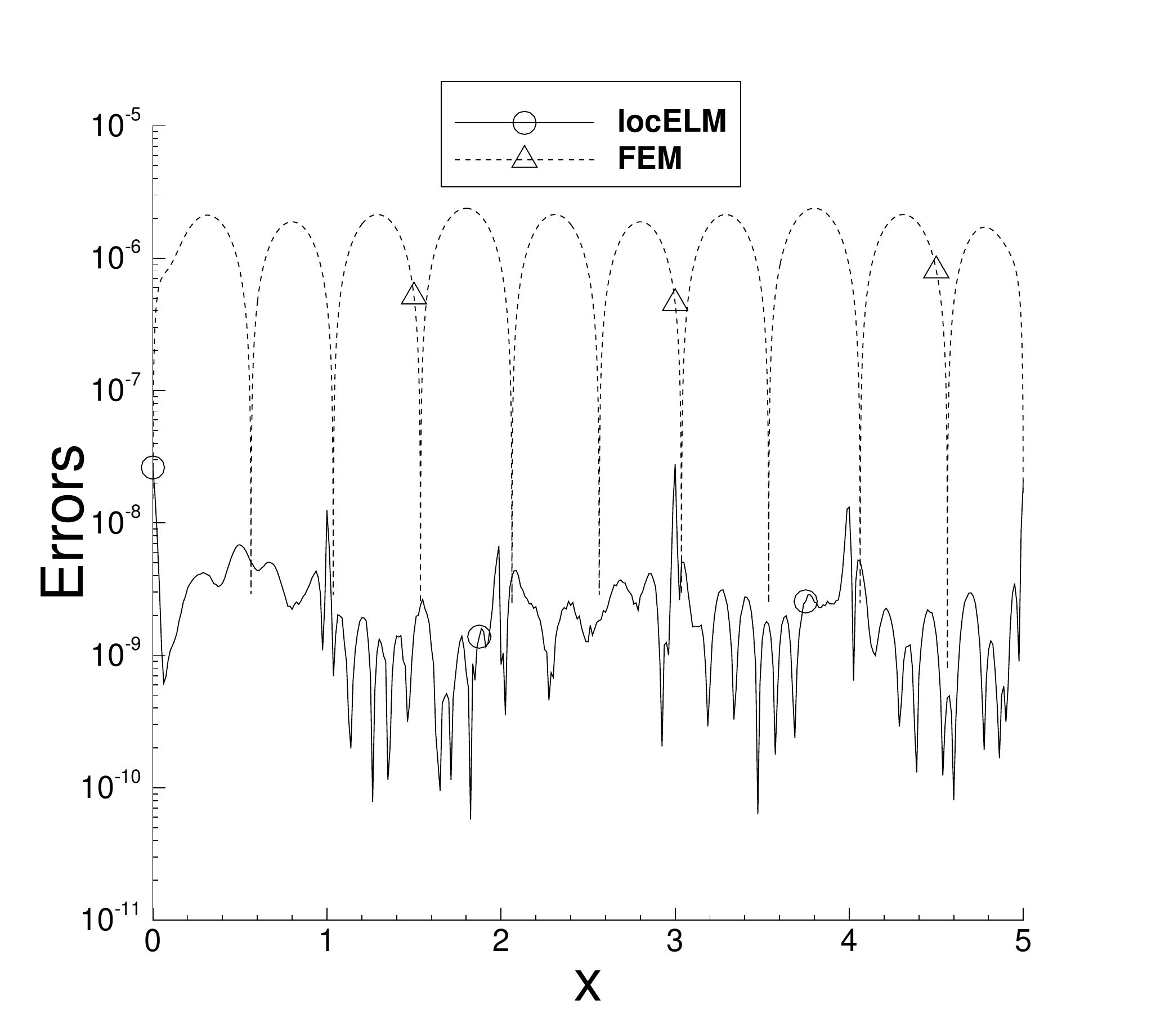}(b)
    \includegraphics[width=2in]{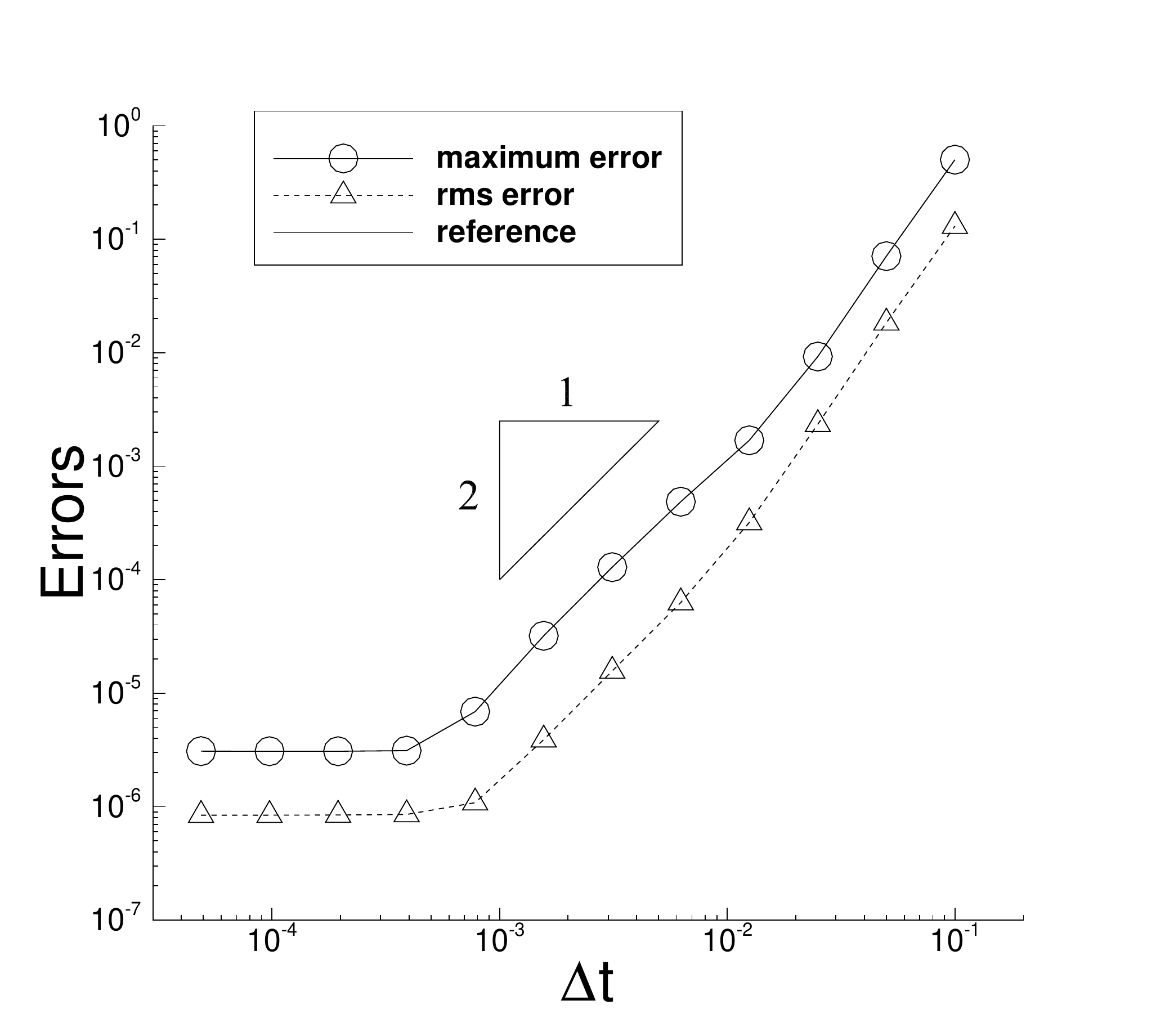}(c)
  }
  \caption{Comparison between locELM and FEM (1D diffusion equation):
    Profiles of (a) the solutions and (b) their absolute errors at $t=1.0$, computed
    using the current locELM method and using the finite element
    method (FEM).
    (c) The FEM maximum and rms
    errors at $t=0.5$ versus $\Delta t$, showing the temporal second-order
    convergence rate of FEM.
  }
  \label{fg_diffu_8}
\end{figure}

\begin{table}
  \centering
  \begin{tabular}{lllllllll}
    \hline
    method & $\Delta t$ & elements & sub-domains & $Q$ & $M$ & maximum & rms  & wall time\\
     & & & & & & error & error & (seconds) \\
    \hline
    locELM & -- & -- & $5$ & $20\times 20$ & $200$ & $2.48e-6$ & $2.23e-7$ & $7.9$ \\
     & -- & -- & $5$ & $20\times 20$ & $250$ & $8.97e-8$ & $2.25e-8$ & $11.3$ \\
     & -- & -- & $5$ & $30\times 30$ & $300$ & $5.82e-8$ & $6.25e-9$ & $28.4$ \\
    \hline
    FEM
    & $0.002$ & $2000$ & -- & -- & -- & $2.42e-4$ & $4.40e-5$ & $5.9$ \\
    & $0.001$ & $2000$ & -- & -- & -- & $9.82e-5$ & $2.01e-5$ & $12.0$ \\
    & $0.0005$ & $2000$ & -- & -- & -- & $1.54e-4$ & $2.61e-5$ & $24.0$ \\
    & $0.00025$ & $2000$ & -- & -- & -- & $1.72e-4$ & $2.85e-5$ & $48.3$ \\
    \cline{2-9}
    & $0.002$ & $5000$ & -- & -- & -- & $3.63e-4$ & $5.98e-5$ & $12.3$ \\
    & $0.001$ & $5000$ & -- & -- & -- & $6.99e-5$ & $1.22e-5$ & $24.6$ \\
    & $0.0005$ & $5000$ & -- & -- & -- & $1.69e-5$ & $3.43e-6$ & $48.8$ \\
    & $0.00025$ & $5000$ & -- & -- & -- & $2.26e-5$ & $3.91e-6$ & $97.9$ \\
    \cline{2-9}
     & $0.002$ & $10000$ & -- & -- & -- & $3.85e-4$ & $6.32e-5$ & $22.2$ \\
     & $0.001$ & $10000$ & -- & -- & -- & $9.11e-5$ & $1.49e-5$ &  $43.9$ \\
     & $0.0005$ & $10000$ & -- & -- & -- & $1.75e-5$ & $3.05e-6$ & $86.9$ \\
     & $0.00025$ & $10000$ & -- & -- & -- & $4.24e-6$ & $8.58e-7$ & $179.0$ \\
    \hline
  \end{tabular}
  \caption{1D diffusion equation: comparison between FEM
    and the current locELM method, in terms of the maximum/rms errors
    in the overall domain and the training/computation time.
    The temporal domain size is $t_f=1$. $R_m=1.0$ in locELM simulations. 
  }
  \label{tab_diffu_9}
\end{table}

Let us next compare the current
locELM method with the classical finite element method
for solving the 1D diffusion equation.
In Figures \ref{fg_diffu_8}(a) and (b) we compare profiles of
the solutions and their absolute errors at $t=1.0$,
obtained using the current locELM method and the finite element
method. The domain and problem settings in these tests correspond to those of
Figures \ref{fg_diffu_6}(e,f), with a temporal domain size $t_f=1$.
The simulation parameters for the locELM computation
also correspond to those of Figures \ref{fg_diffu_6}(e,f).
For the FEM simulation, the diffusion equation~\eqref{eq_diffu_1} is
discretized in time by the second-order backward differentiation formula (BDF2), and
the diffusion term is treated implicitly. We have employed a
time step size $\Delta t=0.00025$ and
$10,000$ uniform linear elements to discretize the spatial domain.  
It is evident from these data that both the FEM and the current method
have produced accurate results.
Figure \ref{fg_diffu_8}(c) shows the maximum and rms errors at $t=0.5$
versus the time step size $\Delta t$ with FEM, showing the
second-order temporal convergence rate. In these tests a fixed mesh
of $10,000$ uniform linear elements has been used, which accounts for
the observed error saturation in Figure \ref{fg_diffu_8}(c) when $\Delta t$
becomes sufficiently small.

Table \ref{tab_diffu_9} provides a comparison of the accuracy and
the computational cost of the locELM method and the finite element method.
In these tests the temporal domain size is set to $t_f=1$.
In the locELM simulations we employ a single time block in the spatial-temporal
domain, $5$ uniform sub-domains in the time block, and several sets of collocation
points/sub-domain and training parameters/sub-domain, a single hidden layer
in the local neural networks, and
$R_m=1.0$ when generating the random coefficients.
In the FEM simulations, we employ several sets of elements and $\Delta t$
values. The maximum error and the rms error in the overall spatial-temporal domain
have been computed, and the wall time for the computation or network training
have been recorded.
In Table \ref{tab_diffu_9} we list these errors and the wall time numbers
corresponding to the different simulation cases with locELM and FEM.
We observe that the current method performs markedly better than FEM.
The current method achieves a considerably better
accuracy with the same computational cost as FEM, and it incurs a lower computational
cost while achieving the same accuracy as FEM.
For example, the locELM case with $(Q,M)=(20\times 20,250)$ has a computational
cost comparable to the FEM cases with $2000$ elements and $\Delta t=0.001$
and with $5000$ elements and $\Delta t=0.002$. But the numerical errors of
locELM are considerably smaller, by around three orders of magnitude,
than those of the FEM cases.
The locELM case with $(Q,M)=(30\times 30,300)$ has a lower computational cost,
by a factor of about three,
and a considerably better accuracy, by a factor of nearly three orders of magnitude,
than the FEM case with $10,000$ elements and $\Delta t=0.0005$.


\subsection{Nonlinear Examples}

\subsubsection{Nonlinear Helmholtz Equation}

As the first nonlinear example,
we test the locELM method using the boundary value problem
with the nonlinear
Helmholtz equation in one dimension.
Consider the domain $[a,b]$ and the following boundary value problem
on this domain,
\begin{subequations}
  \begin{align}
    &
    \frac{d^2u}{d x^2} - \lambda u + \beta\sin(u) = f(x),
    \label{eq_nl_hm_1} \\
    &
    u(a) = h_1, \\
    &
    u(b) = h_2, \label{eq_nl_hm_2}
  \end{align}
\end{subequations}
where $u(x)$ is the function to be solved for, $f(x)$ is a prescribed
source term, $\lambda$ and $\beta$ are constant parameters, 
and $h_1$ and $h_2$ are the boundary values.
We assume the following values for the constant parameters
involved in these equations and domain specification,
\begin{equation*}
  a = 0, \quad
  b = 8, \quad
  \lambda = 50, \quad
  \beta = 10. \quad
\end{equation*}
We choose the source term $f(x)$ and the boundary values $h_1$ and $h_2$
such that the following function satisfies the
equations \eqref{eq_nl_hm_1}--\eqref{eq_nl_hm_2},
\begin{equation}\label{eq_nhm_3}
  u(x) = \sin\left(3\pi x + \frac{3\pi}{20} \right)
  \cos\left(4\pi x - \frac{2\pi}{5} \right) + \frac32 + \frac{x}{10}.
\end{equation}

\begin{figure}
  \centerline{
    \includegraphics[width=1.5in]{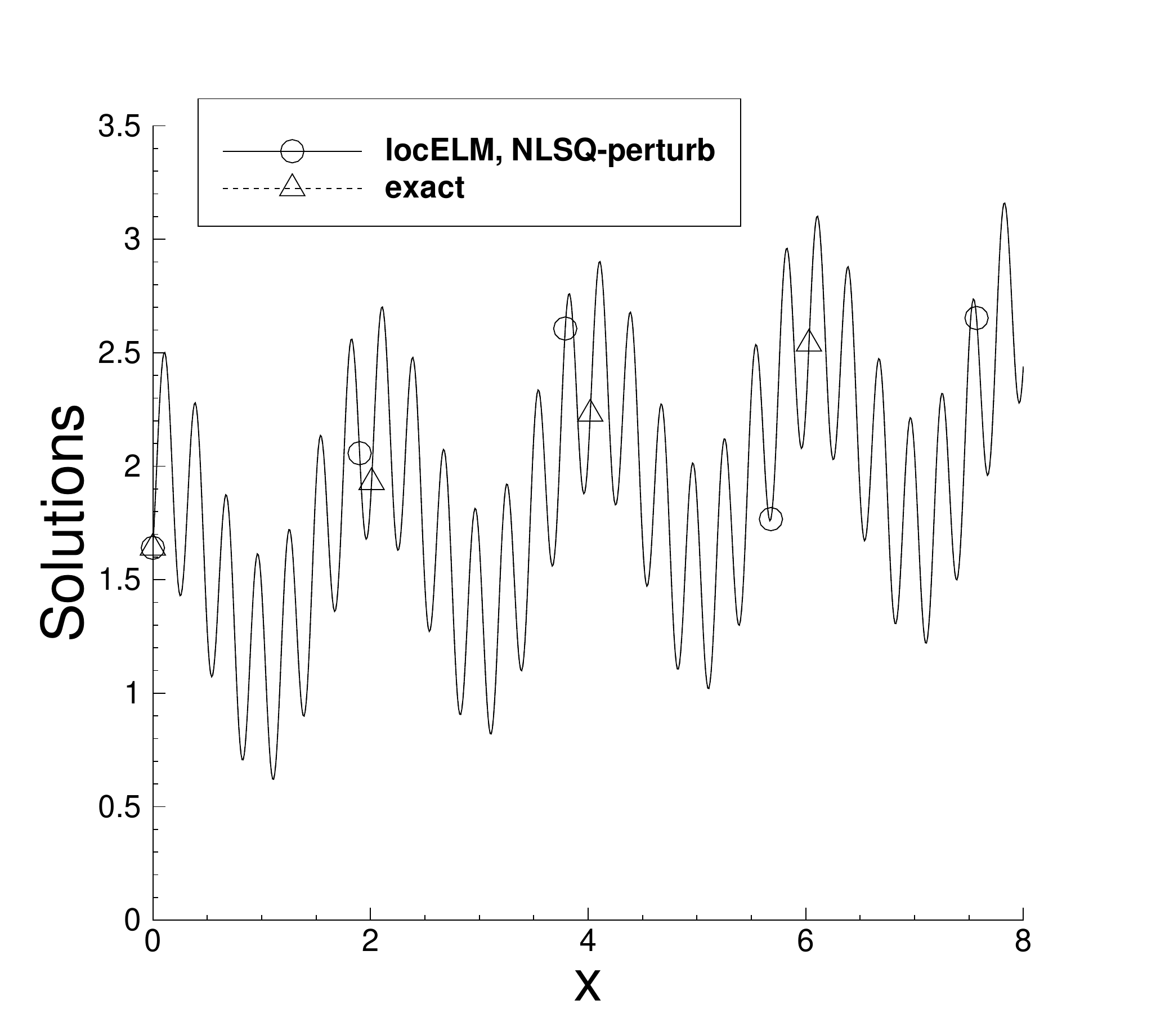}(a)
    \includegraphics[width=1.5in]{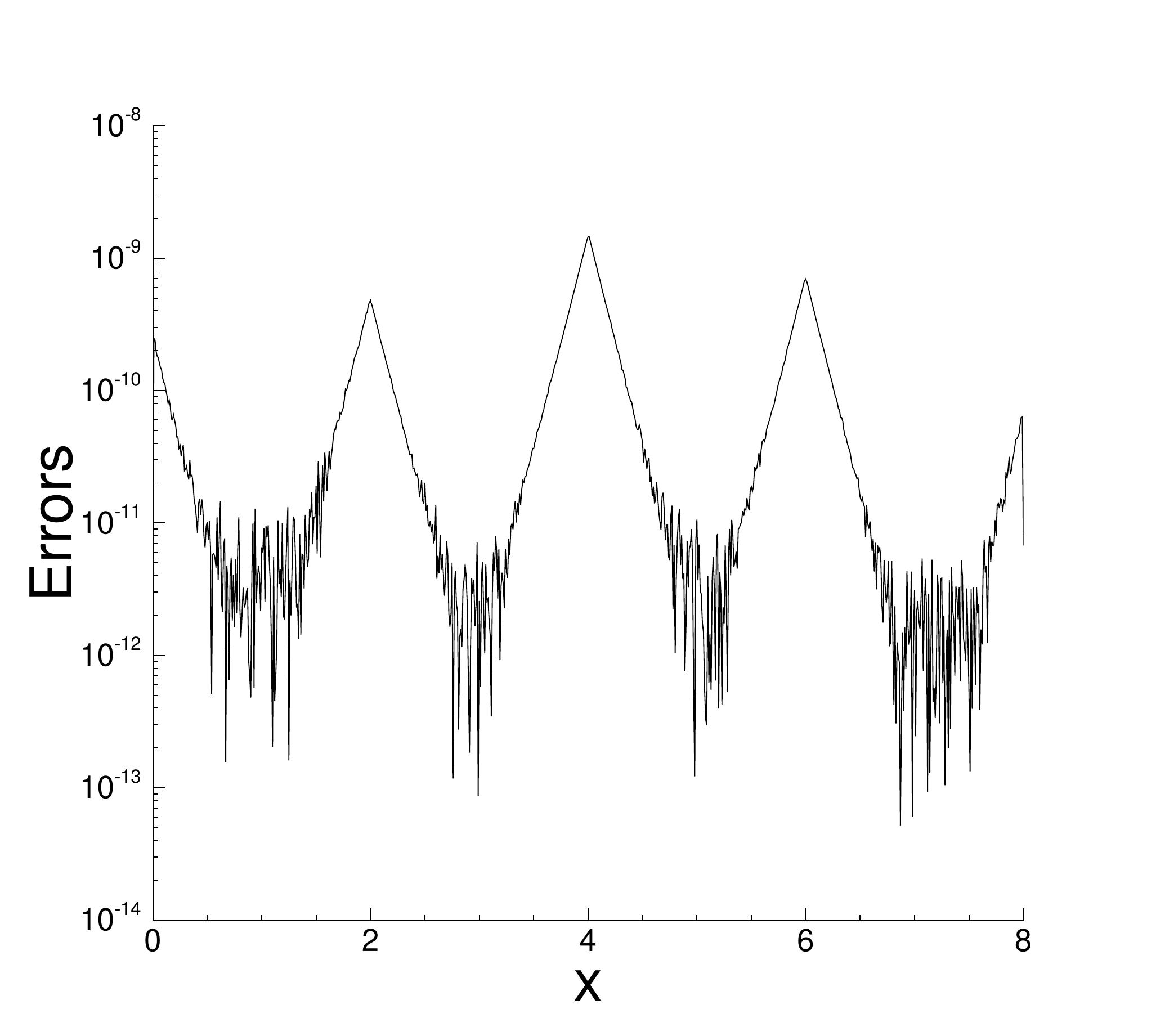}(b)
    \includegraphics[width=1.5in]{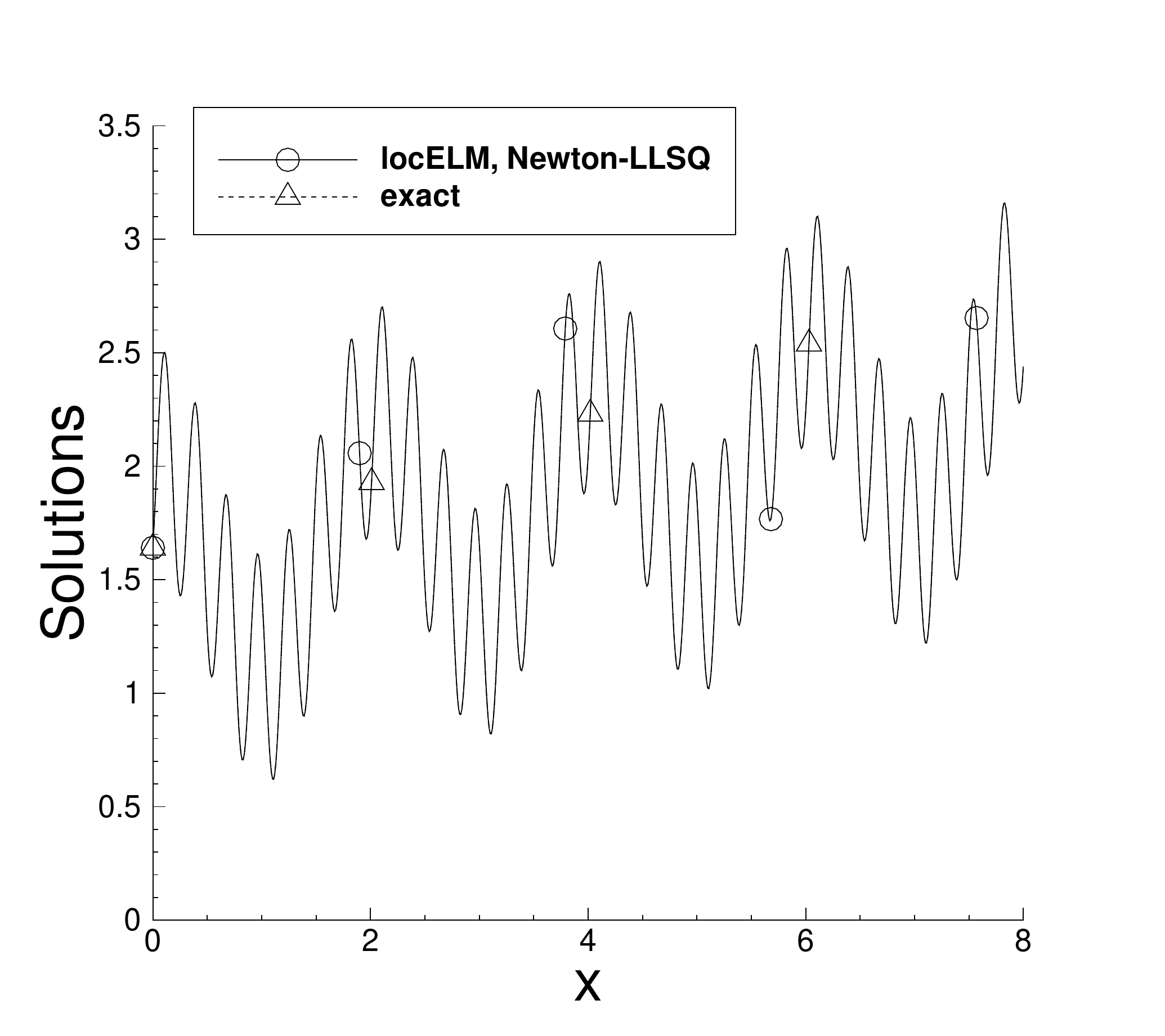}(c)
    \includegraphics[width=1.5in]{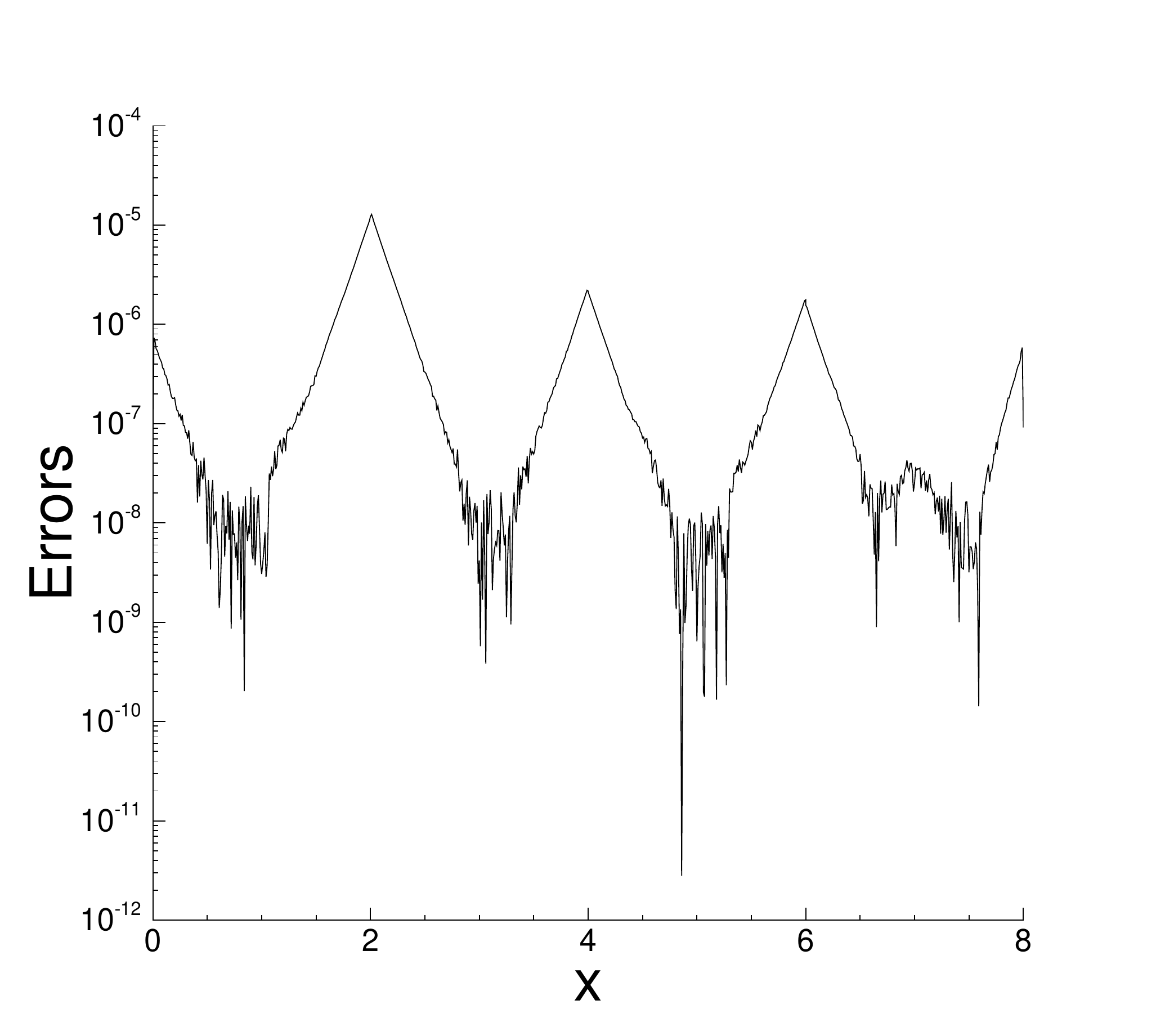}(d)
  }
  \caption{Nonlinear Helmholtz equation: profiles of the locELM
    solutions (a,c) and their absolute errors (b,d), computed
    using NLSQ-perturb (a,b) and Newton-LLSQ (c,d).
  }
  \label{fg_nhm_1}
\end{figure}


We employ the locELM method discussed in Section \ref{sec:nonl_steady}
for solving this problem, by restricting the method to one dimension.
We partition the domain $[a,b]$ into $N_e$ uniform sub-domains (sub-intervals),
and impose the $C^1$ continuity conditions across the sub-domain boundaries.
We employ $Q$ uniform collocation points within each sub-interval.

The local neural network for each sub-domain consists of
an input layer with one node (representing $x$), a single hidden layer
with $M$ nodes and the $\tanh$ activation function, and
an output layer with one node (representing the solution $u$)
and no activation function and no bias.
An additional affine mapping operation normalizing the input $x$ data to
the interval $[-1,1]$ is incorporated into the local neural networks
right behind the input layer for each sub-domain.
The weight and bias coefficients in the hidden layer of the local
neural networks are set to uniform random values generated on the interval
$[-R_m,R_m]$.
We employ a fixed seed value $12$ for the Tensorflow random generator
for all the tests reported in this sub-section.


We employ the nonlinear least squared method with perturbations
(NLSQ-perturb) and the combined Newton/linear least squared method (Newton-LLSQ)
from Section \ref{sec:nonl_steady} for computing the resultant
nonlinear problem. The initial guess to the solution 
is set to zero in all the tests of this subsection.
In the NLSQ-perturb method (see Algorithm \ref{alg:alg_1}),
we have employed $\delta = 0.2$, and $\xi_2=1$ as discussed in Remark \ref{rem_9},
for generating the random perturbations in the following tests.


The locELM simulation parameters include
the number of sub-domains ($N_e$), the number of collocation points per
sub-domain ($Q$), the number of training parameters per sub-domain ($M$),
and the maximum magnitude of the random coefficients 
of the local neural networks ($R_m$).


Figure \ref{fg_nhm_1} illustrates the profiles of
the locELM solutions and their absolute errors computed using
the NLSQ-perturb and Newton-LLSQ methods.
In these simulations, we have employed $N_e=4$ uniform sub-domains,
$Q=100$ uniform collocation points per sub-domain, $M=200$ training
parameters per sub-domain, and $R_m=5.0$ for generating the random
weight/bias coefficients.
The profile of the exact solution given by \eqref{eq_nhm_3}
is also included in these plots.
The solution profiles obtained with the current method exactly overlap
with that of the exact solution.
The error profiles indicate that the NLSQ-perturb method results in more
accurate results than Newton-LLSQ, with error levels
on the order $10^{-12}\sim 10^{-9}$ for NLSQ-perturb
versus $10^{-9}\sim 10^{-5}$ for Newton-LLSQ.


\begin{figure}
  \centerline{
    \includegraphics[width=2.2in]{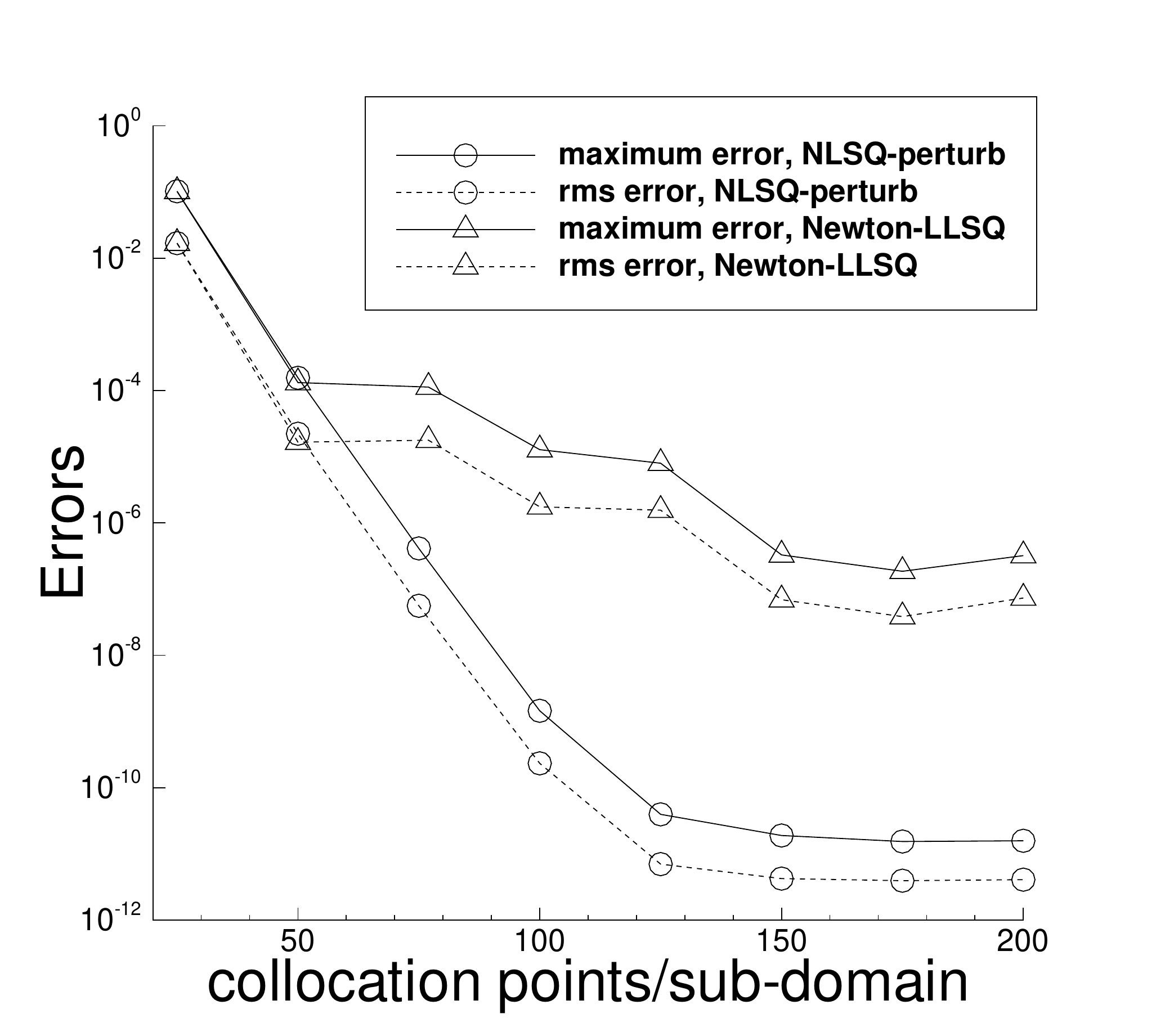}(a)
    \includegraphics[width=2.2in]{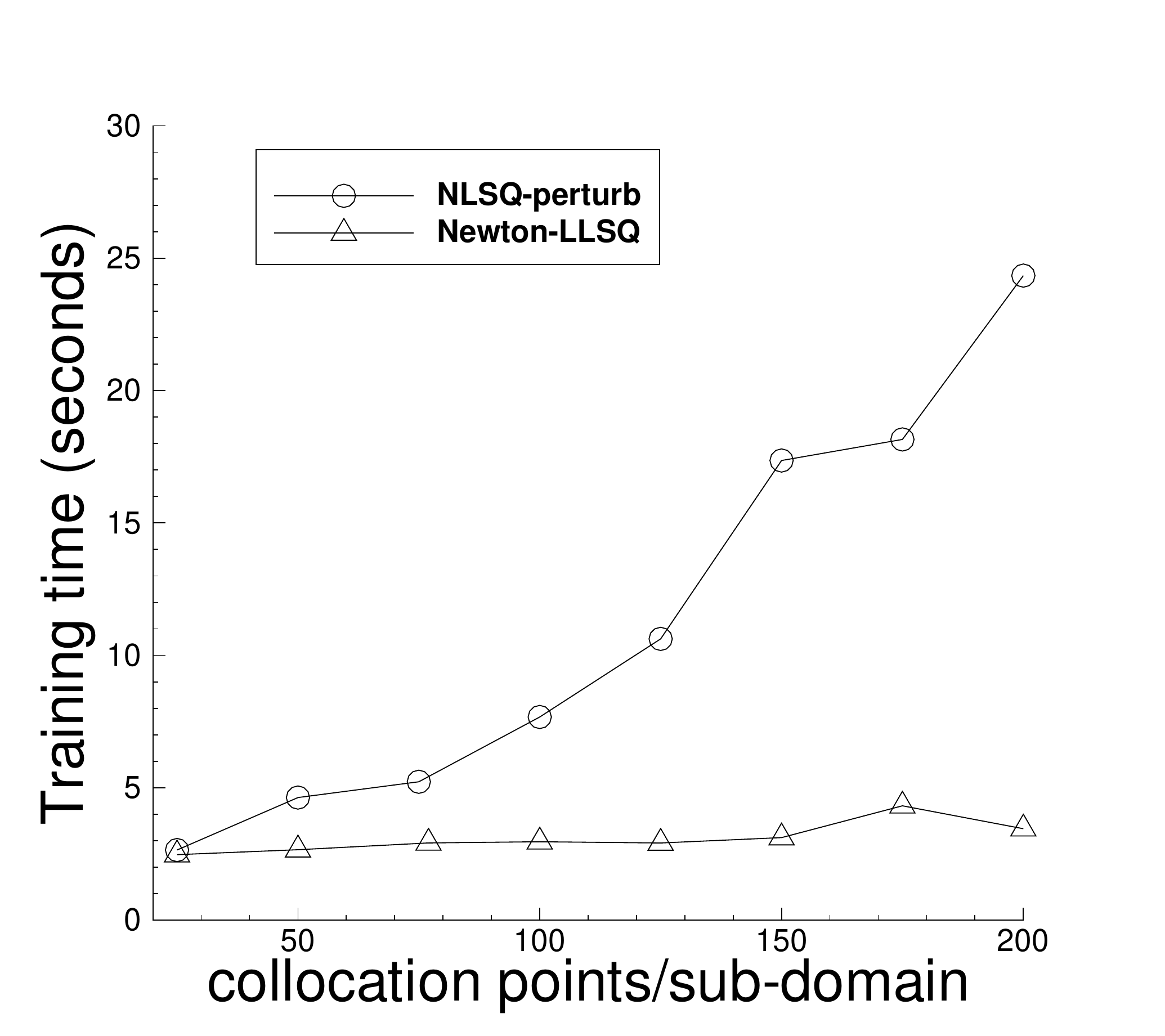}(b)
  }
  \caption{Effect of collocation points (nonlinear Helmholtz equation):
    (a) the maximum and rms errors in the domain,
    and (b) the network training time,
    as a function of the number of collocation points per
    sub-domain, computed using the locELM method
    with NLSQ-perturb and Newton-LLSQ.
  }
  \label{fg_nhm_2}
\end{figure}

Figure \ref{fg_nhm_2} demonstrates the effect of the number of collocation
points per sub-domain on the simulation accuracy and
the computational cost.
In this group of tests, we have employed $N_e=4$ sub-domains,
$M=200$ training parameters per sub-domain, and $R_m=5.0$ when
generating the random coefficients.
The number of uniform collocation points per sub-domain is varied systematically
between $Q=25$ and $Q=200$.
Figure \ref{fg_nhm_2}(a) shows the maximum and rms errors in the domain
as a function of the number of collocation points per sub-domain,
obtained with NLSQ-perturb and Newton-LLSQ.
Figure \ref{fg_nhm_2}(b) shows the corresponding
training time of the overall neural network
versus the number of collocation points per sub-domain.
With the Newton-LLSQ method, the errors are observed to
decrease gradually with increasing
number of collocation points, and appear to stagnate at a level around
$10^{-6}$ when the number of collocation points/sub-domain is
beyond $150$.
With the NLSQ-perturb method, the errors initially decrease exponentially
with increasing number of collocation points (when below $125$),
and then stagnate at a level around $10^{-11}$ when the number of collocation
points/sub-domain increases to $150$ and beyond.
The  NLSQ-perturb results are in general considerably more accurate
than those obtained with Newton-LLSQ.
%
In terms of the training time, the Newton-LLSQ method is consistently faster than
NLSQ-perturb, and the difference becomes larger as the number of
collocation points increases.
With the Newton-LLSQ method, the training time appears not sensitive to
the number of collocation points, and remains nearly the same
with increasing number of collocation points (Figure \ref{fg_nhm_2}(b)).
With the NLSQ-perturb method, the training time increases approximately
linearly with increasing number of collocation points per sub-domain,
and it becomes substantially slower than Newton-LLSQ when the number of
collocation points becomes large.

\begin{figure}
  \centerline{
    \includegraphics[width=2.2in]{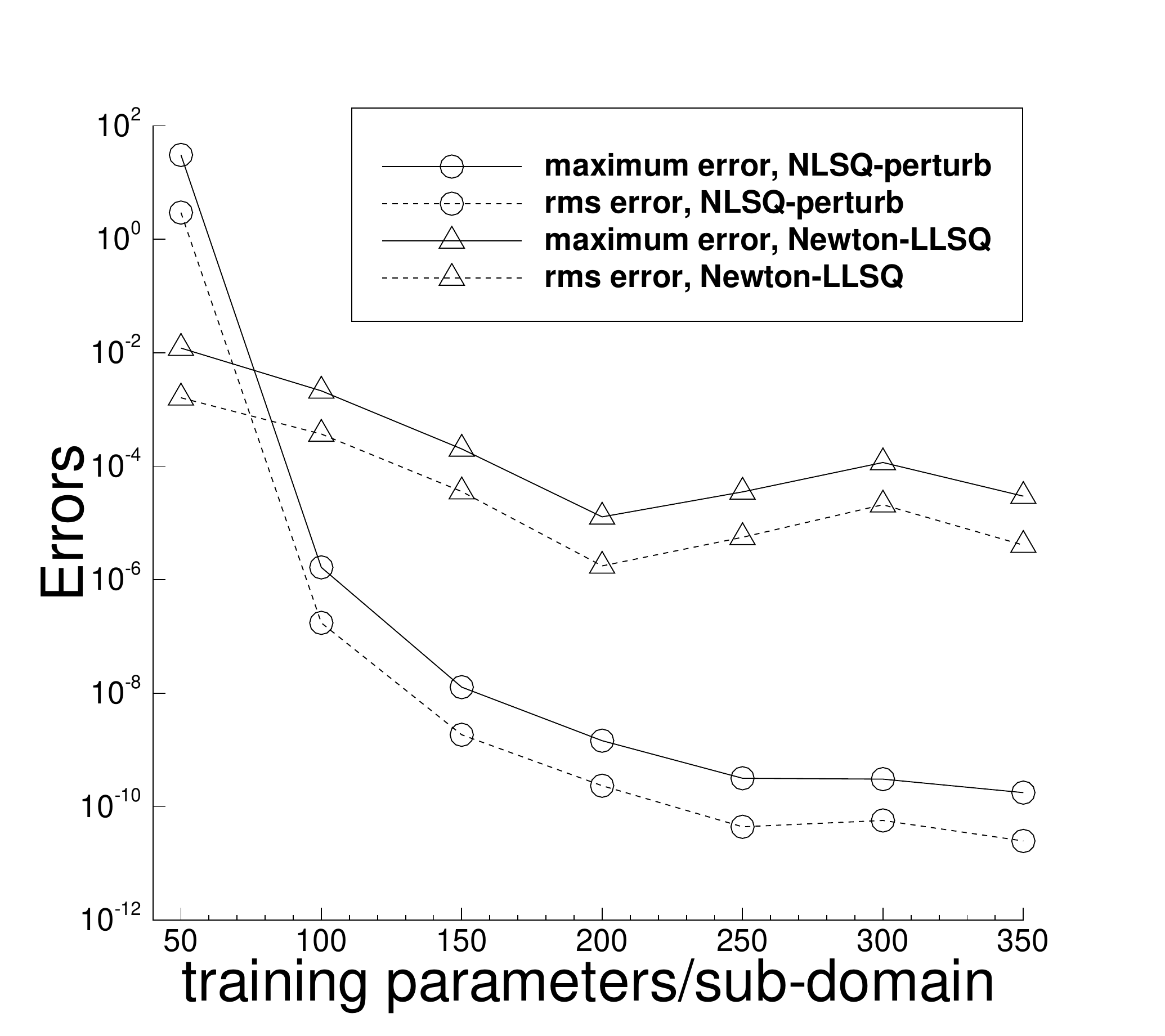}(c)
    \includegraphics[width=2.2in]{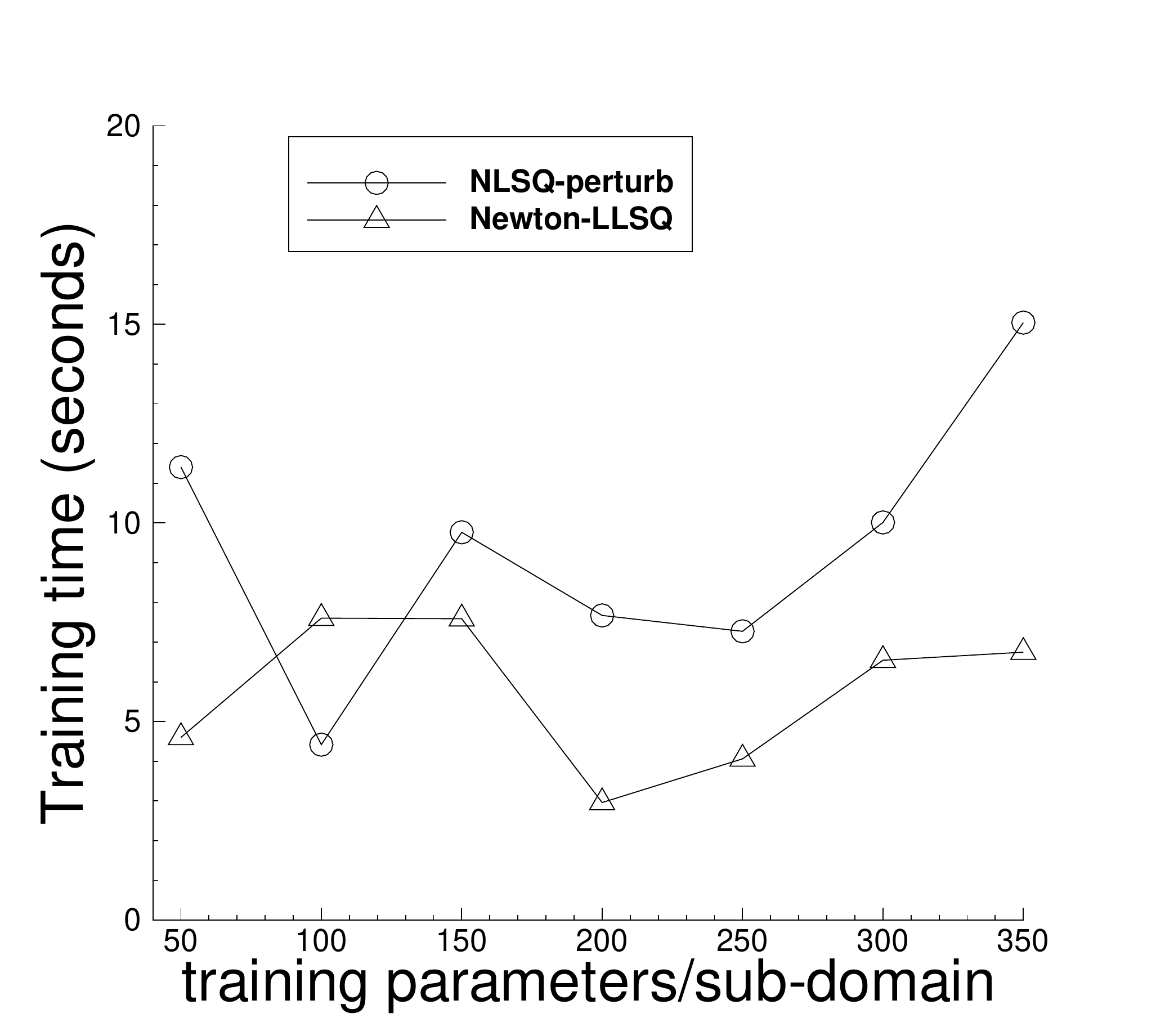}(c)
  }
  \caption{Effect of the number of training parameters (nonlinear Helmholtz equation):
    (a) the maximum and rms errors in the domain,
    and (b) the network training time,
     versus the number of
    training parameters per sub-domain, computed using the locELM method
    with NLSQ-perturb and Newton-LLSQ.
  }
  \label{fg_nhm_3}
\end{figure}

Figure \ref{fg_nhm_3} demonstrates the effect of the number of training parameters
per sub-domain on the simulation accuracy and the computational cost.
In this group of tests, we have employed $N_e=4$ sub-domains,
$Q=100$ uniform collocation points per sub-domain, and $R_m=5.0$
when generating the random coefficients in the hidden layers of the local
neural networks. The number of training parameters per sub-domain
is varied systematically between $50$ and $350$.
Figure \ref{fg_nhm_3}(a) shows the maximum and rms errors of the solutions
as a function of the number of training parameters per sub-domain
obtained with  NLSQ-perturb and Newton-LLSQ.
With NLSQ-perturb, the numerical errors decrease substantially as
the number of training parameters per sub-domain increases,
reaching  a level around
$10^{-10}$ when the number of training parameters increases beyond $200$.
With Newton-LLSQ, one can also observe a decrease in the errors as the number of
training parameters increases. But the error reduction is much slower.
When the number of training parameters per sub-domain exceeds $200$,
the errors with Newton-LLSQ no longer seem to decrease further and remain
at a level around $10^{-5}$. It is evident that
the results from the Newton-LLSQ method are generally much less accurate
than those from the NLSQ-perturb method.
Figure \ref{fg_nhm_3}(b) shows the corresponding network training time
as a function of the number of training parameters per sub-domain.
In the range of training parameters tested here, the training time with both of these
two methods appear to fluctuate around a certain level.
But the training time with the Newton-LLSQ method is generally notably
smaller than that with the NLSQ-perturb method, except for the outlier point
corresponding to $100$ training parameters per sub-domain.
These data suggest
that Newton-LLSQ is generally faster than NLSQ-perturb.


\begin{figure}
  \centerline{
    \includegraphics[width=2in]{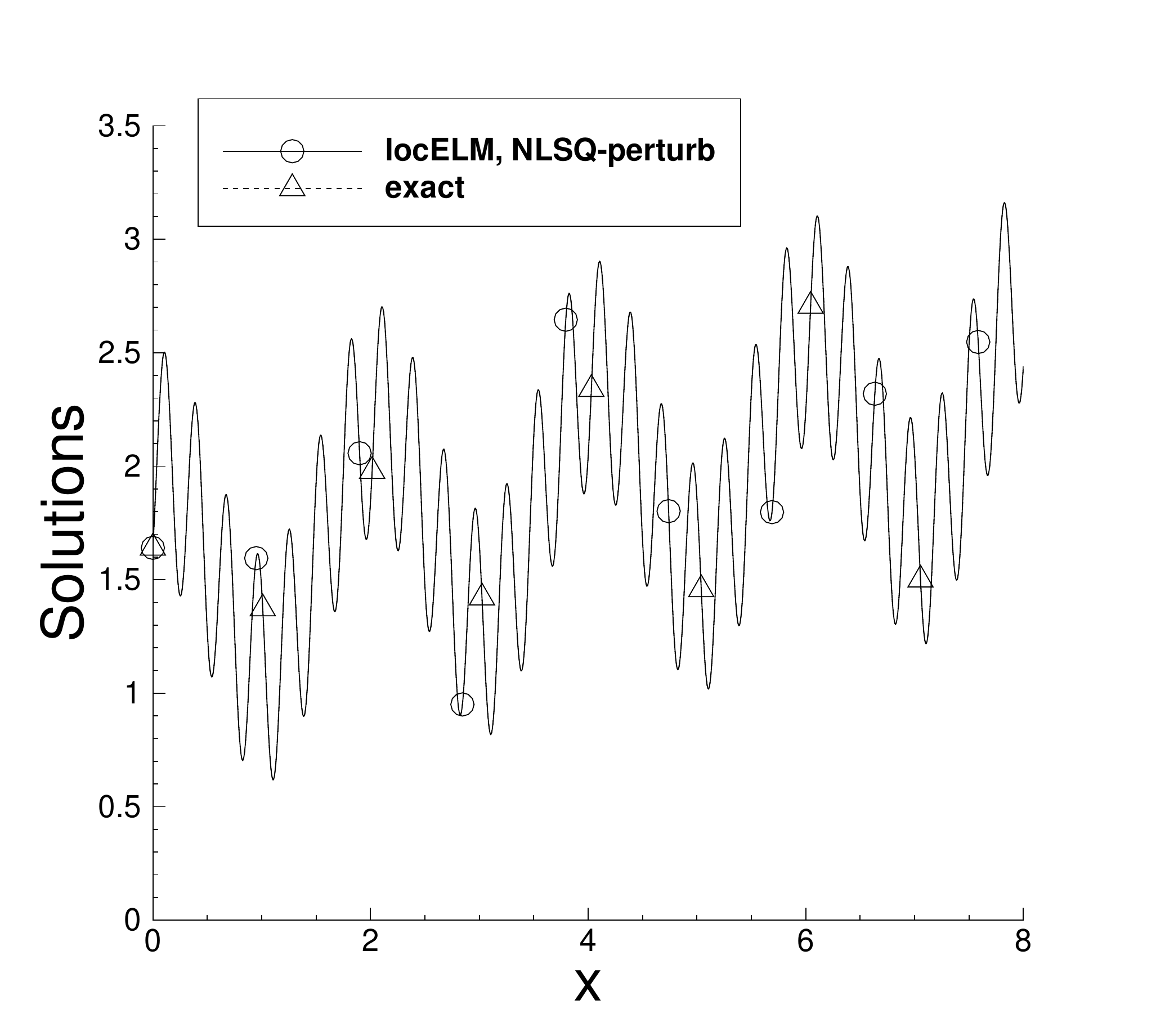}(a)
    \includegraphics[width=2in]{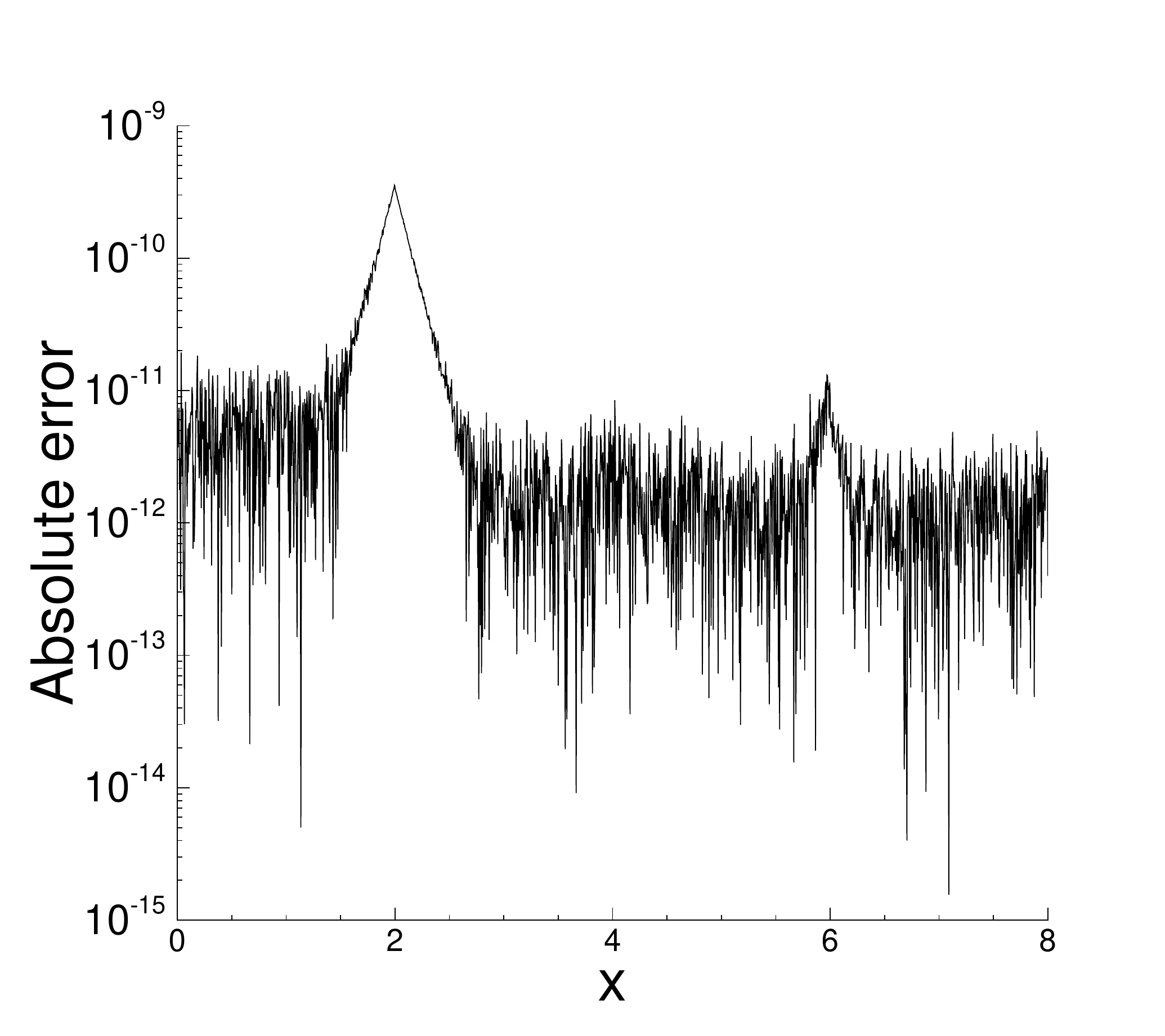}(b)
    \includegraphics[width=2in]{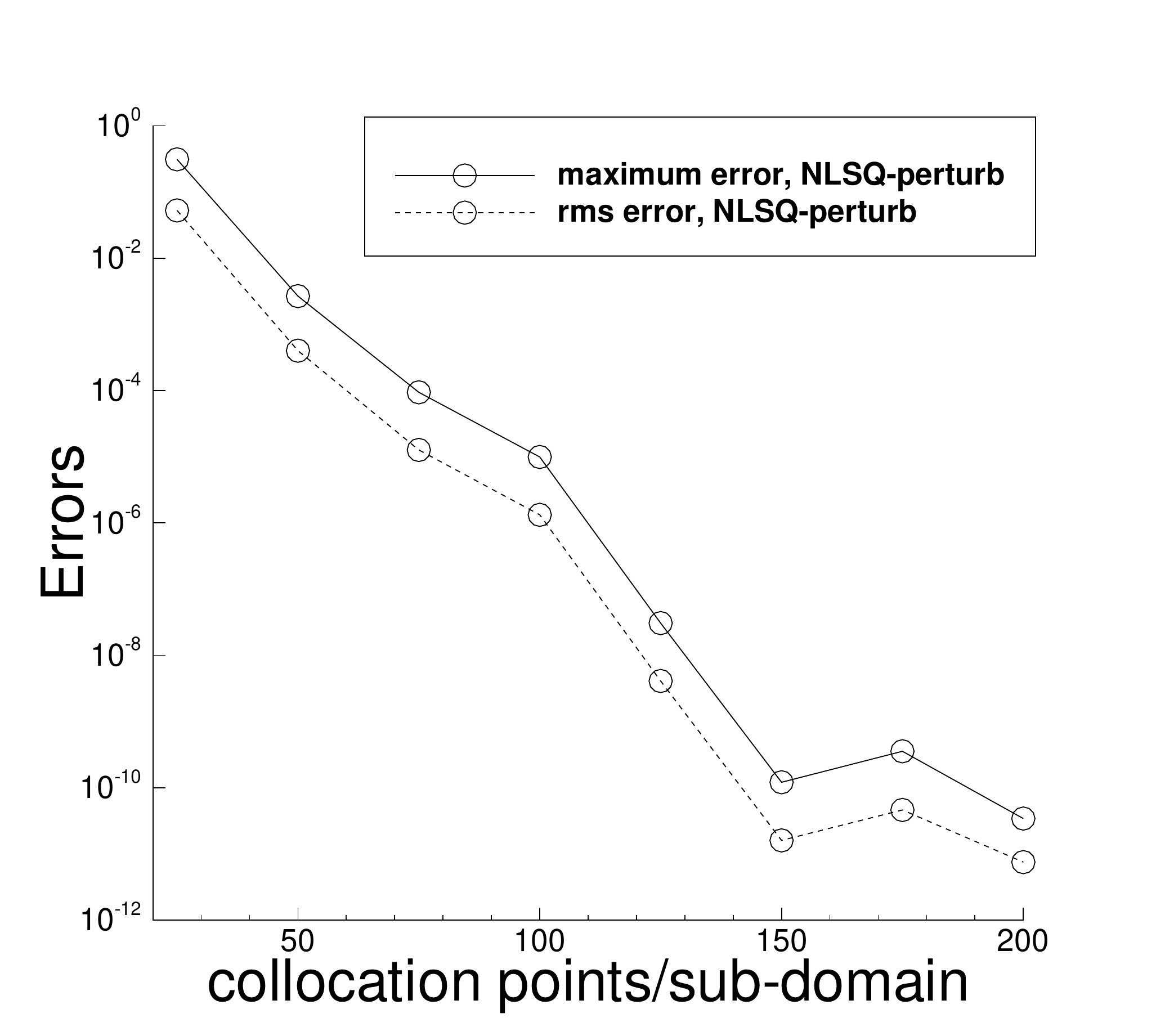}(c)
  }
  \caption{Results obtained with $2$ hidden layers in local neural networks
    (nonlinear Helmholtz equation): profiles of (a) the locELM (NLSQ-perturb) solution
    and (b) its absolute error. (c) The maximum/rms errors in the domain versus
    the number of collocation points per sub-domain.
  }
  \label{fg_nhm_a}
\end{figure}

With the current locELM method, the local neural network
can contain more than one hidden layer. As shown in previous sub-sections,
local neural networks with a small number (more than one) of hidden layers
can also deliver accurate  results using the current method.
Figure \ref{fg_nhm_a} demonstrates again this point with the nonlinear
Helmholtz equation. In this group of tests, we employ $N_e=4$ uniform sub-domains,
$M=250$ training parameters per sub-domain, $2$ hidden layers (with widths
$25$ and $250$, respectively, and the $\tanh$ activation function) in
each local neural network, and $R_m=2.0$ when generating the
random weight/bias coefficients for these hidden layers.
The number of uniform collocation points per sub-domain is varied
systematically in these tests.
Figures \ref{fg_nhm_a}(a) and (b) show the locELM solution and error profiles
obtained with $Q=175$ uniform collocation points per sub-domain
using the NLSQ-perturb method.
Figure  \ref{fg_nhm_a}(c) shows the maximum and rms errors in the domain
as a function of the number of uniform collocation points per sub-domain.
We observe an essentially exponential decrease in the numerical errors
with increasing number of collocation points per sub-domain.

\begin{figure}
  \centerline{
    \includegraphics[width=2.in]{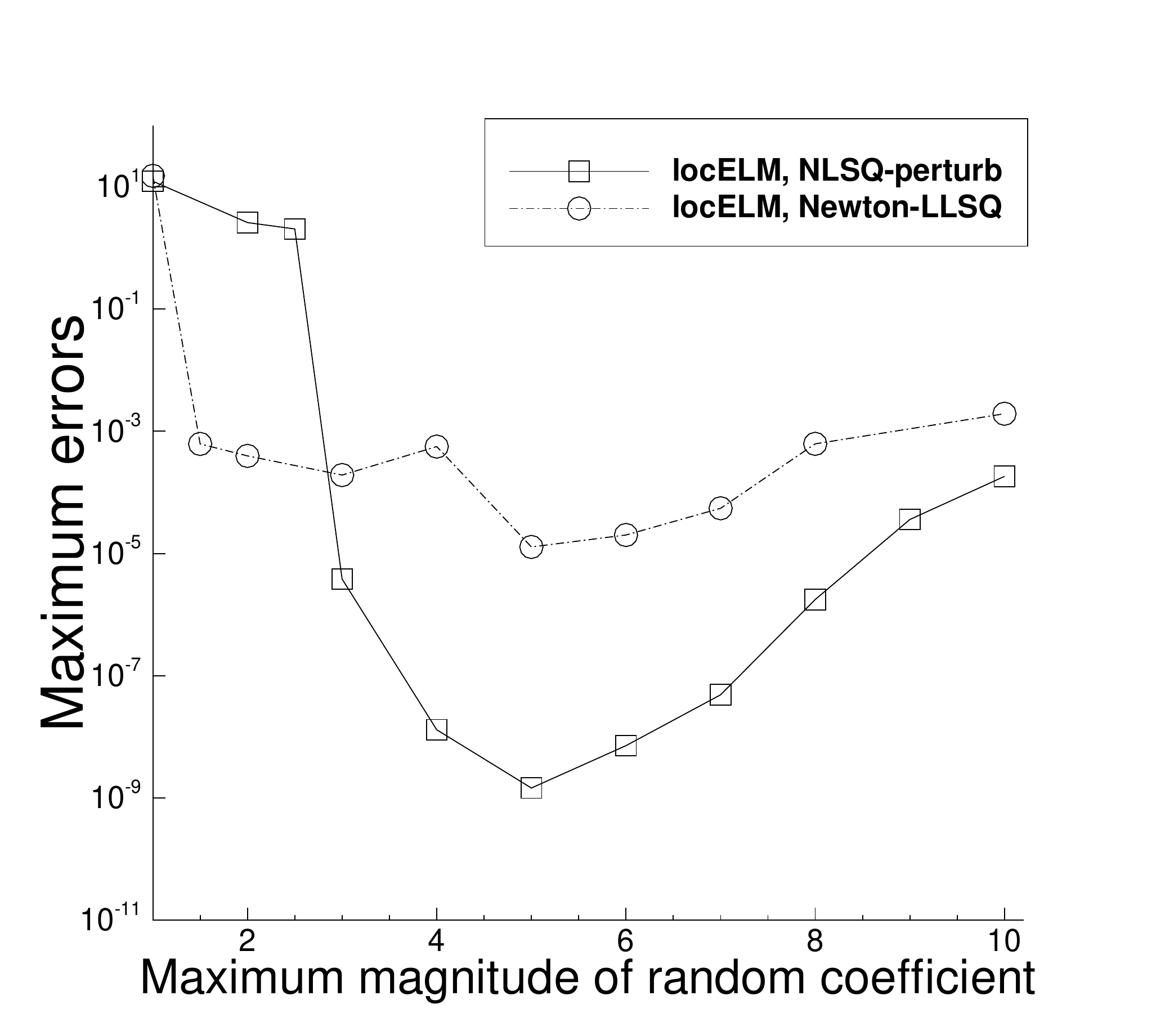}
  }
  \caption{Effect of the random coefficients in 
    local neural networks
    (nonlinear Helmholtz equation): the maximum  error in the domain
    versus $R_m$, 
    obtained with the NLSQ-perturb and Newton-LLSQ methods. 
  }
  \label{fg_nhm_4}
\end{figure}

Figure \ref{fg_nhm_4} illustrates the effect of the random coefficients in
the hidden layers of the local neural networks.
In this group of tests we employ $N_e=4$ sub-domains, $Q=100$ uniform collocation
points per sub-domain, $200$ training parameters per sub-domain, and
a single hidden layer in the local neural networks.
As discussed before, the weight/bias coefficients in the hidden layer of
each local neural network are set to uniform random values generated
on $[-R_m,R_m]$. In these tests, we vary $R_m$ systematically and study its
effect. Figure \ref{fg_nhm_4} shows the maximum  error in the overall
domain as a function of $R_m$, 
obtained with the NLSQ-perturb and the Newton-LLSQ methods.
The  error exhibits a behavior similar to what has been observed
with the linear problems. The methods have a better accuracy with a range
of moderate $R_m$ values, and the results are less accurate
with very large or very small $R_m$ values.
We again observe that the NLSQ-perturb result is significantly more accurate 
than that of Newton-LLSQ, except for a range of small $R_m$ values.

\begin{figure}
  \centerline{
    \includegraphics[width=2in]{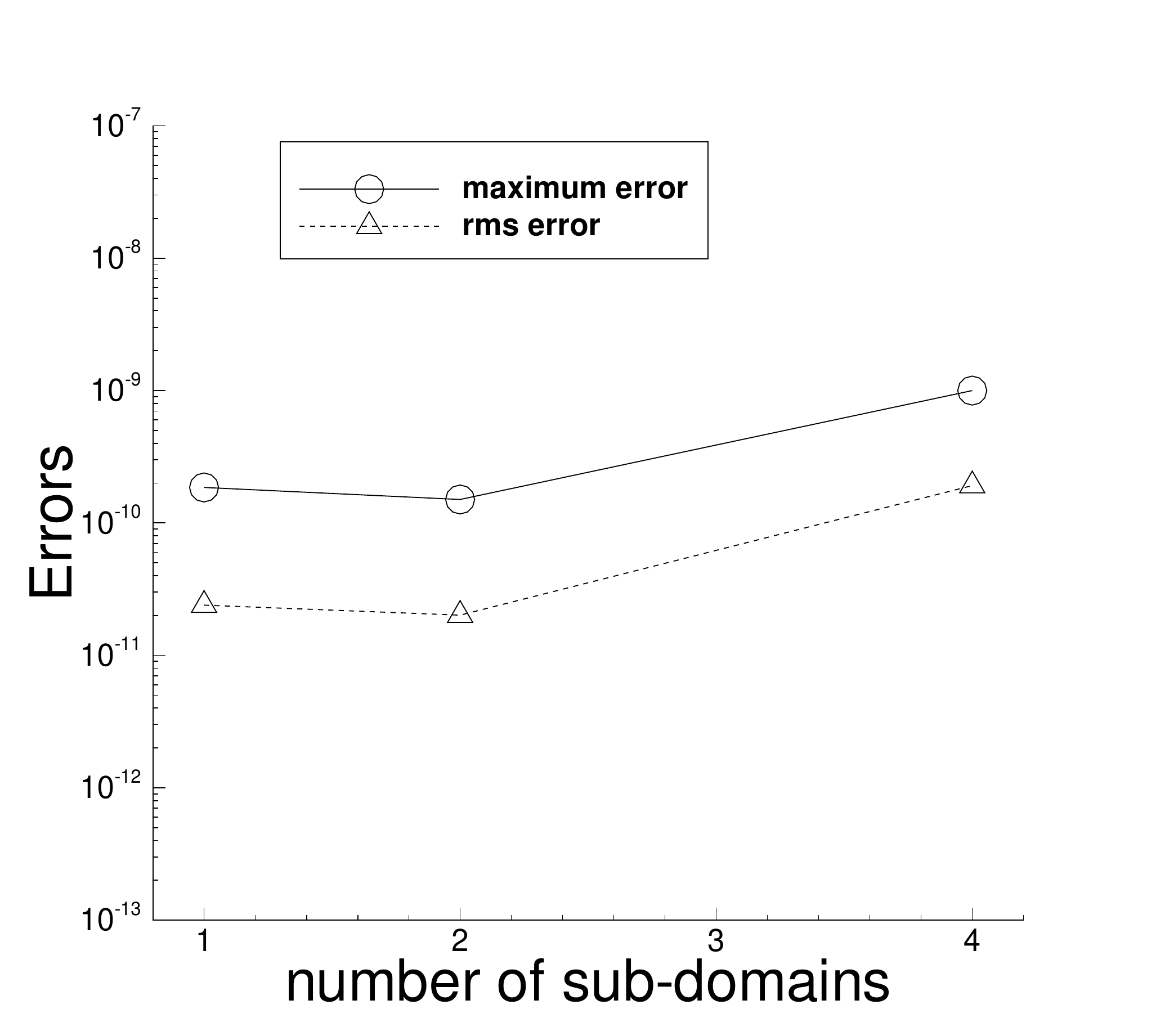}(a)
    \includegraphics[width=2in]{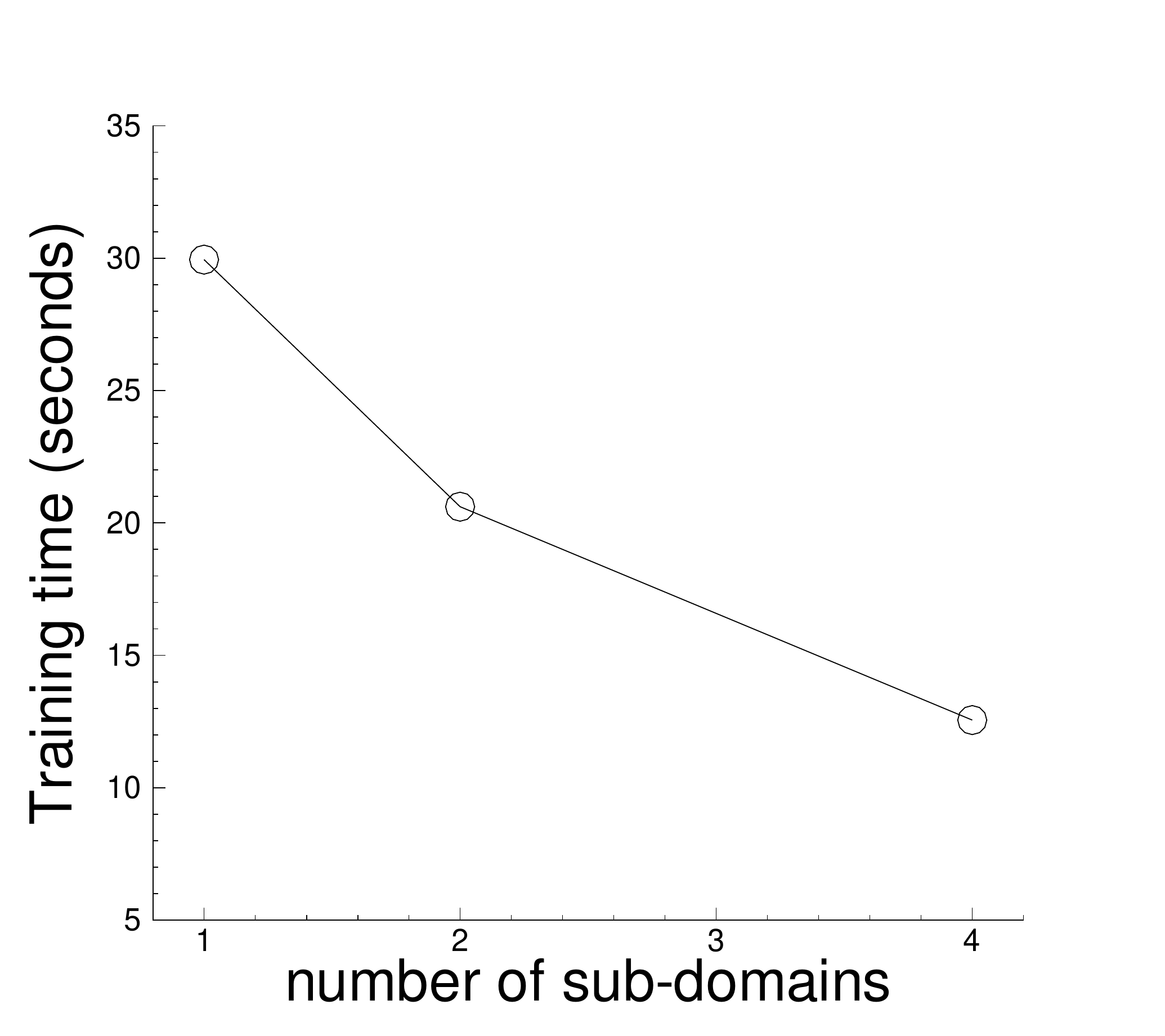}(b)
  }
  \caption{Effect of the number of sub-domains, with fixed total degrees of freedom
    in the domain (nonlinear Helmholtz equation): (a) the maximum and
    rms errors in the domain, and
    (b) the training time, as a function of the number of sub-domains
    in the locELM simulation
    with NLSQ-perturb. 
  }
  \label{fg_nhm_5}
\end{figure}

In Figure \ref{fg_nhm_5} we study the effect of the number of sub-domains
on the simulation accuracy and the computational cost, while the total degrees of
freedom 
in the domain are fixed.
In these tests we vary the number of uniform sub-domains
($N_e$). We choose the number of uniform collocation
points per sub-domain ($Q$) and the training parameters per sub-domain ($M$)
such that the total number of collocation points in the domain
is fixed at $N_eQ=400$ and the total number of training parameters in
the domain is fixed at $N_eM=800$.
We have tested three cases, corresponding to $N_e=1$, $2$ and $4$.
As in the previous sections, the case with one sub-domain ($N_e=1$)
corresponds to use of a global ELM.
Figure \ref{fg_nhm_5}(a) shows the maximum and rms errors in the overall domain
as a function of the number of sub-domains.
Figure \ref{fg_nhm_5}(b) shows the corresponding training time
versus the number of sub-domains.
These results are obtained with the NLSQ-perturb method.
We have employed $R_m=20.0$ when generating the random coefficients with
one sub-domain ($N_e=1$), $R_m=10.0$ with two sub-domains ($N_e=10.0$)
and $R_m=4.5$ with four sub-domains ($N_e=4.5$).
These $R_m$ values approximately reside in the optimal range of $R_m$
values for these cases.
One can observe that the numerical errors obtained with
different number of sub-domains are comparable, with the errors obtained on
four sub-domains a little worse than those of the other cases.
On the other hand, the network training time decreases significantly with
increasing number of sub-domains.

\begin{figure}
  \centerline{
    \includegraphics[width=1.5in]{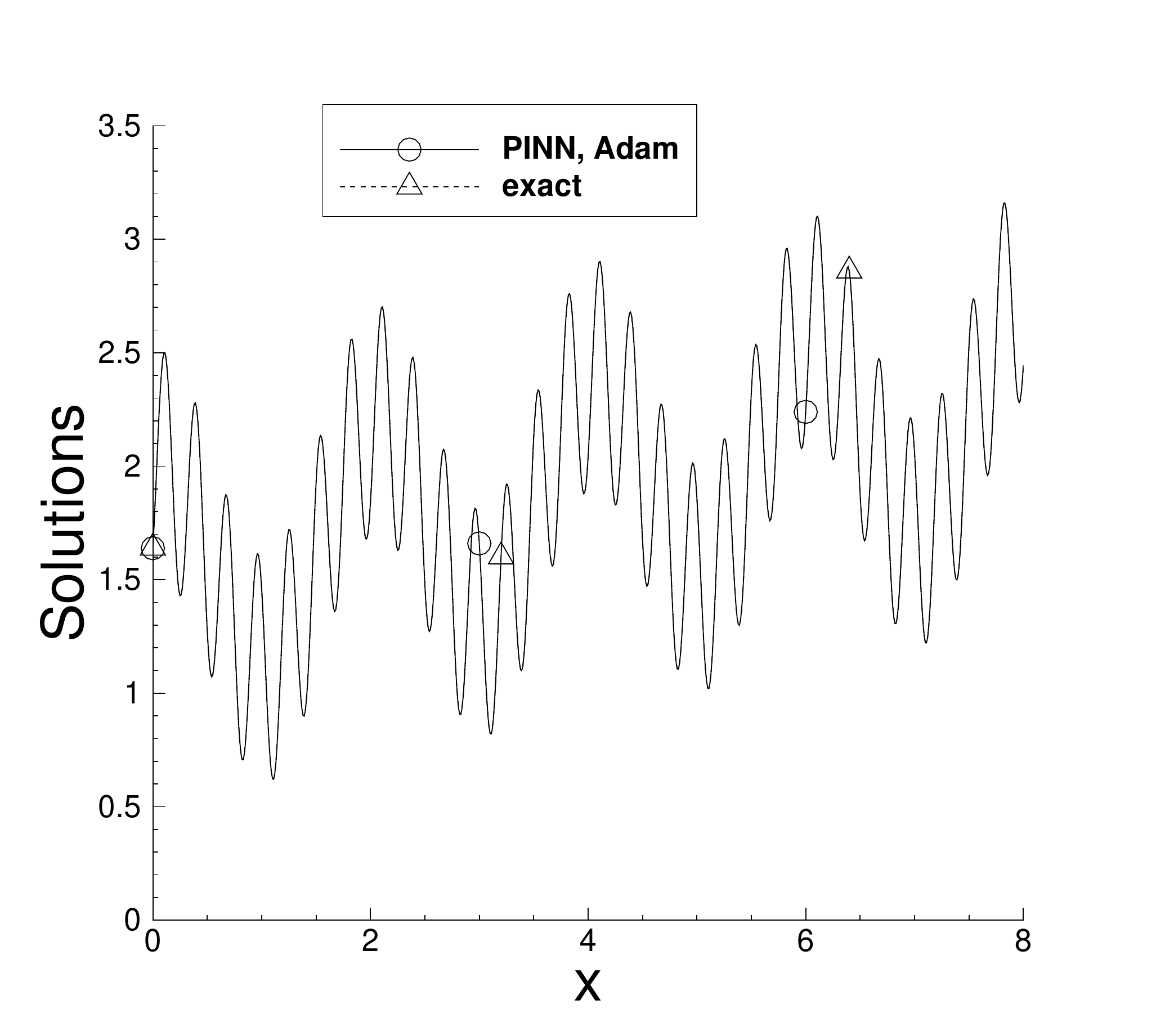}(a)
    \includegraphics[width=1.5in]{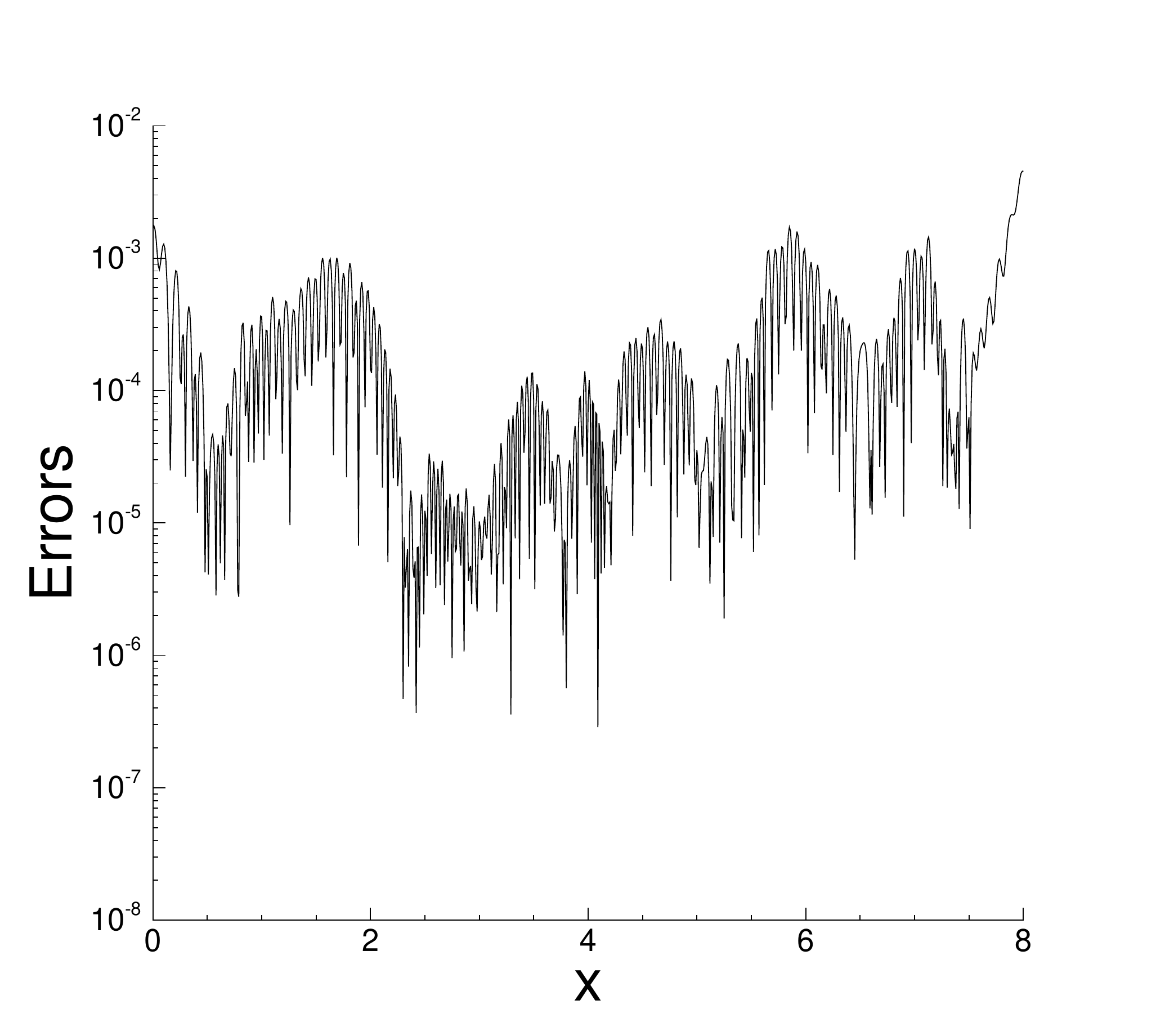}(b)
    \includegraphics[width=1.5in]{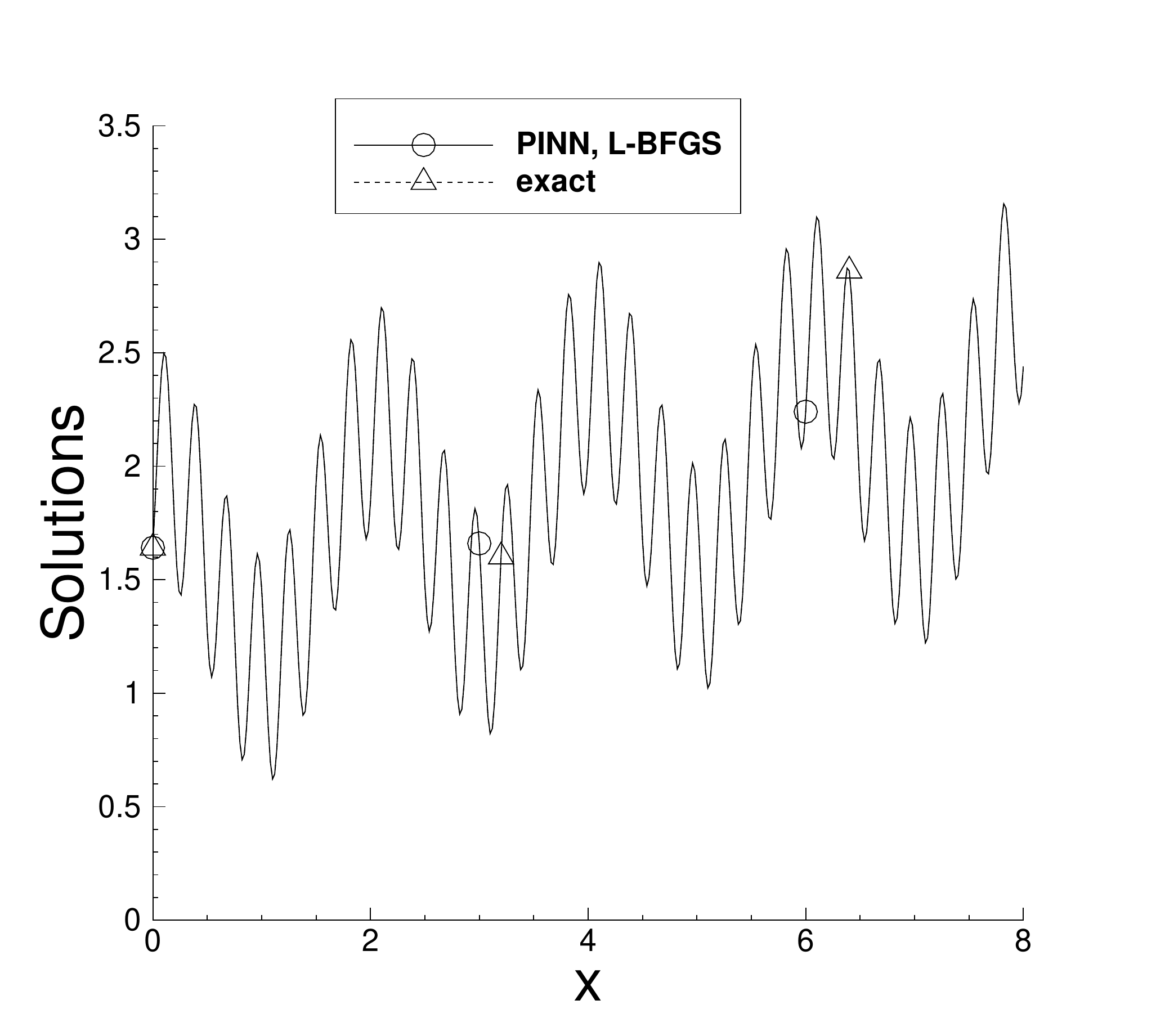}(c)
    \includegraphics[width=1.5in]{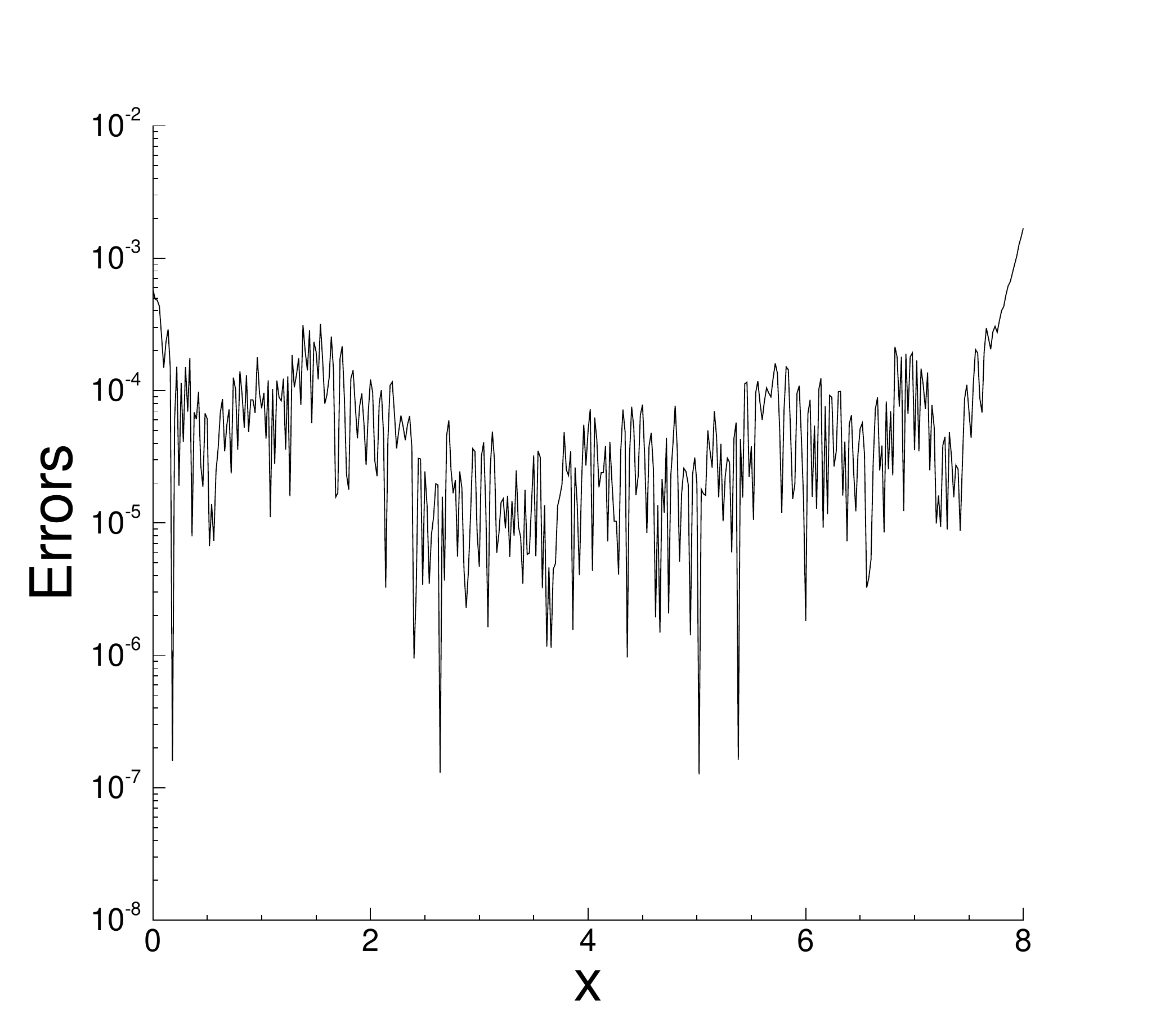}(d)
  }
  \caption{Nonlinear Helmholtz equation:
    Distributions of the solutions (a,c)
    and their absolute errors (b,d) computed using PINN~\cite{RaissiPK2019}
    with the Adam optimizer
    (a,b) and the L-BFGS optimizer (c,d).
    These can be compared with those in Figure \ref{fg_nhm_1}
    computed using locELM.
  }
  \label{fg_nhm_6}
\end{figure}

We next compare the current locELM method with the PINN method~\cite{RaissiPK2019}
for solving the nonlinear Helmholtz equation.
Figure \ref{fg_nhm_6} shows distributions of
the PINN solutions and their absolute errors against the exact solution
given in equation~\eqref{eq_nhm_3}, computed using the Adam optimizer
(Figures \ref{fg_nhm_6}(a,b)) and the L-BFGS optimizer (Figures \ref{fg_nhm_6}(c,d)).
With the Adam optimizer, the neural network consists of $7$ hidden layers,
with a width of $50$ nodes in each layer and the $\tanh$ activation function,
in addition to the input layer of one node (representing $x$) and the output
layer of one node (representing the solution $u$).
The network has been trained on the input data of $400$ uniform collocation points
for $45,000$ epochs, with the learning rate gradually decreasing from $0.001$
at the beginning to $5\times 10^{-6}$ at the end of the training.
With the L-BFGS optimizer, the neural network consists of $4$ hidden layers,
with a width of $50$ nodes in each layer and the $\tanh$ activation function,
apart from the input layer of one node and the output layer of one node.
The network has been trained on the input data of $400$ uniform collocation
points in the domain for $22,000$ L-BFGS iterations.
The results indicate that the PINN method has captured the solution quite accurately,
with the errors on the order $10^{-5}\sim 10^{-3}$ with the Adam optimizer
and on the order $10^{-5}\sim 10^{-4}$ with the L-BFGS optimizer.
Comparing the PINN results in this figure and the locELM  results in Figure~\ref{fg_nhm_1},
we can observe that the locELM method is considerably more accurate than
PINN.

\begin{table}
  \centering
  \begin{tabular}{lllll}
    \hline
    method & maximum error & rms error & epochs/iterations & training time (seconds)\\
    PINN (Adam) & $4.56e-3$ & $5.04e-4$ & $45,000$ & $578.2$ \\
    PINN (L-BFGS) & $1.69e-3$ & $1.69e-4$ & $22,000$ & $806.4$ \\
    locELM (NLSQ-perturb) & $1.45e-9$ & $2.34e-10$ & $71$ & $7.7$ \\
    locELM (Newton-LLSQ) & $1.28e-5$ & $1.75e-6$ & $5$ & $2.7$ \\
    \hline
  \end{tabular}
  \caption{Nonlinear Helmholtz equation: comparison between locELM and PINN
     in terms of the maximum/rms errors in the domain,
    the number of epochs or nonlinear iterations, and the network training time.
    The problem settings and simulation parameters correspond to those of Figures
    \ref{fg_nhm_1} and \ref{fg_nhm_6}.
  }
  \label{tab_nhm_7}
\end{table}

Table~\ref{tab_nhm_7} provides further comparisons between locELM
and PINN in terms of the accuracy and the computational cost.
Here we have listed the maximum and rms errors in the domain, the number
of epochs or nonlinear iterations in the training,
and the network training time, associated with
the PINN (with Adam/L-BFGS optimizers) simulations and the current locELM
simulations. The problem settings and the simulation parameters here
correspond to those in Figure \ref{fg_nhm_1} with locELM
and those in Figure \ref{fg_nhm_6} with PINN.
It is evident that the current locELM method is much more accurate than PINN.
For example, the errors obtained using locELM/NLSQ-perturb
are about six orders of magnitude smaller than those obtained by
PINN/L-BFGS. The errors obtained by locELM/Newton-LLSQ
are about two orders of magnitude smaller than those of PINN/L-BFGS.
Furthermore, the current method is computationally much cheaper than PINN,
with the training time approximately two orders of magnitude smaller
(e.g.~about $8$ seconds with locELM/NLSQ-perturb versus around $806$ seconds
with PINN/L-BFGS).


\begin{figure}
  \centerline{
    \includegraphics[width=2.2in]{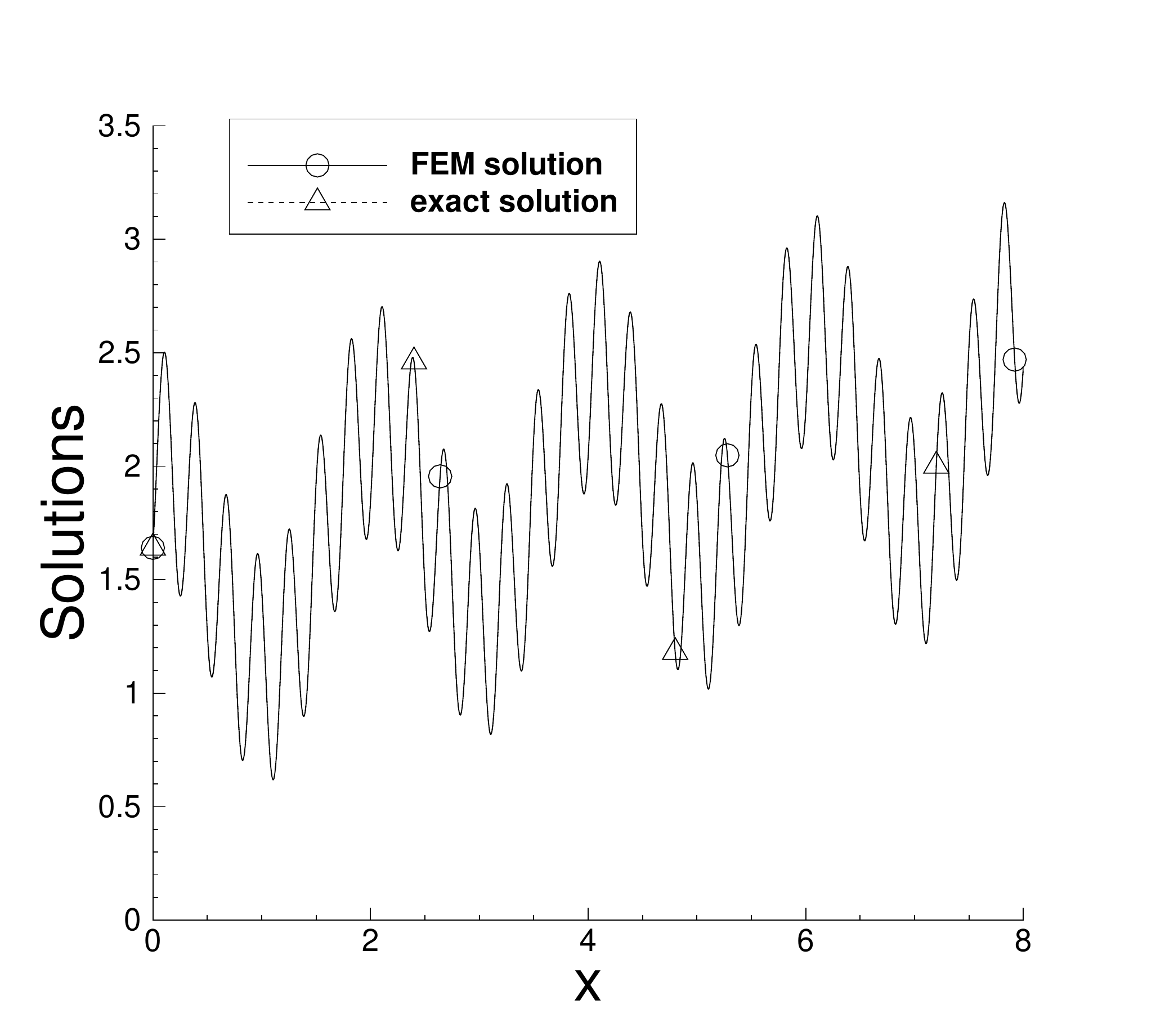}(a)
    \includegraphics[width=2.2in]{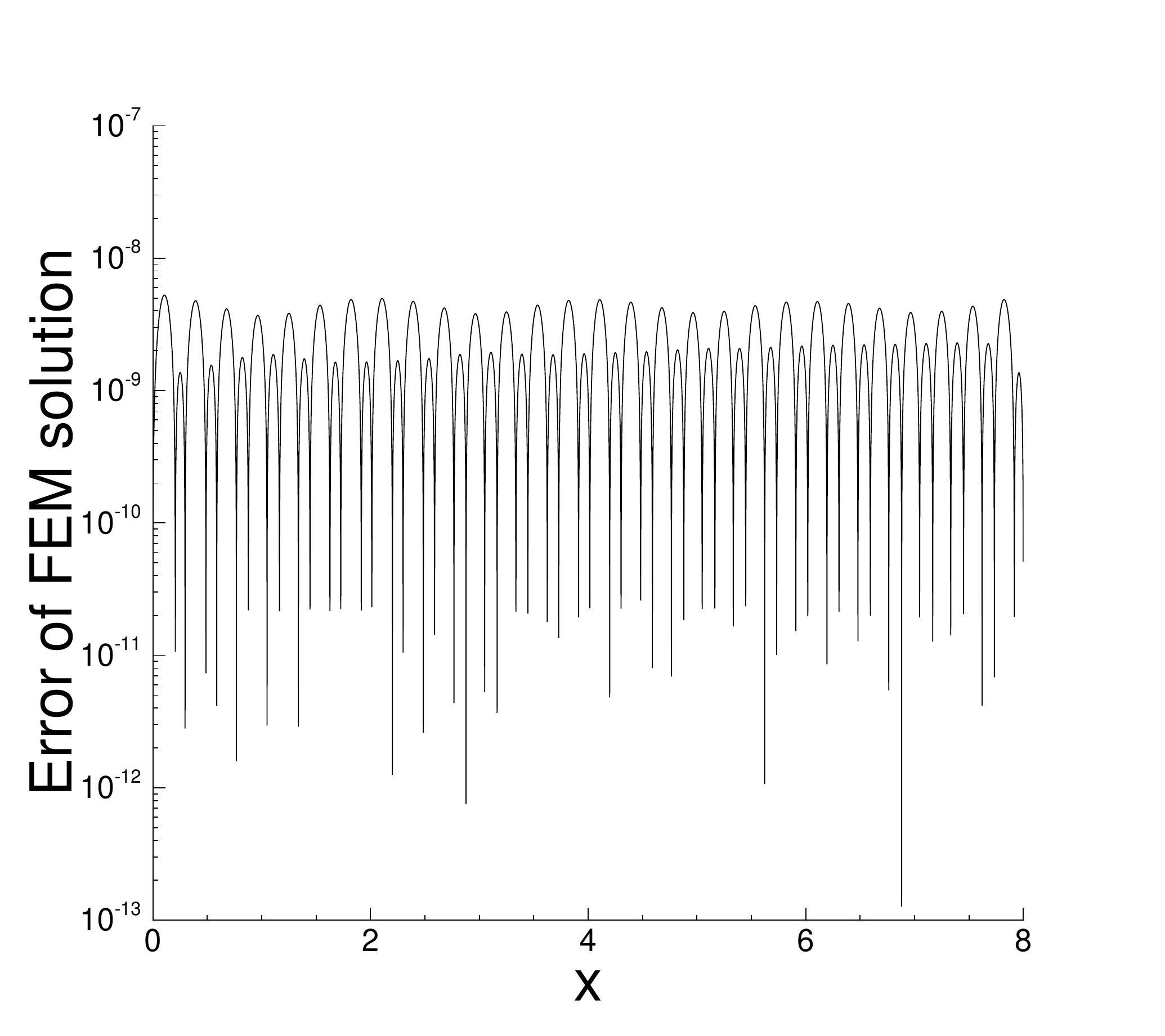}(b)
  }
  \caption{Nonlinear Helmholtz equation: profiles of the solution (a) and its
    absolute error (b) computed using the finite element method (FEM) with
    $200,000$ uniform elements.
  }
  \label{fg_nhm_8}
\end{figure}

\begin{table}
  \centering
  \begin{tabular}{l|lllllll}
    \hline
    method & elements & sub-domains & $Q$ & $M$ & maximum & rms
    & wall-time \\
    & & & & & error & error & (seconds) \\
    \hline
    locELM (NLSQ-perturb) & -- & 4 & $100$ & $200$ & $1.45e-9$ & $2.34e-10$  & $7.7$ \\
    & -- & 4 & $125$ & $200$ & $3.96e-11$ & $7.02e-12$ & $10.6$ \\
    \hline
    FEM & $200,000$ & -- & -- & -- & $5.26e-9$ & $1.37e-9$ & $4.7$ \\
    & $400,000$ & -- & -- & -- & $1.31e-9$ & $3.43e-10$ & $8.8$ \\
    & $800,000$ & -- & -- & -- & $3.29e-10$ & $8.57e-11$ & $18.1$ \\
    \hline
  \end{tabular}
  \caption{Nonlinear Helmholtz equation: comparison between 
    locELM  and FEM,
    in terms of the maximum/rms errors in the domain 
    and the training/computation time.
    The problem settings  correspond to those of Figures
    \ref{fg_nhm_1}(a,b) and \ref{fg_nhm_8}.
  }
  \label{tab_nhm_9}
\end{table}

Let us now compare the current locELM method with the finite element method for
solving the nonlinear Helmholtz equation.
Figure \ref{fg_nhm_8} shows the profiles of the finite element
solution and its absolute error
against the analytic solution, computed on a mesh of $200,000$ uniform elements.
The finite element method is again
implemented using the FEniCS library in Python, and the nonlinear algebraic
equation is solved using a Newton iteration.
The FEM result is observed to be accurate, with an error level
on the order $10^{-9}$.
In Table \ref{tab_nhm_9} we compare the locELM method and the finite element method
with regard to the accuracy and the computational cost.
The table lists the maximum and rms errors in the domain and the wall time
of the training or computation, obtained using locELM
and FEM on several sets of parameters corresponding to different
simulation resolutions.
One can observe that locELM exhibits a comparable, and generally superior,
performance to FEM. For example, the locELM case with $(Q,M)=(100,200)$
has a computational cost comparable to the FEM case with $400,000$ elements, and
their error levels are also comparable.
The locELM case with $(Q,M)=(125,200)$ has a lower cost ($\sim 10$ seconds)
than the FEM case with $800,000$ elements ($\sim 18$ seconds), and
 also has considerably smaller errors, by an order
of magnitude, than the latter.



\subsubsection{Nonlinear Spring Equation}


\begin{figure}
  \centerline{
    \includegraphics[width=2.2in]{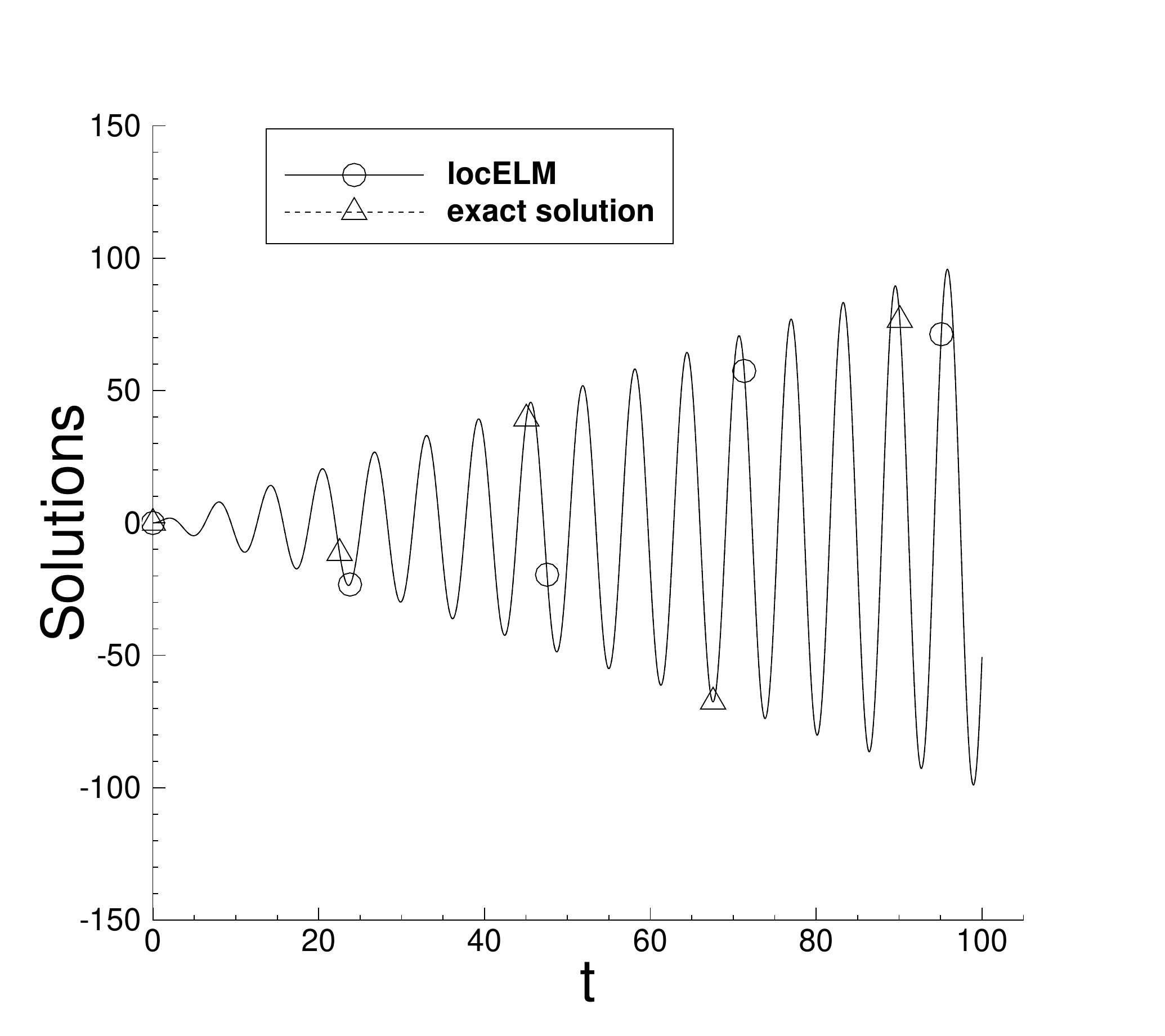}(a)
    \includegraphics[width=2.2in]{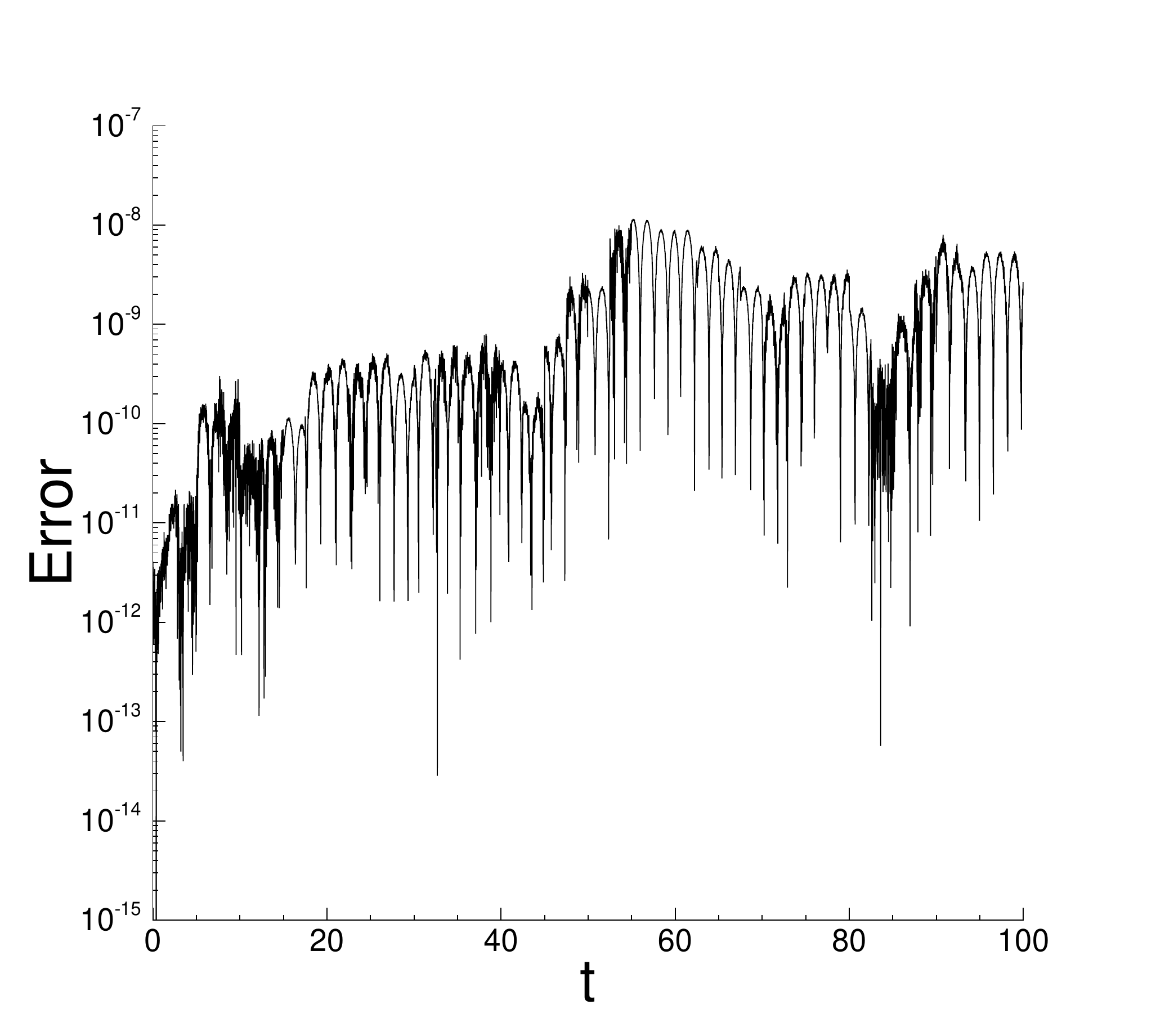}(a)
  }
  \caption{Nonlinear spring: time histories of (a) the locELM solution
  and (b) its absolute error against the exact solution.
  $40$ uniform time blocks are used.
  }
  \label{fg_ns_1}
\end{figure}

In the next example we test the locELM method using an initial value problem,
the nonlinear spring. The goal here is to assess the performance
of the locELM method together with the block time marching scheme,
especially for long-time dynamic simulations.

Consider the temporal domain, $\Omega=[0, t_f]$, and the following
initial value problem on this domain,
\begin{subequations}
\begin{align}
&
\frac{d^2u}{dt^2} + \omega^2 u + \alpha\sin(u) = f(t), \label{eq_ns_1} \\
&
u(0) = u_0, \\
& \left.\frac{du}{dt}\right|_{t=0} = v_0, \label{eq_ns_2} 
\end{align}
\end{subequations}
where $u(t)$ is the displacement, $f(t)$ is an imposed external force,
$\omega$ and $\alpha$ are constant parameters, $u_0$ is the initial
displacement, and $v_0$ is the initial velocity.
The parameters in the above domain and problem specifications assume
the following values in this subsection,
\begin{equation*}
\omega = 2, \quad
\alpha = 0.5, \quad
t_f = 100,\ \text{or}\ 15,\ \text{or}\ 2.5.
\end{equation*}
We choose the external force $f(t)$ such that the following function
satisfies the equation \eqref{eq_ns_1},
\begin{equation}\label{eq_ns_3}
u(t) = t\sin(t).
\end{equation}
We set the initial displacement and the initial velocity both to zero,
i.e.~$u_0=0$ and $v_0=0$.
Under these settings, the initial value problem
consisting of equations \eqref{eq_ns_1}--\eqref{eq_ns_2} has
the solution given by \eqref{eq_ns_3}.


We employ  the locELM method and the block time marching scheme
from Section \ref{sec:tnleq} to solve this initial value problem.
We partition the domain $[0,t_f]$ into $N_b$ uniform
time blocks, and solve this initial value problem on
each time block individually and successively.
For the computation within each time block, we use a single sub-domain
in the simulation,
as the amount of  data involved in is quite small because
the function does not depend on space.
We enforce the equations on $Q$ uniform collocation points
within each time block.
Accordingly, we employ a single neural network within each
time block for this problem.
The neural network consists of an input layer of one
node (representing the time $t$),
a single hidden layer with a width of $M$ nodes and the $\tanh$ activation
function, and an output layer of one node (representing the
solution $u$).
The output layer is assumed to be linear (no activation function) and
contains no bias.
As in previous sections, we incorporate an affine mapping operation
right behind the input layer to normalize the input $t$ data
to the interval $[-1,1]$ for each time block.
The weight and bias coefficients in the hidden layer of the neural network
are pre-set to uniform random values generated on the interval
$[-R_m,R_m]$.
A fixed seed value $1234$ is used for the random number generator.
%
%
We employ the NLSQ-perturb method from Section~\ref{sec:nonl_steady}
for computing the resultant nonlinear algebraic problem.
The initial guess of the solution is set to zero.
In the event the random perturbation is triggered, we employ $\delta=1.0$ and
$\xi_2=1$ (see Algorithm~\ref{alg:alg_1} and Remark~\ref{rem_9})
for generating the random perturbations in the tests of
this subsection.

The locELM simulation parameters  include
the number of time blocks $N_b$,
the number of collocation points per time block $Q$,
the number of training parameters per time block $M$ (i.e.~the number of nodes in
the hidden layer of the neural network), and
the maximum magnitude of the random coefficients $R_m$.

\begin{figure}
  \centerline{
    \includegraphics[width=1.5in]{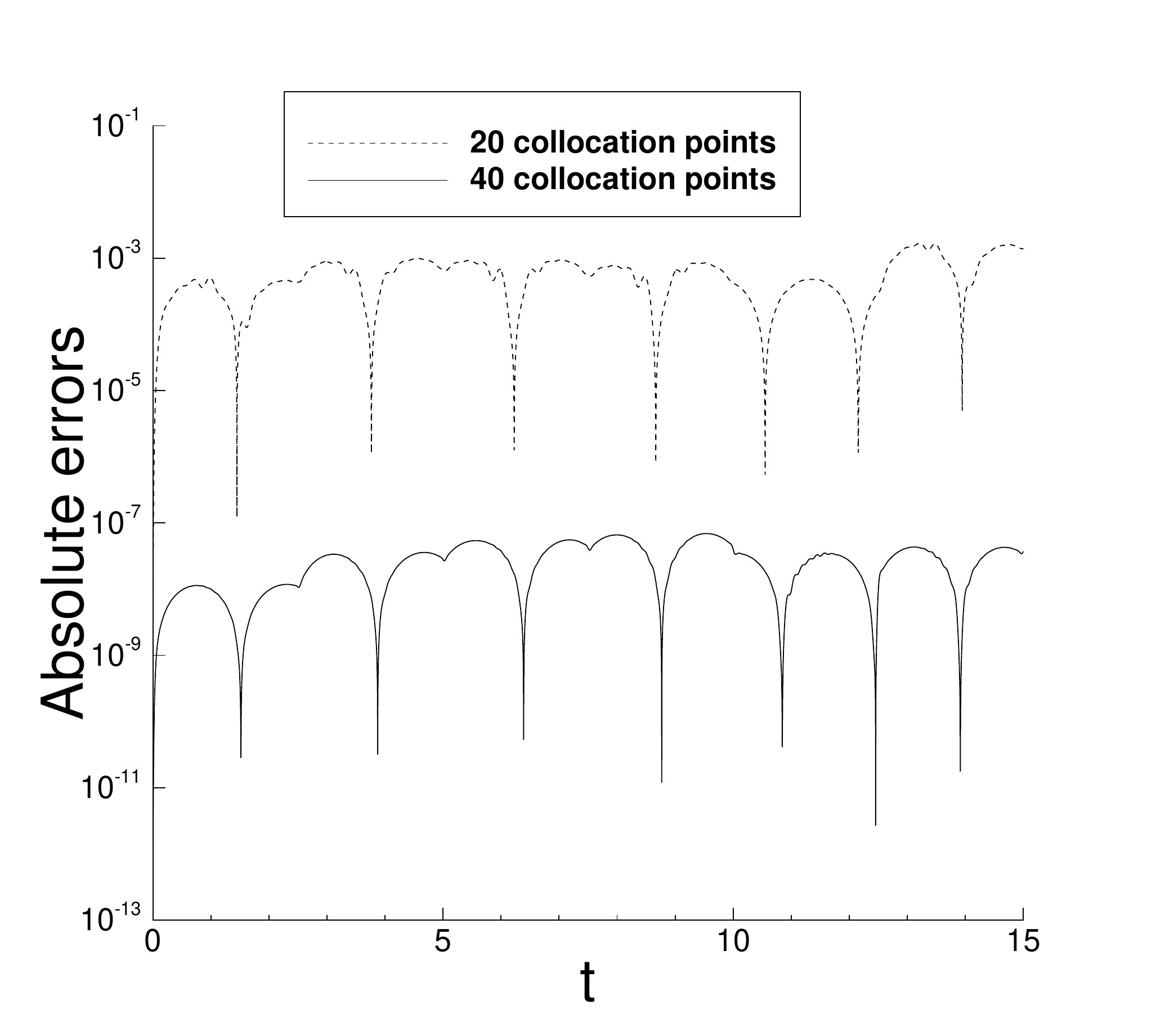}(a)
    \includegraphics[width=1.5in]{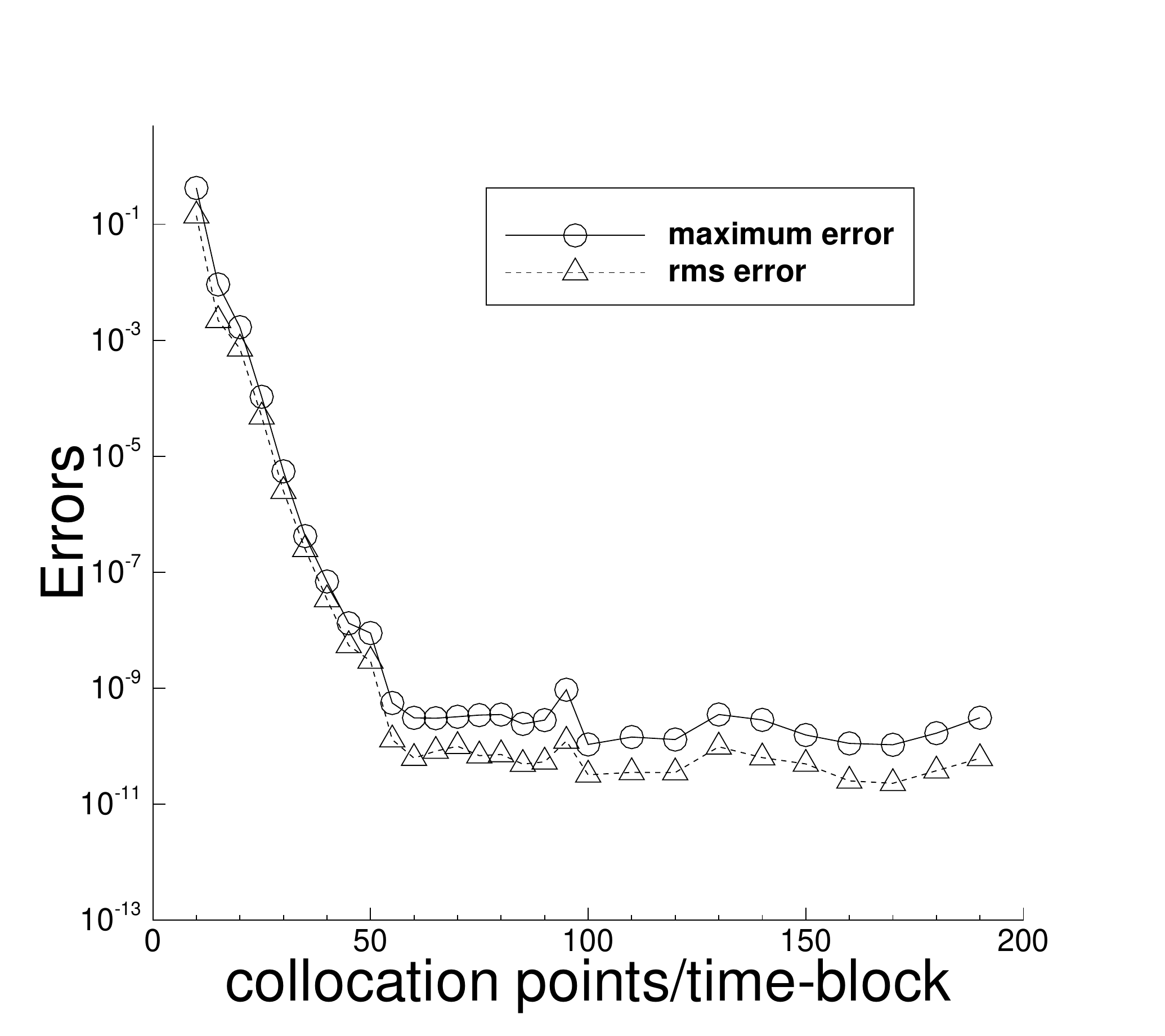}(b)
    \includegraphics[width=1.5in]{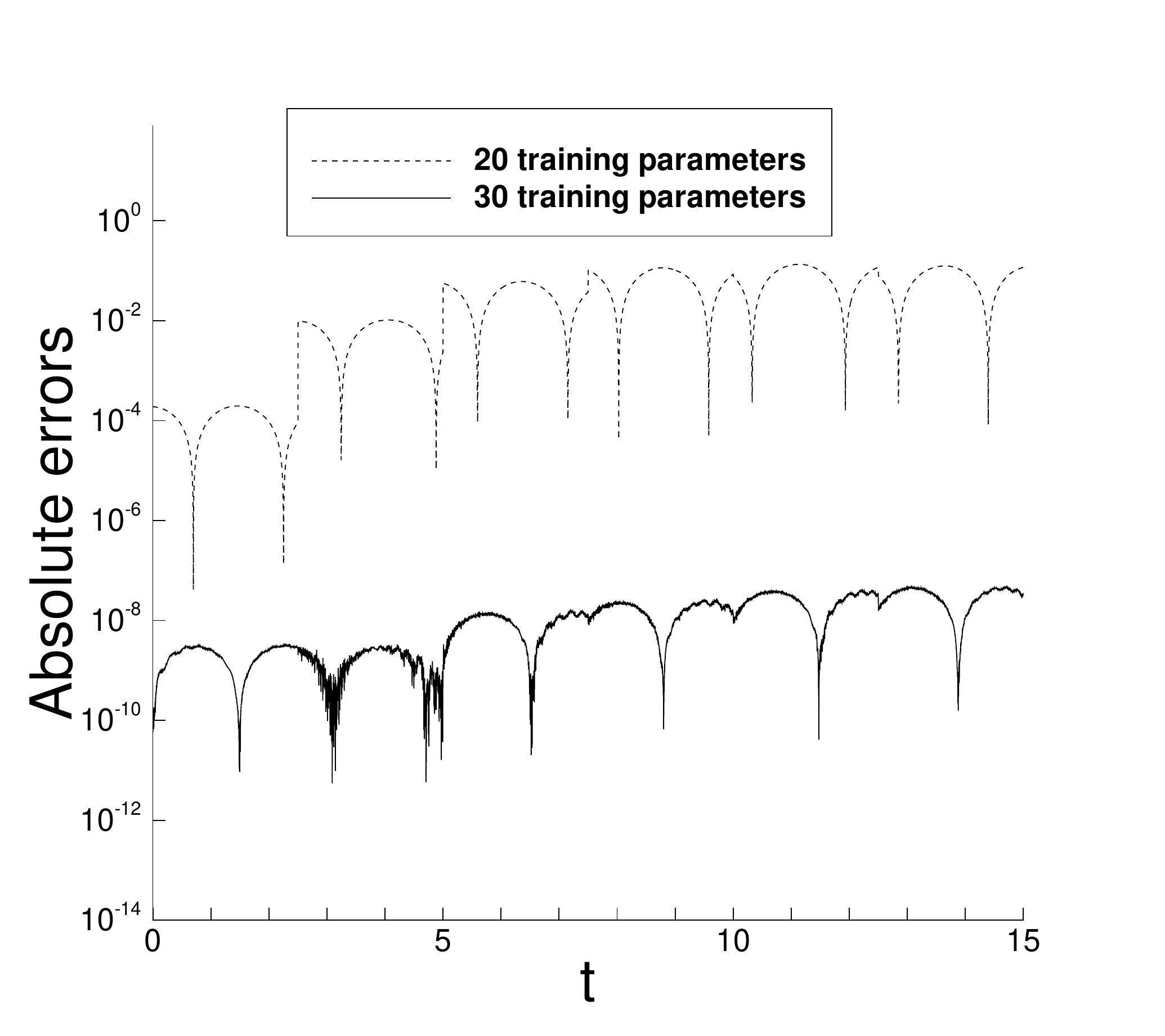}(c)
    \includegraphics[width=1.5in]{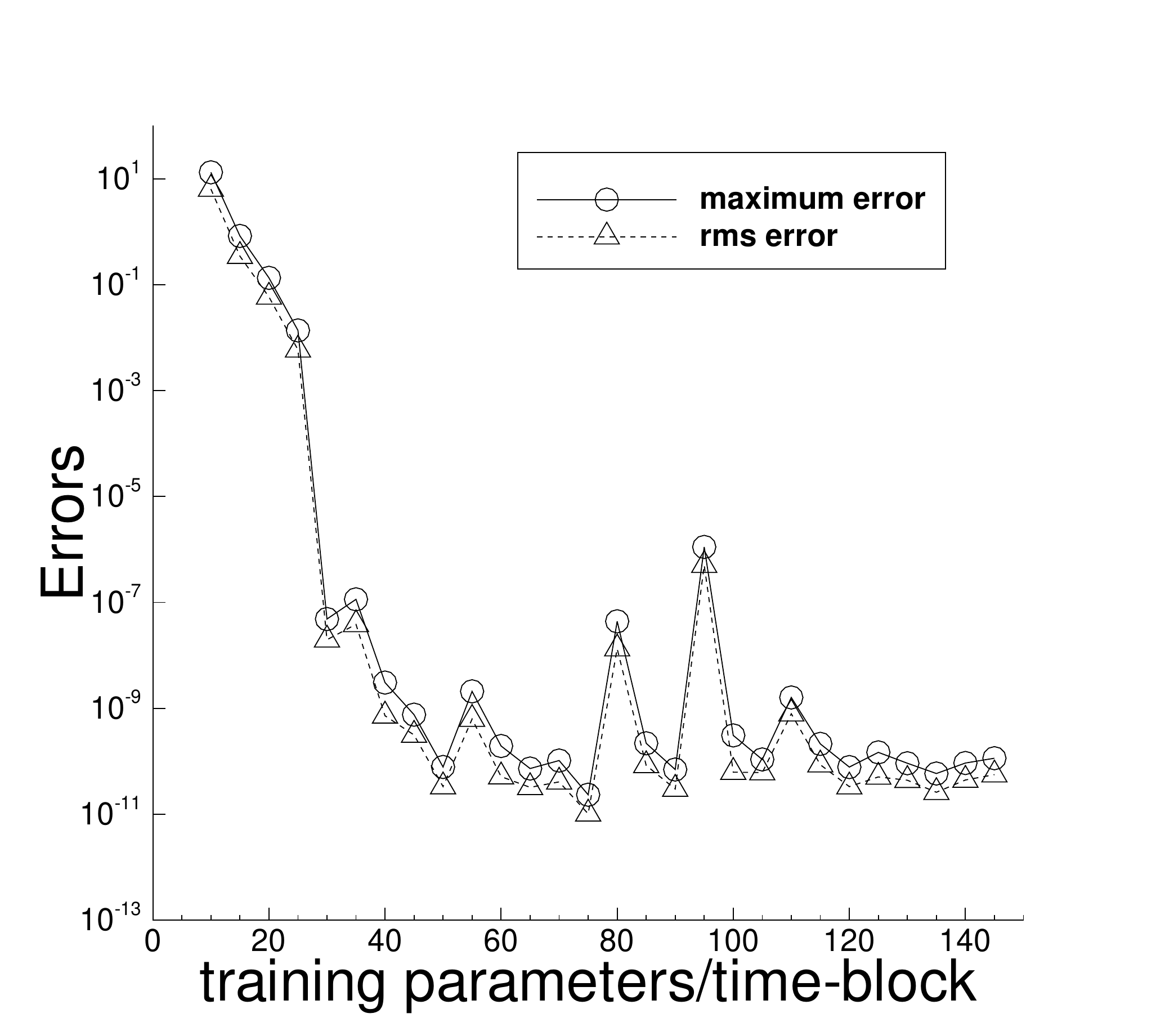}(d)
  }
  \caption{Nonlinear spring: (a) Error histories obtained with $20$ and $40$
    collocation points per time block.
    (b) The maximum/rms errors in the domain versus the number of
    collocation points per time block.
    In (a) and (b), the number of training parameters per time block is fixed at $100$.
    (c) Error histories obtained with $20$ and $30$ training parameters
    per time block.
    (d) The maximum/rms errors versus the number of training parameters per time block.
    In (c) and (d) the number of the collocation points per time block is fixed at $60$.
    }
    \label{fg_ns_2}
\end{figure}


Figure \ref{fg_ns_1} shows the time histories of the displacement
and its absolute error obtained using locELM in a fairly long-time simulation.
The time history of the exact solution given by~\eqref{eq_ns_3} has also
been shown in Figure \ref{fg_ns_1}(a) for comparison.
In this test the domain size is set to $t_f=100$. We have employed
$N_b=40$ uniform time blocks within the domain, $Q=60$ uniform collocation
points per time block, $M=100$ training parameters per time block,
and $R_m=5.0$ when generating the random weight/bias coefficients
for the hidden layer of the neural network.
It is evident that the current locELM method has captured the solution
very accurately, with the maximum level of the absolute error on
the order $10^{-8}$ over the entire domain.

Figure \ref{fg_ns_2} illustrates the effect of the number of
degrees of freedom (collocation points, training parameters)
on the simulation accuracy.
In this group of tests the temporal domain size is set to $t_f=15$, and
we employ $N_b=6$ time blocks within the domain.
Figure \ref{fg_ns_2}(a) shows the absolute-error histories
of the locELM solution against the exact solution, obtained using
$20$ and $40$ collocation points per time block.
Figure \ref{fg_ns_2}(b) shows the maximum and rms errors in
the overall domain obtained with different numbers of collocation
points in the locELM simulation.
The number of training parameters per time block is fixed at $M=100$
with the tests in these two plots.
The errors can be observed to decrease exponentially
as the number of collocation points per time block increases
(when below around $60$), and then become stagnant
as the number of collocation points increases further.
Figure \ref{fg_ns_2}(c) shows time histories of the absolute
errors corresponding to $20$ and $40$ training parameters per time block.
Figure \ref{fg_ns_2}(d) shows the maximum/rms errors in the overall
domain, obtained with different numbers of training parameters per
time block.
In the tests of these two plots, the number of collocation points
per time block has been fixed at $Q=60$.
The convergence  with respect to the training parameters
is not as regular as that for the collocation points.
Nonetheless, one can see that the errors approximately decrease exponentially
with increasing number of training parameters (when below $50$),
and then they essentially stagnate as the number of training parameters
increases further.


\begin{figure}
  \centerline{
    \includegraphics[width=2.2in]{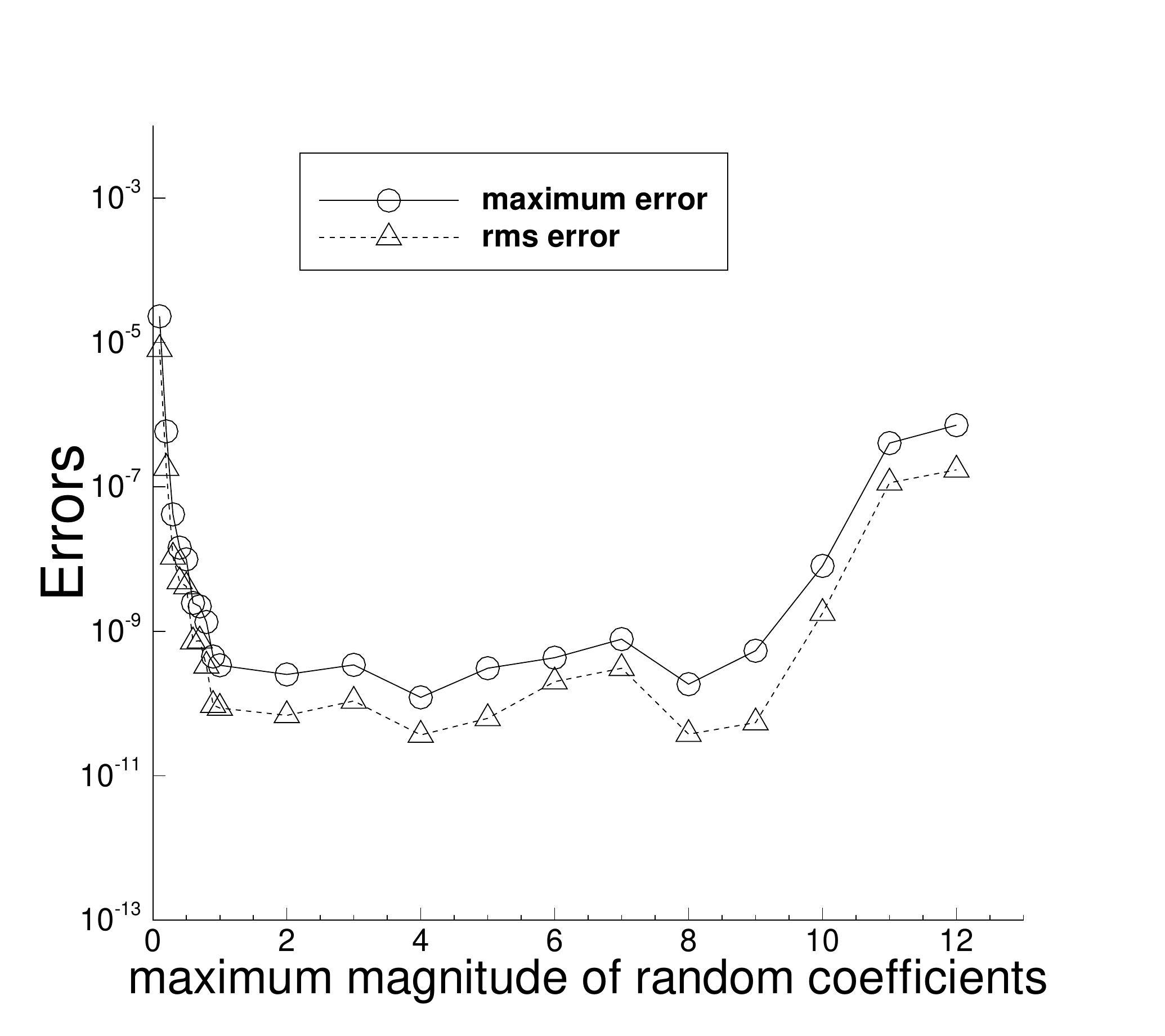}
  }
  \caption{Nonlinear spring: The maximum and rms errors in the overall domain
    as a function of $R_m$, the maximum magnitude of the random coefficients.
  }
  \label{fg_ns_3}
\end{figure}


Figure \ref{fg_ns_3} demonstrates the effect of $R_m$, the maximum magnitude
of the random coefficients, on the simulation accuracy.
In this set of tests, the temporal domain size is $t_f=15$.
We have employed $N_b=6$ uniform time blocks in the domain,
$Q=60$ uniform collocation points in each time block, and
$M=100$ training parameters per time block.
The value of $R_m$ is varied systematically in the tests.
In this figure we plot the maximum and rms errors in the
overall domain corresponding to different $R_m$ values.
The characteristics observed here are consistent with those from previous subsections.
The locELM method has a better accuracy with $R_m$ in a range of moderate
values, and in this case approximately $R_m=1 \sim 9$.
The results are less accurate if $R_m$ is very large or very small.



\begin{figure}
\centerline{
\includegraphics[width=2.2in]{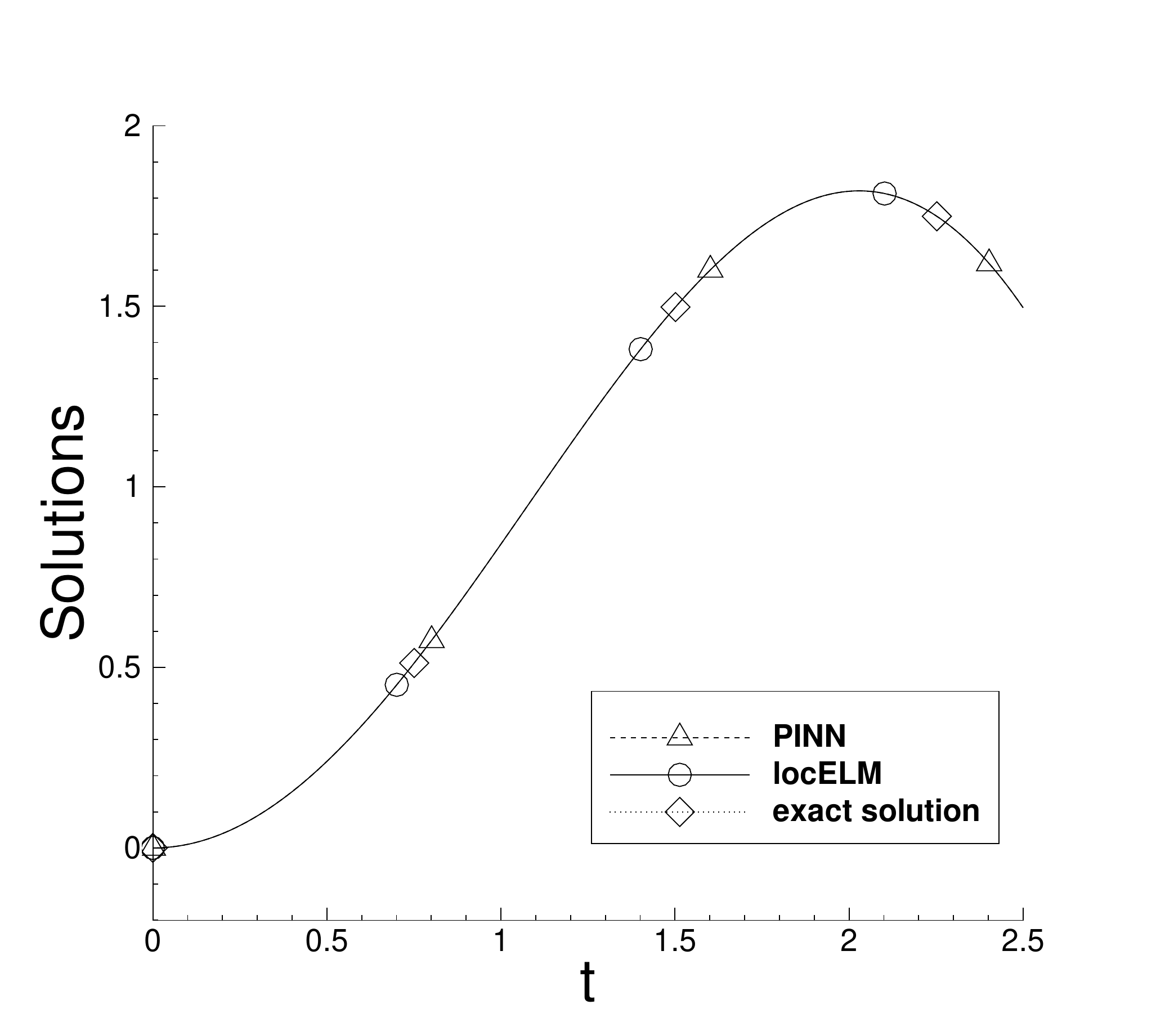}(a)
\includegraphics[width=2.2in]{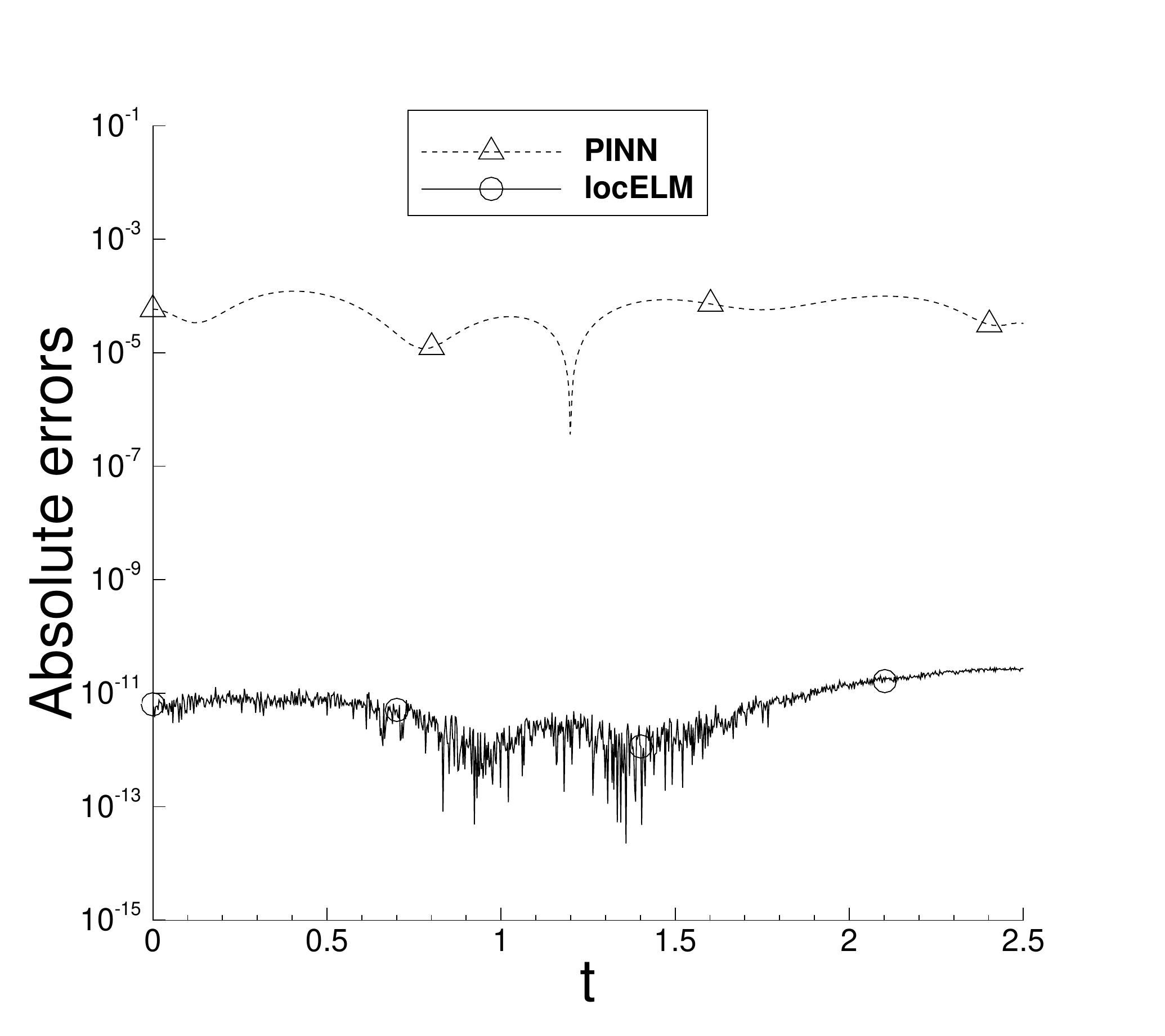}(a)
}
\caption{Comparison between locELM and PINN (nonlinear spring):
Time histories of (a) the solutions and (b) their absolute errors,
computed using PINN~\cite{RaissiPK2019} with the Adam optimizer
and using locELM with the NLSQ-perturb method.
The temporal domain size is $t_f=2.5$.
}
\label{fg_ns_4}
\end{figure}

\begin{table}
\centering
\begin{tabular}{lllll}
\hline
method & maximum error & rms error & epochs/iterations & training time (seconds) \\
PINN (Adam) & $1.21e-4$ & $6.71e-5$ & $20,000$ & $26.3$ \\
locELM (NLSQ-perturb) & $2.82e-11$ & $1.12e-11$ & $48$ & $0.34$ \\
\hline
\end{tabular}
\caption{Nonlinear spring: Comparison between locELM and PINN in terms of
the maximum/rms errors in the domain, the number of epochs or nonlinear iterations in
the training, and the network training time.
The problem settings and the simulation parameters correspond to
those of Figure \ref{fg_ns_4}.
}
\label{tab_ns_5}
\end{table}

Let us next compare the current locELM method with PINN~\cite{RaissiPK2019}
for solving the nonlinear spring equation.
Figure \ref{fg_ns_4} shows a comparison of the time histories of
the solutions and their absolute errors obtained using
PINN with the Adam optimizer and using the current locELM method
with NLSQ-perturb.
In this group of tests, the temporal domain size is set to $t_f=2.5$.
In the PINN simulation, the neural network consists of
an input layer of one node (representing $t$),
three hidden layers with a width of $10$ nodes and
the $\tanh$ activation function in each layer,
and an output layer of one node (representing the solution $u$).
The input data consists of $500$ uniform collocation points
from the domain $[0, t_f]$. The neural network has been
trained using the Adam optimizer for $20,000$ epochs,
with the learning rate decreasing from $0.01$ at the beginning
to $1e-5$ at the end of the training.
In the locELM simulation, we employ a single time block ($N_b=1$)
in the domain, $Q=60$ uniform collocation points within the
time block, $M=100$ training parameters in the time block,
a single hidden layer in the neural network, and $R_m=5.0$
for generating the random weight/bias coefficients.
Figure \ref{fg_ns_4} demonstrates that both PINN and locELM
have captured the solution accurately, but
the error of the locELM result is considerably
smaller than that of PINN.

Table \ref{tab_ns_5} provides a further comparison
between locELM and PINN in terms of their accuracy and
computational cost.
The problem settings and the simulation parameters here correspond to
those of Figure \ref{fg_ns_4}.
We have listed the maximum and rms errors of the PINN and locELM
results in the overall domain, the number of epochs or nonlinear
iterations in the training, and the network training time.
The data demonstrate that the current locELM method is
much more accurate, by six orders of magnitude, than PINN,
and the network training time of locELM is much smaller,
by nearly two orders of magnitude, than that of PINN.

\subsubsection{Viscous Burger's Equation}


In this subsection we further test the locELM method using the
viscous Burger's equation.
Consider the spatial-temporal domain
$\Omega = \{(x,t)\ |\ x\in[a,b],\ t\in[0,t_f] \}$,
and the following initial/boundary value problem with
the Burger's equation,
\begin{subequations}
  \begin{align}
    &
    \frac{\partial u}{\partial t} + u\frac{\partial u}{\partial x}
    =\nu\frac{\partial^2 u}{\partial x^2} + f(x,t),
    \label{eq_bg_1} \\
    &
    u(a,t) = g_1(t), \\
    &
    u(b,t) = g_2(t), \\
    &
    u(x,0) = h(x), \label{eq_bg_2} 
  \end{align}
\end{subequations}
where $u(x,t)$ is the solution to be solved for, the constant $\nu$
denotes the viscosity, $f(x,t)$ is a
prescribed source term, $g_1(t)$ and $g_2(t)$ denote the boundary
distributions, and $h(x)$ is the initial distribution.
We employ the following values for the constant parameters,
\begin{equation*}
  \nu = 0.01, \quad
  a = 0, \quad
  b = 5, \quad
  t_f = 10,\ \text{or}\ 2.5,\ \text{or}\ 0.25.
\end{equation*}
We choose the source term $f(x,t)$ and the boundary/initial distributions
($g_1$, $g_2$ and $h$)
such that the following function
\begin{equation}\label{eq_bg_3}
  \begin{split}
  u(x,t) =& \left(1+\frac{x}{10} \right)\left(1+\frac{t}{10} \right)
  \left[2\cos\left(\pi x+\frac{2\pi}{5} \right)
    + \frac32\cos\left(2\pi x - \frac{3\pi}{5}  \right) \right] 
  \left[2\cos\left(\pi t+\frac{2\pi}{5} \right) \right. \\
    &
    \left.
    + \frac32\cos\left(2\pi t - \frac{3\pi}{5}  \right) \right]
  \end{split}
\end{equation}
satisfies this initial/boundary value problem.


We employ the method presented in Section \ref{sec:tnleq} with block time
marching to solve this problem, by restricting the method to one spatial
dimension. The spatial-temporal domain $\Omega$ is partitioned into
$N_b$ uniform blocks in time, and each time block is computed separately
and successively. Within each time block, we further partition its
spatial-temporal domain into $N_x$ uniform sub-domains in $x$
and $N_t$ uniform sub-domains in time, resulting in $N_e=N_xN_t$ uniform sub-domains
per time block. $C^1$ continuity is imposed on the sub-domain boundaries
in the $x$ direction, and $C^0$ continuity is imposed on the sub-domain
boundaries in the temporal direction. Within each sub-domain
we employ a total of $Q=Q_xQ_t$ uniform collocation points, with $Q_x$
uniform collocation points in the $x$ direction and $Q_t$ uniform collocation
points in time.

\begin{figure}
  \centerline{
    \includegraphics[height=2.3in]{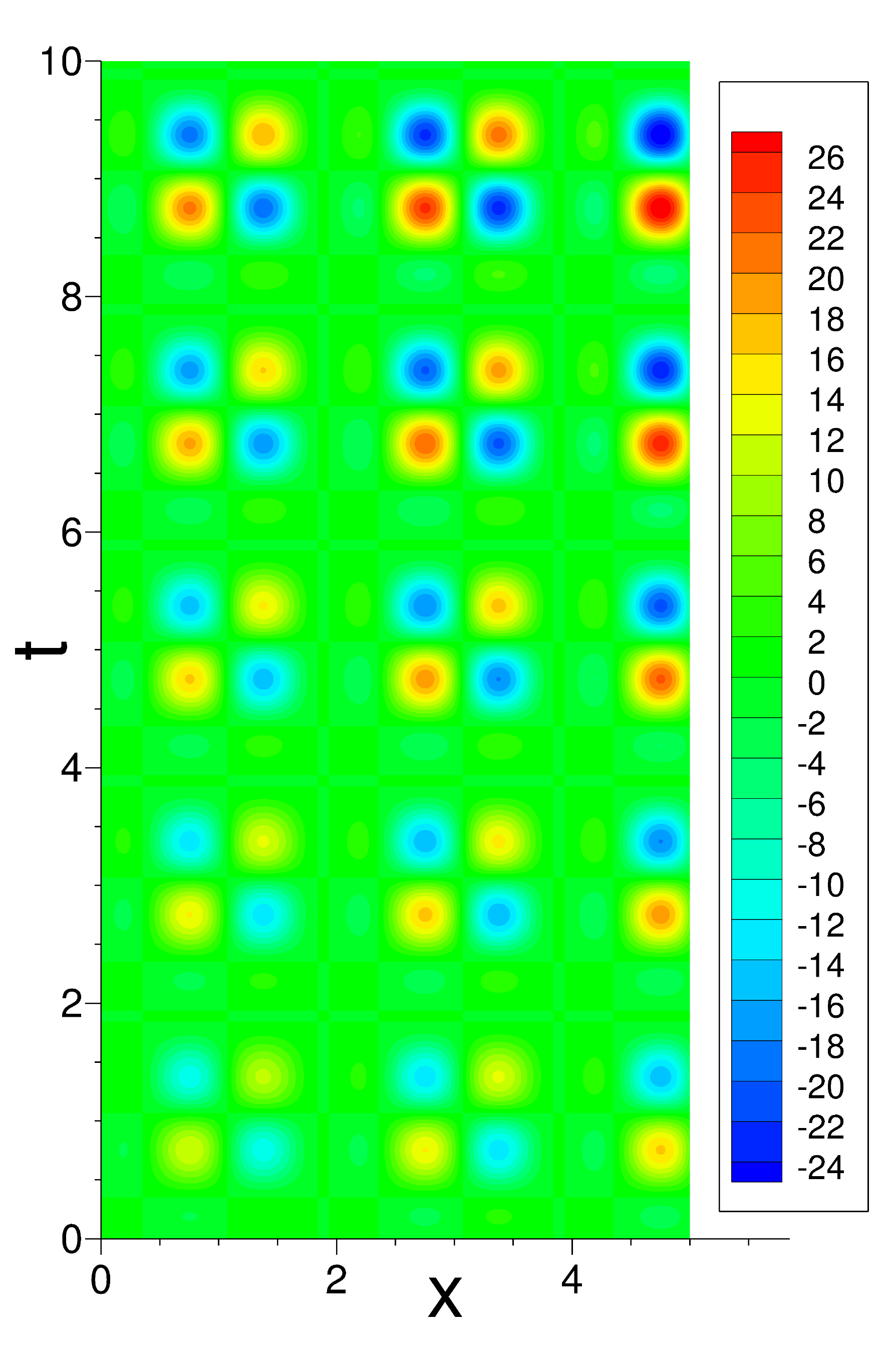}(a)
    \includegraphics[height=2.3in]{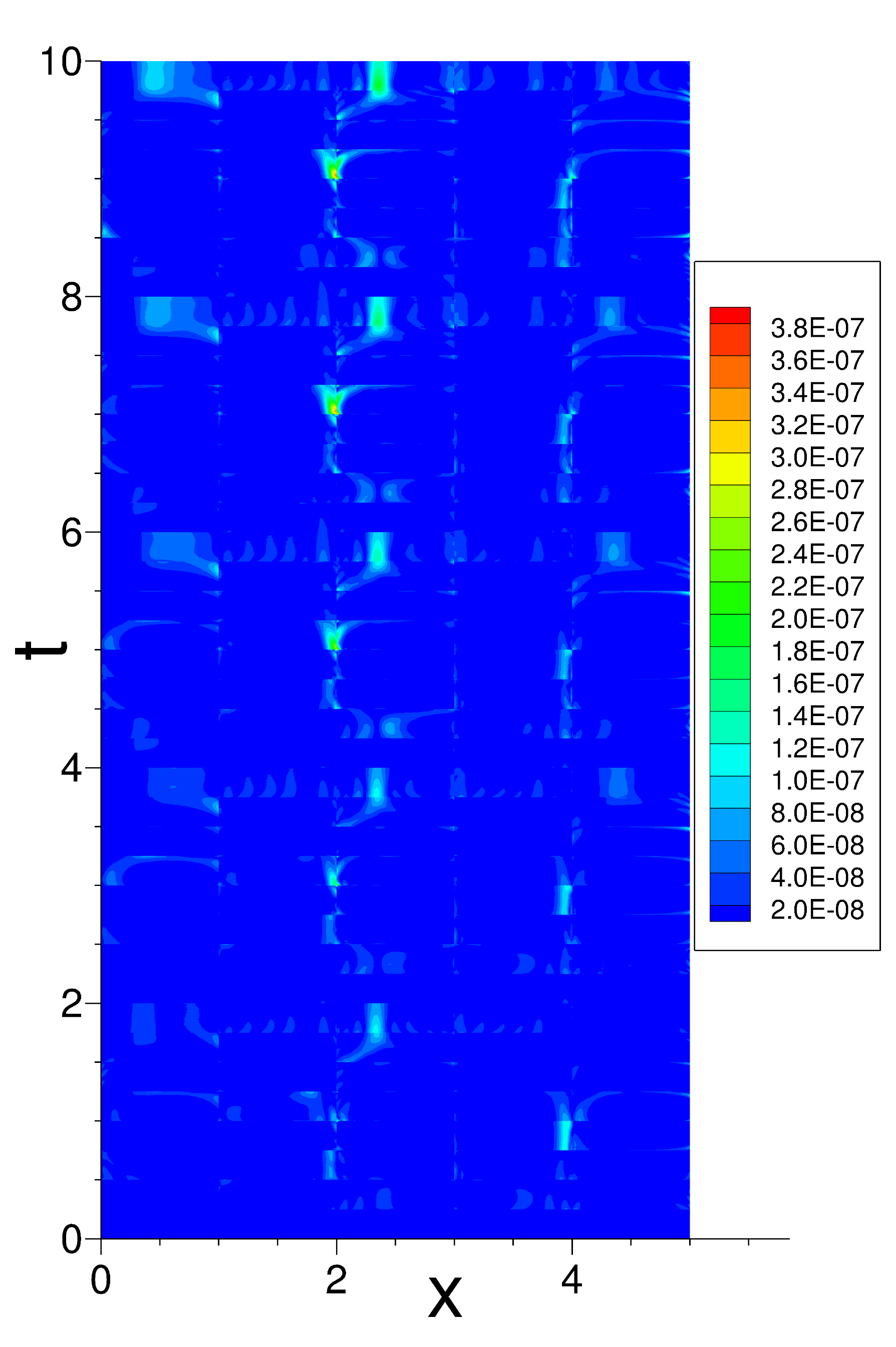}(b)
  }
  \caption{Burgers equation: distributions of (a) the solution, and (b) its
    absolute error in the spatial-temporal plane,
     computed using the current locELM (NLSQ-perturb)
    method.
  }
  \label{fig:burger}
\end{figure}


We employ a local neural network for each sub-domain, leading to a
total of $N_e$ local neural networks in the simulations.
Each local neural network consists of an input layer of two nodes,
representing the $x$ and $t$,
a single hidden layer
with a width of $M$ nodes and the $\tanh$ activation function,
and an output layer of a single node, representing the field solution $u$.
The output layer has no bias and no activation function.
Additionally, an affine mapping operation is incorporated into the network
right behind the input layer to normalize the input $x$ and $t$ data to
the interval $[-1,1]\times[-1,1]$ for each sub-domain.
The weight/bias coefficients in the hidden layer of the local neural networks
are pre-set to uniform random values generated on
the interval $[-R_m,R_m]$, and are fixed during the simulation.
A fixed seed value $22$ is used for the Tensorflow random number generator
for all the tests in this sub-section.

We employ the NLSQ-perturb method from Section \ref{sec:nonl_steady}
for computing the resultant nonlinear algebraic problem in the majority of tests
presented below. The results computed using Newton-LLSQ
are also provided for comparison in some cases. The initial guess
of the solution  is set to zero.
With the NLSQ-perturb method, we employ $\delta=0.5$ and $\xi_2=0$
(see Algorithm \ref{alg:alg_1}
and Remark \ref{rem_9})
for generating the random perturbations in the following tests.

The locELM simulation parameters include
the number of time blocks ($N_b$), the number of sub-domains per time
block ($N_e$, $N_x$, $N_t$), the number of collocation points per
sub-domain ($Q$, $Q_x$, $Q_t$), the number of training parameters
per sub-domain ($M$), and the maximum magnitude of the random coefficients ($R_m$).

\begin{figure}
  \centerline{
    \includegraphics[width=1.5in]{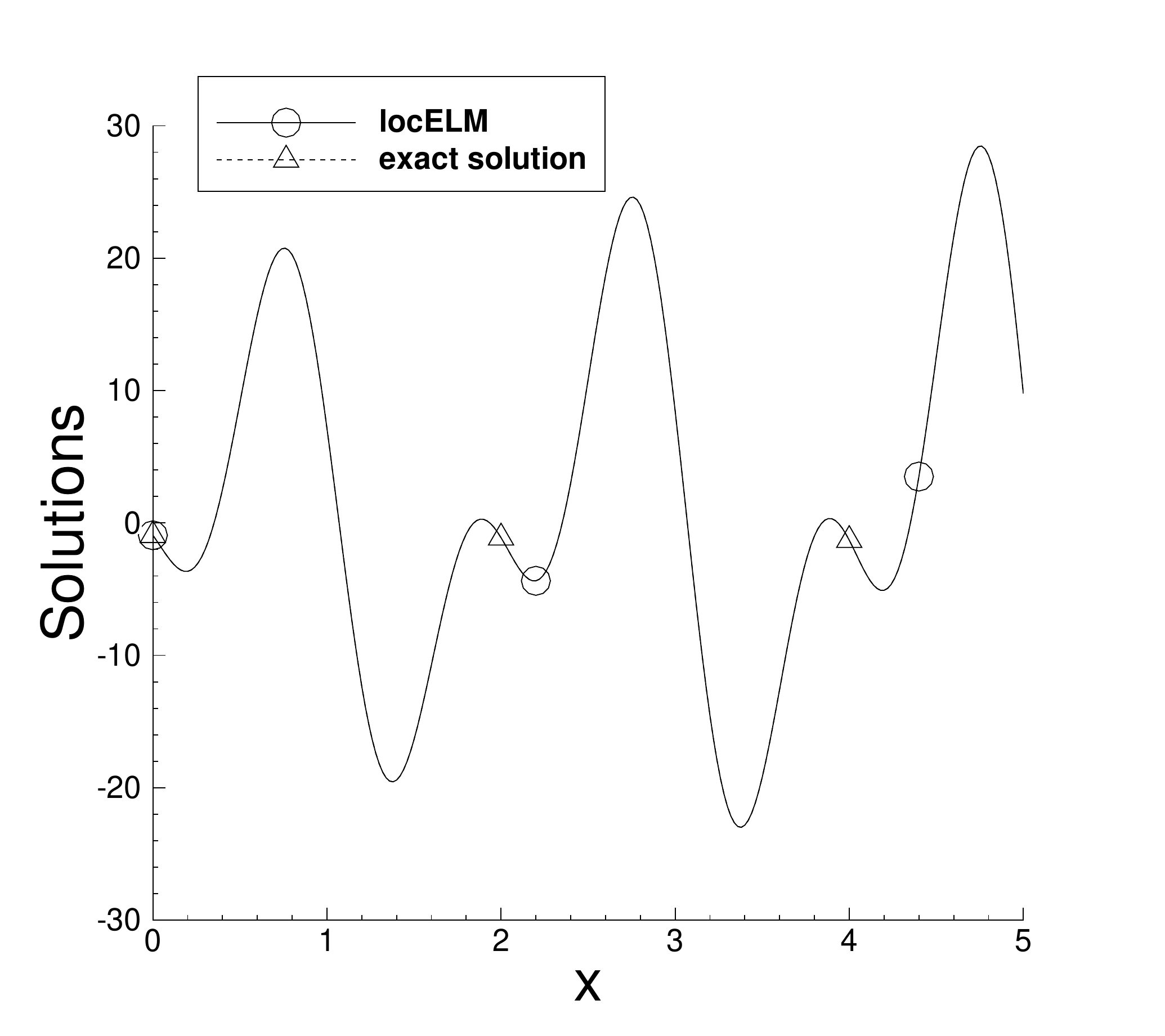}(a)
    \includegraphics[width=1.5in]{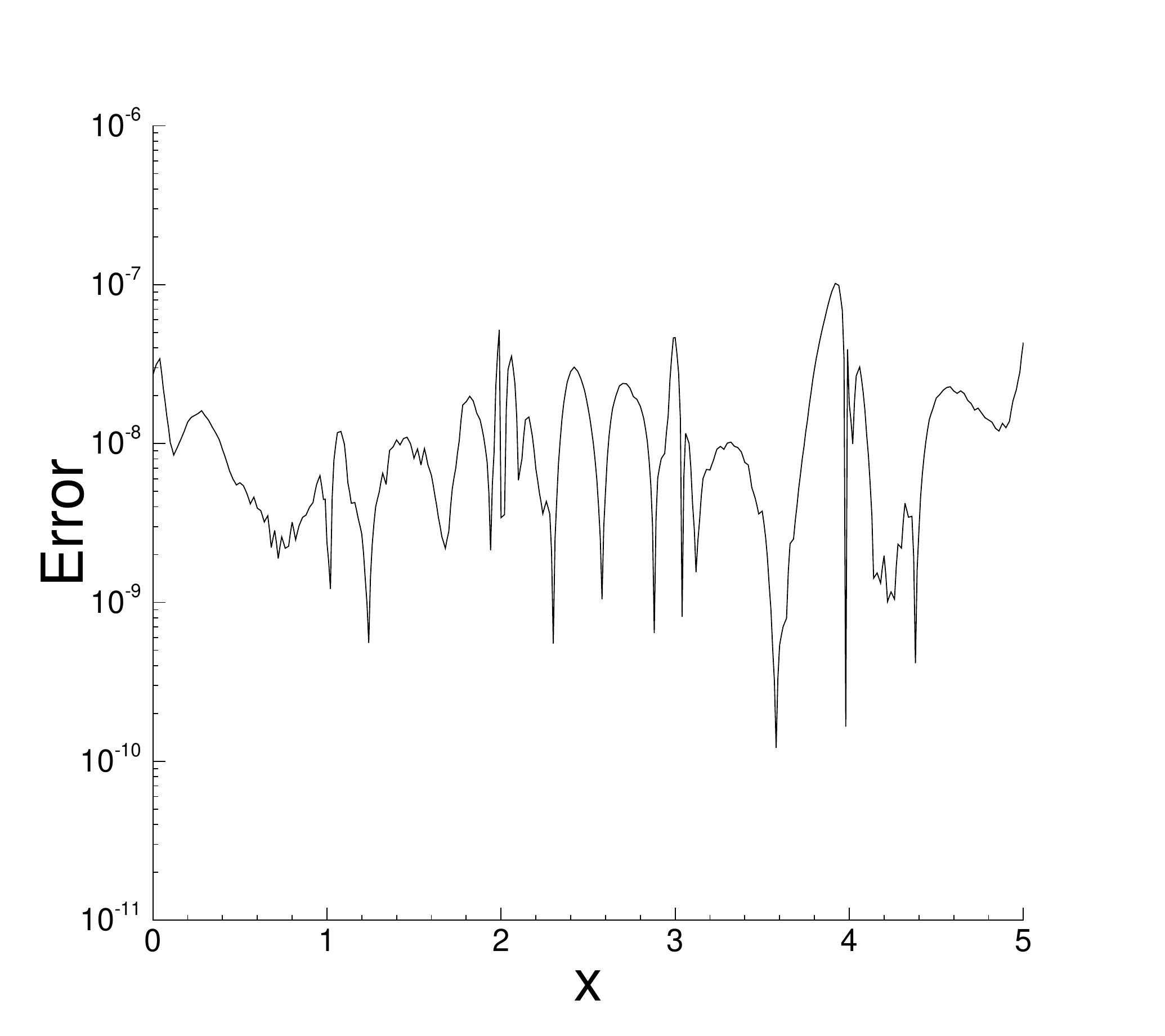}(b)
    \includegraphics[width=1.5in]{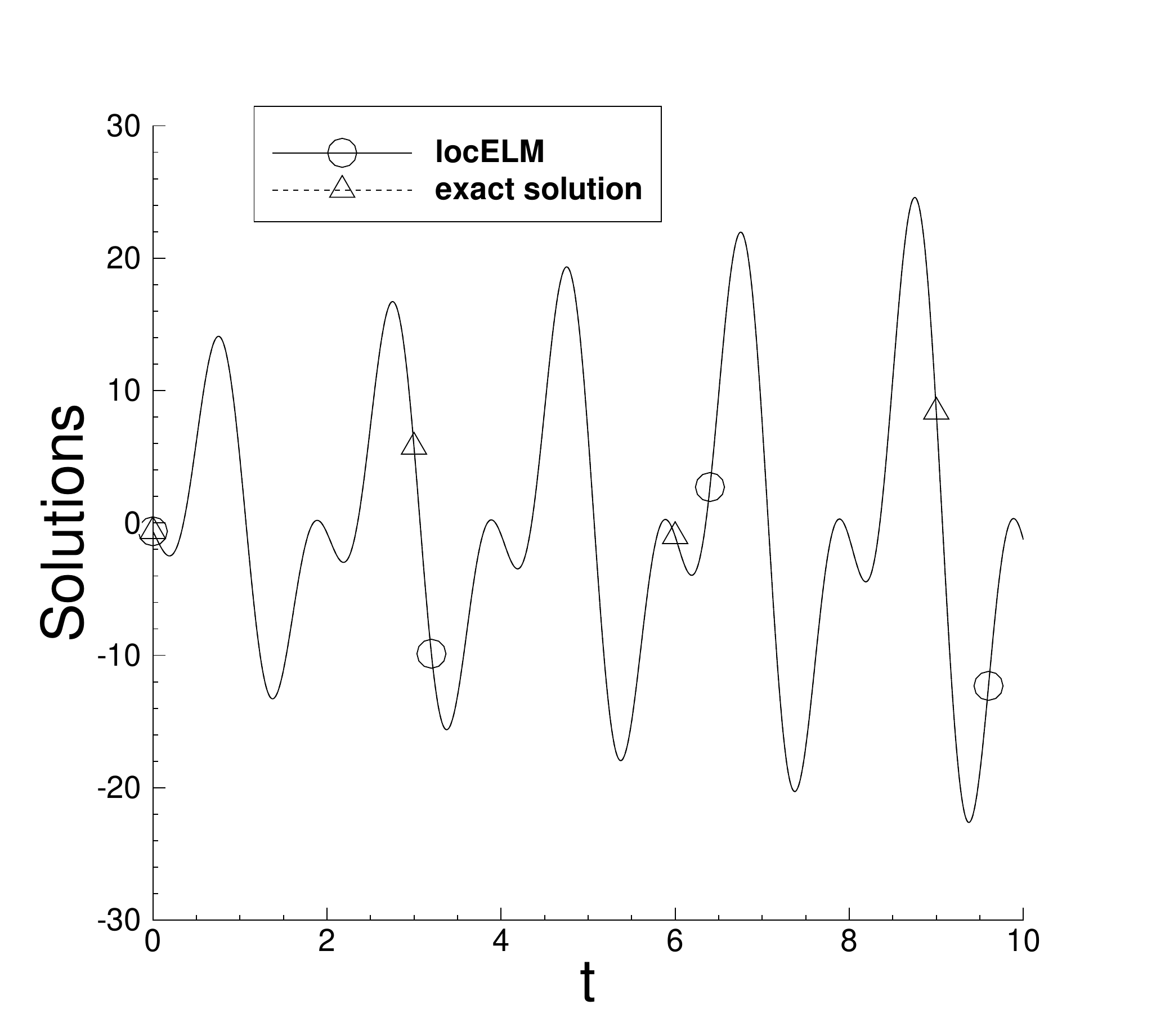}(c)
    \includegraphics[width=1.5in]{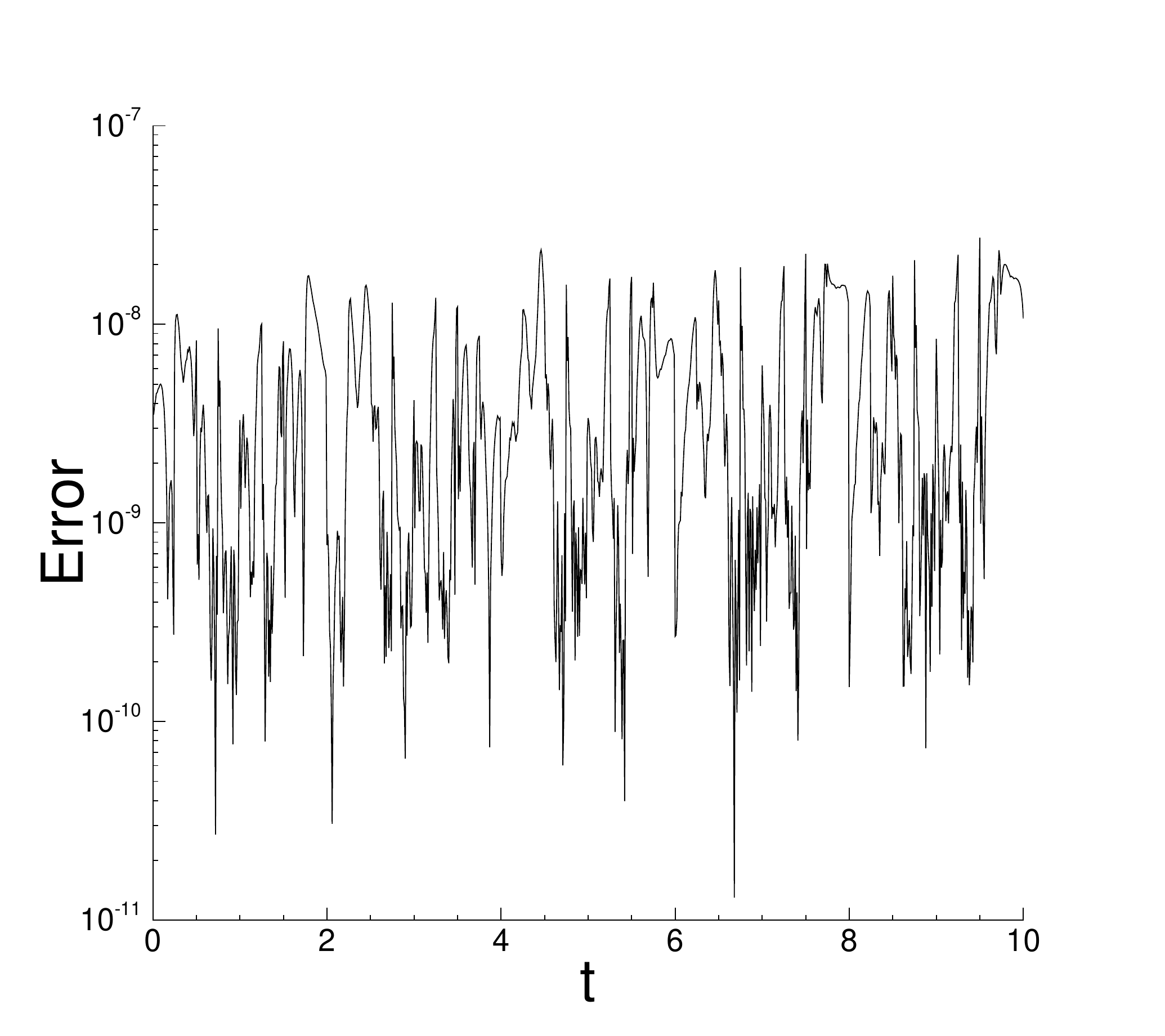}(d)
  }
  \caption{Burger's equation: Profiles of the locELM (NLSQ-perturb) solution (a) and
    its absolute error (b) at $t=8.75$.
    Time histories of the locELM (NLSQ-perturb) solution
    (c) and its absolute error (d) at the
    point $x=2.75$.
    The settings and simulation parameters correspond to those
    of Figure \ref{fig:burger}.
  }
  \label{fg_bg_2}
\end{figure}


Figure \ref{fig:burger} shows distributions of
the solution and its absolute error
in the spatial-temporal plane,
computed using the current locELM method (with NLSQ-perturb).
Here the temporal domain size is set to be $t_f=10$,
and $40$ uniform time blocks ($N_b=40$) are used 
in the spatial-temporal domain.
We have employed  $N_e=5$ uniform sub-domains
with $N_x=5$ and $N_t=1$ within each time block,
$Q=20\times 20$ uniform collocation points per sub-domain
($Q_x=Q_t=20$), $M=200$ training parameters per sub-domain,
and $R_m=0.75$ when generating the random coefficients.
The current method has captured the solution accurately,
with the absolute error on the order $10^{-8}\sim 10^{-7}$
in the overall domain.

Figure \ref{fg_bg_2} further examines  the accuracy of the locELM solution.
The problem settings and the simulation parameters here correspond to those
of Figure \ref{fig:burger}.
Figures \ref{fg_bg_2}(a) and (b) depict the profiles of the locELM (NLSQ-perturb) solution
and its absolute error at the time $t=8.75$. The exact solution profile
at this time instant is also shown in Figure \ref{fg_bg_2}(a).
The locELM solution profile exactly overlaps with
that of the exact solution, and the absolute error is around the level
$10^{-10}\sim 10^{-7}$.
Figures \ref{fg_bg_2}(c) and (d) show the time histories of the locELM (NLSQ-perturb)
solution and its absolute error at the point $x=2.75$.
The time history of the exact solution at this point is also
shown in Figure \ref{fg_bg_2}(c).
The simulated signal overlaps with that of the exact
signal, and the absolute error can be observed to fluctuate around the level
$10^{-10}\sim 10^{-8}$.

\begin{figure}
  \centerline{
    \includegraphics[width=2.2in]{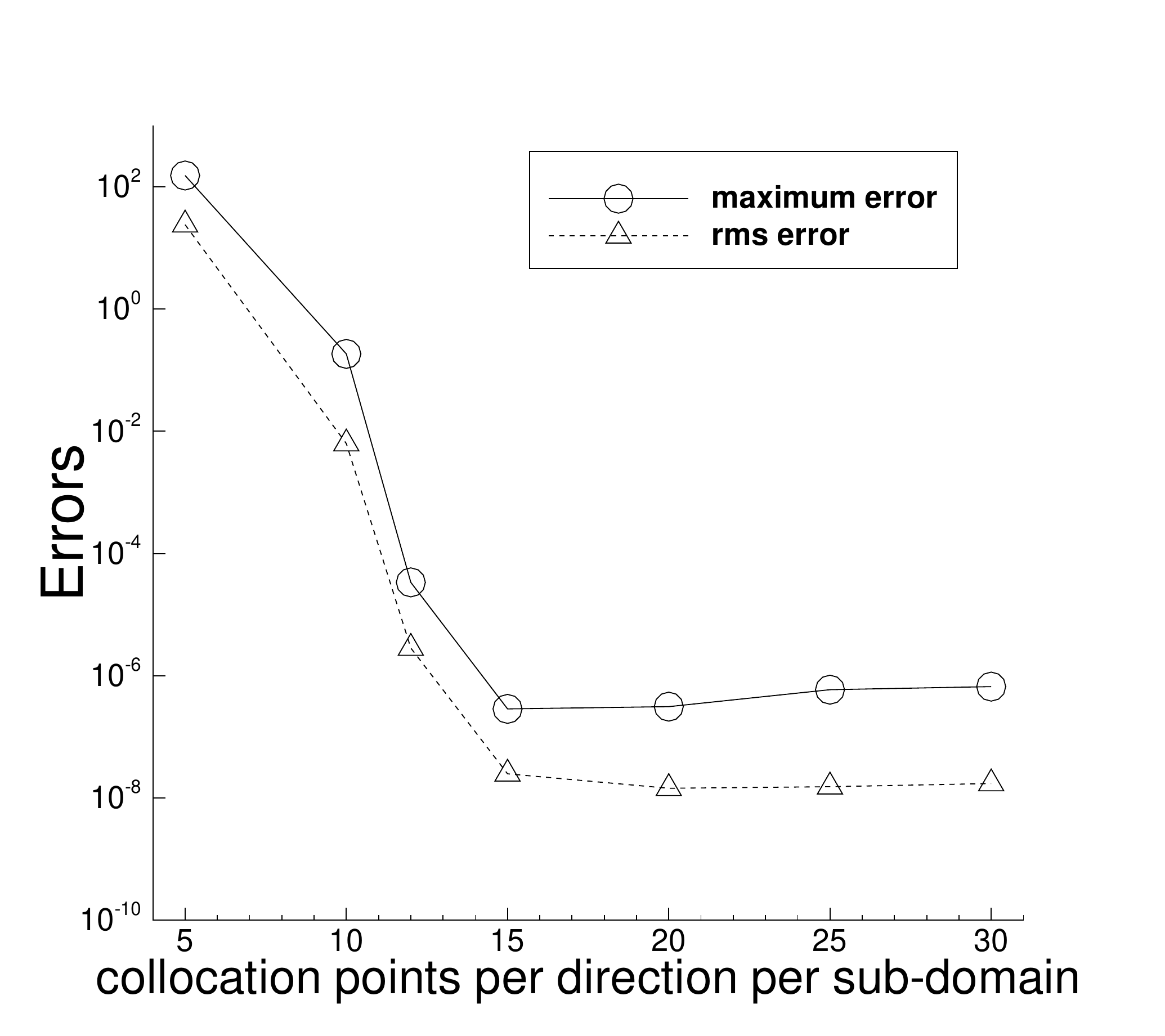}(a)
    \includegraphics[width=2.2in]{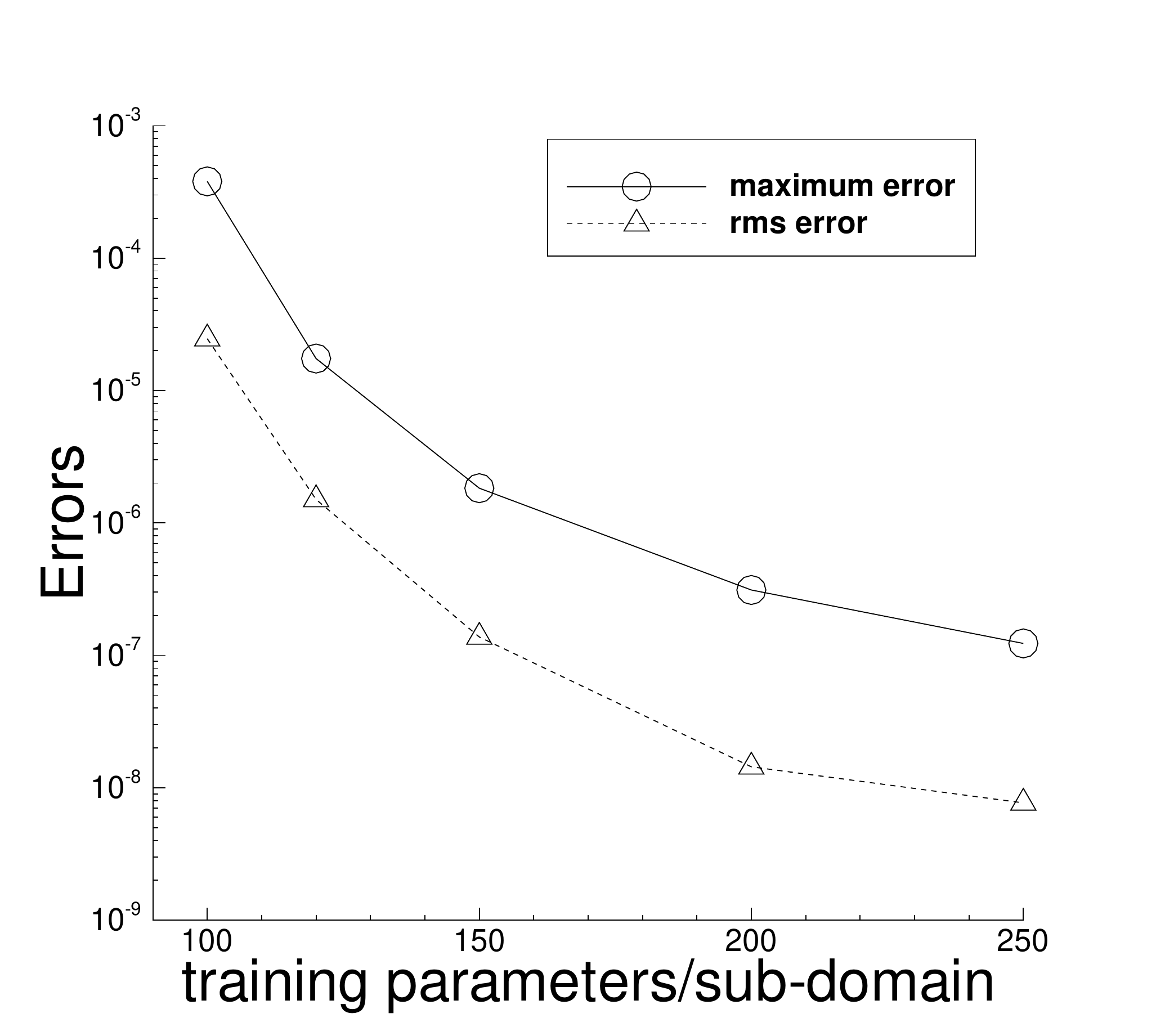}(b)
  }
  \caption{Effect of the degrees of freedom on the accuracy (Burger's equation):
    the maximum and rms errors in the domain as a
    function of (a) the number of collocation points in each direction per
    sub-domain, and
    (b) the number of training parameters per sub-domain.
  }
  \label{fg_bg_3}
\end{figure}

Figure \ref{fg_bg_3} demonstrates the effect of the degrees of freedom
on the simulation accuracy.
In this group of tests the temporal domain size is set to
$t_f=2.5$. We have employed $N_b=10$ uniform time blocks in the overall
spatial-temporal domain, $N_e=5$ uniform sub-domains per time block
(with $N_x=5$ and $N_t=1$), and $R_m=0.5$ when generating random coefficients
for the hidden layers of the local neural networks.
First, we fix the number of training parameters per sub-domain
to $M=200$, and vary the number of (uniform) collocation points
per sub-domain systematically while maintaining $Q_x=Q_t$.
Figure \ref{fg_bg_3}(a) shows the maximum and rms errors in the overall domain
versus the number of collocation points in each direction
per sub-domain.
It is observed that the errors decrease essentially exponentially
with increasing number of collocation points per direction (when
below around $Q_x=Q_t=15$). Then the errors stagnate as the number of
collocation points per direction increases beyond $15$, due to
the saturation associated with the fixed number of training parameters
in the test.
Then, we fix the number of uniform collocation points
to $Q=20\times 20$ per sub-domain, and vary the number of training parameters
per sub-domain systematically in a range of values.
Figure \ref{fg_bg_3}(b) shows the resultant maximum/rms errors in the overall domain
versus the number of training parameters per sub-domain.
As the number of training parameters per sub-domain increases,
the locELM errors can be observed to
decrease substantially.

\begin{figure}
  \centerline{
    \includegraphics[width=2.0in]{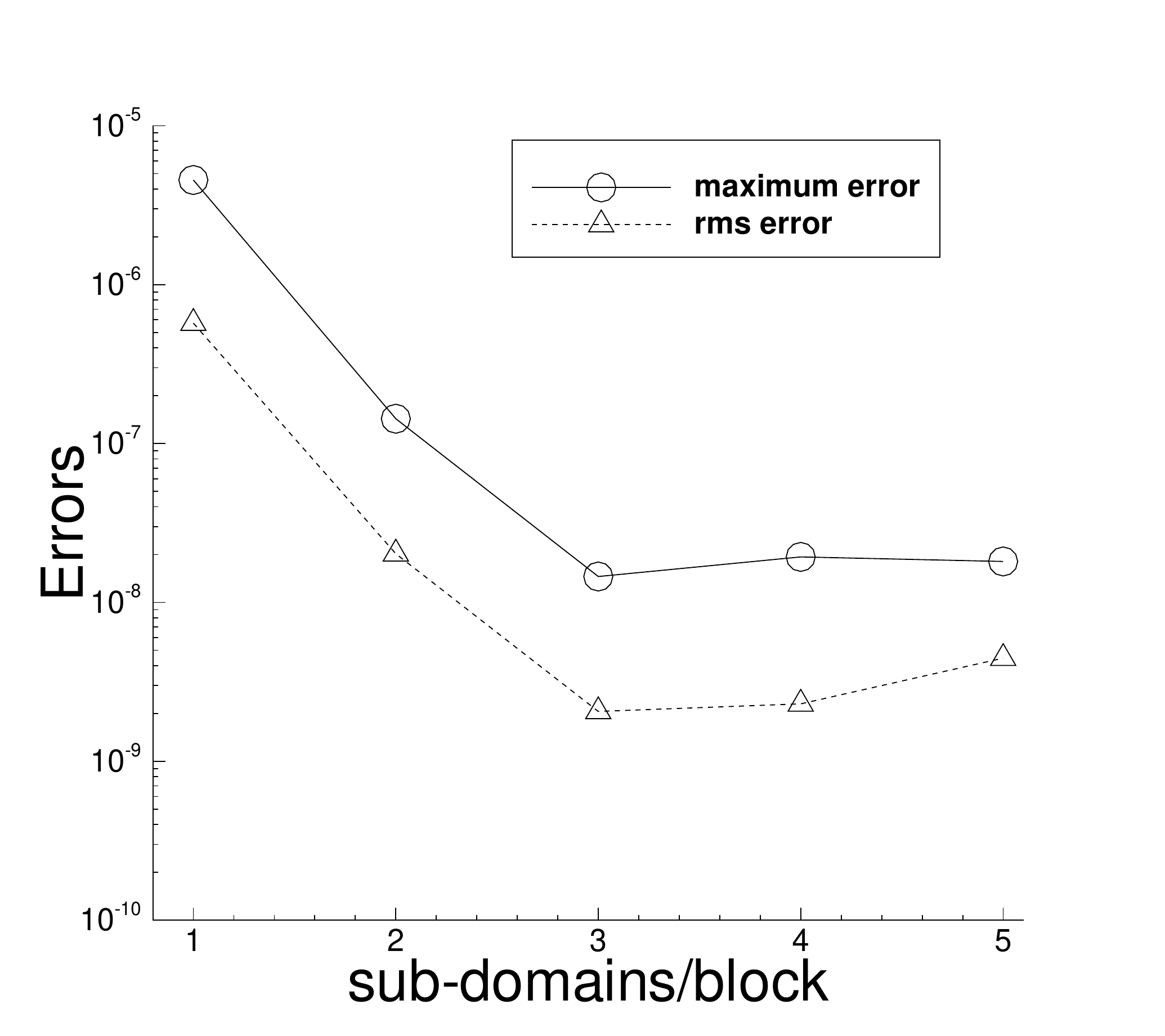}(a)
    \includegraphics[width=2.0in]{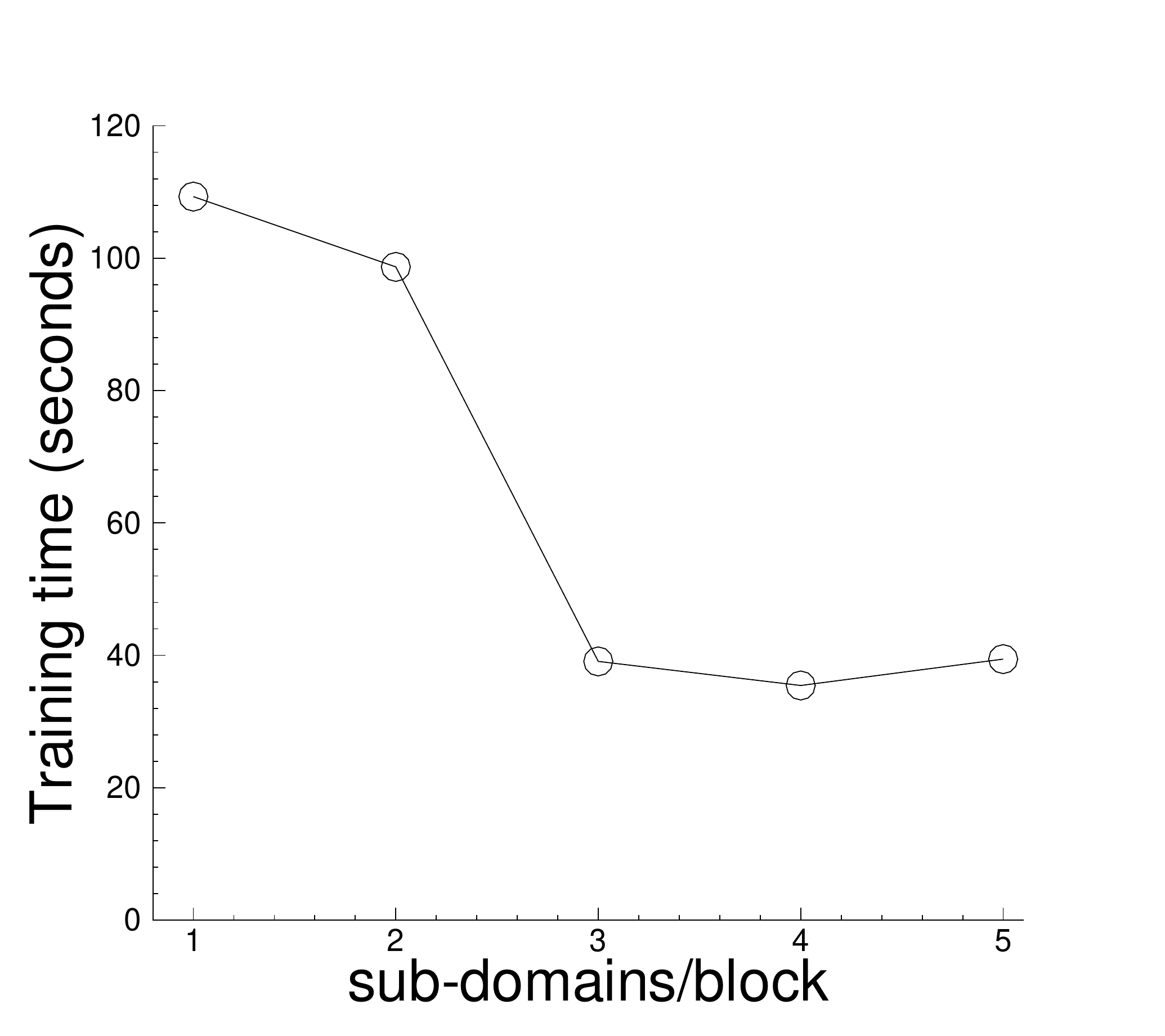}(b)
  }
  \caption{Effect of the number of sub-domains,
    with fixed total degrees of freedom in the domain
    (Burger's equation): (a) the maximum and rms errors in the domain, and (b) the training
    time, as a function of the number of uniform sub-domains per time block.
    Temporal domain size is $t_f=0.25$, and a single time block is used in
    the spatial-temporal domain.
  }
  \label{fg_bg_4}
\end{figure}

Figure \ref{fg_bg_4} demonstrates the effect of the number of sub-domains, with the
total number of degrees of freedom in the domain (approximately) fixed.
In this group of tests, the temporal domain size is set to $t_f=0.25$,
and we employ a single time block in the spatial-temporal domain.
We employ uniform sub-domains, and vary the number of
sub-domains within the time block systematically between $N_e=1$ and $N_e=5$
(with fixed $N_t=1$ and various $N_x$).
The number of (uniform) collocation points per sub-domain and
the number of training parameters per sub-domain are both varied, but
the total number of collocation points and the total number of
training parameters in the time block are fixed
approximately at $N_eQ\approx 2000$ and $N_eM\approx 1000$, respectively.
More specifically, we employ $Q=45\times 45$ collocation points/sub-domain
and $M=1000$ training parameters/sub-domain with $N_e=1$ sub-domain within
the time block, $Q=32\times 32$ collocation points/sub-domain and
$M=500$ training parameters/sub-domain with $N_e=2$ sub-domains,
$Q=26\times 26$ collocation points/sub-domain and $M=333$ training parameters/sub-domain
with $N_e=3$ sub-domains, $Q=22\times 22$ collocation points/sub-domain
and $M=250$ training parameters/sub-domain with $N_e=4$ sub-domains,
and $Q=20\times 20$ collocation points/sub-domain and $M=200$ training parameters/sub-domain
with $N_e=5$ sub-domains within the time block.
When generating the random weight/bias coefficients,
we have employed
$R_m=2.0$ with $N_e=1$ sub-domain in the time block, $R_m=1.0$
with $N_e=2$ and $3$ sub-domains, and $R_m=0.75$ with
$N_e=4$ and $5$ sub-domains within the time block.
These values are approximately in the optimal range of $R_m$ values
for these cases.
Figure \ref{fg_bg_4}(a) shows the maximum and rms errors of the locELM (NLSQ-perturb)
solution in the domain as a function of the number of sub-domains within
the time block.
We observe that the errors decrease quite significantly, by nearly
two orders of magnitude, as the number of sub-domains increases from $N_e=1$
to $N_e=3$. The errors remain approximately at the same level
with three and more sub-domains.
Note that the case with one sub-domain 
corresponds to the global ELM computation.
These results indicate that the local ELM simulation with multiple sub-domains
appears to achieve a better accuracy than the global ELM simulation
for this problem.
Figure \ref{fg_bg_4}(b) shows the training time of the neural network
as a function of the number of sub-domains.
The training time has been reduced substantially as
the number of sub-domains increases from one to three sub-domains
(from around $110$ seconds to about $40$ seconds),
and it remains approximately the same with three and more sub-domains.
These results show that,
compared with the global ELM, the use of
domain decomposition and multiple sub-domains in locELM
can significantly reduce the computational cost 
for the Burger's equation.
This is consistent with the observations with
the other problems in previous sections.

\begin{figure}
  \centerline{
    \includegraphics[width=2.5in]{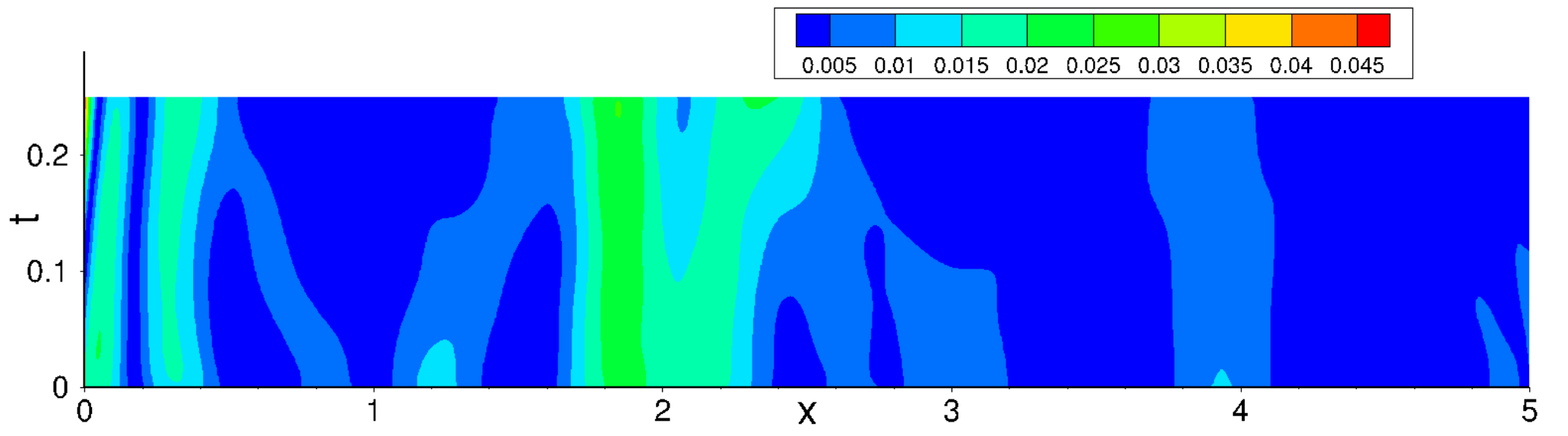}(a)
    \includegraphics[width=2.5in]{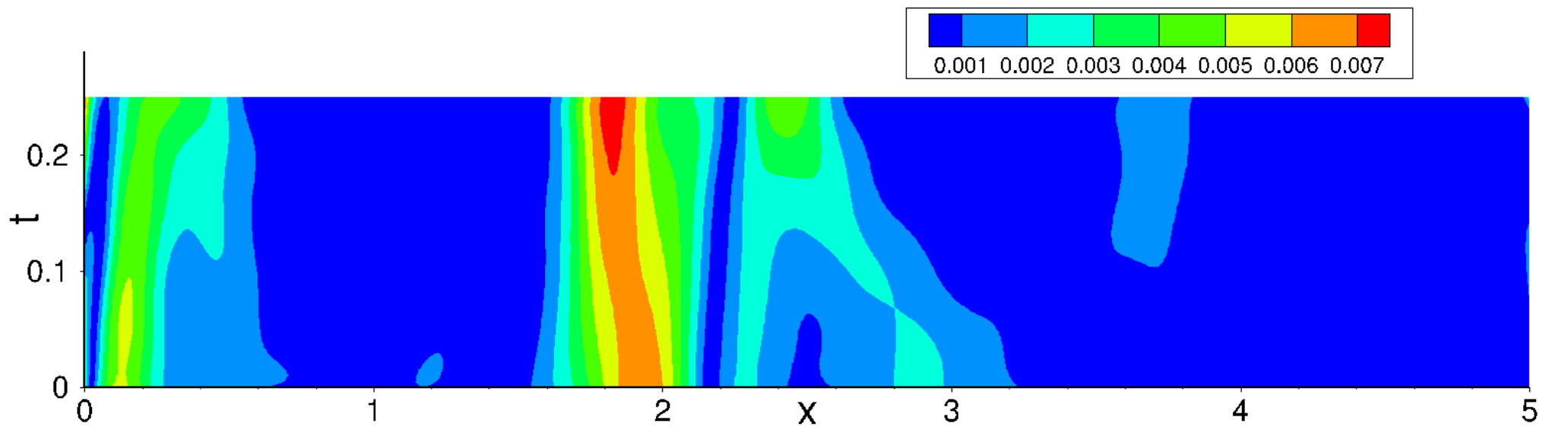}(b)
  }
  \centerline{
    \includegraphics[width=2.5in]{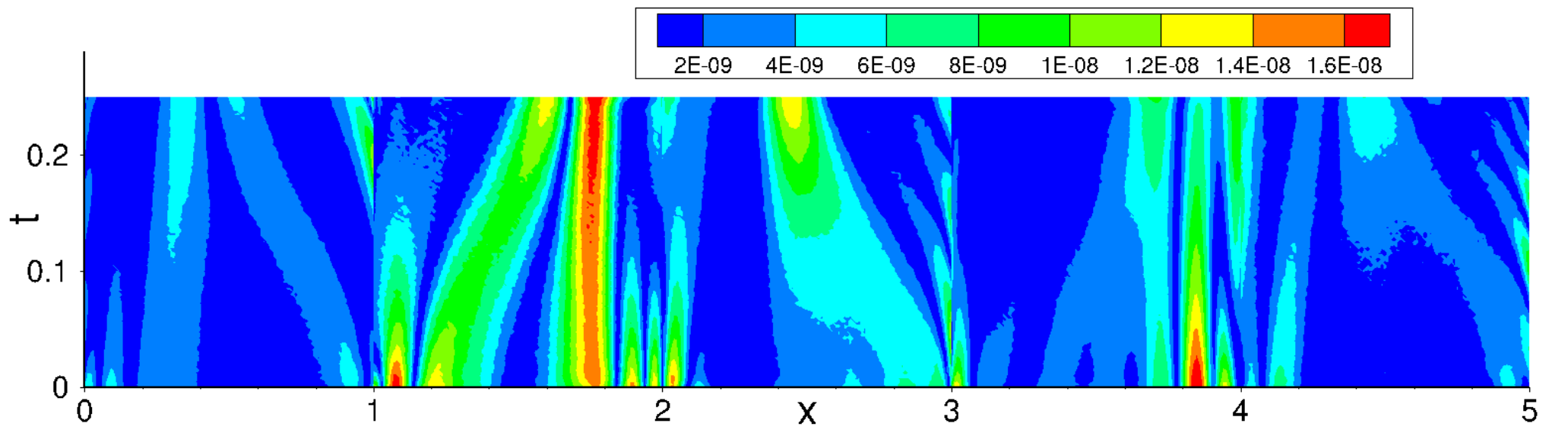}(c)
    \includegraphics[width=2.5in]{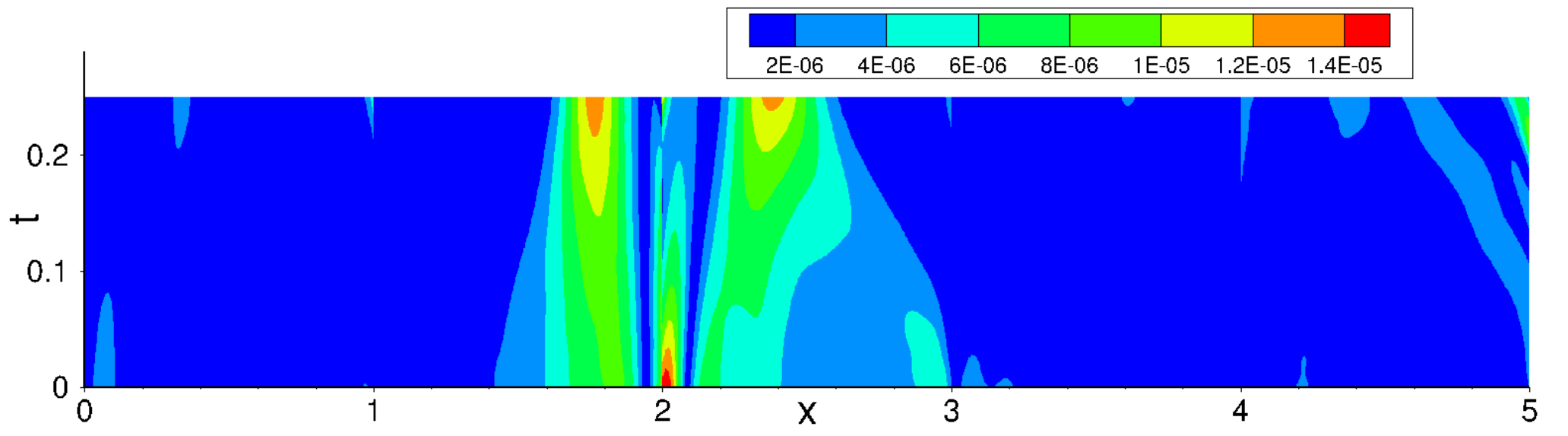}(d)
  }
  \caption{Comparison between locELM and DGM (Burger's equation):
    distributions of the
    absolute errors  computed using DGM with the Adam optimizer
    (a) and L-BFGS optimizer (b), and using locELM
    with NLSQ-perturb (c) and with Newton-LLSQ (d).
  }
  \label{fg_bg_5}
\end{figure}

\begin{figure}
  \centerline{
    \includegraphics[width=2.2in]{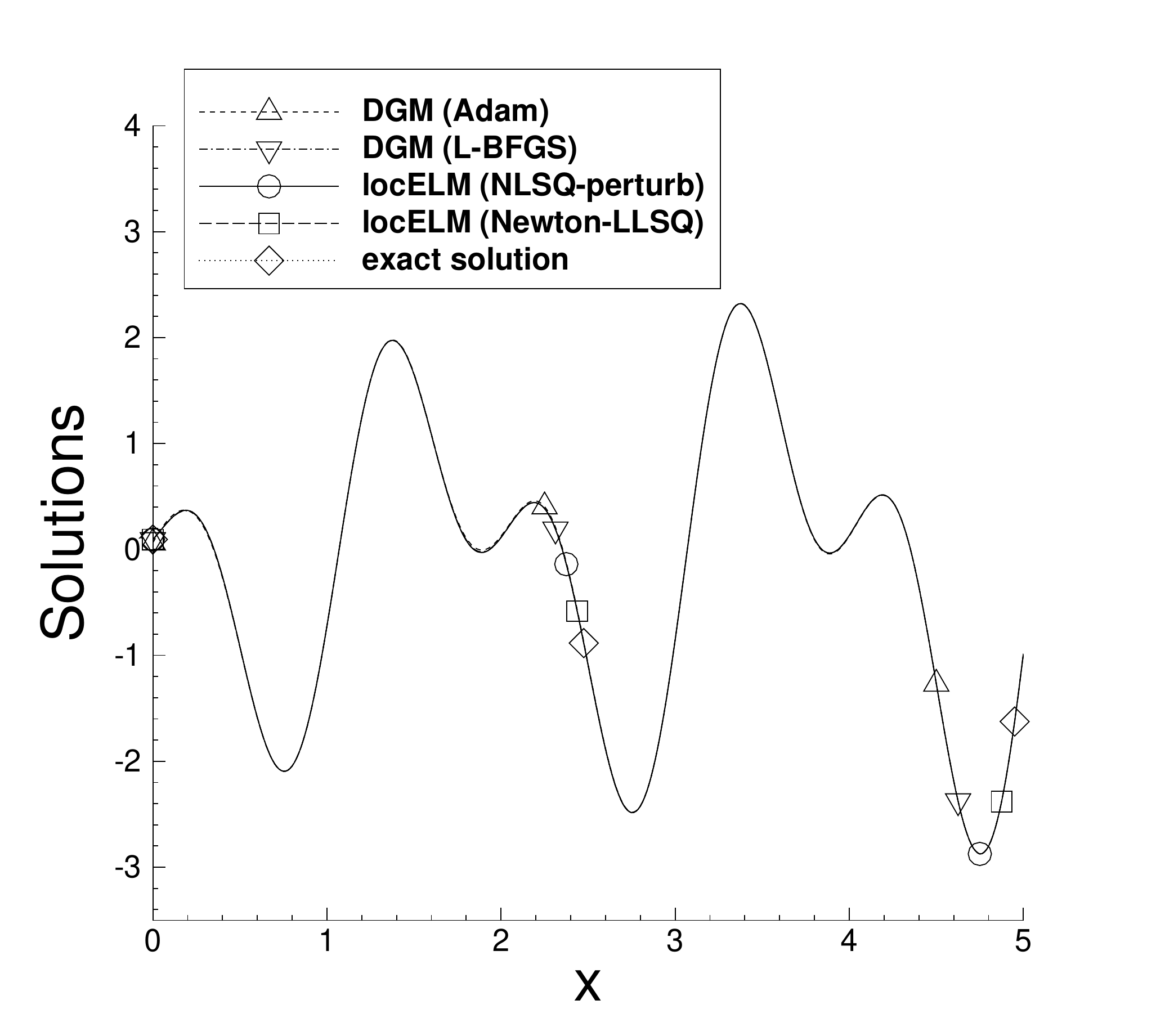}(a)
    \includegraphics[width=2.2in]{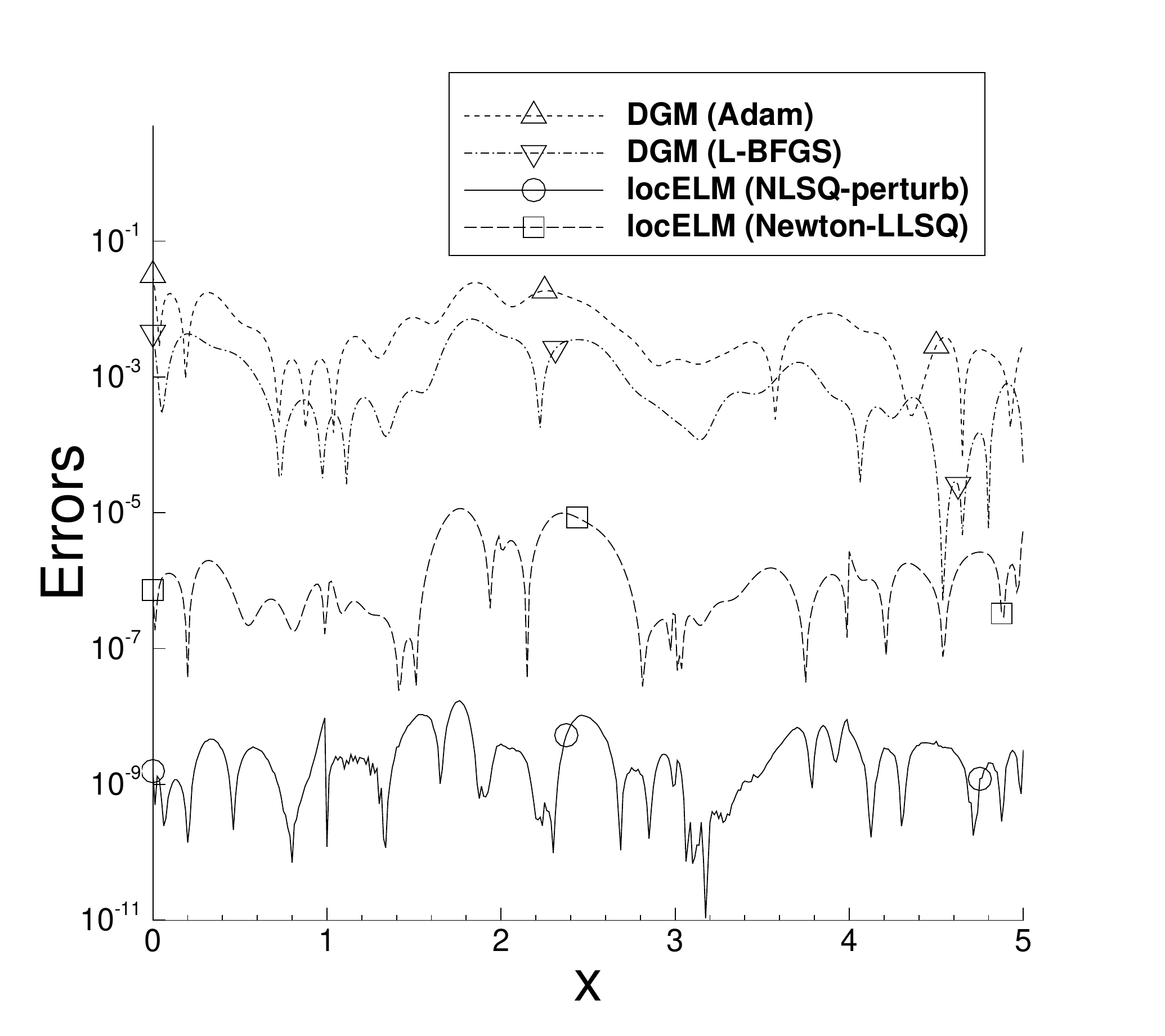}(b)
  }
  \caption{Burger's equation: Profiles of (a) the solutions and (b) their absolute errors
    at $t=0.2$ computed using DGM and locELM.
    The settings and simulation parameters correspond to those of Figure \ref{fg_bg_5}.
  }
  \label{fg_bg_6}
\end{figure}

\begin{table}
  \centering
  \begin{tabular}{lllll}
    \hline
    method & maximum error & rms error & epochs/iterations & training time (seconds)\\
    DGM (Adam) & $4.57e-2$ & $5.76e-3$ & $128,000$ & $1797.8$ \\
    DGM (L-BFGS) & $7.50e-3$ & $1.55e-3$ & $28,000$ & $1813.5$ \\
    locELM (NLSQ-perturb) & $1.85e-8$ & $4.44e-9$ & $27$ & $27.6$ \\
    locELM (Newton-LLSQ) & $1.62e-5$ & $3.11e-6$ & $15$ & $9.1$ \\
    \hline
  \end{tabular}
  \caption{Burger's equation: comparison between locELM and DGM 
    in terms of the maximum/rms errors in the domain, the number of
    epochs or nonlinear iterations, and the network training time.
    The problem settings and simulation parameters correspond to those of
    Figure \ref{fg_bg_5}.
  }
  \label{tab_bg_7}
\end{table}

We next compare the current locELM method with the deep Galerkin method (DGM)
for the Burger's equation.
Figure \ref{fg_bg_5} is a comparison of distributions of the solutions (left column)
and their absolute errors (right column) in the spatial-temporal plane,
obtained using DGM with the Adam and L-BFGS optimizers (top two rows)
and using the current locELM method with NLSQ-perturb and Newton-LLSQ (bottom
two rows). In these tests the temporal domain size is set to $t_f=0.25$.
For DGM, the neural network consists of an input layer of two nodes (representing
$x$ and $t$), $5$ hidden layers with a width of $40$ nodes in each layer
and the $\tanh$ activation function, and an output layer of a single node
(representing $u$) with no bias and no activation function.
When computing the loss function, the spatial-temporal domain has been divided into
$10$ uniform sub-domains along the $x$ direction, and we have used $10\times 10$
Gauss-Lobatto-Legendre quadrature points in each sub-domain
for computing the residual norms.
With the Adam optimizer, the neural network has been trained for $128,000$
epochs, with the learning rate gradually decreasing from $0.001$ at the beginning
to $10^{-5}$ at the end of the training.
With the L-BFGS optimizer, the neural network has been trained for $28,000$ iterations.
For the current locELM method, we have employed a single time block
in the spatial-temporal domain and $N_e=5$ uniform sub-domains along the $x$ direction
within this time block.
With NLSQ-perturb, we have employed $Q=20\times 20$ uniform collocation points
per sub-domain, $M=200$ training parameters per sub-domain, and $R_m=0.75$
when generating the random coefficients.
With Newton-LLSQ, we have employed $Q=20\times 20$ uniform collocation points
per sub-domain, $M=150$ training parameters per sub-domain,
and $R_m=1.0$ when generating the random coefficients.
The results in Figure \ref{fg_bg_5} indicate that the current locELM method is
considerably more accurate than DGM for the Burger's equation.
The errors of the current method is generally several orders of magnitude smaller
than those of DGM. The locELM method with NLSQ-perturb provides the best
accuracy, with the errors on the order $10^{-9}\sim 10^{-8}$.
Then it is the locELM method with Newton-LLSQ, with the errors on the level
$10^{-6}\sim 10^{-5}$. In contrast, the errors of the DGM with Adam and L-BFGS
are generally on the levels $10^{-3}\sim 10^{-2}$ and $10^{-3}$, respectively.

Figure \ref{fg_bg_6} compares the profiles of the DGM and locELM solutions
(plot (a))
and their errors (plot (b)) at the time instant $t=0.2$.
The profile of the exact solution at this instant is also included
in Figure \ref{fg_bg_6}(a) for comparison.
The problem settings and the simulation parameters here correspond to
those of Figure \ref{fg_bg_5}.
The solution profiles from DGM and locELM simulations are in good agreement
with that of the exact solution.
The error profiles, on the other hand,
reveal disperate accuracies in the results obtained using
these methods. They confirm the ordering of these methods,
from the most to the least accurate, to be locELM/NLSQ-perturb,
locELM/Newton-LLSQ, DGM/L-BFGS, and DGM/Adam.

Table \ref{tab_bg_7} provides a further comparison between locELM and DGM
for the Burger's equation,
in terms of their accuracy and computational cost.
We have listed the maximum and rms errors in the overall spatial-temporal domain,
the number of epochs or nonlinear iterations in the training,
and the training time of the neural network corresponding to DGM with
the Adam/L-BFGS optimizers and the current locELM method with NLSQ-perturb
and Newton-LLSQ.
The observations here are consistent with those of previous sections.
The locELM method is orders of magnitude more accurate than DGM
(e.g.~$10^{-8}$ with locELM/NLSQ-perturb versus $10^{-3}$ with DGM/L-BFGS),
and its training time is orders of magnitude smaller than that of
DGM (e.g.~around $28$ seconds with locELM/NLSQ-perturb versus
around $1800$ seconds with DGM/L-BFGS).


\begin{figure}
  \centerline{
    \includegraphics[width=2in]{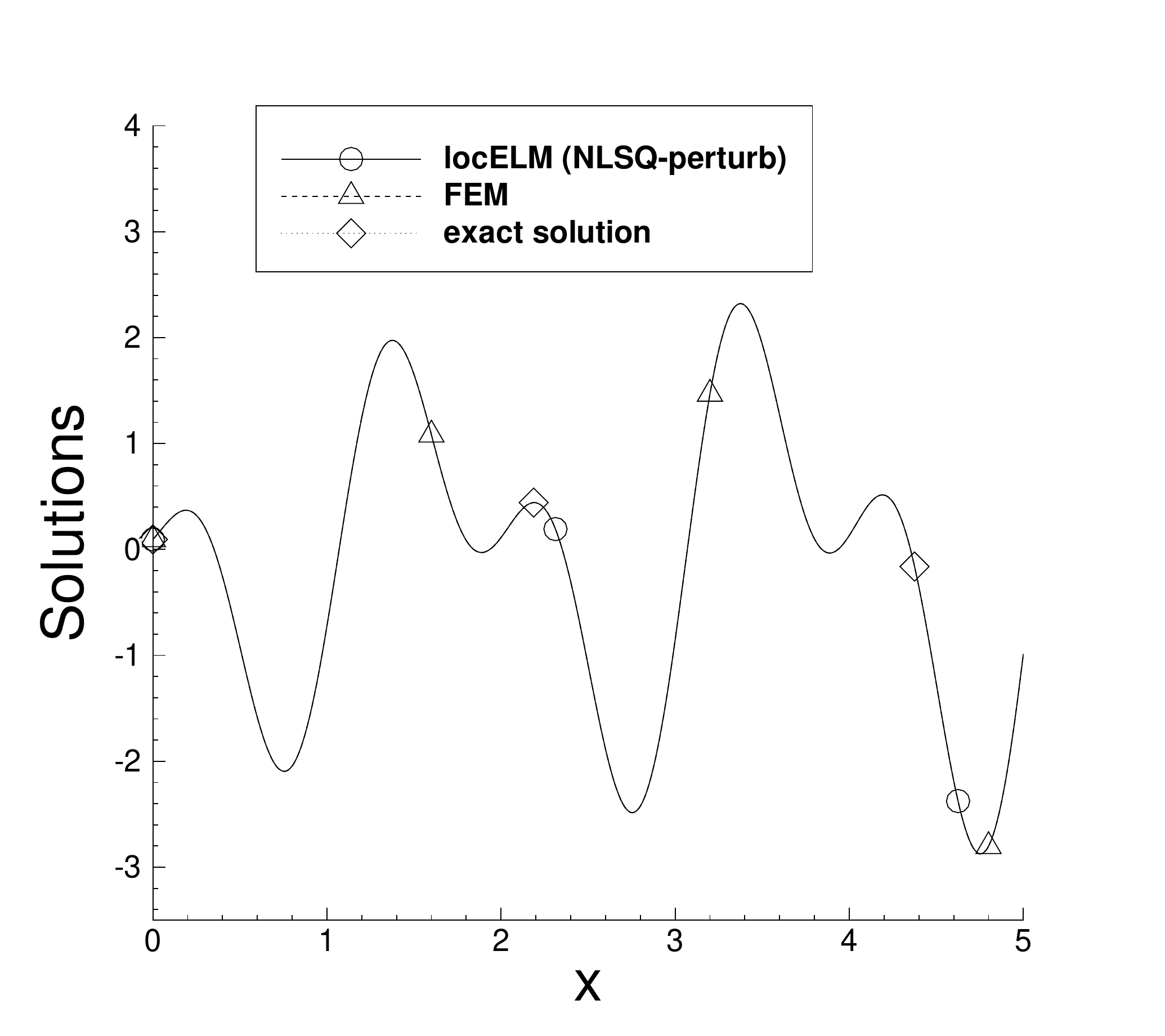}(a)
    \includegraphics[width=2in]{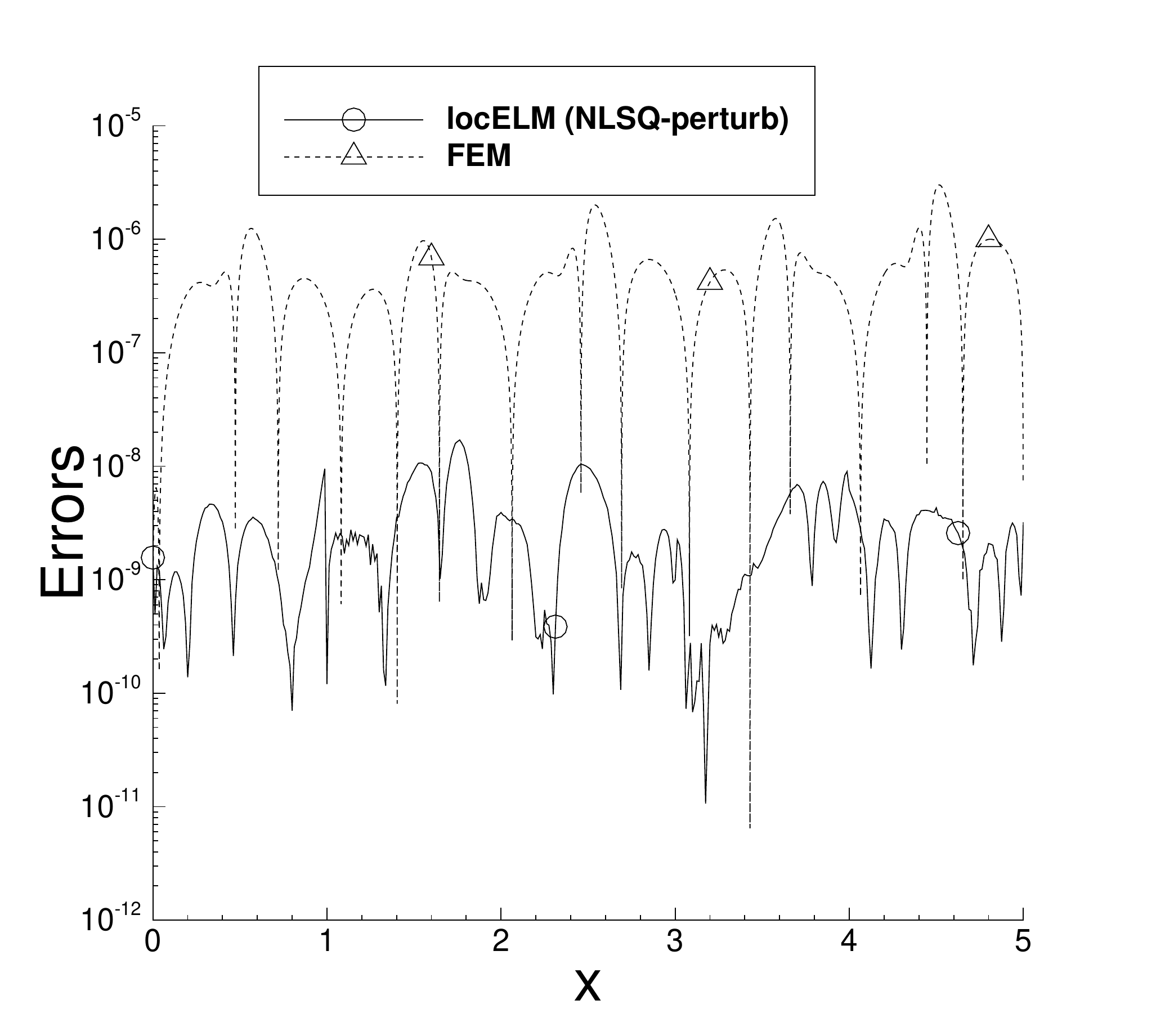}(b)
    \includegraphics[width=2in]{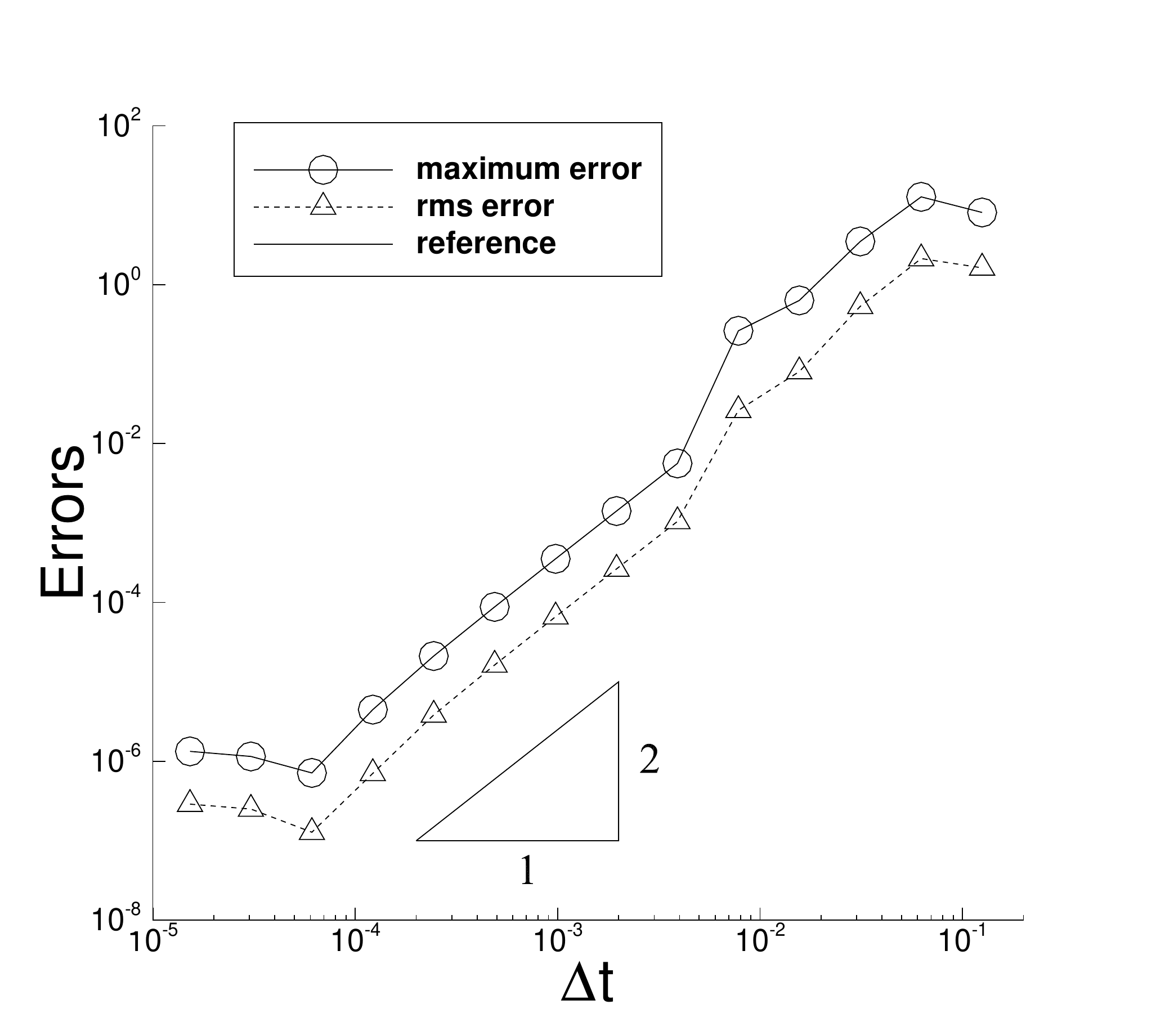}(c)
  }
  \caption{Comparison between locELM and FEM (Burger's equation): Profiles
    of (a) the solutions and (b) their absolute errors at $t=0.2$,
    computed using locELM (with NLSQ-perturb) and using FEM.
    (c) The maximum and rms errors at $t=0.25$ versus
    $\Delta t$ computed using FEM (with a mesh of $10,000$ uniform elements),
    showing its second-order convergence rate in time.
  }
  \label{fg_bg_7}
\end{figure}

\begin{table}
  \centering
  \begin{tabular}{l|llllllll}
    \hline
    method & elements  & $\Delta t$ & sub- & $Q$ & $M$ & maximum & rms & wall-time \\
      & & & domains & & & error & error & (seconds) \\
    \hline
    locELM 
     & -- & -- & $5$ & $15\times 15$ & $150$ & $2.10e-6$ & $4.35e-7$ & $14.7$ \\
    (NLSQ-perturb)  & -- & -- & $5$ & $20\times 20$ & $200$ & $1.85e-8$ & $4.44e-9$  & $27.6$ \\
    \hline
    locELM 
     & -- & -- &  $5$ & $15\times 15$ & $150$ & $1.25e-5$ & $2.71e-6$ & $6.8$ \\
    (Newton-LLSQ) &  -- & -- & $5$ & $20\times 20$ & $150$  & $1.62e-5$ & $3.11e-6$  & $9.1$ \\
    \hline
    FEM
    & 2000 & $0.001$ & -- & -- & -- & $2.64e-5$ & $5.15e-6$ & $12.5$ \\
    & 2000 & $0.0005$ & -- & -- & -- & $3.07e-5$ & $5.76e-6$ & $25.4$ \\ \cline{2-9}
    & 5000 & $0.001$ & -- & -- & -- & $1.89e-5$ & $1.74e-6$ & $26.0$ \\
    & 5000 & $0.0005$ & -- & -- & -- & $4.13e-6$ & $7.90e-7$ & $50.8$ \\ \cline{2-9}
    & 10000 & $0.001$ & -- & -- & -- & $2.22e-5$ & $1.99e-6$ & $47.7$ \\
    & 10000 & $0.0005$ & -- & -- & -- & $4.74e-6$ & $4.36e-7$ & $92.6$ \\
    \hline
  \end{tabular}
  \caption{Burger's equation: comparison between locELM and FEM
    in terms of the maximum/rms errors in the domain and the training/computation time.
    $Q$ and $M$ denote the number of collocation points per sub-domain and
    the number of training parameters per sub-domain, respectively.
  }
  \label{tab_bg_8}
\end{table}

Finally, we compare the current locELM method with the classical finite element
method for solving the Burger's equation.
In the FEM simulation, we discretize the Burger's equation~\eqref{eq_bg_1}
in time using a semi-implicit scheme.
We treat the nonlinear term explicitly and the viscous term implicitly,
and discretize the time derivative by the second-order backward differentiation
formula (BDF2). The method is again implemented using the FEniCS library
in Python.
Figures \ref{fg_bg_7}(a) and (b) show a comparison of the solution and error profiles
at $t=0.2$ obtained using the current locELM (NLSQ-perturb)
method and using the finite element method.
Figure \ref{fg_bg_7}(c) shows the numerical errors
at $t=0.25$ as a function of the time step
size $\Delta t$ computed using the finite element method.
In these simulations the temporal domain size is $t_f=0.25$.
In Figures \ref{fg_bg_7}(a,b),
the FEM simulation is conducted with $\Delta t=1.25e-4$ on a mesh of
$10,000$ uniform elements, and the locELM simulation is conducted with
a single time block in the domain and $N_e=5$ sub-domains in the time block, with
$(Q,M)=(20\times 20, 200)$ and $R_m=0.75$.
In Figure \ref{fg_bg_7}(c), the simulations are performed
with a fixed mesh of $10,000$ uniform elements.
It can be observed that both locELM and FEM have produced accurate solutions,
and that the FEM exhibits a second-order convergence rate in time
before the error saturation when $\Delta t$ becomes very small.

Table \ref{tab_bg_8} provides a comparison between locELM and FEM
in terms of their accuracy and computational cost for the Burger's equation.
The temporal domain size is $t_f=0.25$ in these tests.
We solve the problem using locELM and FEM
on several sets of simulation parameters with different numerical resolutions.
The maximum and rms errors in the spatial-temporal domain are computed,
and we also record the training time of locELM and the computation time
of FEM in these simulations.
We list in this table the maximum and rms errors, as well
as the training/computation time, corresponding to
different simulation parameters for the locELM method with NLSQ-perturb
and Newton-LLSQ and for the finite element method.
A single time block has been used in the spatial-temporal domain 
for the locELM simulations, and we employ $R_m=0.75$ with
locELM/NLSQ-perturb and $R_m=1.0$ with locELM/Newton-LLSQ for generating
the random coefficients.
%
It is observed that the current locELM method with both NLSQ-perturb and
Newton-LLSQ shows a superior performance to the FEM.
For example,
the two cases with locELM/Newton-LLSQ have numerical errors comparable to
the FEM cases with $2000$ elements (for both $\Delta t$), $5000$ elements ($\Delta t=0.001$)
and $10000$ elements ($\Delta t=0.001$), but the computational cost
of locELM/Newton-LLSQ is notably smaller than the cost of these FEM cases.
The locELM/NLSQ-perturb case with $(Q,M)=(15\times 15, 150)$ has numerical
errors comparable to the FEM cases with $5000$ elements ($\Delta t=0.0005$)
and $10000$ elements ($\Delta t=0.0005$), but the computational cost of 
this locELM/NLSQ-perturb case is only a fraction of those of these two FEM cases.
The locELM/NLSQ-perturb case with $(Q,M)=(20\times 20, 200)$ has a computational
cost comparable to the FEM cases with $2000$ elements ($\Delta t=0.0005$)
and $5000$ elements ($\Delta t=0.001$), but the errors of this
locELM/NLSQ-perturb case are nearly three orders of magnitude smaller than
those of these two FEM cases.


\section{Concluding Remarks}
\label{sec:summary}

%
%

In this paper we have developed an efficient method
based on domain decomposition and local extreme learning machines (termed locELM) 
for solving linear and nonlinear  partial differential equations.
The problem domain is partitioned into sub-domains,
and the field solution on each sub-domain is represented
by a local shallow feed-forward neural network,
consisting of a small number (one or more) of
hidden layers.
$C^k$ continuity, with $k$ determined by the order of the PDE,
is imposed on the sub-domain boundaries.
The overall neural network constitutes a multi-input multi-output model
consisting of the local neural networks.
The weight/bias coefficients in the hidden layers
of all the local neural networks are pre-set to random values,
and are fixed throughout the computation.
The  training parameters are composed of
the weight coefficients in the output layers of
the local neural networks.


We employ a set of collocation points within each sub-domain, the collection of
which constitutes
the input data to the neural network.
%
The PDE is enforced on the collocation points in
each sub-domain, and the  derivatives involved therein are computed
by auto-differentiation. The boundary and initial conditions are enforced on
those collocation points that reside on the spatial
and temporal boundaries of the spatial-temporal domain.
The $C^k$ continuity conditions are enforced on those collocation
points that reside on the corresponding sub-domain boundaries.
%
These operations result in a system of linear or nonlinear algebraic equations about
the set of training parameters.
%
We seek a least squares solution to this system, and 
 compute the solution by a linear least squares routine or a nonlinear
least squares method.
Training the overall neural network
consists of the linear or nonlinear least squares computations.
It should be noted that this training method is different from
the back propagation-type algorithms.

%

For  longer-time simulations of time-dependent PDEs,
we have developed a block time-marching scheme together with the locELM method.
The spatial-temporal domain is first divided into a number of windows in time,
referred to as time blocks, and we solve the PDE on each time block
separately and successively. The locELM method is then applied to the spatial-temporal
domain of each time block to find the solution,
as discussed in foregoing paragraphs.
We observe that when the temporal dimension of the domain is large, if without block
time marching, the
neural network can become very difficult to train.
On the other hand, with block time marching and using a moderate time block size,
the problem is more manageable and much easier to solve. 
Block time marching requires re-training of the overall neural network on
different time blocks, and so all network trainings become online operations.
This is feasible with the current locELM method thanks to its high accuracy
and low computational cost. 
We have demonstrated the capability of the current method for long-time dynamic
simulations
with the advection equation.


We have performed extensive numerical experiments to test the locELM method,
and compared extensively the locELM method with DGM, PINN, global ELM, and
the classical finite element method (FEM).
We have the following observations:
\begin{itemize}

\item
  The locELM method exhibits a clear sense of
  convergence with increasing number of degrees of freedom.
  Its errors typically decrease exponentially or nearly exponentially
  as the number of sub-domains, or the number of collocation points/sub-domain,
  or the number of training parameters/sub-domain increases.

\item
  The random weight/bias coefficients in the hidden layers of local neural networks
  influence the simulation accuracy.
  In the current work, these weight/bias coefficients are set to uniform random
  values generated on $[-R_m,R_m]$.
  The simulation accuracy tends to decrease
  with very large or very small $R_m$ values. Higher accuracy is generally associated
  with a range of moderate $R_m$ values.
  This range of optimal $R_m$ values tends to expand when the number of
  collocation points/sub-domain or the number of training parameters/sub-domain
  increases.

\item
  The network training time  generally increases
  linearly (or super-linearly for some problems) with respect to the number of
  sub-domains in the simulation. It also tends to increase with respect to
  the number of collocation points/sub-domain and to
  the number of training parameters/sub-domain, but the relation is
  not quite regular.

\item
  When the total degrees of freedom (total collocation points,
  total training parameters) in the system are fixed,
  increasing the number of sub-domains in the simulation, hence with
  the number of collocation points/training parameters per sub-domain correspondingly reduced,
  generally leads to simulation results with comparable accuracy,
  but it can dramatically reduce the network training time.
  Compared with global ELM, which corresponds to the
  locELM configuration with a single sub-domain, the use of
  domain decomposition and multiple sub-domains in locELM
  can significantly reduce the network training time, and
  produce results with comparable accuracy.

\item
  The current locELM method shows a clear superiority to DGM and PINN, which are some of
  the commonly-used DNN-based PDE solvers,
  in terms of both accuracy and computational cost.
  The numerical errors and the network training time of locELM are considerably smaller,
  typically by orders of magnitude, than those of DGM and PINN.

\item
  The current locELM method exhibits a computational performance that is comparable,
  and oftentimes superior, to the classical finite element method.
  With the same
  computational cost, the locELM errors
  are comparable to, and oftentimes considerably smaller than, the FEM errors.
  To achieve the same accuracy,
  the training time of locELM 
  is comparable to, and oftentimes markedly smaller than, the FEM computation time.

\end{itemize}


We would like to make some further comments with regard to $R_m$,
the maximum magnitude of the random weight/bias coefficients in the hidden layers of
the local neural networks.
As discussed above, the simulation results have a better accuracy
if $R_m$ falls into a range of moderate values for a given problem.
Let us consider the following question:
given a new problem (e.g.~a new PDE), how do we know
what this range is and how do we
find this range of optimal $R_m$ values in practice? 
The approximate range of these optimal $R_m$ values can be estimated readily
by preliminary numerical experiments. Here is the basic idea.
Given a new problem, one can always add some source terms to the PDE or to
the boundary/initial conditions, and then manufacture a solution to
the given problem, with the augmented source terms.
Then one can use the manufactured solution
to evaluate the accuracy of a set of preliminary simulations by varying the
$R_m$ systematically.
This will provide a reasonable estimate for the range of optimal $R_m$ values.
After that, one can conduct actual simulations of the given problem,
without the added source term, by using an $R_m$ value from
the estimated range.


Some further comments are also in order concerning the numerical tests
with fixed total degrees of
freedom in the domain, while the number of sub-domains is varied.
Because the total degrees of freedom in the domain is fixed,
the degrees of freedom
(number of collocation points, number of training parameters) per sub-domain 
decrease as the number of sub-domains increases.
One can anticipate that, when the number of sub-domains becomes sufficiently large, 
the number of degrees of freedom per sub-domain can become very small.
This will be bound to adversely affect the simulation accuracy,
because the solution is represented locally by these degrees of
freedom on each sub-domain.
Therefore, if the total degrees of freedom in the domain  are fixed,
when the number of sub-domains in domain decomposition
increases beyond a certain point, the simulation accuracy will
start to deteriorate.
The afore-mentioned observation about comparable accuracy with increasing number of
sub-domains, with fixed  total degrees of freedom in the domain,
is for the cases where the number of sub-domains is below that point.

As demonstrated by ample examples in this paper,
the computational performance of the current locELM method is on par with,
and oftentimes exceeds, that of
the classical finite element method.
The importance of this point cannot be overstated. To the best of the authors'
knowledge, this seems to be the first time when a neural network-based method delivers
the same performance as, or a better performance than, a traditional
numerical method for the commonly-encountered computational problems in
low dimensions.
The current method demonstrates the great potential, and perhaps points toward
a path forward, for neural network-based
methods to be truly competitive, and excel, in computational science and engineering
simulations.



\section*{Acknowledgement}
This work was partially supported by
NSF (DMS-2012415, DMS-1522537). 

\section{Appendix. Additional Numerical Tests}

\subsection{Two-Dimensional Helmholtz Equation}

We consider the  boundary value problem with
the two-dimensional (2D) Helmholtz equation on
a rectangular domain, $\Omega=[a_1,b_1]\times [a_2,b_2]$, as follows,
\begin{subequations}
  \begin{align}
    &
    \frac{\partial^2u}{\partial x^2} + \frac{\partial^2u}{\partial y^2}
    - \lambda u = f(x,y), \label{equ_hm2_1} \\
    &
    u(a_1,y) = h_1(y), \\
    & u(b_1,y) = h_2(y), \\
    & u(x,a_2) = h_3(x), \\
    & u(x,b_2) = h_4(x),
  \end{align}
\end{subequations}
where $u(x,y)$ is the field function to be solved for, $f(x,y)$ is
a prescribed source term, and $h_i$ ($1\leqslant i\leqslant 4$) denote
the boundary distributions. The constant parameters
are given by
\begin{equation*}
  \lambda = 10, \quad
  a_1 = a_2 = 0, \quad
  b_1 = b_2 = 3.6.
\end{equation*}
We choose the source term $f(x,y)$ such that the following function
satisfies equation \eqref{equ_hm2_1},
\begin{equation}\label{equ_hm2_2}
  u = -\left[\frac32\cos\left(\pi x+\frac{2\pi}{5}\right)
    +2\cos\left(2\pi x-\frac{\pi}{5}\right) \right]
  \left[\frac32\cos\left(\pi y+\frac{2\pi}{5}\right)
    +2\cos\left(2\pi y-\frac{\pi}{5}\right) \right].
\end{equation}
We choose the boundary distributions $h_1(y)$, $h_2(y)$, $h_3(x)$
and $h_4(x)$ in accordance with \eqref{equ_hm2_2}, by setting
\eqref{equ_hm2_2} on the
corresponding boundaries.
Consequently, \eqref{equ_hm2_2} provides the solution to
this boundary value problem.


We employ the locELM method from Section \ref{sec:steady} to solve this problem.
$\Omega$ is partitioned into $N_x$ and $N_y$ uniform sub-domains along
the $x$ and $y$ directions, respectively, leading to a total of
$N_e=N_xN_y$ uniform sub-domains. We impose $C^1$
continuity conditions on the sub-domain boundaries.
In each sub-domain, we employ $Q_x$
and $Q_y$ uniform collocation points  in
the $x$ and $y$ directions, respectively, 
leading to a total of $Q=Q_xQ_y$ collocation points per sub-domain.

\begin{figure}
  \centerline{
    \includegraphics[width=1.5in]{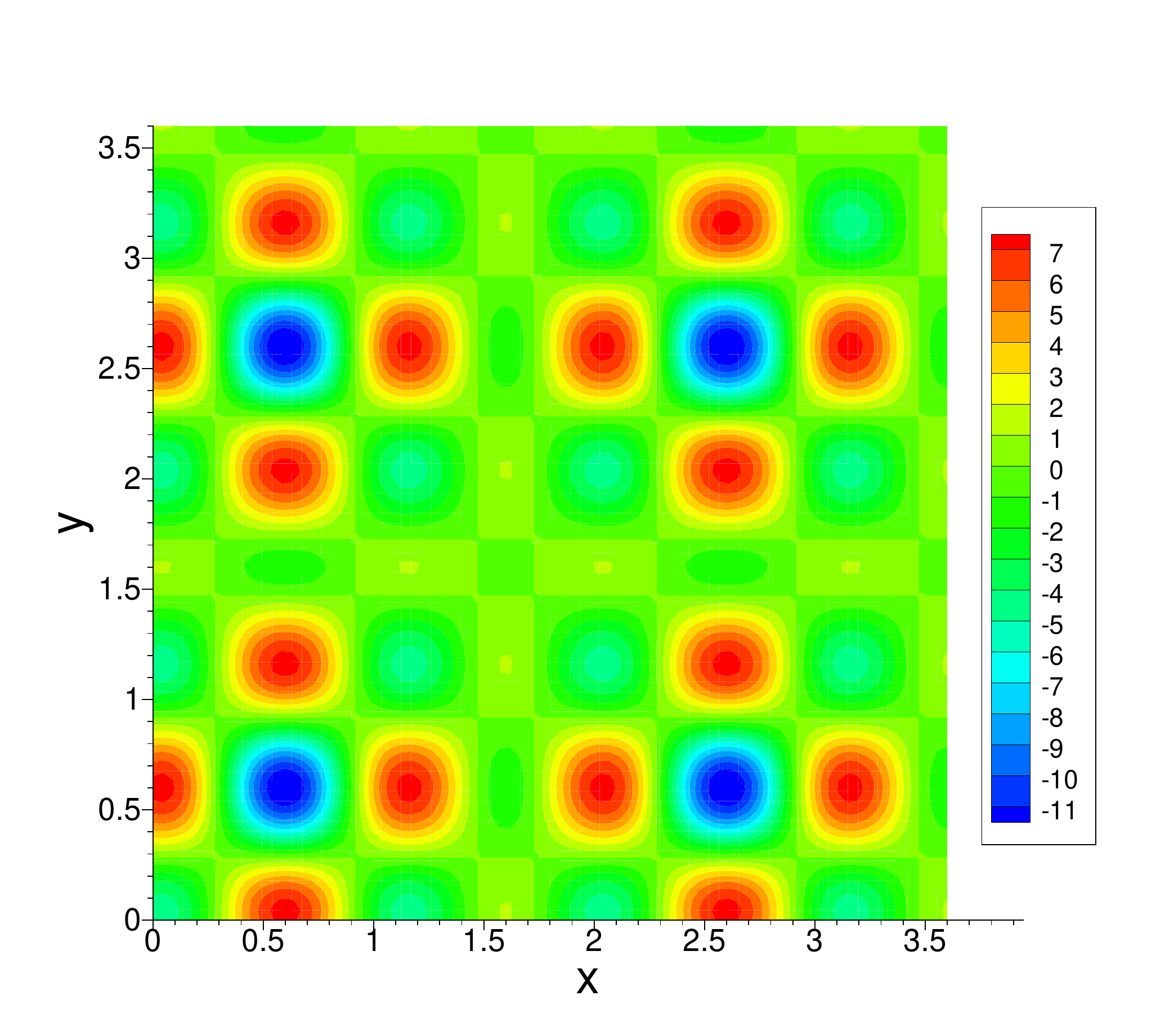}(a)
    \includegraphics[width=1.5in]{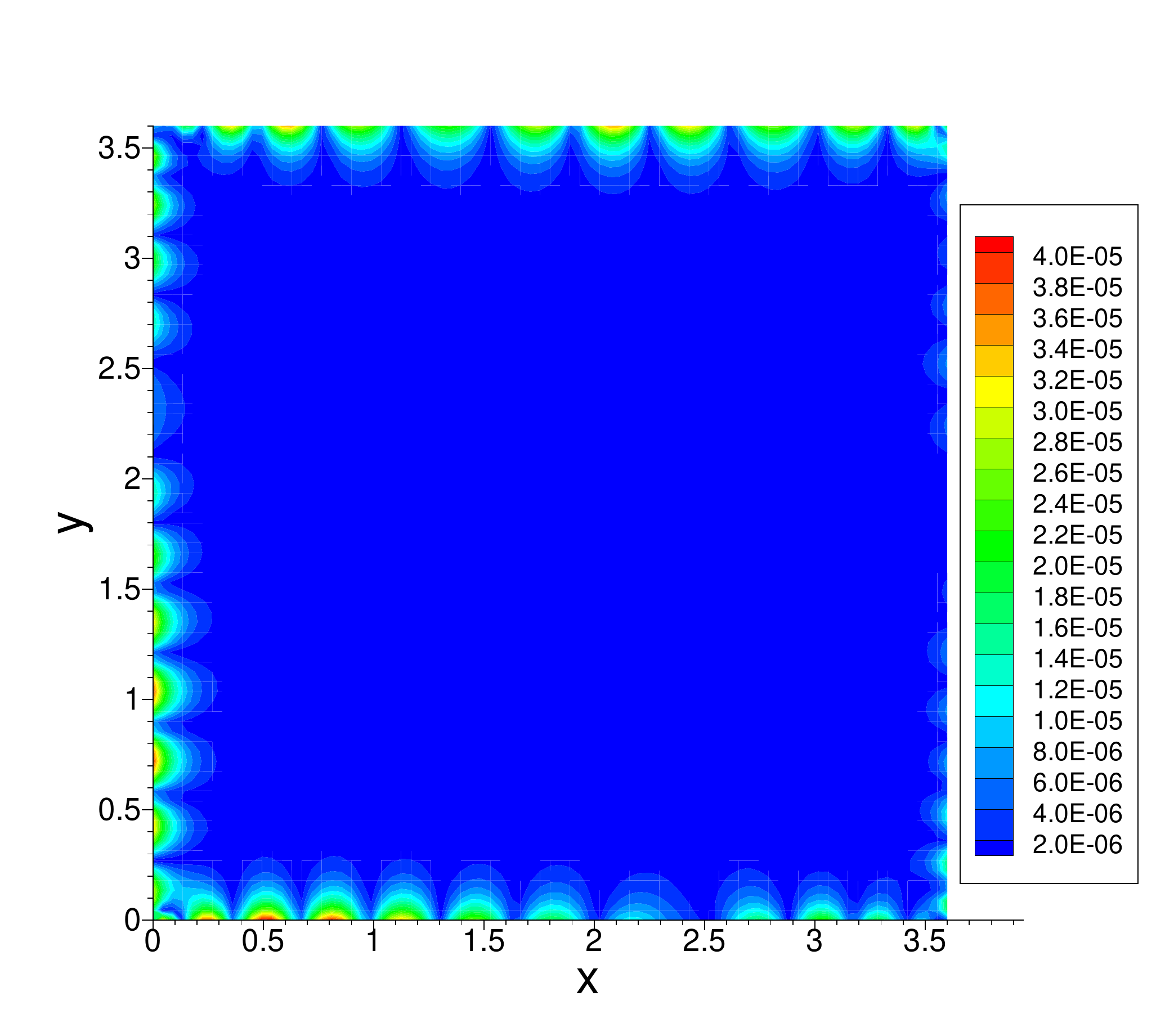}(b)
    \includegraphics[width=1.5in]{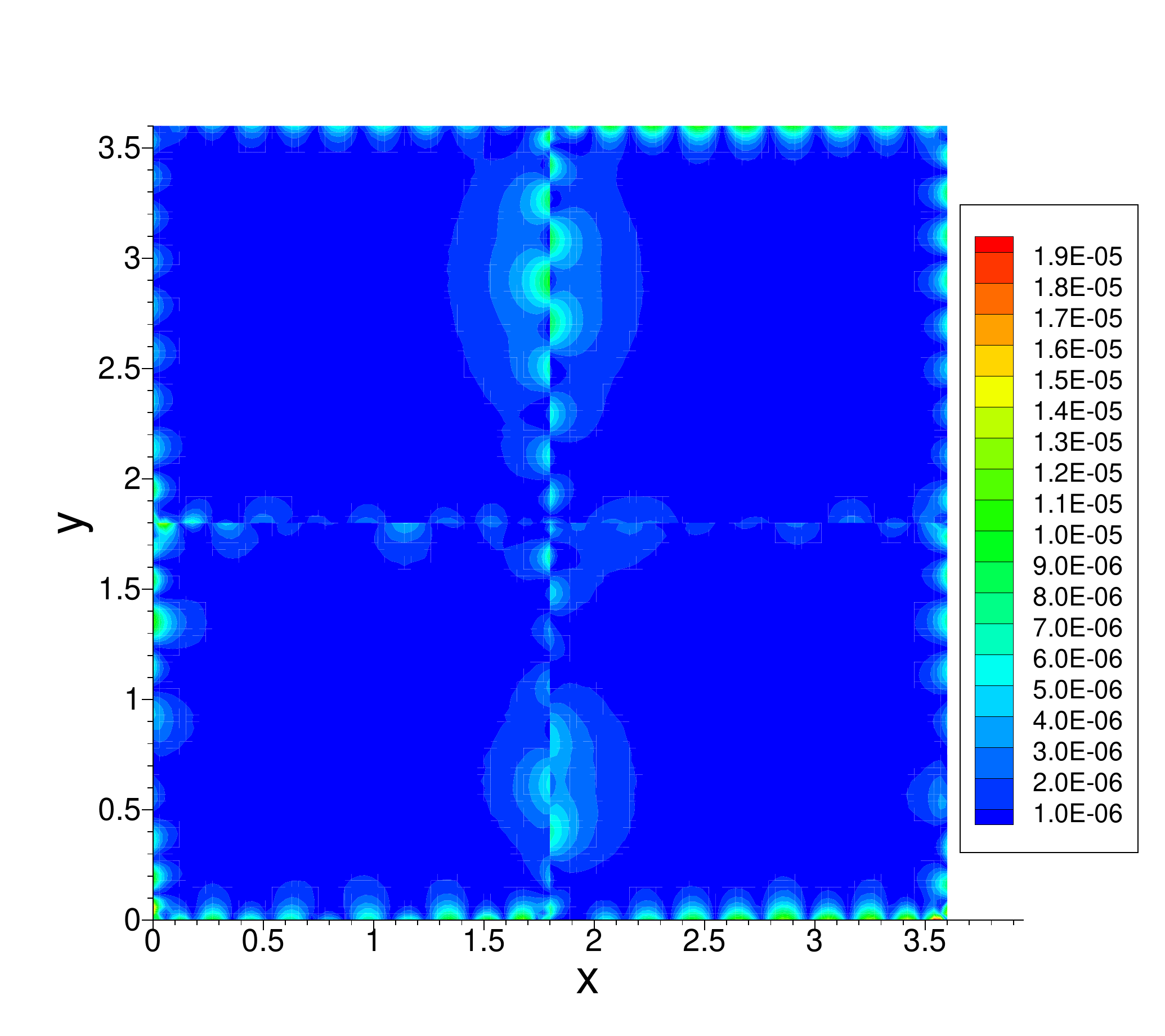}(d)
  }
  \caption{
    2D Helmholtz equation: (a) Solution distribution.
    Distributions of the absolute error obtained using one sub-domain (b)
    and using $4$ sub-domains (c) in locELM simulation.
  }
  \label{fig:helm2d_1}
\end{figure}

On each sub-domain, we employ a local neural network consisting of
an input layer with two nodes (representing
$x$ and $y$),
a single hidden layer with $M$ nodes and the $\tanh$ activation
function, and an output layer with one node (representing $u$).
The output layer is linear,
with no bias and no activation function.
An additional affine mapping normalizing the
input $x$ and $y$ data to the interval $[-1,1]\times[-1,1]$
has been incorporated
right behind the input layer for all sub-domains.
The weight and bias coefficients in the hidden layers
are set to uniform random values generated
on the interval $[-R_m,R_m]$.

The simulation parameters  include the number of sub-domains
($N_x$, $N_y$, $N_e$), the number of collocation points per sub-domain
($Q_x$, $Q_y$, $Q$), the number of training parameters per sub-domain ($M$),
and the maximum magnitude of the random coefficients ($R_m$).
We employ a fixed seed value $1$ for the Tensorflow random number generator
for all the tests in this subsection.

Figure \ref{fig:helm2d_1} depicts the field distributions of the
locELM solution (a), and the absolute errors of the solution computed
using one sub-domain (b) and $4$ uniform sub-domains (c).
The case with one sub-domain in locELM computation
is equivalent to the configuration
of a global ELM.
In this case, we have employed a total of $Q=50\times 50$
uniform collocation points in the domain, $M=1600$ training parameters
in the neural network, and $R_m=2.0$ when generating the random
weight/bias coefficients for the hidden layers.
For the case with $4$ sub-domains, 
we have partitioned the domain into $2$ sub-domains in each direction
($N_x=N_y=2$), and 
employed $Q=25\times 25$ uniform collocation points in each sub-domain
(i.e.~$Q_x=Q_y=25$), $M=400$ training parameters per sub-domain,
and $R_m=1.5$ when generating the random weight/bias coefficients.
Therefore the total degrees of freedom for these two cases,
i.e.~the number of collocation points and the number of training parameters
in the domain, are the same.
These simulations have both captured the solution accurately.
The resultant errors  are also
comparable, with the case
of $4$ sub-domains slightly better.

\begin{figure}
  \centerline{
    \includegraphics[width=2in]{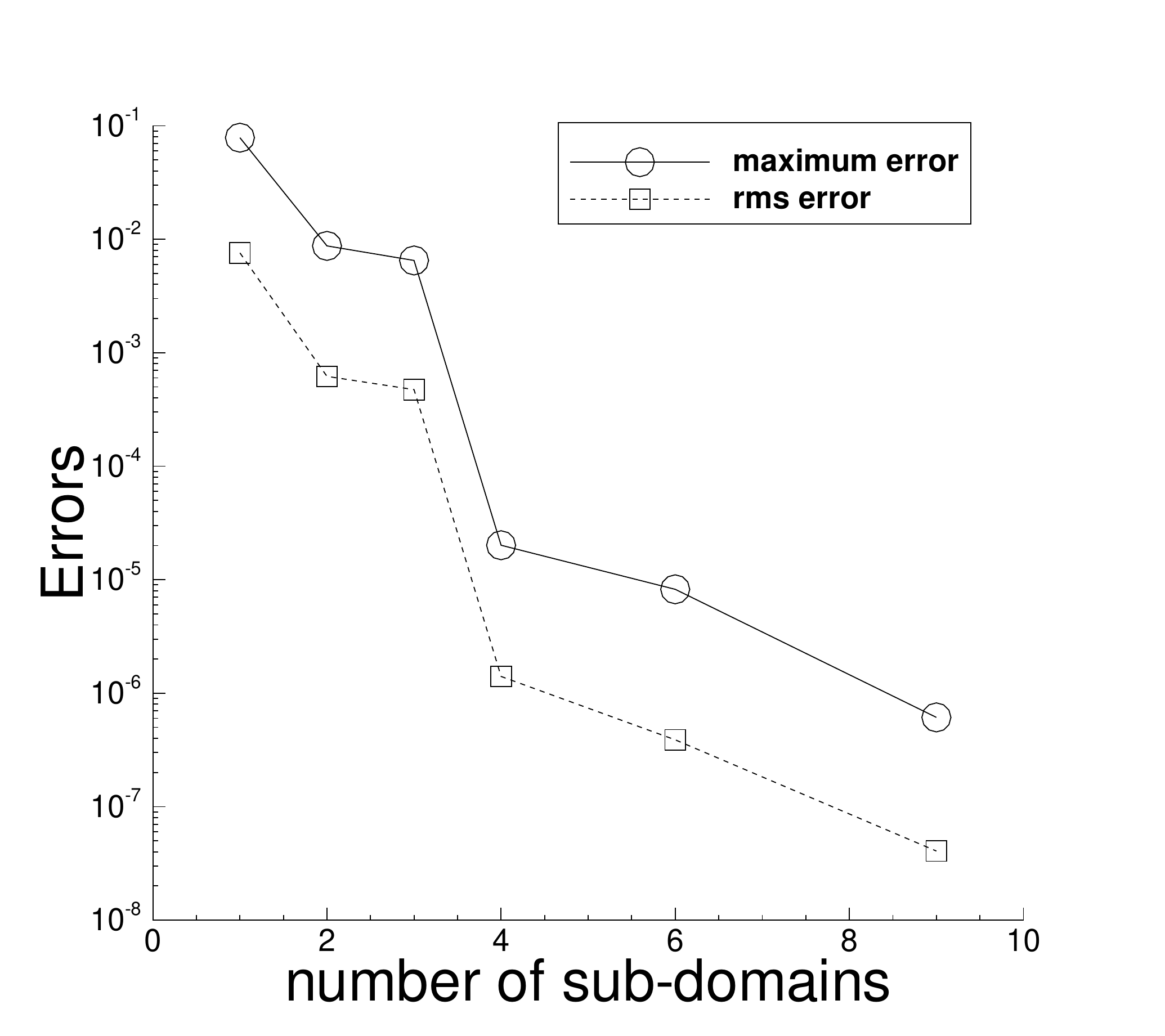}(a)
    \includegraphics[width=2in]{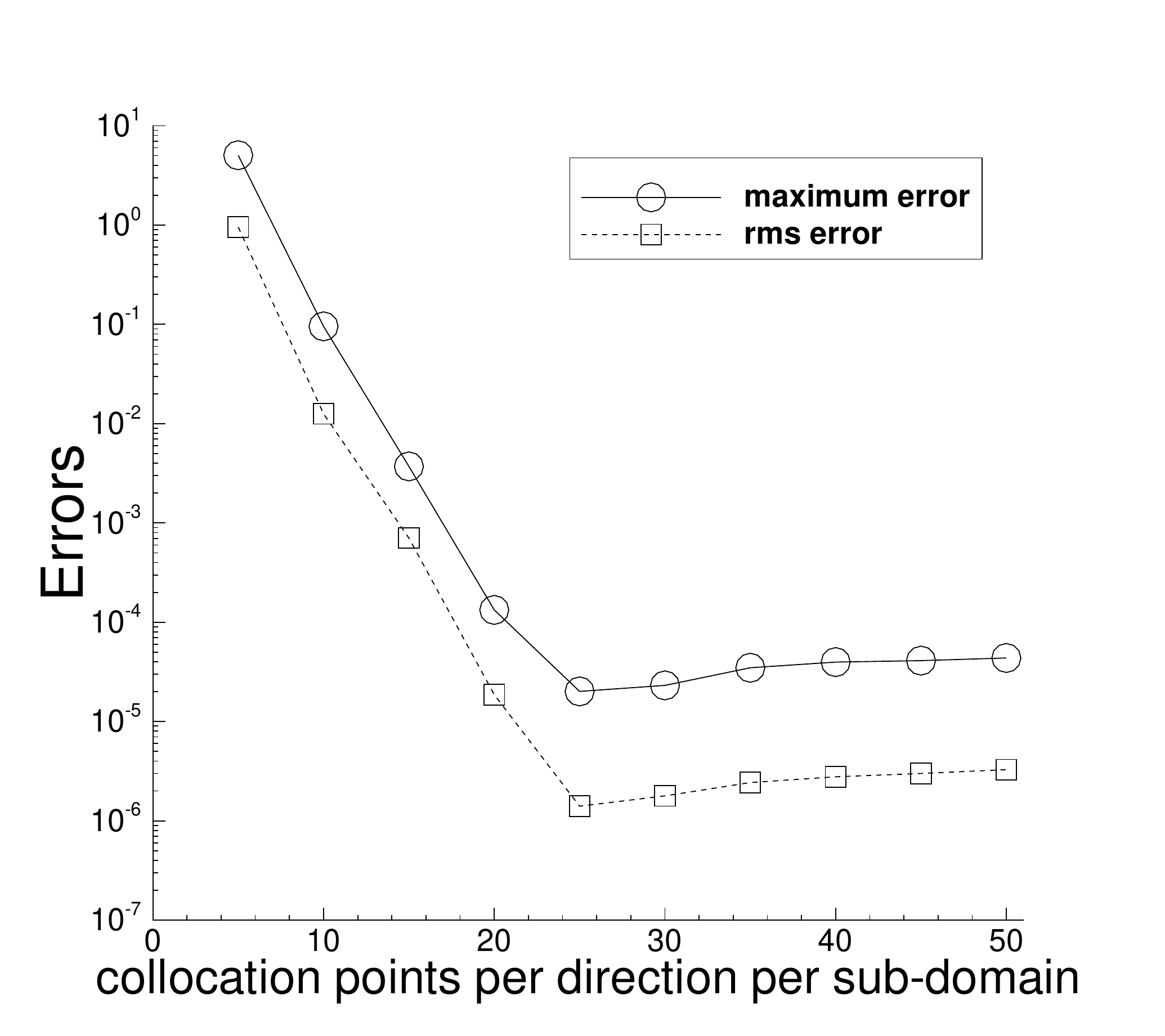}(b)
    \includegraphics[width=2in]{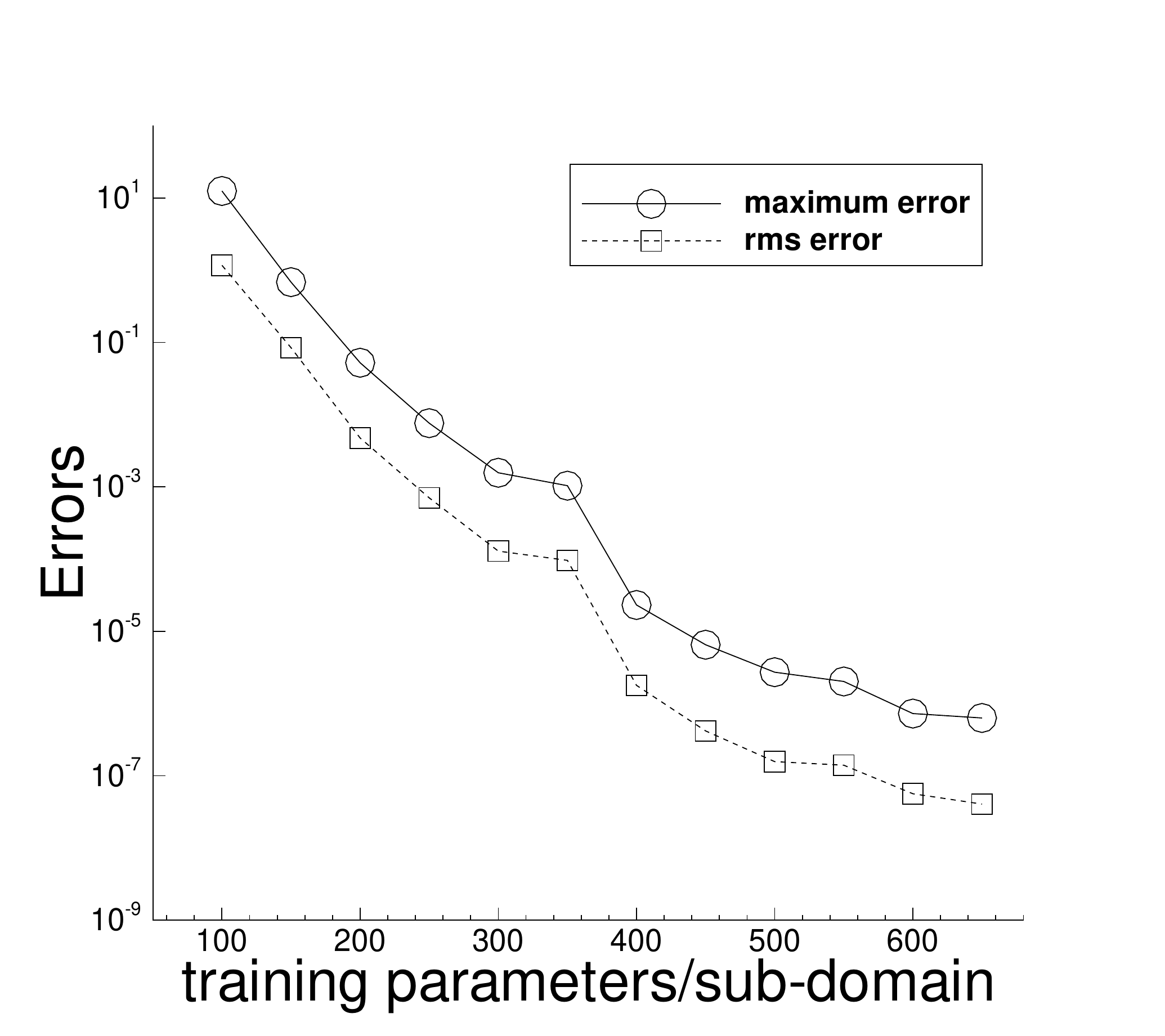}(c)
  }
  \caption{Effect of degrees of freedom
    (2D Helmholtz equation): the maximum/rms errors in the domain as a function
    of (a) the number of sub-domains,
    (b) the number of collocation points per direction per sub-domain,
    and (c) the number of
    training parameters per sub-domain.
    In (b,c) four sub-domains are used.
  }
  \label{fig:helm2d_3}
\end{figure}

Figure \ref{fig:helm2d_3} shows the effect of the degrees of freedom on
the locELM accuracy.
The effect of varying the number of sub-domains,
with the degrees of freedom per
sub-domain fixed, is illustrated in Figure \ref{fig:helm2d_3}(a).
In this group of tests, we employ $Q=25\times 25$ uniform collocation
points/sub-domain ($Q_x=Q_y=25$),
$M=400$ training parameters/sub-domain,
and $R_m=1.5$ for generating the random weight/bias coefficients.
The domain is partitioned into one to nine sub-domains.
Figure \ref{fig:helm2d_3}(a) shows the maximum and rms errors
in the domain versus the number of sub-domains in the locELM simulation.
We can observe that the  errors decrease
essentially exponentially
with increasing number of sub-domains.

The effect of varying the number of collocation points or
the number of training parameters per sub-domain is illustrated in
Figures \ref{fig:helm2d_3}(b) and (c).
In this group of tests we have employed four uniform sub-domains ($N_x=N_y=2$),
and $R_m=1.5$ when generating the random weight/bias
coefficients,
while the number of collocation points/sub-domain or the number of
training parameters/sub-domain is varied systematically.
Figure \ref{fig:helm2d_3}(b) depicts the maximum and rms errors in the domain
as a function of the number of uniform collocation points in each direction ($Q_x=Q_y$)
per sub-domain, where the number of training parameters per sub-domain
is fixed at $M=400$.
It can be observed that the errors initially decrease exponentially with
increasing collocation points, until the number of collocation points reaches
a certain level. Then the errors stagnate, and remains essentially
the same as the number of collocation points
further increases.
Figure \ref{fig:helm2d_3}(c) shows the maximum and rms errors in the domain
as a function of the training parameters per sub-domain, where the number of
uniform collocation points per sub-domain is fixed at $30\times 30$
($Q_x=Q_y=30$). It is observed that the errors decrease dramatically
with increasing training parameters per sub-domain. The error reduction rate
is nearly exponential initially, and the reduction slows down
as the number of training parameters/sub-domain becomes large.
These results show that
the locELM method exhibits a clear sense of convergence with increasing
degrees of freedom.
The locELM  errors decrease
essentially exponentially as the number of sub-domains, or the number of
collocation points in each direction
per sub-domain, or the number training parameters per sub-domain increases.

\begin{figure}
  \centerline{
    \includegraphics[width=1.5in]{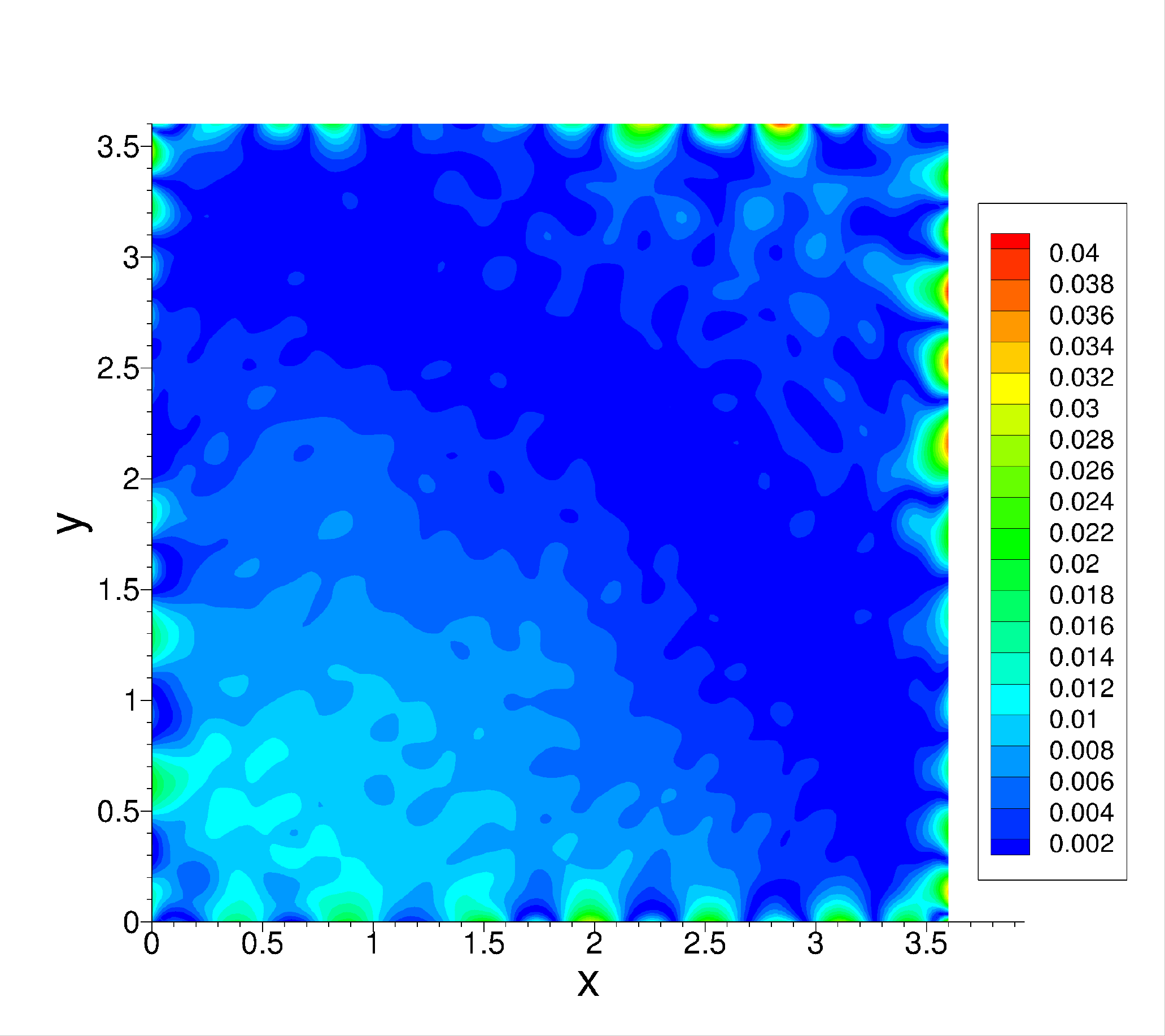}(a)
    \includegraphics[width=1.5in]{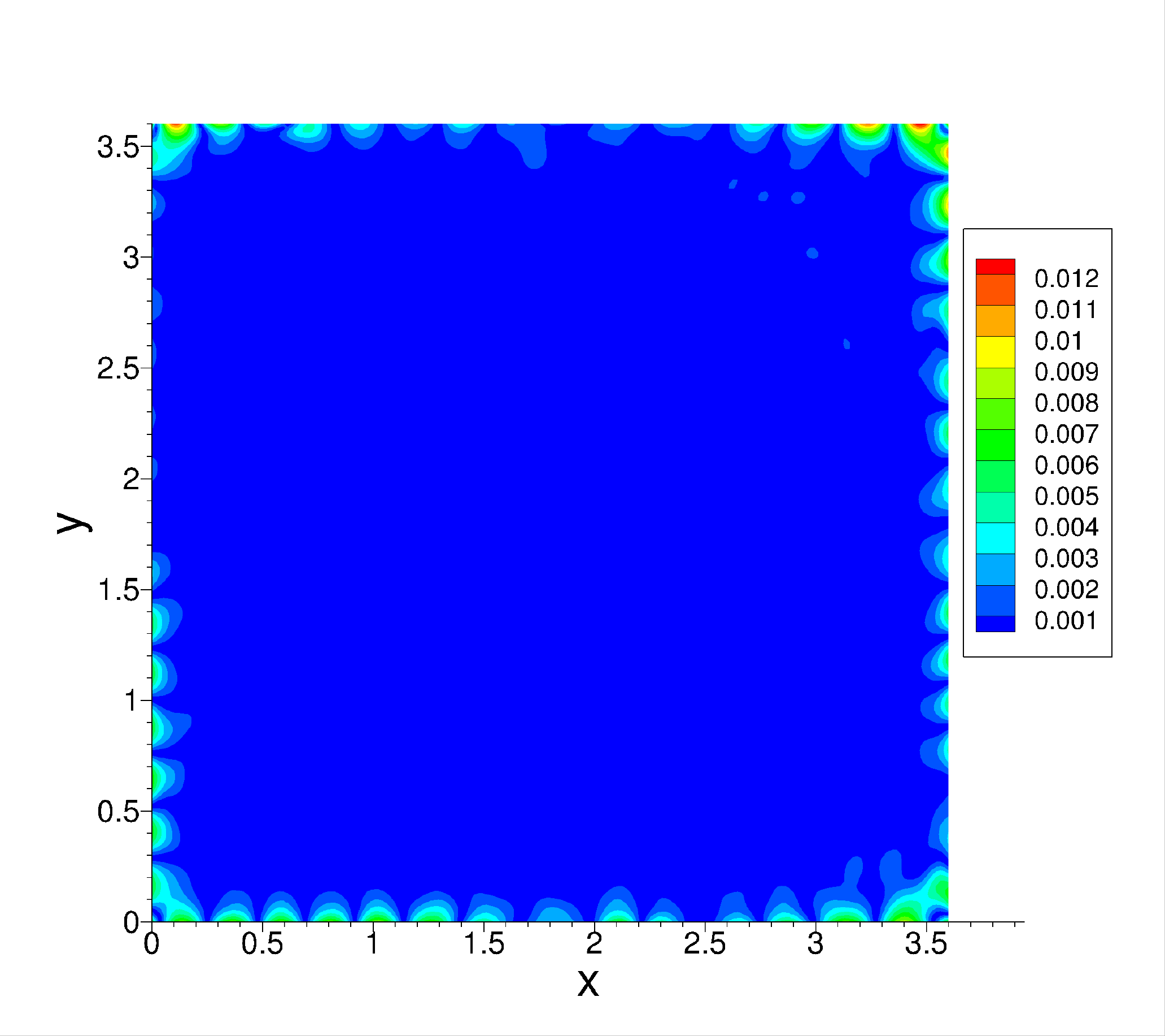}(b)
  }
  \caption{PINN solution of
    2D Helmholtz equation: Field distributions of 
    and the absolute errors obtained using
    PINN~\cite{RaissiPK2019} with the Adam optimizer (a)
    and the L-BFGS optimizer (b).
    This figure can be compared with Figure \ref{fig:helm2d_1}, which
    is obtained using the current locELM method.
  }
  \label{fig:helm2d_5}
\end{figure}

\begin{figure}
  \centerline{
    \includegraphics[width=1.5in]{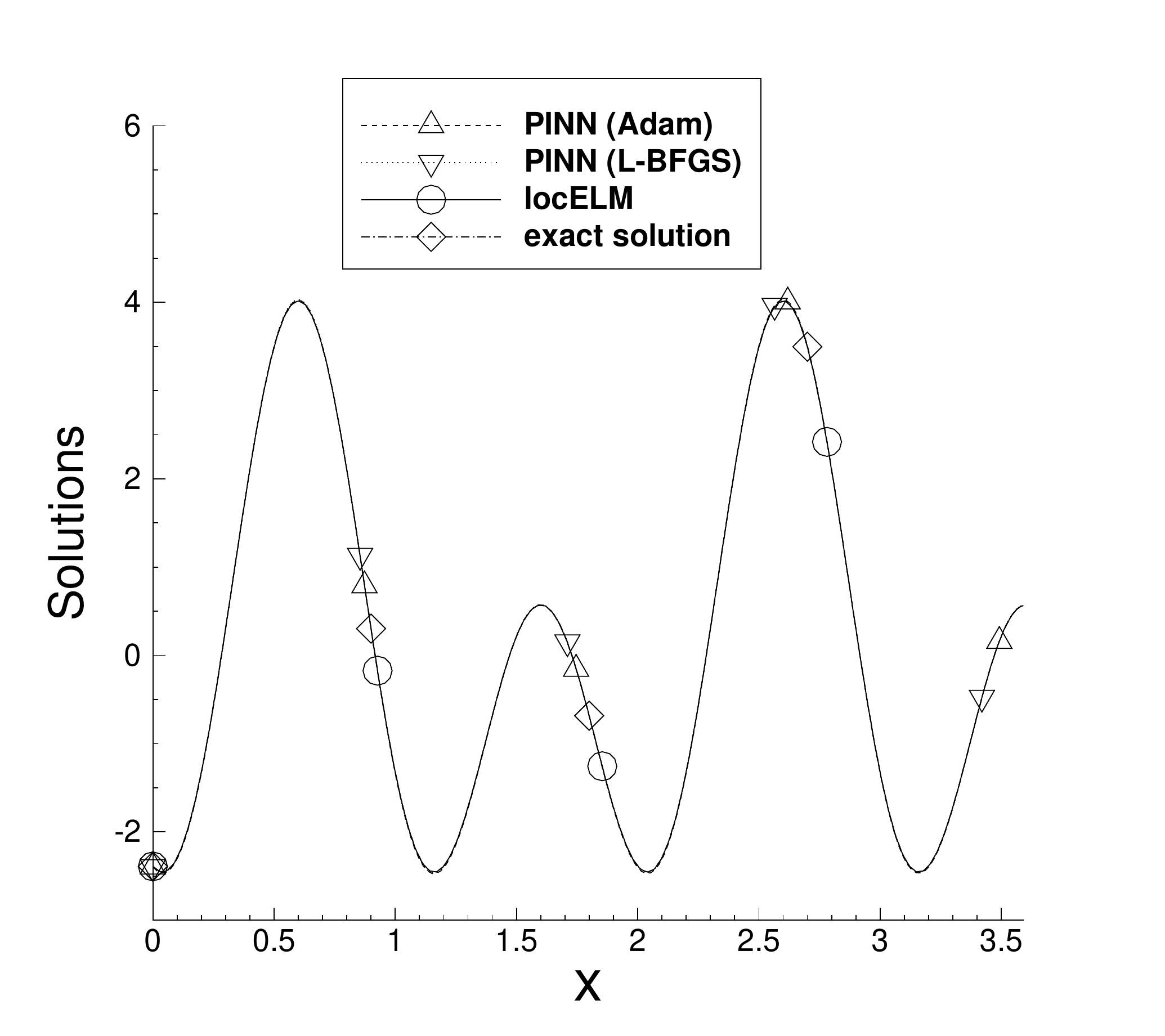}(a)
    \includegraphics[width=1.5in]{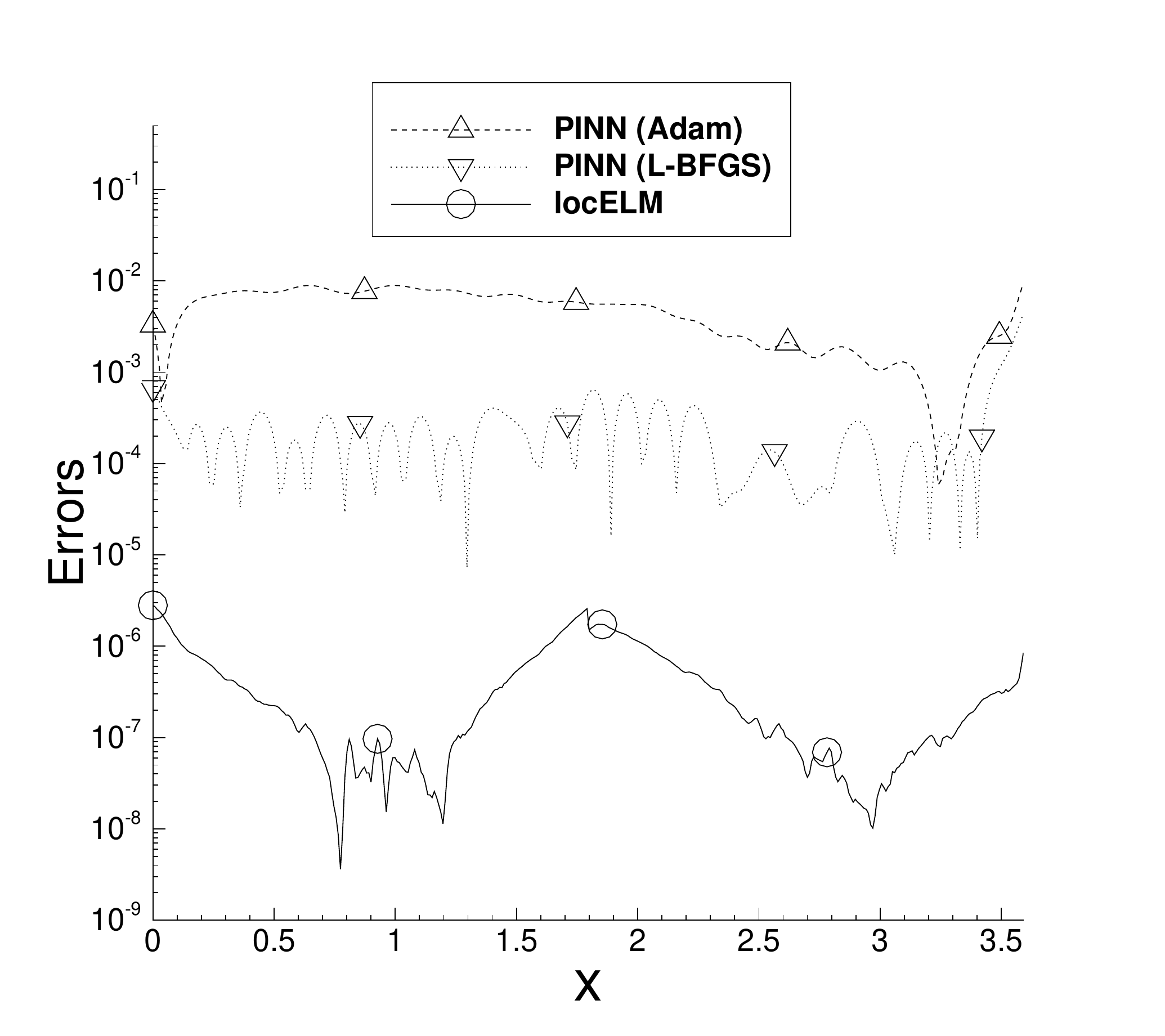}(b)
    \includegraphics[width=1.5in]{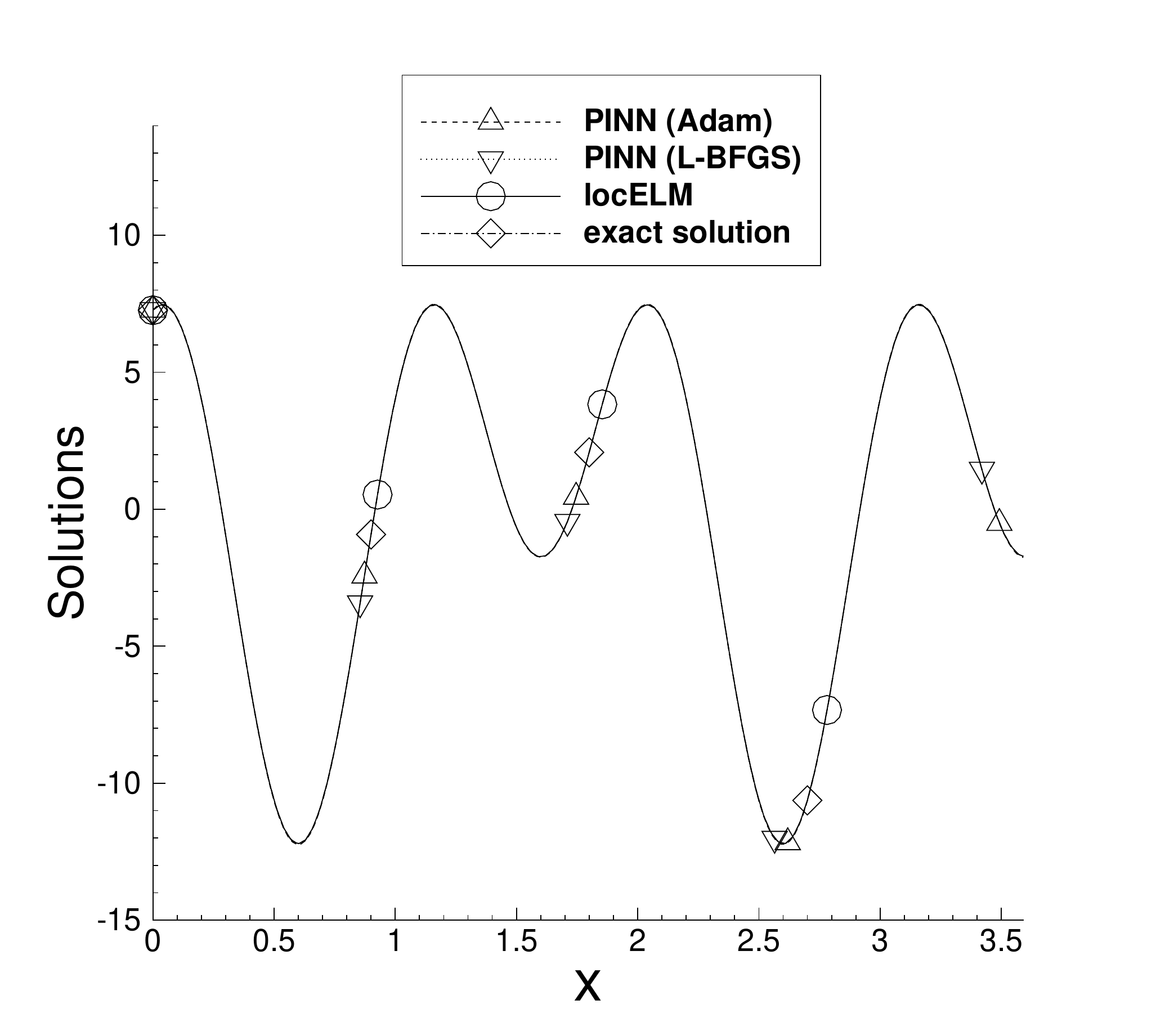}(c)
    \includegraphics[width=1.5in]{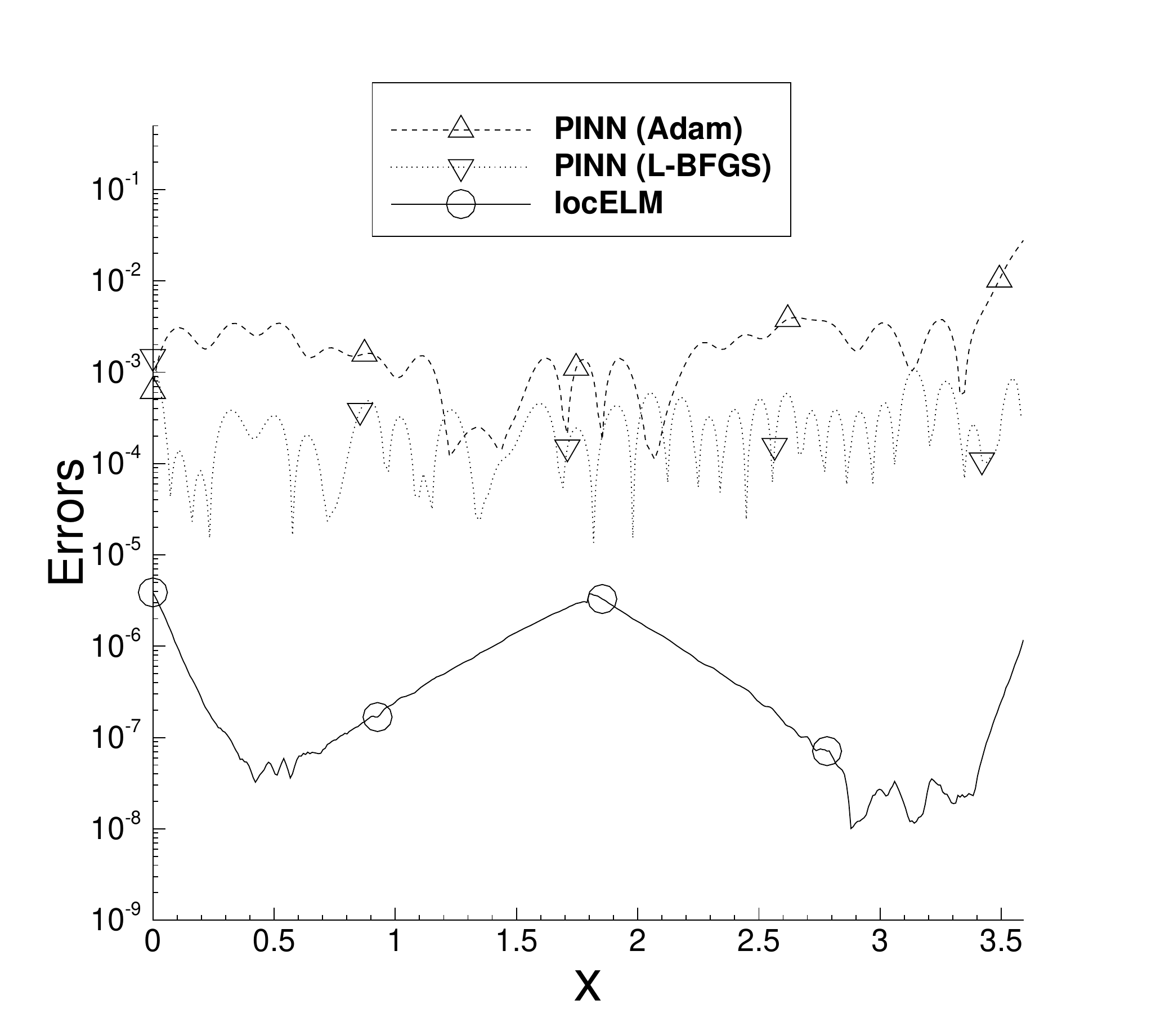}(d)
  }
  \caption{Comparison between locELM and PINN 
    (2D Helmholtz equation): profiles of the solutions (a,c)
    and their absolute errors (b,d) along two horizontal lines across
    the domain at $y=1.0$ (a,b) and $y=2.6$ (c,d),
    obtained using PINN~\cite{RaissiPK2019} with the Adam/L-BFGS optimizers and
    using locELM with 4 sub-domains.
    The PINN results correspond to Figure \ref{fig:helm2d_5}.
    The locELM results correspond to Figure \ref{fig:helm2d_1}(c,d).
  }
  \label{fig:helm2d_6}
\end{figure}

We next compare the locELM method  with
PINN~\cite{RaissiPK2019}
for solving the 2D Helmholtz equation.
Figure \ref{fig:helm2d_5} shows the field distributions of
the solutions and their errors obtained using PINN together with
the Adam and L-BFGS optimizers.
The input data consist of $50\times 50$ uniform points in the domain. 
For the Adam optimizer, the neural network contains $5$ hidden layers,
with a width of $40$ nodes for each layer and the $\tanh$ activation function.
The neural network has been trained on the input data for $107,000$ epochs, with
the learning rate gradually decreasing from $0.001$ at the beginning
to $10^{-5}$ at the end of training.
For the L-BFGS optimizer, the neural network contains $4$ hidden layers,
with a width of $50$ nodes for each layer and the $\tanh$ activation function.
The network has been trained on the input data for $24,500$ L-BFGS
iterations.
The L-BFGS result is observed to be generally more accurate
than the Adam result.
The results in this figure can be compared with those of Figure \ref{fig:helm2d_1},
which are obtained using the current locELM method.
The field distributions indicate that the locELM results 
are considerably more accurate than the PINN results.

Figure \ref{fig:helm2d_6} shows a comparison of the solution and 
error profiles along two horizontal lines across the domain,
at $y=1.0$ and $y=2.6$, obtained using PINN and using the current locELM
method (with $4$ sub-domains).
The superior accuracy of the current method is evident.

\begin{table}
  \centering
  \begin{tabular}{lllll}
    \hline
    method  & maximum error & rms error & epochs/iterations & training time (seconds) \\
    PINN (Adam)  & $4.13e-2$ & $3.91e-3$ & $107,000$ & $7213.7$ \\
    PINN (L-BFGS)  & $1.27e-2$ & $9.25e-4$ & $24,500$ & $5051.2$ \\
    global ELM  & $4.17e-5$ & $4.54e-6$ & $0$ & $410.6$ \\
    locELM (4 sub-domains)  & $2.01e-5$ & $1.41e-6$ & $0$ & $33.6$ \\
    \hline
  \end{tabular}
  \caption{2D Helmholtz equation: comparison between locELM and PINN 
    in terms of the accuracy (maximum/rms errors in domain) and the
    computational cost (epochs/iterations and the training time).
    The PINN results correspond to those in Figure \ref{fig:helm2d_5}.
    The locELM results correspond to those in Figure
    \ref{fig:helm2d_1}.
  }
  \label{tab:helm2d_1}
\end{table}

Table \ref{tab:helm2d_1} is a further comparison between the current locELM method
and PINN in terms of the maximum and rms errors in the domain,
the number of epochs or iterations during the training, and the network training time.
The PINN results with the Adam/L-BFGS optimizers correspond to those
from Figure \ref{fig:helm2d_5}.
The global ELM results correspond to those in Figure \ref{fig:helm2d_1}(a,b), and
the locELM results with 
$4$ sub-domains correspond to those in Figure \ref{fig:helm2d_1}(c,d).
Note that the input data for all these cases (PINN and locELM)
are the same (total $50\times 50$ uniform collocation points in the domain).
It can be observed that the locELM results are two to
three orders of magnitude more accurate than the PINN results.
In terms of the training time, the current locELM method is much faster than
PINN (by over two orders of magnitude).
The global ELM method, which is equivalent to the locELM method
with one sub-domain, is more than an order of magnitude faster than
PINN (around $410$ seconds
versus over $5000$ seconds). With $4$ sub-domains, the current locELM method
is also much faster than the global ELM (around $33$ seconds versus $410$ seconds).
These results clearly signify the superiority of the current
method to PINN, in term of both accuracy and the computational cost.


\begin{figure}
  \centerline{
    \includegraphics[height=2.1in]{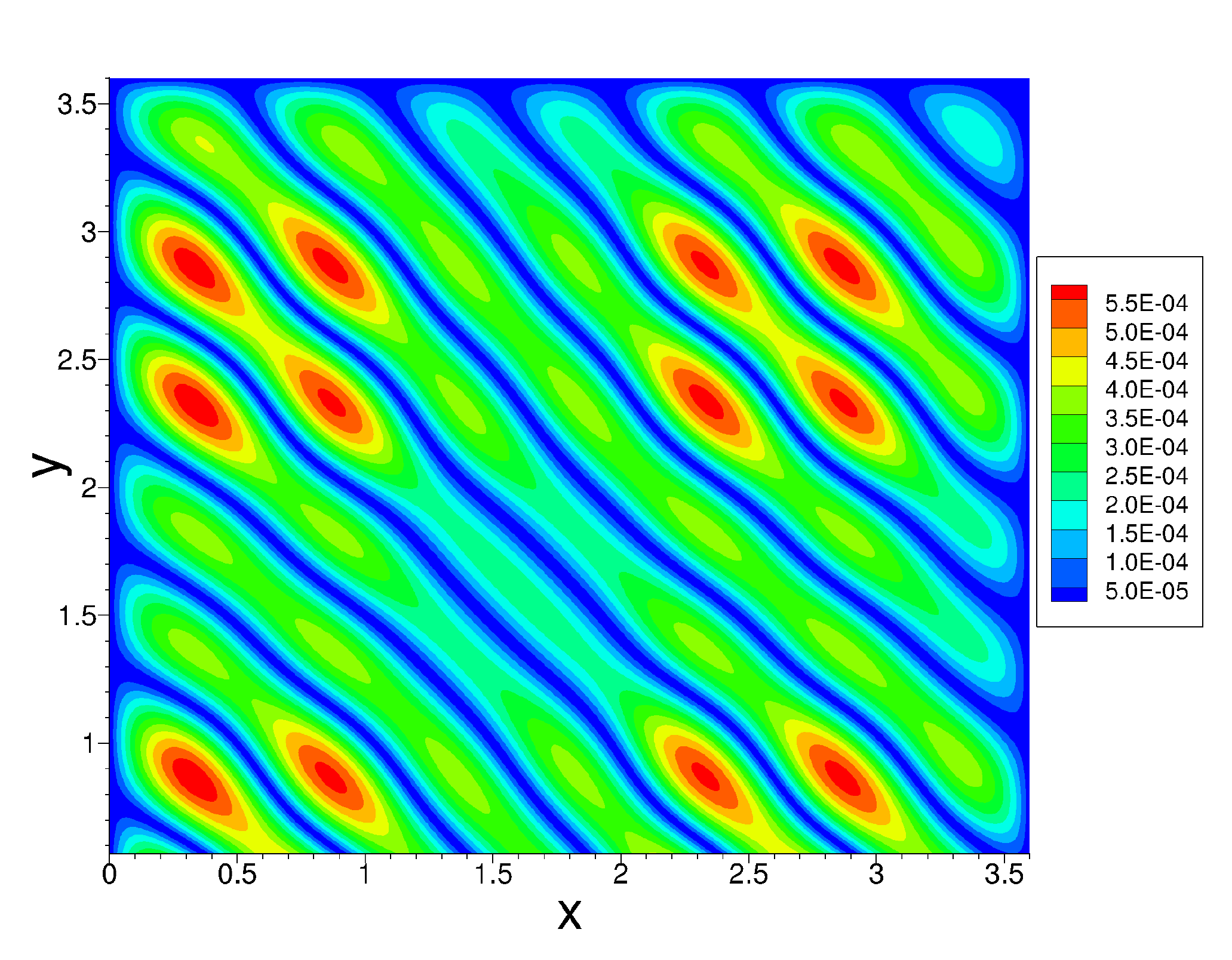}(a)
    \includegraphics[height=2.2in]{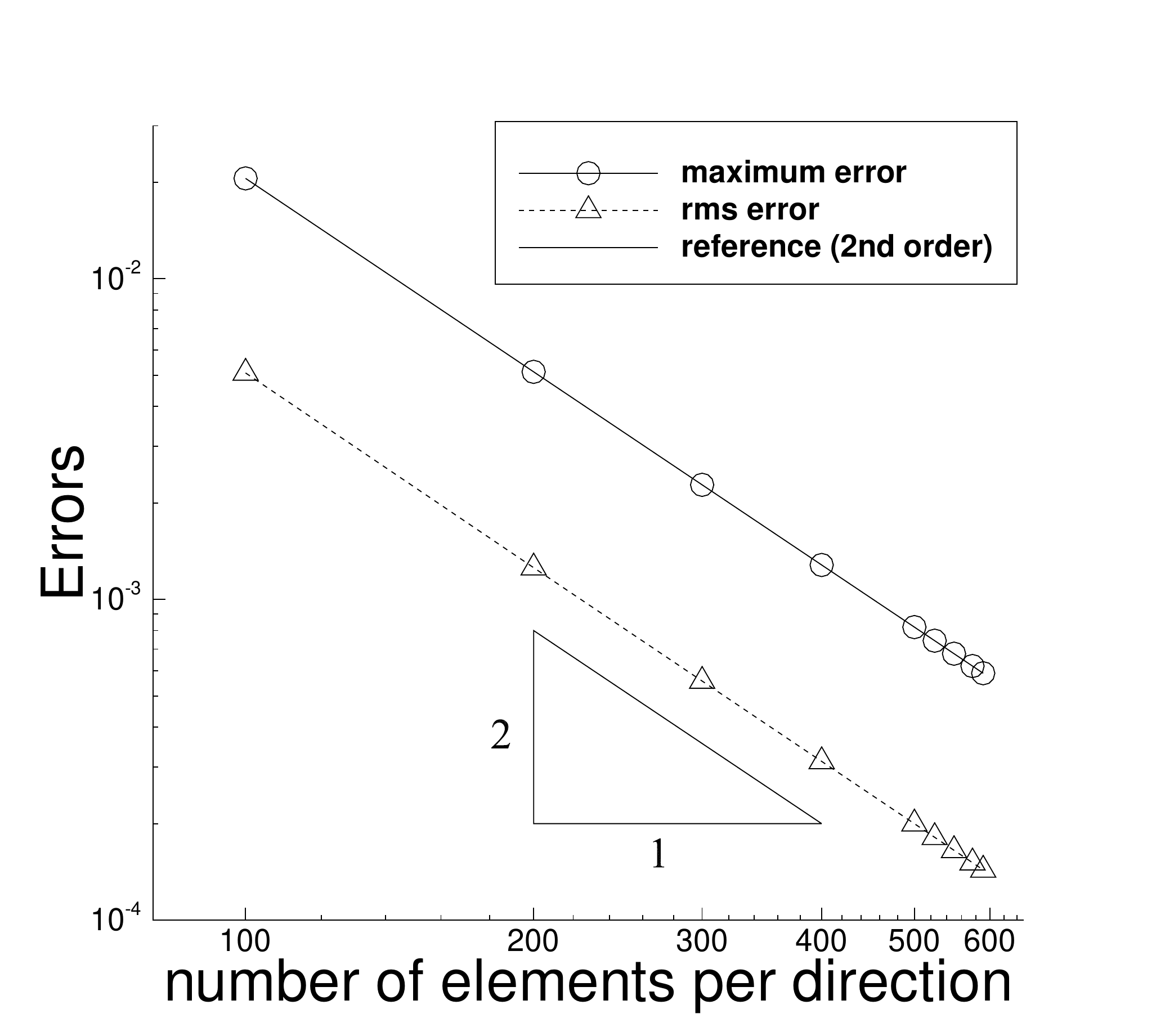}(b)
  }
  \caption{FEM solution of the 2D Helmholtz equation: (a) FEM error distribution
    computed on a $590\times 590$ mesh.
    (b) The FEM maximum/rms errors in the domain 
    versus the number of elements in each direction,
    showing the second-order convergence rate.
    On a $N\times N$ mesh, the number of triangular elements is $2N^2$.
  }
  \label{fg_helm2d_7}
\end{figure}

\begin{table}
  \centering
  \begin{tabular}{l|lllllll}
    \hline
    method & mesh & sub-domains & $Q$ & $M$  & max-error & rms-error &  wall-time (seconds) \\ \hline
    FEM & $500\times 500$ & -- & -- & -- & $8.20e-4$ & $2.00e-4$ & $18.5$ \\
    & $590\times 590$ & -- & -- & -- & $5.89e-4$ & $1.51e-4$ & $35.4$ \\
    \hline
    locELM & --  & $4$ & $20\times 20$ & $300$ & $7.28e-4$ & $5.28e-5$ & $17.1$ \\
     & --  & $4$ & $25\times 25$ & $400$ & $2.01e-5$ & $1.41e-6$ & $33.6$ \\
    \hline
  \end{tabular}
  \caption{
    2D Helmholtz equation: comparison between locELM and
    the finite element method (FEM)
    in terms of the maximum/rms errors in the domain and the
    training or computation time.
    The FEM results correspond to those in Figure \ref{fg_helm2d_7}.
    The locELM results correspond to those in Figure
    \ref{fig:helm2d_1}.
  }
  \label{tb_helm2d_2}
\end{table}

Finally we compare the performance of the current locELM method with
the classical finite element method for the 2D Helmholtz equation.
Figure \ref{fg_helm2d_7} illustrates the FEM solution and its
second-order convergence rate. 
Figure \ref{fg_helm2d_7}(a) shows the distribution of
the absolute error of
the FEM solution computed on a $590\times 590$ uniform rectangular mesh.
Note that
each rectangle in the mesh is further divided along its diagonal into two triangular
linear elements, as stipulated by the FEniCS library.
Therefore, in the current work,
an $N_1\times N_2$ rectangular mesh
contains a total of $2N_1N_2$ triangular elements
for the FEM simulations.
Figure \ref{fg_helm2d_7}(b) shows the maximum/rms errors
of the FEM solution in the domain
versus the number of rectangles in each direction
in the rectangular mesh, demonstrating the second-order convergence
rate of the method.

In Table \ref{tb_helm2d_2} we compare the locELM method
and the finite element method with regard to their
accuracy and the computational cost.
Here we list the maximum and rms errors in the domain,
and the wall time for the training or computation, corresponding to a set of
different meshes or simulation parameters obtained using locELM and
FEM. We observe that the locELM performance is on par with or better than
that of the FEM. For example,
the FEM case with the $500\times 500$ mesh and the locELM
case with $Q=20\times 20$ and $M=300$
have a comparable computational cost and also a similar accuracy.
The FEM case with the $590\times 590$ mesh has a comparable computational
cost to the locELM case with $Q=25\times 25$ and $M=400$,
but its errors are  over an order of
magnitude larger than those of the latter.



\bibliographystyle{plain}
\bibliography{elm,mypub,dnn,sem,obc}

\end{document}